\documentclass[11pt, reqno]{amsart}
\makeatletter
\usepackage[utf8]{inputenc}
\usepackage[T1]{fontenc}
\usepackage{lmodern}

\usepackage{textcmds}  
\usepackage{amsmath, amssymb, amsfonts, amstext, verbatim, amsthm, mathrsfs}
\usepackage[mathcal]{eucal}
\usepackage{microtype}
\usepackage{graphicx}
\usepackage[all]{xy}
\usepackage{tikz} 
\usetikzlibrary{shapes,arrows}
\usepackage{stmaryrd}

\usepackage{enumerate}
\usepackage{enumitem}
\usepackage{xspace}

\usepackage[colorlinks=true,linkcolor=blue,citecolor=blue,urlcolor=blue,citebordercolor={0 0 1},urlbordercolor={0 0 1},linkbordercolor={0 0 1}]{hyperref} 
\usepackage[shortalphabetic]{amsrefs} 
\usepackage[nameinlink]{cleveref}

\usepackage[nopostdot,sort=use,automake]{glossaries}
\makeglossaries
\usepackage{imakeidx}
\makeindex

\renewcommand*{\glossarymark}[1]{}

\newglossaryentry{filtered_objects}
{
    name={stack of filtered objects},
    description={The stack parameterizing algebraic families of maps $\Theta \to \X$.},
    symbol={$\Filt(\X)$},
    sort=f
}

\newglossaryentry{graded_objects}
{
    name={stack of graded objects},
    description={The stack parameterizing algebraic families of maps $B\Gm \to \X$.},
    symbol={$\Grad(\X)$},
    sort=g
}

\newglossaryentry{flag_space}
{
    name={flag space},
    description={An algebraic space parameterizing algebraic families of filtrations of a given point $\xi : \Spec(R) \to \X$. It typically has many connected components.},
    symbol={$\Flag(\xi)$},
    sort=0
}

\newglossaryentry{formal_fan}
{
    name={formal fan},
    description={A structure analogous to the fans studied in toric geometry, except cones are not embedded in $\bR^n$, and they can meet along arbitrary rational polyhedral subcones rather than just faces. The definition is analogous to that of a semisimplicial set.},
    symbol={$F_\bullet, G_\bullet, \ldots$},
    sort=a
}

\newglossaryentry{degeneration_fan}
{
    name={degeneration fan},
    description={Formal fan which encodes all possible filtrations of an object in a moduli problem.},
    symbol={$\Deg(\X,p)$ or $\Deg(p)$ for $p \in \X(k)$},
    sort=b
}

\newglossaryentry{degeneration_space}
{
    name={degeneration space},
    description={A topological space parameterizing filtrations of a given point in a moduli problem. It is the ``projective geometric realization" of the degeneration fan.},
    symbol={$\iDeg(\X,p)$ or $\iDeg(p)$ for $p \in \X(k)$},
    sort=c
}

\newglossaryentry{component_fan}
{
    name={component fan},
    description={Formal fan which encodes filtered objects in a stack $\X$ up to algebraic equivalence, meaning two filtrations are equivalent if they occur as points in an algebraic family of filtrations over a connected base scheme.},
    symbol={$\Comp(\X)$},
    sort=d
}

\newglossaryentry{component_space}
{
    name={component space},
    description={A topological space parameterizing filtered objects in a moduli problem $\X$ up to algebraic equivalence. It is the projective geometric realization of the component fan.},
    symbol={$\iComp(\X)$},
    sort=e
}

\newglossaryentry{theta_reductive}
{
    name={$\Theta$-reductive},
    description={A property of a stack $\X$. It means the flag spaces of $\X$ satisfy the valuative criterion for properness.},
    sort=f
}

\newglossaryentry{numerical_invariant}
{
    name={numerical invariant},
    description={The data provided for a stack $\X$ which allows one to formulate the question of existence and uniqueness of Harder-Narasimhan filtrations.},
    symbol={$\mu : \cU \to \bR$ with $\cU \subset \iComp(\X)$},
    sort=g
}

\newglossaryentry{standard_numerical_invariant}
{
	name={standard numerical invariant},
	description={A numerical invariant which satisfies several additional hypotheses which are useful in practice and hold for the numerical invariants arising in GIT.},
	sort=h
}

\def\makeCal#1{%
\expandafter\newcommand\csname c#1\endcsname{\mathcal{#1}}}
\def\makeBB#1{%
\expandafter\newcommand\csname b#1\endcsname{\mathbb{#1}}}
\def\makeFrak#1{%
\expandafter\newcommand\csname f#1\endcsname{\mathfrak{#1}}}

\count@=0
\loop
\advance\count@ 1
\edef\y{\@Alph\count@}%
\expandafter\makeCal\y
\expandafter\makeBB\y
\expandafter\makeFrak\y
\ifnum\count@<26
\repeat

\theoremstyle{plain}
\newtheorem*{thm*}{Theorem}
\newtheorem{thm}{Theorem}[subsection]
\newtheorem{cor}[thm]{Corollary}
\newtheorem{lem}[thm]{Lemma}

\newtheorem{prop}[thm]{Proposition}

\newtheorem{simplification}[thm]{Simplification}

\newtheorem{thmx}{Theorem}

\theoremstyle{definition}
\newtheorem{rem}[thm]{Remark}
\newtheorem{defn}[thm]{Definition}
\newtheorem{notn}[thm]{Notation}
\newtheorem{const}[thm]{Construction}
\newtheorem{ex}[thm]{Example}

\newtheorem{problem}[thm]{Problem}
\newtheorem{warning}[thm]{Warning}

\newcommand{\id}{\mathop{{\rm id}}\nolimits}
\newcommand{\inner}[1]{{\underline{#1}}}
\newcommand{\dual}{\vee}
\newcommand{\colim}{\mathop{\operatorname{colim}}}
\newcommand{\Map}{\op{Map}}
\newcommand{\iMap}{\inner{\op{Map}}}

\newcommand{\op}[1]{\!\!\mathop{\rm ~#1}\nolimits}

\newcommand{\pt}{\mathrm{pt}}
\newcommand{\stack}[1]{\op{St}_{#1}}
\DeclareMathOperator{\Spec}{Spec}
\DeclareMathOperator{\Isom}{Isom}

\newcommand{\Gm}{\mathbb{G}_m}
\newcommand{\GL}{\mathrm{GL}}

\newcommand{\SL}{\mathrm{SL}}
\newcommand{\fg}{\mathfrak{g}}
\newcommand{\fp}{\mathfrak{p}}
\DeclareMathOperator{\Stab}{Stab}

\DeclareMathOperator{\D}{D}

\DeclareMathOperator{\Coh}{Coh}
\DeclareMathOperator{\QCoh}{QCoh}
\DeclareMathOperator{\Perf}{Perf}
\DeclareMathOperator{\APerf}{APerf}
\newcommand{\Mod}{\operatorname{-Mod}}
\DeclareMathOperator{\Hom}{Hom}
\newcommand{\twist}[1]{\langle #1 \rangle}
\DeclareMathOperator{\Aut}{Aut}
\DeclareMathOperator{\End}{End}

\DeclareMathOperator{\coker}{coker}

\newcommand{\ev}{\mathop{{\rm ev}}\nolimits}
\DeclareMathOperator{\agr}{gr}

\DeclareMathOperator{\ch}{ch}
\DeclareMathOperator{\gr}{gr}

\newcommand{\Knum}{\cN}
\newcommand{\simp}[1]{\Delta_{#1}}
\newcommand{\ray}[1][\bullet]{R_{#1}}
\newcommand{\X}{\fX}
\newcommand{\Y}{\fY}
\newcommand{\Z}{\fZ}
\renewcommand{\S}{\fS}
\newcommand{\gitmod}{/\!/}
\newcommand{\canon}{\mathfrak{c}}
\newcommand{\one}{\mathbf{1}}

\newcommand{\Cones}{\mathfrak{Cone}}

\newcommand{\iComp}{\mathop{{\mathscr Comp}}\nolimits}
\newcommand{\Comp}{\mathop{\mathbf{CF}}\nolimits}
\newcommand{\Deg}{\mathop{\mathbf{DF}}\nolimits}
\newcommand{\iDeg}{\mathop{{\mathscr Deg}}\nolimits}

\newcommand{\ST}{\overline{\mathop{ST}}}

\newcommand{\filt}[2][]{\iMap(\Theta^{#1}, {#2})}
\newcommand{\grad}[2][]{\iMap(\pt/\Gm^{#1}, {#2})}
\DeclareMathOperator{\Filt}{Filt}
\DeclareMathOperator{\Grad}{Grad}
\DeclareMathOperator{\Flag}{Flag}

\DeclareRobustCommand{\SkipTocEntry}[5]{}

\begin{document}


\title{On the structure of instability in moduli theory}
\author{Daniel Halpern-Leistner}
\date{\today}

\begin{abstract}
We formulate a theory of instability and Harder-Narasimhan filtrations for an arbitrary moduli problem in algebraic geometry. We introduce the notion of a $\Theta$-stratification of a moduli problem, which generalizes the Kempf-Ness stratification in GIT as well as the Harder-Narasimhan stratification of the moduli of coherent sheaves on a projective scheme. Our main theorems establish necessary and sufficient conditions for the existence of these stratifications. We define a structure on an algebraic stack called a numerical invariant, and we show that in many situations a numerical invariant defines a $\Theta$-stratification on the stack, assuming a certain ``HN boundedness'' condition holds. We also discuss criteria under which the semistable locus has a moduli space. We apply our methods to an example that lies beyond the reach of classical methods: the stratification of the stack of objects in the heart of a Bridgeland stability condition.
\end{abstract}

\maketitle


\tableofcontents

One fundamental observation in algebraic geometry is that small deformations of an algebro-geometric object are often controlled by a finite dimensional parameter space that is itself an algebraic variety. A moduli problem seeks to globalize this observation, i.e., find a geometric object that parameterizes all isomorphism classes of a particular type of algebro-geometric object. The pursuit of a general framework for studying moduli problems has introduced quite a bit of machinery and has ultimately lead to the theory of algebraic stacks.

The result is that many moduli problems have a tautological solution: when formulated correctly, the moduli problem \emph{is} and algebraic stack, which is a geometric object. Unfortunately, the stacks that arise in nature are rarely representable by algebraic varieties or even algebraic spaces, which limits their applicability. Many natural examples, such as the moduli stack of vector bundles of a fixed rank and degree on a smooth curve, are highly non-separated and are not bounded, i.e., the algebraic stack that represents the moduli problem is not quasi-compact.

Sometimes these issue are mild, and the stack admits a moduli space. The first general result in this direction was the Keel-Mori theorem \cite{KeelMori}, which constructs coarse moduli spaces under mild hypotheses for stacks in which objects have finite automorphism groups. The notion of a good moduli space for an algebraic stack was introduced in \cite{alper2013good} as a suitable notion of moduli space for moduli problems with infinite automorphism groups, and generalizes the kinds of moduli spaces constructed using geometric invariant theory in characteristic $0$. Recently, necessary and sufficient conditions were discovered for an algebraic stack to have a good moduli space \cite{AHLH}.

In this paper, we focus on moduli problems that do not necessarily admit moduli spaces. We propose a new standard for what it means to ``solve" a moduli problem in algebraic geometry: in addition to defining an algebraic stack $\X$ that represents the moduli problem, one should construct a certain type of stratification of $\X$, which we call a \emph{$\Theta$-stratification}. Often the open piece of the stratification admits as good moduli space, as do the ``centers'' of the other strata, so this paper is complementary to \cite{AHLH}, and together these papers provide a very general framework for studying moduli problems. See \Cref{fig:flow_chart} below.

\begin{ex} \label{E:vector_bundles}
The prototypical example is the moduli stack $\cM_{r,d}$ of vector bundles of rank $r$ and degree $d$ on a smooth Riemann surface $C$. $\cM_{r,d}$ is not quasi-compact or separated. There is an open substack $\cM_{r,d}^{\rm{ss}} \subset \cM_{r,d}$, parameterizing \emph{semistable} vector bundles, that admits a projective good moduli space $\cM_{r,d}^{\rm{ss}} \to M$, which is a coarse moduli space if $r$ and $d$ are coprime. On the other hand, every unstable vector bundle (meaning one which is not semistable) has a canonical Harder-Narasimhan filtration \cite{HN}, and the unstable locus in the moduli stack admits a stratification $\cM_{r,d}^{us} = \bigcup \S_\alpha$, where $\S_\alpha$ is the locally closed substack of bundles for which the associated graded pieces of the Harder-Narasimhan filtration have particular rank and degree recorded by the index $\alpha$ \cite{Sh77}. The $\S_\alpha$ fiber over moduli stacks $\Z_\alpha^{\rm{ss}}$ of semistable graded vector bundles, and the $\Z_\alpha^{\rm{ss}}$ again have projective good moduli spaces.
\end{ex}

Our definition of a $\Theta$-stratification (\Cref{defn:theta_stratification}) gives an intrinsic description of this structure. The stack $\Theta := \bA^1 / \bG_m$ plays a central role. If $k$ is a field, then we define a \emph{filtration} of a $k$-point $p \in \X(k)$ to be a map $f : \Theta_k \to \X$ along with an isomorphism $f(1) \simeq p$. More generally, for any map $\xi : \Spec(R)  \to \X$, we define a filtration of $\xi$ to be a map $f : \Theta_R \to \X$ along with an isomorphism $f|_{\{1\} \times \Spec(R)} \cong \xi$.

\begin{ex} \label{E:vector_bundles_2}
A map $\Theta_k \to \cM_{r,d}$ is a $(\bG_m)_k$-equivariant vector bundle on $\bA^1 \times C$. Using the Rees construction, one can show that this is equivalent to a $\bZ$-indexed diagram $\cdots \subset E_{w+1} \subset E_w \subset \cdots$ of vector bundles on $C$ such that $E_w / E_{w+1}$ is locally free, $E_w = 0$ for $w \gg 0$, and $E_w$ stabilizes for $w \ll 0$. This justifies the terminology ``filtration.''
\end{ex}

A weak $\Theta$-stratification (\Cref{defn:theta_stratification})\footnote{The distinction between a weak $\Theta$-stratification and the slightly stronger notion of a $\Theta$-stratification is not relevant in this introduction. In characteristic $0$, any weak $\Theta$-stratification is a $\Theta$-stratification by \Cref{cor:char_0_theta_strat}.} is a stratification of $\X$ into locally closed substacks that parameterize points in $\X$ along with a canonical filtration. In order to generalize the canonical Harder-Narasimhan filtrations of \Cref{E:vector_bundles}, we introduce the following:

\begin{defn} \label{defn:numerical_invariant_simple}
A \emph{numerical invariant}\footnote{We will give a slightly more general definition in \Cref{defn:numerical_invariant} that specializes to this one in most cases of interest. See \Cref{rem:alternate_numerical}.} $\mu$ with values in a totally ordered set $\Gamma$ consists of an assignment of a scaling-invariant function $\mu_\gamma : \bR^n \setminus 0 \to \Gamma$ for any $p \in \X(k)$ and any homomorphism of $k$-groups $\gamma: (\bG_m^n)_k \to \Aut(p)$ with finite kernel. This data must satisfy:
\begin{enumerate}
\item If $\gamma' : (\bG_m^n)_{k'} \to \Aut(p_{k'})$ is the extension of $\gamma$ to a field extension $k \subset k'$, then $\mu_\gamma = \mu_{\gamma'}$.
\item Given a homomorphism with finite kernel $\phi : (\bG_m^q)_k \to (\bG_m^n)_k$, $\mu_{\gamma \circ \phi}$ is the restriction of $\mu_\gamma$ along the inclusion $\bR^q \hookrightarrow \bR^n$ induced by $\phi$; and
\item Given a scheme $S$, an $S$-point $\xi : S \to \X$, and a homomorphism of $S$-group-schemes $\gamma : (\bG_m^n)_S \to \Aut(\xi)$ with finite kernel, if one considers the maps $\gamma_s : (\bG_m^n)_{k(s)} \to \Aut(\xi(s))$ indexed by points in $s \in S$, then the function $\mu_{\gamma_s}$ on $\bR^n \setminus 0$ is locally constant on $S$.
\end{enumerate}
\end{defn}

A numerical invariant defines a $\Gamma$-valued function on the set of non-trivial filtrations in $\X$ by assigning $\mu(f) = \mu_\gamma(1)$, where $\gamma : (\bG_m)_k \to \Aut(f(0))$ is the canonical homomorphism induced by $f : \Theta_k \to \X$. We will always assume that $\Gamma$ comes with a marked point $0 \in \Gamma$.
\begin{defn}[\Cref{defn:HN_filtration}]
A point $p \in |\X|$ is said to be unstable if there is a filtration $f$ such that $\mu(f)>0$ and $f(1)=p$. A filtration $f: \Theta_k \to \X$ is defined to be an \emph{HN filtration} of the point $p=f(1)\in |\X|$ if
\[
\mu(f) = \sup \{ \mu(f') | f' \text{ is a filtration with } f'(1)=p\text{ in } |\X|\},
\]
\end{defn}

\begin{ex} \label{E:GIT}
Geometric invariant theory (GIT) is the special case of our theory in which $\X = X/G$ is a quotient stack for the action of a reductive group $G$ on a scheme $X$ that is projective over its affinization. In this context, the notion of a numerical invariant generalizes the normalized Hilbert-Mumford numerical invariant \cite{MFK94}*{Chap.~2.2}. The function $\mu_\gamma(x)$ associated to a homomorphism with finite kernel $\gamma : (\Gm^n)_k \to \Aut_\X(p) = \Stab_G(p)$ has the form $\ell(x) / \sqrt{b(x)}$, where $\ell(x)$ is a linear function on $\bR^n$ and $b(x)$ is a positive definite rational quadratic form. $\ell$ is the character of the action of $(\Gm^n)_k$ on the fiber $L_p$ of some $G$-ample line bundle on $X$, and $b$ comes from a choice of Weyl-invariant inner product on the coweight lattice of $G$. The HN filtration of an unstable point corresponds to Kempf's optimally destabilizing one-parameter-subgroup \cite{Ke78}, and the $\Theta$-stratifications that we construct agree with the Hesselink-Kempf-Kirwan-Ness stratification \cites{Ness84,Ki84,He79}. See \Cref{ex:git_numerical_invariant} and \Cref{defn:induced_invariant}.
\end{ex}

Given a numerical invariant $\mu$ on a stack $\X$, one can ask if it defines a (weak) $\Theta$-stratification on $\X$. By this we mean there is a (weak) $\Theta$-stratification of $\X$ such that the canonical filtration $f$ of any unstable point $p \in |\X|$ is an HN filtration, and the HN filtration is unique up to scaling, i.e., precomposing $f : \Theta_k \to \X$ with a ramified covering $(-)^n : \Theta_k \to \Theta_k$ for $n>0$.

Our first main result pertains to numerical invariants satisfying two hypotheses:

We say that $\mu$ is \emph{standard} if 1) $\mu_\gamma(x)$ and $\mu_\gamma(-x)$ can not both be positive, and 2) each function $\mu_\gamma$ is strictly quasi-concave in the sense that for any $t \in (0,1)$ and any linearly independent pair $x_0,x_1 \in \bR^{n} \setminus \{0\}$ such that $\mu_\gamma(x_i) \geq 0$ for $i=0,1$,
\[
\mu_\gamma(t x_0 + (1-t) x_1 ) \geq \max\{\mu_\gamma(x_0), \mu_\gamma(x_1)\},
\]
with strict inequality if either $\mu_\gamma(x_0) > 0$ or $\mu_\gamma(x_1)>0$.

We say that $\mu$ satisfies condition \ref{princ:R} if for any $\gamma$ for which $\mu_\gamma$ attains a positive value, the $\mu_\gamma$ obtains a maximum at a rational point. In practice, this is straightforward to verify when it holds. The numerical invariants of \Cref{E:GIT} are standard and satisfy condition \ref{princ:R}.

\begin{thmx}[\Cref{thm:main_improved}] \label{T:main_summary}
Let $\X \to B$ be a stack satisfying \ref{hyp3} over a locally noetherian base stack $B$, and let $\mu$ be a standard numerical invariant on $\X$ satisfying \ref{princ:R}. Then $\mu$ defines a weak $\Theta$-stratification if and only if it satisfies the following:
\begin{enumerate}
\item[\textbf{\ref{princ:S}}] \textbf{HN-Specialization:} Let $R$ be a discrete valuation ring essentially of finite type over $B$, with fraction field $K$ and residue field $k$. Let $\xi : \Spec(R) \to \X$ be a map over $B$ whose generic point is unstable, and let $f_K$ be a HN filtration of $\xi_K$. Then one has $$\mu(f_K) \leq \sup \left\{\mu(f') \left| f' \text{ is a filtration of } \xi|_{\Spec(k)} \right.\right\},$$
and when equality holds there is an extension of valuation rings $R' \supset R$ with fraction field $K'$ such that $f_K|_{K'}$ extends to a filtration $f_{R'}$ of $\xi_{R'}$. \\
\item[\textbf{\ref{princ:B2}}] \textbf{HN-Boundedness:} For any map from a finite type affine scheme $\xi : T \to \X$, $\exists$ a quasi-compact substack $\X' \subset \X$ such that for any finite type point $p \in T(k)$ and any filtration $f$ of $\xi(p)$ for which $\mu(f)>0$, there is another filtration $f'$ of $\xi(p)$ with $\mu(f') \geq \mu(f)$ and whose associated graded point $f'(0)$ lies in $\X'$.
\end{enumerate}
\end{thmx}

Condition \ref{princ:B2} holds automatically if the stack $\X$ is quasi-compact, but for non-quasi-compact algebraic stacks this condition can be quite subtle. On the other hand, condition \ref{princ:S} holds automatically if the stack $\X$ has the property that \emph{any} filtration of $\xi_K$ extends uniquely to a filtration of $\xi$ over $R$. We call such a stack $\Theta$-reductive (\Cref{defn:reductive_stack}).

\addtocontents{toc}{\SkipTocEntry}
\subsection{Codimension-two filling conditions}

Many stacks in nature are not $\Theta$-reductive, so we will introduce another condition that guarantees that \ref{princ:S} holds, the condition that $\mu$ is ``strictly $\Theta$-monotone'' (\Cref{defn:theta_monotone}). Any numerical invariant on a $\Theta$-reductive stack is $\Theta$-monotone, and the numerical invariants of \Cref{E:GIT} are $\Theta$-monotone.

Roughly speaking, $\mu$ is (strictly) $\Theta$-monotone if given any DVR $R$ and a $\Gm$-equivariant map $f : \bA^1_R \setminus \{(0,0)\} \to \X$, there is a proper birational $\Gm$-equivariant morphism $\Sigma \to \bA^1_R$ that is an isomorphism over $\bA^1_R \setminus \{(0,0)\}$ and an extension of $f$ to a $\Gm$-equivariant map $\Sigma \to \X$ such that the value of $\mu$ on the $\Gm$-fixed points of $\Sigma$ is (strictly) monotone increasing.

We also discuss a closely related notion called (strict) {\textsf{S}}-monotonicity (\Cref{D:S_monotone}), which combined with $\Theta$-monotonicity can be used to establish the existence of a good moduli space for the semistable locus. These two monotonicity conditions are independent, but in practice the same methods can often be used to verify both conditions.

The monotonicity conditions lead to two theorems, which are the most directly applicable versions of our results. Here we state a simplified combined version, and refer the reader to the full statements below for slightly stronger results.

\begin{thmx}[Intrinsic GIT, \Cref{thm:beyond_git}] \label{T:monotone_summary}
Let $\X \to B$ be an algebraic stack satisfying \ref{hyp3} over a locally noetherian base stack $B$, and let $\mu$ be a numerical invariant on $\X$ that is standard and satisfies \ref{princ:R}.
\begin{enumerate}
\item If $\mu$ is strictly $\Theta$-monotone, then it defines a weak $\Theta$-stratification if and only if it satisfies condition \ref{princ:B2} above. If $\X$ has characteristic $0$, then this is a $\Theta$-stratification. (\Cref{T:monotone_stratifications})\\
\item If furthermore $B$ has characteristic $0$, $\mu$ is strictly {\textsf{S}}-monotone, the connected components of $\X^{\rm ss}$ are quasi-compact relative to $B$, and for any $\gamma$, $\mu_\gamma(x)>0$ if and only if $\mu_\gamma(-x)<0$, then $\X^{\rm ss}$ has a good moduli space relative to $B$ that is separated and locally of finite type. The connected components of the moduli space are proper over $B$ if $\X$ satisfies the existence part of the valuative criterion for properness and the $\Theta$-stratification of $\X$ is well-ordered. (\Cref{thm:monotone_moduli_spaces})
\end{enumerate}
\end{thmx}

This result suggests a general program for analyzing a moduli problem, illustrated in \Cref{fig:flow_chart}.

\tikzstyle{decision} = [diamond, draw, fill=blue!20,
    text width=6em, text badly centered, node distance=4cm, inner sep=1pt]
\tikzstyle{block} = [rectangle, draw, fill=blue!20,
    text width=10em, text centered, rounded corners, minimum height=4em]
\tikzstyle{line} = [draw, very thick, color=black!70, -latex']
\tikzstyle{cloud} = [draw, ellipse,fill=red!20, node distance=3cm, text width=5em, node distance=4cm, minimum height=2em]

\begin{center}
\begin{figure}
\begin{tikzpicture}[node distance = 2cm, auto]
    \node [block] (init) {Identify moduli problem of interest.};
    \node [block, below of=init] (stack) {Represent it by an algebraic stack $\X$.};
    \node [block, below of=stack] (num) {Identify numerical invariant $\mu$ on $\X$ that satisfies condition $(R)$.};
    \node [decision, below of=num] (monotone) {$\mu$ satisfies strict $\Theta$/{\textsf{S}}-monotonicity?};
    \node [block, left of=monotone, node distance=5cm] (update) {Modify the algebraic stack $\X$.};
    \node [decision, below of=monotone] (bounded) {$\mu$ satisfies HN boundedness?};
    \node [cloud, right of=monotone] (reductive) {Automatic if $\X$ is $\Theta$-reductive/{\textsf{S}}-complete.};
    \node [cloud, right of=bounded] (bnd) {Automatic if $\X$ is bounded.};
    \node [block, below of=bounded, node distance=4cm, text width = 20em] (win) {\textbf{Main theorem:} \\ \begin{enumerate} \item $\mu$ defines a $\Theta$-stratification of $\X$, \item $\X^{\rm ss}$ admits a separated good moduli space if it is bounded, and \item the moduli space is proper if $\X$ satisfies the existence part of the valuative criterion for properness. \end{enumerate}};
    \path [line] (init) -- (stack);
    \path [line] (stack) -- (num);
    \path [line] (num) -- (monotone);
    \path [line] (monotone) -- node {yes}(bounded);
    \path [line] (bounded) -- node {yes}(win);
    \path [line] (monotone) -- node {no}(update);
    \path [line] (bounded) -| node {no}(update);
    \path [line] (update) |- (num);
    \path [line,dashed] (reductive) -- (monotone);
    \path [line,dashed] (bnd) -- (bounded);
\end{tikzpicture}
\caption{Flow chart for analyzing a moduli problem in algebraic geometry. The conclusion of the main theorem is stated over a field of characteristic $0$, for simplicity.}\label{fig:flow_chart}
\end{figure}
\end{center}

\addtocontents{toc}{\SkipTocEntry}
\subsection{Foundations of the theory}

In \Cref{sect:filtrations} we discuss properties of the stack of filtered points of an algebraic stack $\X$, which we denote $\Filt(\X)$. In \Cref{sect:stratifications} we define the notion of a (weak) $\Theta$-stratification. This amounts to specifying an open substack of $\Filt(\X)$ whose points have properties analogous to Harder-Narasimhan filtrations in \Cref{E:vector_bundles}. The data of the stratification is encoded completely by a collection of irreducible components of $\Filt(\X)$ along with an ordering of those components.

Our first main result, whose proof is relatively straightforward, is \Cref{thm:main_stratification}. It establishes a list of five necessary and sufficient conditions for an ordered collection of irreducible components of $\Filt(\X)$ to define a $\Theta$-stratification. \Cref{T:main_summary} and part (1) of \Cref{T:monotone_summary} amount to showing that some of the necessary and sufficient conditions of \Cref{thm:main_stratification} hold automatically when the irreducible components of $\Filt(\X)$ are chosen to contain maximizers of a numerical invariant.

The theory of numerical invariants, and many of the ideas leading to our main results, is based on a new combinatorial object, studied in \Cref{sect:structures}. We call it the \emph{degeneration fan}, $\Deg(\X,p)_\bullet$, of a point $p$ in a stack $\X$, see \Cref{defn:degeneration_space}. We believe it is of independent interest. For instance, a normal toric variety $X$ is encoded by a fan in the space of cocharacters of the torus $T$ acting on $X$, and $\Deg(X/T,p)_\bullet$ encodes this fan. More generally, $\Deg(\X,p)_\bullet$ plays the same role in our theory that the spherical building of a reductive group plays in geometric invariant theory \cite{MFK94}*{Chap.~2}.

\addtocontents{toc}{\SkipTocEntry}
\subsection{Summary of contents}

In Thomason's terminology (see \cite{thomason1985algebraic}*{Introduction}) the ``thrill seeker" might want to look at the list of key terms at the end of the introduction and then skip immediately to the statements of the main theorems:
\begin{itemize}
\item First and second characterization theorems for $\Theta$-stratifications, \Cref{thm:main_stratification} and \Cref{thm:main_improved},
\item Existence of $\Theta$-stratifications for strictly $\Theta$-monotone numerical invariants, \Cref{T:monotone_stratifications},
\item The ``Intrinsic GIT Theorem'' \Cref{thm:beyond_git},
\item The ``Recognition Theorem'' for HN filtrations, \Cref{thm:recognition_revisit},
\item The ``Perturbation Theorem'' for filtrations, \Cref{thm:perturbation}.
\end{itemize}
As mentioned above, sections 1 through 3 are mostly foundational, sections 4 and 5 are the conceptual core of the paper, and section 6 contains a completely worked example, a certain moduli of complexes on a projective scheme over an arbitrary field associated to a Bridgeland stability condition (See \Cref{thm:theta_stratification_torsion_free}). We refer the reader to the beginning of each section for a more detailed summary of the contents.

\addtocontents{toc}{\SkipTocEntry}
\subsection{Applications}

Our main application in this paper is in \Cref{sect:moduli_derived}, where we construct $\Theta$-stratifications on the moduli of ``torsion free'' objects in the heart of certain $t$-structures on the derived category of coherent sheaves $\D^b(X)$ on a projective scheme $X$. This example illustrates how the main theorems are applied in practice. Furthermore, the moduli stacks in this example are not known to be local quotient stacks, and thus are hard to study using geometric invariant theory.

Since the original version of this paper was released, the methods here have been used to construct new examples of $\Theta$-stratifications in several additional moduli problems of interest \cites{blum2020properness, halpernleistner2021moduli, gomez2021moduli, GHLFH}.

We conclude this introduction by describing various applications that $\Theta$-stratifications have had, which fall into two categories.

\medskip
\noindent \textit{Comparing the geometry of $\X$ with that of $\X^{\rm{ss}}$:}
\medskip

Atiyah \& Bott \cite{atiyah1983yang} used the Harder-Narasimhan stratification to derive explicit formulas for the Betti numbers of the stack $\cM^{\rm ss}$ of semistable $G$-bundles on a Riemann surface. A closely related phenomenon appears in Teleman \& Woodward's proof \cite{teleman2009index} of a vast generalization of the Verlinde formula for the dimension of the sheaf cohomology of certain tautological vector bundles on $\cM^{\rm ss}$ -- cohomology is first computed on the stack of all bundles, and the HN stratification is used to compare the cohomology on the stack with the cohomology on the semistable locus.

Likewise for a smooth projective variety $X$ over $\bC$ with an action of a reductive group $G$, Kirwan \cite{Ki84} used the stratification of the unstable locus to derived explicit formulas for the Betti numbers of $X^{\rm{ss}}/G$. The theme in all of these examples is that sometimes the geometry of the stack $\X$ is simpler than the geometry of the semistable locus $\X^{\rm{ss}} \subset \X$, and the presence of a $\Theta$-stratification allows one to make precise comparisons between the two.

More recent examples of applications of $\Theta$-stratifications of this flavor include the proof of properness of the moduli space of $K$-semistable klt Fano varieties \cite{blum2020properness}, and the virtual non-abelian localization formula \cite{halpernleistner2015remarks}, which has been applied to generalize the Verlinde formula to moduli spaces of Higgs bundles \cite{halpernleistner2016equivariant}.

\medskip
\noindent \textit{Variation of stability and wall-crossing:}
\medskip

Numerical invariants typically depend on a rational cohomology class $\ell \in H^2(\X;\bQ)$ (see \Cref{defn:induced_invariant}). Classically, the theory of variation of GIT quotient \cites{DH98,TH96} studies how the semistable locus $\X^{\rm ss}$ varies as $\ell$ changes when $\X = X/G$ is a quotient stack. We formulate and prove an intrinsic version of variation of GIT in \Cref{thm:main_GIT}. The $\Theta$-stratification of the unstable locus allows one to identify unstable points that become semistable as $\ell$ varies (those strata for which the value of the numerical invariant tends to $0$), and to compare the geometry of $\X^{\rm ss}$ for different values of $\ell$.

We expect many interesting examples of variation of $\Theta$-stratification to arise via the following meta-principle: Whenever one has a birational isomorphism of algebraic spaces $M_1 \dashrightarrow M_2$, if $M_1$ and $M_2$ are good moduli spaces for algebraic stacks $\X_1$ and $\X_2$ parameterizing objects of geometric interest, then there should be a third stack $\X$ that is $\Theta$-reductive and also parameterizes objects of geometric interest, in which $\X_i \subset \X, i=1,2$ is the open substack of semistable points for two different numerical invariants $\mu_i,i=1,2$. Ideally $\X$ also has a good moduli space $\X \to M$, and the resulting maps $M_1 \to M \leftarrow M_2$ are projective and birational. In this case the $\Theta$-stratifications induced by $\mu_i$ would allow one to compare the geometry of $\X_1$ and $\X_2$.

Recent applications of this kind include the proof of the $D$-equivalence conjecture for Calabi-Yau manifolds that are birationally equivalent to a moduli space of sheaves on a K3-surface \cite{halpernleistner2020derived}, and other results on how wall crossings affect the derived categories of moduli spaces \cites{toda2021semiorthogonal, toda2021window}.


\addtocontents{toc}{\SkipTocEntry}
\subsection{Context}

We will typically denote stacks with fraktur font ($\X$,$\Y$, etc.), with the exception of our fixed base stack, which we denote $B$. All of our stacks $\X$ will be stacks on the big \'{e}tale site of affine schemes over $B$, or equivalently stacks $\X$ on the big \'{e}tale site of affine schemes with a fixed map $\X \to B$. We will denote the quotient stack for a group $G$ acting on a space $X$ as $X/G$. Although we will often not mention this explicitly, all conditions on $\X$ will be relative to $B$. For instance, we might omit $B$ from discussions of the \emph{inertia of $\X$}, which is the first projection morphism $I_{\X/B}:=\X \times_{\X \times_B \X} \X \to \X$. If $S$ is a $B$-scheme, we will use the notation $\X_S$ to denote the $S$-stack $S \times_B \X$.

We will often make the following technical hypotheses on $\X$, whose main use is to guarantee that the stack of filtered objects is algebraic (see \Cref{prop:existence_mapping}). The condition holds, for example,  if $\X$ admits a Zariski open cover by stacks of the form $X/G$ where $X$ is a quasi-separated, quasi-compact algebraic space and $G$ is a smooth affine group scheme.

\begin{enumerate}[label=$(\dagger)$]
\item \label{hyp2} $B$ is a quasi-separated algebraic stack, and $\X \to B$ is a morphism of algebraic stacks that is locally of finite presentation, quasi-separated, and such that points of $\X$ have affine automorphism groups in their fiber over $B$.
\end{enumerate}

We will sometimes refer to the following stronger hypotheses, which guarantee that filtrations in $\X$ have no automorphisms that act trivially on the underlying point of $\X$ (see \Cref{prop:representable_flags}).

\begin{enumerate}[label=$(\dagger\dagger)$]
\item \label{hyp3} \ref{hyp2} holds, and the relative inertia morphism $I_{\X/B} = \X \times_{\X \times_B \X} \X \to \X$ is separated.
\end{enumerate}

\addtocontents{toc}{\SkipTocEntry}
\subsection{Author's note}

This project began as the third chapter of my PhD thesis, which contained much weaker versions of many of the results here. As I came to understand the relative ubiquity of examples of $\Theta$-stratifications (known by other names), I began to appreciate the value of a text that would unify these examples and establish tools for studying these structures in new moduli problems. This has led me to expand, strengthen, and rewrite this paper several times over the intervening years.

As a result, I owe a large debt to the many formal and informal mentors I have had over this time: Constantin Teleman, Andrei Okounkov, Davesh Maulik, Johan de Jong, and Michael Thaddeus. This work has benefited from useful conversations with many other researchers: Jarod Alper, Arend Bayer, Gergely B{\'e}rczi, Tristan Collins, Brian Conrad, Anand Deopurkar, Galyna Dobrovolska, Maksym Fedorchuk, Simion Filip, Joshua Jackson, Frances Kirwan, S\'andor Kov\'acs, Jochen Heinloth, Victoria Hoskins, Jacob Lurie, Johan Martens , Ian Morrison, Yuji Odaka, Martin Olsson, Alexander Polishchuk, Simon Schieder, Carlos Simpson, Pablo Solis, David Swinarski, Richard Thomas, Yukinobu Toda, Roman Travkin, Xiaowei Wang, Chris Woodward, and Chenyang Xu.

This work was partially supported by NSF grants DMS-1945478 (CAREER), DMS-1601967, and DMS-1303960 (MSPRF), and a Simons Collaboration Grant. I would like to thank Columbia University, the Institute for Advanced Study, Cornell University, and the Mathematical Sciences Research Institute for providing the stimulating environments in which I completed this work.

\label{S:list}
\printglossary[title={List of key terms},style=index, nogroupskip]
\bigskip


\section{Filtrations} \label{sect:filtrations}

In this section we lay the foundational framework for the theory of $\Theta$-stability in an algebraic stack $\X$. As mentioned in the introduction, the quotient stack $\Theta := \bA^1 / \Gm$ plays a key role. Motivated by the concept of the Rees module associated to a filtered vector space (discussed in \Cref{E:vector_bundles_2} and recalled in \Cref{prop:quasicoh_theta} below), we interpret a map $f : \Theta \to \X$ as a ``filtration'' of the point $f(1) \in \X$.  Similarly, a graded point of $\X$ is a map $(B\Gm)_k \to \X$.

This leads us to define an algebraic stack $\Filt(\X) := \filt{\X}$ parameterizing families of filtered points of $\X$, and an algebraic stack $\Grad(\X) := \grad{\X}$ parameterizing families of graded points of $\X$ (\Cref{prop:existence_mapping}). These are related via several universal maps:

\begin{figure}[htbp]
\centering
\fbox{
\begin{minipage}{.95\textwidth}
\begin{equation} \label{eqn:universal_maps}
\xymatrix{ \Grad(\X) \ar `u[rr] `[rr]^u [rr] \ar@/_10pt/[r]_\sigma & \Filt(\X) \ar[l]_{\agr} \ar@/^/[r]^{\ev_1} \ar@/_/[r]_{\ev_0} & \X}
\end{equation}
\caption[Canonical morphisms]{Summary of important maps from the stack of graded and filtered points of $\X$. They are:\vskip1em
\begin{minipage}{.8\textwidth}
\begin{itemize}
\item[$\ev_1$:] ``forgetting the filtration'' map restricts a filtration $\Theta_k \to \X$ to the point $\{1\} \in \Theta_k$.\\
\item[$\sigma$:] ``split filtration'' map pulls back a graded point $\pt / (\Gm)_k \to \X$ along the projection $\Theta_k \to \pt / (\Gm)_k$.\\
\item[$\agr$:] ``associated graded" map restricts a filtration $f : \Theta_k \to \X$ to the closed substack $\{0\}/(\Gm)_k \hookrightarrow \Theta_k$. $\agr \circ \sigma \simeq \id_{\Grad(\X)}$, and $\agr$ is an algebraic deformation retract (\Cref{lem:retract}).\\
\item[$\ev_0$:] ``associated graded'' map assigns a filtration $f : \Theta_k \to \X$ to $f(0)$ without its $\Gm$-action; $u \circ \agr \simeq \ev_0$ canonically.\\
\item[$u$:] ``forgetful'' map restricts a map $\pt / (\Gm)_k \to \X$ along the map $\pt \to \pt/(\Gm)_k$; $u \simeq \ev_1 \circ \sigma \simeq \ev_0 \circ \sigma$ canonically.\\
\end{itemize}\end{minipage}}
\label{F:main_maps}
\end{minipage}
}
\end{figure}

For any scheme $T$ and map $\xi : T \to \X$, we also define (\Cref{defn:flag_space}) an algebraic space $\Flag(\xi)$ over $T$, the generalized flag space, that parameterizes families of filtrations of the point $\xi$. Formally, $\Flag(\xi)$ is the fiber of $\ev_1$ over $\xi : T \to \X$.

This section establishes the basic properties of these objects. For instance, we study a natural action of the monoid $\bN^\times$ acts on the stack $\Filt(\X)$, we give an explicit description of $\Grad(\X)$ and $\Filt(\X)$ when $\X = X/G$ is a quotient stack (\Cref{thm:describe_strata_quotient_GLN}), and we show that $\Flag(\xi)$ is separated if $\X$ has affine diagonal (\Cref{prop:representable_flags}).

\subsubsection{Motivating example}

Maps $S \times \Theta \to B \GL_N$ classify locally free sheaves on $S \times \Theta$. Recall the Rees construction:

\begin{prop}\label{prop:quasicoh_theta}
Let $S$ be a scheme. The category of quasi-coherent sheaves on $S \times \Theta$ is equivalent to the category of diagrams of quasi-coherent sheaves on $S$ of the form $\cdots \to F_i \to F_{i-1} \to \cdots$.
\end{prop}
\begin{proof}
Using descent one sees that quasi-coherent sheaves on $S \times \Theta$ are the same as graded $\cO_S[t]$ modules, where $t$ has degree $-1$. The equivalence assigns a diagram $\cdots \to F_i \to F_{i-1} \to \cdots$ to the module $\bigoplus F_i$ with $F_i$ in degree $i$ with the maps $F_i \to F_{i-1}$ corresponding to multiplication by $t$.
\end{proof}

Under this equivalence, locally free sheaves on $S \times \Theta$ correspond to diagrams such that each $F_i$ is locally free on $S$, $F_i \to F_{i-1}$ is injective and $\op{gr}_i (F_\bullet) = F_i / F_{i+1}$ is locally free for each $i$, $F_i$ stabilizes for $i\ll 0$, and $F_i = 0$ for $i\gg 0$.  In other words $\GL_N$ bundles on $S \times \Theta$ correspond to locally free sheaves with decreasing, weighted filtrations, and the shape of this filtration must be constant along connected schemes.


\subsection{Graded objects, filtered objects, and flag spaces}

Given two stacks $\X$ and $\Y$ over a site, one can form the mapping stack $\iMap(\Y,\X)$ as the presheaf of groupoids
\begin{equation}\label{eqn:mapping_stack}
\iMap(\Y,\X) : T \mapsto \Map(\Y \times T, \fX),
\end{equation}
where $\Map$ denotes groupoid of $1$-morphisms between stacks. If $\Y = \Theta := \bA^1 / \Gm$ and $\X$ is an algebraic stack, then one can apply smooth descent \cite{Vi05} to describe $\filt{\X}$ more explicitly. We consider the first three levels of the simplicial scheme determined by the action of $\Gm$ on $\bA^1 \times T$
\begin{equation} \label{eqn:descent_diagram}
\xymatrix{ \Gm \times \Gm \times \bA^1 \times T \ar@<1ex>[r]^-\mu \ar[r]|-\sigma \ar@<-1ex>[r]_-{a} & \Gm \times \bA^1 \times T \ar@<.5ex>[r]^-\sigma \ar@<-.5ex>[r]_-a & \bA^1 \times T }
\end{equation}
Where $\mu$ denotes group multiplication, $\sigma$ denotes the action of $\Gm$ on $\bA^1$, and $a$ forgets the leftmost group element. Then the category $\filt{\X}(T)$ has
\begin{itemize}
\item \emph{objects:} $\eta \in \X(\bA^1 \times T)$ along with a morphism $\phi : a^\ast \eta \to \sigma^\ast \eta$ satisfying the cocycle condition $\sigma^\ast \phi \circ a^\ast \phi = \mu^\ast \phi$
\item \emph{morphisms:} $f : \eta_1 \to \eta_2$ such that $\phi_2 \circ a^\ast(f) = \sigma^\ast(f) \circ \phi_1 : a^\ast \eta_1 \to \sigma^\ast \eta_2$
\end{itemize}
Informally, if $\X(T)$ is the groupoid of $T$-families of geometric objects of a certain kind, then $\filt{\X}(T)$ is the groupoid of $\Gm$-equivariant families over $\bA^1 \times T$.

As with any mapping stack, one has a universal evaluation $1$-morphism
\begin{equation}\label{eqn:ev_map}
\ev : \Theta \times \filt{\X} \to \X.
\end{equation}
$\Theta$ has two canonical $B$-points, the generic point corresponding to $1 \in \bA^1$ and the special point $0 \in \bA^1$, and the restriction to the open (respectively closed) substack $\{1\} \subset \Theta$ (respectively $\{0\}/\Gm \hookrightarrow \Theta$) define the maps $\ev_1$ (respectively $\agr$) of \Cref{eqn:universal_maps}. More generally, we consider the $B$-point of $\Theta^n$ determined by $(1,\ldots,1) \in \bA^n$, and the map $\{0\} / \Gm^n \to \Theta^n$ determined by the point $(0,\ldots,0) \in \bA^n$. We denote the resulting restriction maps $\ev_1 : \filt[n]{\X} \to \X$ and $\agr : \filt[n]{\X} \to \grad[n]{\X}$ as well.

Note that in our context we implicitly work relative to $B$, so all test schemes $S$ are $B$-schemes, $\bA^1$ and $\Gm$ refer to $\bA^1_B$ and $(\Gm)_B$ unless otherwise specified, and the mapping stack is formed relative to $B$. Concretely, a $T$-point of $\filt[n]{\cX}$ is a map of stacks $\Theta^n_T \to \cX$ along with a factorization of the composition $\Theta^n_T \to B$ through the projection $\Theta^n_T \to T$. If $B$ is a quasi-separated algebraic stack, then any factorization of a morphism $\Theta^n_T \to B$ through $\Theta^n_T \to T$ is unique up to unique isomorphism, so the factorization is a condition rather than additional data. This follows from the universal property of good moduli space morphisms \cite{ahr2}*{Thm.~17.2}, applied to the good moduli space morphism $\Theta^n_T \to T$.

The formation of mapping stacks behaves well with respect to base change, so we note the following purely formal observation:
\begin{lem} \label{lem:base_change_mapping}
Let $T$ be an algebraic stack over $B$, and let $S$ be a $T$-scheme. The groupoid of $S$-points of $\iMap_T(\Theta^n_T,\X_T)$ is canonically equivalent to the groupoid of maps $\Theta^n_S \to \X$ relative to $B$, and likewise for $\iMap_T((\pt/\Gm^n)_T,\X_T)$. This gives canonical equivalences of $T$-stacks
\begin{equation} \label{eqn:base_changing_mapping}
\begin{array}{c}
\iMap_T(\Theta^n_T,\X_T) \simeq \filt[n]{\X} \times_B T \\ \iMap_T((\pt / \Gm^n)_T,\X_T) \simeq \grad[n]{\X} \times_B T
\end{array}.
\end{equation}
Furthermore this equivalence identifies the canonical evaluation map $\Theta^n_T \times_T \iMap_T(\Theta^n_T,\X_T) \to \X_T$ with the base change of \eqref{eqn:ev_map}, and likewise for the maps $\agr$, $\ev_1$, and $\sigma$ appearing in \eqref{eqn:universal_maps}.
\end{lem}

The mapping stack $\iMap(\bA^1,\X)$ is rarely algebraic, because $\bA^1$ is not proper. Hence, it is not obvious from the description of $\filt{\X}$ in terms of descent data that the latter is algebraic, but $\Gm$-equivariance fixes this non-compactness issue.

\begin{prop}\label{prop:existence_mapping}
Let $\X \to B$ satisfy \ref{hyp2}. Then for any $n \geq 1$, $\filt[n]{\X}$ and $\grad[n]{\X}$ are algebraic stacks, locally of finite presentation and quasi-separated over $B$. If $\X \to B$ has affine (resp. quasi-affine, resp. separated) diagonal then so do $\filt[n]{\X}$ and $\grad[n]{\X}$.
\end{prop}
\begin{proof}
The definition of the mapping stack commutes with base change along any morphism $B' \to B$, so it suffices to prove the claim when $B = \Spec(R)$ is affine. The stacks $\pt/\Gm^n$ and $\Theta^n$ are cohomologically projective and flat over $B$, so this is an immediate application of \cite{HLPreygel}*{Thm.~5.1.1} when $\X \to B$ has quasi-affine diagonal. The more general version stated here is based on \cite{ahr2}*{Thm.~14.9}, which implies that because $\Theta^n_B \to B$ and $(B\Gm)^n_B \to B$ are flat good moduli space morphisms, $\filt[n]{\X}$ is an algebraic stack locally of finite presentation and with quasi-separated diagonal over $B$, and that it has affine/quasi-affine/separated diagonal whenever $\X \to B$ does.

We claim, in addition, that $\filt[n]{\X} \to B$ has quasi-compact diagonal, hence is quasi-separated. This amounts to the claim that for any two maps from an affine scheme $\xi_1,\xi_2 : T \to \filt[n]{\X}$ the $T$-space classifying isomorphisms $\underline{\Isom}_{\filt[n]{\X}}(\xi_1,\xi_2)$ is quasi-compact. Say $f_1,f_2:\Theta^n_T \to \X$ are the maps classified by $\xi_1$ and $\xi_2$. $\underline{\Isom}_{\X}(f_1,f_2)$ corresponds to a quasi-separated and quasi-compact algebraic space $X$ with a $(\Gm^n)_T$-action along with an equivariant map $X \to \bA^n_T$, and $\underline{\Isom}_{\filt[n]{\X}}(\xi_1,\xi_2)$ is the Weil restriction of $X/\Gm^n \to \Theta^n_T$ along the projection $\Theta_T \to T$. So, we have a cartesian diagram
\[
\xymatrix{\underline{\Isom}_{\filt[n]{\X}}(\xi_1,\xi_2) \ar[r] \ar[d] & \underline{\Map}_T(\Theta_T^n, X/\Gm^n) \ar[d] \\
T \ar[r]^{\rm id} & \underline{\Map}_T(\Theta_T^n, \Theta_T^n) }.
\]
The stack $X/\Gm^n$ has quasi-affine diagonal because $X$ is quasi-separated \cite{stacks-project}*{\href{https://stacks.math.columbia.edu/tag/02LR}{Tag 02LR}}, so we already know the theorem for $\underline{\Map}_T(\Theta_T^n,X/\Gm^n)$. In fact, we will see in the proof of \Cref{thm:describe_strata_quotient_GLN} below that
\[
\filt[n]{X/\Gm^n} = \bigsqcup_{\psi : \Gm^n \to \Gm^n} X^{\psi,+}/\Gm^n,
\]
where the $\psi$ are group homomorphisms, $X^{\psi,+}$ is the algebraic space representing the functor of $\Gm^n$-equivariant maps $\bA^n \to X$ for the action of $\Gm^n$ on $X$ via $\psi$ (see \Cref{prop:Hesselink_concentration}). Each space $X^{\psi,+}$ admits a quasi-compact morphism to the closed algebraic space $X^{\psi,0} \hookrightarrow X$ corresponding to the $\psi(\Gm^n)$-fixed locus, and hence is quasi-compact. Applying this to the cartesian square above, we see that $\underline{\Isom}_{\filt[n]{\X}}(\xi_1,\xi_2)$ is the fiber over $(1,\ldots,1) : T \to \bA^n_T$ of the projection $X^{\id,+} \to \bA^n_T = (\bA^n_T)^{\id,+}$, hence it is quasi-compact. The proof that $\grad[n]{\X}$ has quasi-compact diagonal is similar, so we omit it.

\end{proof}

We will be using these mapping stacks so often that we introduce more concise and intuitive notation.

\subsubsection{The \gls{graded_objects}}

\begin{ex} \label{ex:main1}
Assume that $B$ is a scheme, and let $X$ be a projective scheme over $B$. Let $\X$ denote the stack of flat families of coherent sheaves on $X$. Then $\grad[n]{\X}$ parameterizes flat families of graded coherent sheaves, graded by the abelian group $\bZ^n$.
\end{ex}

Motivated by this example, we make the following
\begin{defn} \label{defn:graded_objects}
Given a $B$-stack $T$, we define the stack of \emph{$\bZ^n$-graded objects} of $\X$ over $T$ to be
\[
\Grad^n_T(\X) := \grad[n]{\X} \times_B T \simeq \iMap_T((\pt/\Gm^n)_T,\X_T).
\]
We simplify notation by writing $\Grad^n$ for $\Grad^n_B$.
\end{defn}

There is a canonical ``forgetful" map $u : \Grad^n(\X) = \grad[n]{\X} \to \X$ that restricts along the map $B \to (\pt / \Gm)_B$. Concretely for any $B$-scheme $S$ this assigns a map $S\times \pt / \Gm \to \X$ to its restriction $\xi : S \to \X$. In fact, lifting an $S$ point $\xi \in \X(S)$ along the map $u$ is equivalent to giving a homomorphism from $\Gm^n$ to the automorphism group 
\[
\Aut_{\X/B}(\xi) := \ker(\Aut_\X(\xi) \to \Aut_B(\pi(\xi))),
\]
where $\pi : \X \to B$ is the structure map. More precisely we have

\begin{lem} \label{lem:graded_stack}
For any $B$-scheme $S$ and any $S$-point $\xi : S \to \X$ the following diagram is cartesian
\[
\xymatrix{\underline{\op{Hom}}_{S\rm{-gp}}((\Gm^n)_S,\underline{\Aut}_{\X/B}(\xi)) \ar[r] \ar[d] & \grad[n]{\X} \ar[d]^{\text{restriction along } B \to (\pt/\Gm^n)_B} \\ S \ar[r]^\xi & \X }.
\]
Where $\underline{\Aut}_{\X/B}(\xi)$ denotes the sheaf of groups on $S$ that assigns to any map $U \to S$ the group $\Aut_{\X/B}(\xi|_U)$. In particular if $\X$ satisfies \ref{hyp2}, then the vertical maps are representable by algebraic spaces.
\end{lem}

\begin{proof}
It suffices to show the non-sheafified version of the claim for any $S$-scheme $T$. This is basically a formal consequence of the description of the groupoid $\Map_B((\pt / \Gm^n)_T,\X)$ via descent. This groupoid is canonically equivalent to the groupoid of maps $\xi : T \to \X$ relative to $B$ along with an automorphism of the composition $(\Gm^n)_T \to T \to \X$ relative to $B$ satisfying a cocycle condition. The data of an automorphism of the map $(\Gm^n)_T \to T \to \X$ relative to $B$ is equivalent to specifying a section of the pullback of the $T$-sheaf $\underline{\Aut}_{\X/B}(\xi)$ along $(\Gm^n)_T \to T$, which in turn is equivalent to a map of sheaves of sets $(\Gm^n)_T \to \underline{\Aut}_{\X/B}(\xi)$ over $T$. The cocycle condition translates to the condition that this map of sheaves respects the group structure. Thus we see that $\Map_B((\pt/\Gm^n)_T,\X)$ is equivalent to the groupoid consisting of a pairs $(\xi : T \to \X / B, \phi \in \Hom_{T\rm{-gp}}((\Gm^n)_T,\underline{\Aut}_{\X/B}(\xi)))$ and isomorphisms $(\xi_1,\phi_1) \simeq (\xi_2,\phi_2)$ are isomorphisms $\eta : \xi_1 \simeq \xi_2$ relative to $B$ for which the induced bijection $\Hom_{T\rm{-gp}}((\Gm^n)_T,\underline{\Aut}_{\X/B}(\xi_1)) \simeq \Hom_{T\rm{-gp}}((\Gm^n)_T,\underline{\Aut}_{\X/B}(\xi_2))$ maps $\phi_1 \mapsto \phi_2$.
\end{proof}

\begin{rem}
There are general existence theorems for schemes representing the functor $\underline{\op{Hom}}_{S\rm{-gp}}(\Gm^n,G)$ for an $S$ group scheme $G$, but they typically require $G$ to be flat. \Cref{lem:graded_stack} extends these results to any $G$ that is \'etale locally modeled on the inertia $I_{\X/B} \to \X$ for some stack satisfying \ref{hyp2}.
\end{rem}

Let $I_{\X/B} := \X \times_{\X \times_B \X} \X$ denote the inertia stack of $\X$ relative to $B$. It is a group scheme over $\X$ whose $R$ points are $\xi \in \X(R)$ along with an automorphism in $\Aut_{\X/B}(\xi)$. The previous lemma can be interpreted as providing a canonical equivalence of stacks over $\X$
\[
\Grad^n(\X) \simeq \underline{\Hom}_{gp}((\Gm^n)_\X, I_{\X/B}),
\]
where the latter denotes the Hom sheaf between group sheaves over $\X$.

\begin{cor} \label{cor:inertia_preserving}
Let $\phi : \X \to \Y$ be a map of stacks such that the canonical map on inertia groups fits into a short exact sequence of group sheaves over $\X$
\[
\{1\} \to I_{\X/B} \to \X \times_\Y I_{\Y/B} \to \cQ \to \{1\},
\]
where the fibers of $\cQ$ over any field-valued point of $\cX$ are representable and quasi-finite. Then the canonical map
\[
\Grad^n(\X) \to \X \times_\Y \Grad^n(\Y)
\]
induces an equivalence of groupoids of $k$-points for any field $k$ over $B$. Furthermore, if either $\cQ = \{1\}$ (i.e. $\phi$ is inertia preserving), or $\X$ and $\Y$ satisfy \ref{hyp2} and $\phi$ is representable and \'etale, then $\Grad^n(\X) \to \X \times_\Y \Grad^n(\Y)$ is an isomorphism of stacks.
\end{cor}
\begin{proof}
This follows from \Cref{lem:graded_stack} and the fact that the Hom sheaf between group schemes commutes with pullback, so $\phi^{-1} (\underline{\Hom}((\Gm^n)_\Y,I_{\Y/B})) \simeq \underline{\Hom}((\Gm^n)_\X,I_{\X/B})$. The claim when $\X$ and $\Y$ satisfy \ref{hyp2} and $\phi$ is representable and \'etale follows from the observation that $\Grad^n(\X) \to \X \times_\Y \Grad^n(\Y)$ is a universally bijective map of algebraic stacks that induces an isomorphism on tangent complexes by \Cref{prop:open_closed_filtrations}(4).
\end{proof}

Applying the previous corollary gives

\begin{cor} \label{cor:closed_open_graded}
If $\X \to \Y$ is a closed immersion (resp. surjective closed immersion, resp. open immersion) relative to $B$, then so is the map $\Grad^n(\X) \to \Grad^n(\Y)$.
\end{cor}

We can apply this as well to establish a base change result.

\begin{cor} \label{cor:base_change_grad}
Let $\X$ be a stack over an algebraic base stack $B$ and consider a map $\Spec(R) \to B$. Then the canonical maps define equivalences
\[
\Grad^n(\X_R) \simeq \Spec(R) \times_B \Grad^n(\X) \simeq \iMap_{\Spec(R)}((\pt/\Gm)^n,\X_R),
\]
where the stack of graded objects is formed relative to $B$.
\end{cor}
\begin{proof}
The second equivalence is \Cref{lem:base_change_mapping}. For the first equivalence, consider the composition $\X_R \to \Spec(R) \to B$. This leads to a left-exact sequence of relative inertia group sheaves over $\X$
\[
\{1\} \to I_{\X_R/\Spec(R)} \to I_{\X_R/B} \to \pi^{-1}(I_{\Spec(R)/B}),
\]
where $\pi : \X_R \to \Spec(R)$ is the projection. $I_{\Spec(R)/B} = \{1\}$ because $\Spec(R) \to B$ is representable, and it follows that $I_{\X_R/R} \to I_{\X_R/B}$ is an equivalence. On the other hand if $p : \X_R \to \X$ is the projection, then the composition $I_{\X_R/R} \to I_{\X_R/B} \to p^{-1}(I_{\X/B})$ is an equivalence, and it follows that $\X_R \to \X$ is inertia preserving relative to $B$. We can therefore apply \Cref{cor:inertia_preserving}.
\end{proof}

A closely related observation is the following:
\begin{cor} \label{cor:grad_over_affine}
Let $\X$ be a stack over a base stack $B$. Let $B' \to B$ be a map of stacks with trivial relative inertia $I_{B'/B} = \{1\}$, and let $\X' = \X \times_B B'$. Then the natural map $\iMap_{B'}((\pt/\Gm)^n,\X') \to \iMap_{B}((\pt/\Gm)^n,\X')$ as stacks over $\X'$ is an equivalence.
\end{cor}
\begin{proof}
Under \Cref{lem:graded_stack}, we can identify this map of stacks with the map
\[
\underline{\Hom}_{gp}((\Gm^n)_{\X'},I_{\X'/B'}) \to \underline{\Hom}_{gp}((\Gm^n)_{\X'},I_{\X'/B})
\]
of group sheaves over $\X'$. The argument of \Cref{cor:base_change_grad} shows that $I_{\X'/B'} \to I_{\X'/B}$ is an isomorphism, so the map above is an isomorphism.
\end{proof}

\subsubsection{The \gls{filtered_objects}}

\begin{defn} \label{defn:stack_filtered}
Given a $B$-stack $T$, we define the stack of \emph{$\bZ^n$-filtered objects} of $\X$ over $T$ to be
\[
\Filt^n_T(\X) := \filt[n]{\X} \times_B T \simeq \iMap_T(\Theta^n_T,\X_T),
\]
We simplify notation by writing $\Filt^n$ for $\Filt^n_B$, and we say that a point of $\Filt^n(\X)$, a filtered point of $\X$, is \emph{split} if it lies in the image of $\sigma : \Grad^n(\X) \to \Filt^n(\X)$.
\end{defn}

\begin{ex} \label{ex:main2}
Continuing \Cref{ex:main1}, we will see in \Cref{lem:filtrations_derived} that $\filt[n]{\X}$ parameterizes flat families of coherent sheaves along with a flat family of filtrations indexed by the partially ordered abelian group $\bZ^n$. The map $\ev_1 : \filt[n]{\X} \to \X$ is the map that forgets the filtration, and the map $\agr : \filt[n]{\X} \to \grad[n]{\X}$ maps a flat family of filtered coherent sheaves to its associated graded family of coherent sheaves. The map $\sigma$ is the canonical map regarding a graded coherent sheaf as a filtered coherent sheaf, where the filtered pieces are $\cE_{\geq w} = \bigoplus_{i \geq w} \cE_i$.
\end{ex}

\subsubsection{Flag spaces}

Finally, we can define an analog of the classical flag scheme for general moduli problems, using the following:

\begin{prop} \label{prop:representable_flags}
If $\X$ satisfies \ref{hyp3}, then $\ev_1 : \filt[n]{\X} \to \X$ is representable by algebraic spaces, locally of finite presentation, and quasi-separated. Furthermore:
\begin{enumerate}
\item If $\X \to B$ has quasi-affine (resp. separated) diagonal, then the same is true for $\ev_1$; and
\item If $\X \to B$ has affine diagonal, then $\ev_1$ is separated.
\end{enumerate}
\end{prop}
\begin{proof}
\Cref{prop:existence_mapping} implies that $\ev_1$ is relatively representable by algebraic stacks, locally of finite presentation, and quasi-separated. To show that $\ev_1$ is representable by algebraic spaces, it thus suffices to show that $\filt[n]{\X}$ is a category fibered in sets over $\X$.

Let $S$ be a $B$-scheme and let $f : \Theta^n_S \to \X$ be a morphism, an element of $\filt[n]{\X}(S)$. Then $\Aut(f)$ is equivalent to the group of sections of $\Y : = \Theta^n_S \times_{\X \times \X} \X \to \Theta^n_S$, where $\Theta^n_S \to \X \times \X$ in this fiber product is classified by $(f,f)$ and $\X \to \X \times \X$ is the diagonal. $\ev_1(f) \in \X(S)$ is the restriction of $f$ to $\{1\} \times S$, and automorphisms of $f$ that induce the identity on $\ev_1(f)$ correspond to those sections of $\Y \to \Theta^n_S$ that agree with identity on the open substack $(\bA^1 - \{0\})^n \times S / \Gm^n \subset \Theta^n_S$. By hypothesis $\Y \to \Theta^n$ is representable by separated algebraic spaces, so a section is uniquely determined by its restriction to $(\bA^1 -\{0\})^n \times S$. Hence $\Aut(f) \to \Aut_{\X} (\ev_1(f))$ has trivial kernel.

\Cref{prop:existence_mapping} also implies the stronger conclusions about the diagonal in (1), and when $\X \to B$ has affine diagonal, then so does $\ev_1$.

We now prove (2). For any $n \geq 1$, we use the identity $\Filt^n(\X) \simeq \Filt(\Filt^{n-1}(\X))$ to obtain a sequence of maps
\[
\Filt^n(\X) \xrightarrow{\ev_1} \Filt^{n-1}(\X) \xrightarrow{\ev_1} \cdots \xrightarrow{\ev_1} \Filt(\X) \xrightarrow{\ev_1} \X,
\]
and we have already shown that each is representable by locally finitely presented algebraic spaces with affine diagonal. It therefore suffices to verify that $\ev_1 : \Filt(\X) \to \X$ satisfies the valuative criterion for separatedness with respect to all discrete valuation rings $R$. 

Let $R$ be a discrete valuation ring over $B$, and consider two morphisms $f_1,f_2 : \Theta_R \to \X$. We must show that any isomorphism $f_1|_\cU \simeq f_2|_{\cU}$, where $\cU \subset \Theta_R$ is the open compliment of the unique codimension 2 point, extends uniquely to an isomorphism $f_1 \simeq f_2$. By hypothesis the sheaf of isomorphisms $\underline{\Isom}(f_1,f_2)$ is representable by an algebraic stack $\cI$ that is affine over $\Theta_R$, and we must show that any section over $\cU \subset \Theta_R$ extends uniquely to a section over $\Theta_R$. It suffices to do this after base change along the morphism $\bA^1_R \to \Theta_R$, after which one has an affine morphism $\pi : I \to \bA^1_R$, and a section of $\pi$ over the open subset $U := \bA^1_R \setminus \{(0,0)\}$. This section must extend uniquely because $\bA^1_R$ is the affinization of $U$, and $I$ is affine.
\end{proof}

\begin{rem} \label{rem:theta_stack_fibers}
The fact that for any $B$-scheme $S$ the fiber of $\ev_1$ over an $S$ point of $\X$ is equivalent to a set does not depend on the representability of $\filt[n]{\X}$ or even the representability of $\X$. It only relies on the fact that the inertia stack $I_{\X/B} \to \X$ is representable by separated algebraic spaces.
\end{rem}

\begin{defn} \label{defn:flag_space}
Given a map $\xi : T \to \X$ over $B$, we define the \emph{$\bZ^n$-\gls{flag_space}} of $\xi$ to be
\[
\Flag^n(\xi) := \filt[n]{\X} \times_{\ev_1,\X,\xi} T \simeq \Filt^n_T(\X) \times_{\ev_1,\X_T,\xi_T} T.
\]
When $n=1$ we drop ``$n$" from the notation and we drop ``$\bZ^n$-'' from the terminology.
\end{defn}

Under the hypotheses of \Cref{prop:representable_flags}, $\Flag^n(\xi)$ is an algebraic space locally of finite presentation and quasi-separated over $T$, and it is separated if $\X \to B$ has affine diagonal.

\begin{ex} \label{ex:main3}
Continuing \Cref{ex:main2}, if $\xi : T \to \X$ parameterizes a $T$-flat coherent sheaf $\cE$ on $T \times X$, then $\Flag(\xi)$ consists of countably many connected components. Each component is isomorphic to a classical flag scheme parameterizing $T$-flat filtrations of $\cE$ whose associated graded coherent sheaves have specified Hilbert polynomials. Each flag scheme of this kind, however, appears infinitely many times in $\Flag(\xi)$, corresponding to all of the ways to assign integer weights to the associated graded pieces.
\end{ex}


\subsection{Deformation theory and the spectral mapping stack}
\label{S:deformation_theory}

We will often consider the deformation theory of the stacks $\grad[n]{\X}$ and $\filt[n]{\X}$, but the easiest way to access that will be to use a bit of spectral algebraic geometry. Spectral algebraic geometry will not play a very central role, so in this section we will simply summarize the short list of results we will need and provide references to thorough treatments of the subject. We suggest the reader skip this section on a first reading of this paper.

In the context of spectral algebraic geometry, we will consider stacks as sheaves of spaces over the $\infty$-category $\op{CAlg}^{cn}$ of connective $E_\infty$-algebras with its \'{e}tale topology \cite{lurie2012higher}*{Sect.~7.5}. We say that $\X : \op{CAlg}^{cn} \to \cS$ is a $1$-stack if for any $A \in \op{CAlg}^{cn}$ with $\pi_0(A) \simeq A$, the space $\X(A)$ is $1$-truncated, i.e. weakly equivalent to the classifying space of a groupoid. For any $1$-stack $\X$, we define the underlying classical stack $\X^{\rm{cl}}$ as the restriction of $\X$ to the category of classical rings $\op{Ring} \subset \op{CAlg}^{cn}$. On the other hand given a classical stack $\X$, we can define $\X^{\rm{sp}}$ to be the left Kan extension of $\X$ along this same inclusion, and refer to this as the stack $\X$ regarded as a spectral stack.

More concretely if $\X$ is an algebraic classical stack, then one can choose a simplicial scheme $U_\bullet$ presenting $\X$ such that $U_n = \op{Spec}(R_n)$ is affine for every $n$, and $\X^{\rm{sp}}$ is the spectral stack that is the colimit of the simplicial spectral scheme $U_\bullet^{\rm{sp}}$ obtained by regarding each $R_n \in \op{Ring}$ as an $E_\infty$-algebra. The functor $(-)^{\rm{sp}}$ is fully faithful. Furthermore these operations are adjoint to one another, in the sense that if $\X$ is classical and $\Y$ spectral, then
\[
\Map(\X^{\rm{sp}},\Y) \simeq \Map(\X,\Y^{\rm{cl}})
\]
where the left hand side is the space of maps of spectral stacks, and the right hand side is the classifying space of the groupoid of maps of classical stacks. In particular the unit of adjunction $\X \to (\X^{\rm{sp}})^{\rm{cl}}$ is an equivalence of classical stacks.

In the context of spectral algebraic geometry, one defines the mapping stack via the same functor of points \eqref{eqn:mapping_stack}, but where $T$ is an arbitrary affine spectral scheme, i.e. connective $E_\infty$-algebra.
\begin{lem}\label{lem:spectral_maps}
If $\X$ is a spectral $1$-stack and $\Y$ is a classical stack, then the adjunction between $(-)^{\rm{sp}}$ and $(-)^{\rm{cl}}$ provides a canonical equivalence of classical stacks
\[
\iMap(\Y^{\rm{sp}},\X)^{\rm{cl}} \simeq \iMap(\Y,\X^{\rm{cl}}),
\]
where the right hand side denotes the classical mapping stack.
\end{lem}

If $\X$ is a spectral algebraic stack over a classical base $B$ such that the underlying classical stack $\X^{\rm cl}$ satisfies \ref{hyp2}, then the spectral mapping stack $\iMap((\Theta^n)^{\rm{sp}},\X)$ is algebraic, locally almost of finite presentation and quasi-separated over $B$. Indeed, because $\X^{\rm cl}$ is algebraic this is equivalent to verifying that $\X$ admits an almost perfect cotangent complex, which follows from \cite{HLPreygel}*{Prop.~5.1.10}. We will be concerned mostly with the situation where we start with a classical algebraic stack $\X$ satisfying \ref{hyp2}, and we will consider the \emph{spectral stack of $\bZ^n$-filtered objects}
\[
\Filt^n(\X^{\rm{sp}}) := \iMap((\Theta^n)^{\rm{sp}},\X^{\rm{sp}}).
\]
Note that \Cref{lem:spectral_maps} provides a canonical equivalence $\Filt^n(\X^{\rm{sp}})^{\rm{cl}} \simeq \Filt^n(\X)$.

Crucially, however, the spectral stack of $\bZ^n$-filtered objects $\Filt^n(\X^{\rm{sp}})$ is \emph{not} equivalent to $\Filt^n(\X)^{\rm{sp}}$ as a spectral stack. The difference lies in their deformation theory. Any spectral algebraic $1$-stack $\X$ over $B$ has a canonical cotangent complex $\bL_{\X} \in \D_{qc}(\X)$, which is almost perfect if $\X$ is locally almost finitely presented over $B$.\footnote{Technically the absolute cotangent complex $\bL_\X$ in this context refers to the relative cotangent complex $\bL_{\X/B}$, but as is customary we will suppress $B$ from the notation.} The cotangent complex of $\Filt^n(\X)^{\rm{sp}}$ in the context of spectral algebraic geometry, as well as the cotangent complex of the classical stack $\Filt^n(\X)$ defined using simplicial commutative rings as in \cite{laumon2018champs}, are hard to compute in practice. In contrast, we have

\begin{lem}[\cite{HLPreygel}*{Prop.~4.13}] \label{lem:cotangent_complex_mapping_1}
At a $k$-point of $\Filt^n(\X^{\rm{sp}})$, corresponding to a map $f : \Theta_k^n \to \X$, we have a canonical quasi-isomorphism
\[
(\bL_{\Filt^n(\X^{\rm{sp}})})_f \simeq R\Gamma(\Theta_k^n,(f^\ast \bL_\X)^\dual)^\dual.
\]
Likewise for a graded point $g : \Spec(k) / (\Gm)^{n} \to \X$, we have
\[
(\bL_{\Grad^n(\X^{\rm{sp}})})_g \simeq (g^\ast (\bL_\X))^{0},
\]
where $(-)^0$ denotes the invariant piece of an object of $\D_{qc}(\Spec(k)/(\Gm)^n)$ regarded as a $\bZ^n$-graded complex of vector spaces.
\end{lem}

From this formula we can explicitly compute the cotangent complex at a point $f \in \Filt(\X^{\rm{sp}})(k)$ classifying a split filtration. For a complex of graded vector spaces $E$, we let $E^{\geq 0}$ and $E^{<0}$ denote the summand with nonnegative (resp. negative) weights.
\begin{lem} \label{lem:cotangent_complex_mapping}
Let $g : (\pt/\Gm)_k \to \X$ be a map and let $f = \sigma(g) : \Theta_k \to \X$ be the corresponding split filtration. Then we have a canonical quasi-isomorphism of exact triangles
\[
\xymatrix{(\bL_{\Filt(\X^{\rm{sp}})/\X^{\rm{sp}}}[-1])_f \ar[r] \ar[d]^\simeq &  (\bL_{\X^{\rm{sp}}})_{f(1)} \ar[r]^{\ev_1^\ast} \ar[d]^\simeq & (\bL_{\Filt(\X^{\rm{sp}})})_f \ar[r] \ar[d]^\simeq & \\ (g^\ast \bL_\X)^{>0} \ar[r] & g^\ast \bL_\X \ar[r] & (g^\ast \bL_\X)^{\leq 0} \ar[r] & }
\]
\end{lem}
\begin{proof}
The top row is the exact triangle of cotangent complexes induced by the map $\ev_1 : \Filt(\X^{\rm{sp}}) \to \X^{\rm{sp}}$. Using \Cref{lem:cotangent_complex_mapping_1} we can identify the canonical map $(\bL_{\X^{\rm{sp}}})_{f(1)} \to (\bL_{\Filt(\X^{\rm{sp}})})_f$ with the linear dual of the map of restriction to the point $\one \in \Theta_k$:
\[
((\bL_{\X^{\rm{sp}}}|_{f(1)})^\dual)^\dual \simeq \bL_{\X^{\rm{sp}}}|_{f(1)} \to R\Gamma(\Theta_k,(f^\ast \bL_{\X^{\rm{sp}}})^\dual)^\dual.
\]
Let $p : \Theta_k \to \Spec(k)/\Gm$ denote the projection. The fact that $f = \sigma(g)$ means, by definition, that $f = g \circ p$, which gives a canonical isomorphism $f^\ast(\bL_{\X^{\rm{sp}}}) \simeq p^\ast (g^\ast \bL_{\X^{\rm{sp}}})$. In particular $\bL_{\X^{\rm{sp}}}|_{f(1)} \simeq g^\ast (\bL_{\X^{\rm{sp}}})$ after forgetting the $\Gm$-action on the latter.

After chasing through a few dualizations, the claim follows from the following more concrete claim: for any graded complex of vector spaces, regarded as an object $E \in \D_{qc}((\pt/\Gm)_k)$, the map $R\Gamma(\Theta_k, p^\ast(E)) \to p^\ast(E)|_1 \simeq E$ is isomorphic to the inclusion of the direct summand $E^{\geq 0} \to E$. One can verify this by identifying $R\Gamma(\Theta_k, p^\ast(E))$ with $(E \otimes_k k[t])^{\Gm}$, where $t$ has weight $-1$.
\end{proof}

For most of this paper we will be concerned instead with classical algebraic stacks. As mentioned above, the deformation theory of the spectral stack of filtrations differs from the classical stack of filtrations, so it will be useful to have a comparison theorem between the two. We introduce some temporary notation. For a spectral algebraic 1-stack $\X$, we let $\bL^{E_\infty}_{\X} \in \D_{qc}(\X)$ denote its cotangent complex in the context of spectral algebraic geometry, and by a slight abuse of notation let $\bL^{E_\infty}_{\X^{\rm{cl}}} \in \D_{qc}(\X^{\rm{cl}})$ denote the cotangent complex of the spectral algebraic stack $(\X^{\rm{cl}})^{\rm{sp}}$. Here we are making use of the canonical equivalence of stable $\infty$-categories $\D_{qc}(\X^{\rm{cl}}) \to \D_{qc}((\X^{\rm{cl}})^{\rm{sp}})$, where the former denotes the derived category of complexes with quasi-coherent homology on the classical stack. Finally, we let $\bL_{\X^{\rm{cl}}}^{\rm{cl}} \in \D_{qc}(\X^{\rm{cl}})$ denote the cotangent complex computed in the context of classical algebraic geometry.

There are two canonical maps in $\D_{qc}(\X^{\rm{cl}})$
\[
\xymatrix{ \bL_{\X}^{E_\infty}|_{\X^{\rm{cl}}} \ar[r]^{\phi^{1}_\X} & \bL^{E_\infty}_{\X^{\rm{cl}}} \ar[r]^{\phi^2_\X} & \bL_{\X^{\rm{cl}}}^{\rm{cl}}},
\]
where the first is induced by the inclusion of spectral stacks $(\X^{\rm{cl}})^{\rm{sp}} \hookrightarrow \X$, and the second is induced by a direct comparison between the functors on $\D_{qc}(\X^{\rm{cl}})$ that are corepresented by $\bL^{\rm{cl}}$ and $\bL^{E_\infty}$ respectively. The canonical cofiber sequence of cotangent complexes associated to a map of spectral stacks $\X \to \Y$ are compatible with these comparison maps, which induces comparison maps
\[
\xymatrix{\bL_{\X/\Y}^{E_\infty}|_{\X^{\rm{cl}}} \ar[r]^{\phi^{1}_{\X/\Y}} & \bL^{E_\infty}_{\X^{\rm{cl}}/\Y^{\rm{cl}}} \ar[r]^{\phi^2_\X} & \bL_{\X^{\rm{cl}}/\Y^{\rm{cl}}}^{\rm{cl}} }.
\]
These maps are not equivalences, but they induce isomorphisms in high cohomological degree.
\begin{lem} \label{lem:spectral_cotangent}
Given a map of quasi-separated spectral algebraic $1$-stacks $\X \to \Y$, the canonical map on homology sheaves
\begin{enumerate}
\item $H^n(\phi^1_{\X/\Y}) : H^n (\bL_{\X/\Y}^{E_\infty}|_{\X^{\rm{cl}}}) \to H^n (\bL^{E_\infty}_{\X^{\rm{cl}}/\Y^{\rm{cl}}})$ is surjective for $n \geq -1$ and an isomorphism for $n \geq 0$, and
\item $H^n(\phi^2_{\X/\Y}) : H^n(\bL^{E_\infty}_{\X^{\rm{cl}}/\Y^{\rm{cl}}}) \to H^n(\bL_{\X^{\rm{cl}}/\Y^{\rm{cl}}}^{\rm{cl}})$ is surjective for $n \geq -2$ and an isomorphism for $n\geq -1$.
\end{enumerate}
\end{lem}
\begin{proof}
The claims are equivalent to the claim that $H^n(\op{cofib}(\phi^i_{\X/\Y}))=0$ for $n \geq -1$ in the case $i=1$ and for $n\geq -2$ in the case $i=2$. The canonical exact triangle of cotangent complexes associated to the map $\pi : \X \to \Y$ induces an exact triangle
\[
\pi^\ast \op{cofib}(\phi^i_{\Y}) \to \op{cofib}(\phi^i_\X) \to \op{cofib}(\phi^i_{\X/\Y}) \to,
\]
which shows that the claim of the lemma can be deduced from the absolute case, i.e. the claim that $H^n(\op{cofib}(\phi^i_\X)) = 0$ for $n$ in the appropriate range and for any spectral algebraic $1$-stack $\X$.

Let $\X$ be a spectral stack and consider a smooth surjective map $X \to \X$ from a spectral scheme $X$. The exact triangle of cotangent complexes associated to the map $X \to \X$ allows one to reduce both claims to the corresponding claim for $\bL_X$ and $\bL_{X/\X}$. All three versions of the cotangent complex $\bL^{\rm{cl}}_{X^{\rm{cl}}/\X^{\rm{cl}}}$, $\bL^{E_\infty}_{X^{\rm{cl}}/\X^{\rm{cl}}}$, and $\bL^{E_\infty}_{X/\X}|_{\X^{\rm{cl}}}$ are compatible with base change over $\X$, which allows us to reduce the general case to the case for spectral schemes, and ultimately to showing the relevant vanishing of $H^n(\op{cofib}(\phi^i_{X})$ where $X$ is an affine spectral scheme. Claim (1) now follows from \cite{lurie2012higher}*{Thm.~7.4.3.1} and Claim (2) follows from the identification of topological Andr\'e-Quillen cohomology and classical Andr\'e-Quillen cohomology of commutative rings in low homological degree (See \cite{lurie2012representability}*{Warn.~1.0.7}).
\end{proof}

In the remainder of the article, we will not make use of the intermediate object $\bL^{E_{\infty}}_{\X^{\rm{cl}}/\Y^{\rm{cl}}}$, so we can safely drop the superscript from our notation: for spectral stacks $\bL_{\X/\Y}$ will always denote the cotangent complex in the spectral context, and for classical stacks $\bL_{\X/\Y}$ will denote the cotangent complex in the classical context.


\subsection{Some general properties of \texorpdfstring{$\Filt^n(\X)$}{Filtn(X)}}

\begin{prop} \label{prop:open_closed_filtrations}
Let $\X$ and $\Y$ be stacks satisfying \ref{hyp2}, and let $\phi: \Y \to \X$ be a morphism that is representable by algebraic spaces. Then so is the induced morphism $\Filt^n(\phi) : \Filt^n(\Y) \to \Filt^n(\X)$. Furthermore:
\begin{enumerate}
\item If $\phi$ is a monomorphism, then so is $\Filt^n(\phi)$.
\item If $\phi$ is a closed immersion, then so is $\Filt^n(\phi)$, and $\Filt^n(\phi)$ identifies $\Filt^n{\Y}$ with the closed substack $\ev_1^{-1} \Y \subset \Filt^n(\X)$.
\item If $\phi$ is an open immersion, then so is $\Filt^n(\phi)$, and $\Filt^n(\phi)$ identifies $\Filt^n(\Y)$ with the preimage of $\Y \subset \X$ under the composition $$\Filt^n(\X) \xrightarrow{\agr} \Grad^n(\X) \to \X.$$
\item If $\phi$ is smooth (respectively \'etale), then so are $\Filt^n(\phi)$ and $\Grad^n(\phi) : \Grad^n(\Y) \to \Grad^n(\X)$, and this holds even if $\phi$ is not representable.
\end{enumerate}
\end{prop}
\begin{proof}

Let $S \to \Filt^n(\X) =  \filt[n]{\X}$ be an $S$-point defined by a morphism $\Theta^n_S \to \X$. Then the fiber product $\Theta^n_S \times_{\X} \Y \to \Theta^n_S$ is representable and is thus isomorphic to $E/\Gm^n$ for some algebraic space $E$ with a $\Gm^n$-equivariant map $E \to \bA^n_S$. The fiber of $\Filt^n(\Y) \to \Filt^n(\X)$ over the given $S$-point of $\Filt^n(\X)$ corresponds to the groupoid of sections of $E/\Gm^n \to \bA^n_S/\Gm^n$, which form a set. Thus $\Filt^n(\Y)$ is equivalent to a sheaf of sets as a category fibered in groupoids over $\Filt^n(\X)$. Both $\Filt^n(\X)$ and $\Filt^n(\Y)$ are algebraic by \Cref{prop:existence_mapping}, so $\Filt^n(\phi)$ is representable by algebraic spaces.

Say $\phi$ is a monomorphism, meaning it induces a fully faithful embedding on groupoids of $S$-points for any $B$-scheme $S$. Then smooth descent implies implies that $\Map(\Theta^n_S,\X) \to \Map(\Theta^n_S,\Y)$ is fully faithful as well, so $\Filt^n(\Y) \to \Filt^n(\X)$ is a monomorphism. To identify the full subfunctor $\Filt^n(\Y) \subset \Filt^n(\X)$, one only needs to identify which maps $f : \Theta^n_S \to \X$ factor through $\Y$ for each $B$-scheme $S$.

The smallest closed substack of $\Theta^n_S$ containing the $S$-point determined by the point $(1,\ldots,1) \in \bA^n$ is $\Theta^n_S$ itself, so if $\phi$ is a closed immersion then $f$ factors through $\Y$ if and only if the composition $\{1\} \times S \to \Theta^n_S \to \X$ factors through $\Y$, which shows (2). Likewise any open substack of $\Theta_S^n$ containing $(\{0\}/\Gm^n)_S$ contains all of $\Theta^n_S$, so if $\phi$ is an open immersion then $f$ factors through $\Y$ if and only if the induced map $(\{0\}/\Gm^n)_S  \to \X$ factors through $\Y$. Finally the map $(\pt/\Gm^n)_S \to S$ induces a bijection between posets of open substacks, which shows (3).

It suffices to show (4) in the case $n=1$, because $\Filt^n(-) \simeq \Filt(\Filt(\cdots))$ and likewise for $\Grad^n(-)$. The computation of the cotangent complex of the spectral mapping stack in \Cref{lem:cotangent_complex_mapping_1} and \Cref{lem:cotangent_complex_mapping}, combined with the fact that the projection functors onto weight spaces $(-)^0$,$(-)^{\leq 0}$, and $(-)^{>0}$ are exact, implies that $\Grad(\Y^{sp}) \to \Grad(\X^{sp})$ is smooth, and $\Filt(\Y^{sp}) \to \Filt(\X^{sp})$ is smooth at every split filtration in $\Filt(\Y)$. The comparison result \Cref{lem:spectral_cotangent} shows that the same is true for the classical mapping stacks. Finally, we will see below in \Cref{lem:retract} that every point of $\Filt(\Y)$ specializes to split point, so because smoothness is an open condition on the source of map that is locally of finite presentation, it follows that $\Filt(\Y) \to \Filt(\X)$ is smooth.
\end{proof}

Another important property of the induced map $\Filt^n(\X) \to \Filt^n(\Y)$ is the following:

\begin{prop}\label{prop:affine_map}
Let $\pi : \X \to \Y$ be an affine (respectively finite) representable morphism of stacks. Then for any $n$ the map $\Filt^n(\X) \to \Filt^n(\Y)$ is an affine (respectively finite) representable morphism, and the canonical map $\Filt^n(\X) \to \Filt^n(\Y) \times_\Y \X$ is a closed immersion (respectively a surjective closed immersion), where the fiber product is taken with respect to the morphism $\ev_1 : \Filt^n(\Y) \to \Y$.
\end{prop}

\begin{lem}
Let $T$ be a scheme, and let $\cA = \bigoplus_{n \in \bZ} \cA_n$ be a graded $\cO_T[t]$-algebra, where $t$ has degree $-1$. Then a graded map of $\cO_T[t^{\pm}]$-algebras $\psi : \cA[t^{-1}] \to \cO[t^{\pm}]$ is the localization of a map of graded $\cO_T[t]$-algebras $\tilde{\psi} : \cA \to \cO_T[t]$ if and only if the composition
\[\cA_1 \xrightarrow{\times t} \cA[t^{-1}] \xrightarrow{\psi} \cO_T[t^\pm]\]
vanishes. If such a $\tilde{\psi}$ exists then it is unique.
\end{lem}

\begin{proof}
$\psi$ is uniquely determined by its restriction to $\cA$, by the universal property of the localization. For degree reasons, graded maps $\cA \to \cO_T[t]$ factor uniquely through the quotient $\cA \to \bigoplus_{n \leq 0} \cA_n / t^{n+1} \cdot \cA_1$, so it follows that the restriction of $\psi$ to $\cA$ factors through $\cO_T[t] \subset \cO_T[t^\pm]$ if and only if $\psi(t \cdot \cA_1)=0$.

If $\psi$ is the localization of a map $\tilde{\psi} : \cA \to \cO_T[t]$, then the restriction of $\psi$ to $\cA$ factors through $\cO_T[t]$ and thus annihilates $t \cdot \cA_1$. Conversely if $\psi(t \cdot \cA_1) = 0$ then we have a map of algebras $\tilde{\psi}:\cA \to \cO_T[t]$ such that $\tilde{\psi}[t^{-1}]$ agrees with $\psi$ after restricting to $\cA$, and thus $\psi = \tilde{\psi}[t^{-1}]$.
\end{proof}

\begin{proof}[Proof of \Cref{prop:affine_map}]
We first prove the claim when $n=1$. Consider a map $T \to \Filt(\Y) \times_\Y \X$. This corresponds to a map $f:T \times \Theta \to \Y$, and a cosection of the sheaf of algebras $f^\ast \pi_\ast \cO_\X$ over $T \times (\bA^1-\{0\}) \times T / \Gm$. Under the identification between quasi-coherent sheaves of $T \times \Theta$ and graded $\cO_T[t]$-modules \Cref{prop:quasicoh_theta}, we can identify $f^\ast \pi_\ast \cO_\X$ with a graded $\cO_T[t]$-algebra, $\cA$, and the section is the same as a map of $\cO_T[t^\pm]$-algebras $\psi : \cA[t^{-1}] \to \cO_T[t^\pm]$.

On the other hand, the set of lifts of $T \to \Filt(\Y) \times_\Y \X$ to $\X$ corresponds bijectively to the set of maps $\tilde{\psi} : \cA[t] \to \cO_T[t]$ for which $\psi$ is the localization. By the previous lemma this is unique and exists if and only if the ideal $\psi(t \cdot \cA_1) \subset \cO_T$ vanishes. Because the formation of this ideal is compatible with pulling back along a map $(g,\id_\Theta) : T' \times \Theta \to T \times \Theta$, it follows that the fiber product $T \times_{\Filt(\Y) \times_\Y \X} \Filt(\X)$ is represented by the closed subscheme of $T$ defined by this ideal.

We have shown the morphism $\Filt(\X) \to \Filt(\Y) \times_\Y \X$ is a closed immersion and hence affine. It follows that the composition
\[
\Filt(\X) \to \Filt(\Y) \times_\Y \X \to \Filt(\Y)
\]
is affine as well, and it is finite if the map $\pi$ is finite. Furthermore if $\pi$ is finite, then the valuative criterion for properness for the map $\pi$ implies that $\Filt(\X) \to \Filt(\Y) \times_\Y \X$ is surjective, because the fiber for any $k$-point of $\Filt(\Y) \times_\Y \X$ consists of the set of lifts of the map $\Theta_k \to \Y$ to a map $\Theta_k \to \X$ given a lift at the generic point. We can identify $\Filt^n(\X) \simeq \Filt(\Filt^{n-1}(\X))$, so the morphism $\Filt^n(\X) \to \Filt^n(\Y)$ is affine (respectively finite) for any $n \geq 1$ by induction.

Finally, we prove that $\Filt^n(\X) \to \Filt^n(\Y) \times_\Y \X$ is a closed immersion for $n > 1$. Again using the identification of $\Filt^n(\X)$ as an iterated mapping stack, we can factor the map $\ev_1 : \Filt^n(\X) \to \X$ into a sequence of maps $\Filt^n(\X) \to \Filt^{n-1}(\X) \to \cdots \to \X$. We consider the following commutative diagram in which every square is Cartesian
$$\xymatrix{
\Filt^n(\X) \ar[r]_<>(.5){\op{cl. imm.}}|{\cdots} & \Z_{n,2} \ar[r]_<>(.5){\op{cl. imm.}} \ar[d]|{\vdots} & \Z_{n,1} \ar[r]_<>(.5){\op{cl. imm.}} \ar[d]|{\vdots} & \Filt^n(\Y) \times_\Y \X \ar[r]_<>(.5){\op{affine}} \ar[d]|{\vdots} & \Filt^n(\Y) \ar[d]|{\vdots} \\
& \Filt^2(\X) \ar[r]_<>(.5){\op{cl. imm.}} & \Z_{2,1} \ar[r]_<>(.5){\op{cl. imm.}} \ar[d] & \Filt^2(\Y) \times_\Y \X \ar[r]_<>(.5){\op{affine}} \ar[d] & \Filt^2(\Y) \ar[d] \\
& & \Filt(\X) \ar[r]_<>(.5){\op{cl. imm.}} \ar[r] & \Filt(\Y) \times_\Y \X \ar[r]_<>(.5){\op{affine}} \ar[d] & \Filt(\Y) \ar[d] \\
& & & \X \ar[r]^\pi_<>(.5){\op{affine}} & \Y}$$
In each row, we know that the left-most arrow is a closed immersion because the composition of the maps in the row below is affine and we have already shown the claim when $n=1$. The claim thus follows for all $n$ by induction. The same argument shows that $\Filt^n(\X) \to \Filt^n(\Y) \times_\Y \X$ is a surjective closed immersion if $\pi$ is finite, by replacing the label ``affine" with ``finite" and replacing the label ``closed immersion" with ``surjective closed immersion" in the diagram above.
\end{proof}

\begin{cor} \label{cor:embedding_flag_space}
Let $\pi : \Y \to \X$ be a finite map of stacks that both satisfy \ref{hyp2}. Then for any $B$-algebra $R$ and $\xi : \Spec(R) \to \Y$ the canonical map $\Flag^n(\xi) \to \Flag^n(\pi \circ \xi)$ is a surjective closed immersion (and hence a universal homeomorphism) of algebraic stacks and an isomorphism if $\pi$ is a closed immersion.
\end{cor}
\begin{proof}
The case of a closed immersion follows immediately from part (2) of  \Cref{prop:open_closed_filtrations}, which implies that in the diagram
\[
\xymatrix{\Flag^n(\xi) \ar[r] \ar[d] & \Filt^n(\Y) \ar[r] \ar[d]^{\ev_1} & \Filt^n(\X) \ar[d]^{\ev_1} \\ \Spec(R) \ar[r]^\xi & \Y \ar[r]^\pi & \X},
\]
the right square is cartesian and thus so is the outermost rectangle. For the more general case where $\pi$ is finite, one argues not that the squares are cartesian, but that the canonical map from the left corner to the fiber product is a surjective closed immersion, using \Cref{prop:affine_map}. A diagram chase shows that for this property of $2$-commutative squares, if the rightmost square has the property, then the left square has the property if and only if the combined square does.
\end{proof}

\subsubsection{Deformation retract of $\Filt^n(\X)$.}

We shall say that a \emph{$\Theta$-deformation retract of a stack $\Y$ onto a stack $\Z$} consists of morphisms $\pi : \Y \to \Z$, $\sigma : \Z \to \Y$, and $r : \Theta \times \Y \to \Y$ such that: 1) $\pi \circ \sigma \simeq \id_\Z$, 2) $r|_{\{1\} \times \Y} \simeq \id_\Y$, 3) $r|_{\{0\} \times \Y} \simeq \sigma \circ \pi$, and 4) the two compositions in the square
\[
\xymatrix@R=10pt{\Theta \times \Y \ar[r]^{r} \ar[d]_{\rm{project}_\Y} & \Y \ar[d]^\pi \\ \Y \ar[r]^\pi & \Z}
\]
are isomorphic. We shall see in \Cref{thm:describe_strata_global_quotient} that $\sigma :\Z \to \Y$ need not be a monomorphism. Note that restricting $r$ to a map $\bA^1 \times\Y \to \Y$ gives an $\bA^1$-deformation retract, so a $\Theta$-deformation retract is a stronger notion.

\begin{lem} \label{lem:full_opens}
Given a $\Theta$-deformation retract of an algebraic stack $\Y$ onto $\Z$, if $\fU \subset \Y$ is an open substack such that $\sigma^{-1}(\fU) = \Z$, then $\fU = \Y$.
\end{lem}

\begin{proof}
Consider the open substack $r^{-1}(\fU) \subset \Theta \times \Y$. Because the composition $\{0\} \times \Y \to \Theta \times \Y \to \Y$ factors through $\sigma$ and $\sigma^{-1}(\fU) = \Z$, we have
\[
(\{0\}/\Gm) \times \Y \subset r^{-1}(\fU).
\]
We choose an atlas $Y \to \Y$ and note that the restriction of $r^{-1}(\fU)$ along the map $\Theta \times Y \to \Theta \times \Y$ is an open substack of $\Theta \times Y$ that contains $(\{0\}/\Gm) \times Y$. Such an open substack is all of $\Theta \times Y$ by necessity, so $r^{-1}(\fU) = \Theta \times \Y$. This implies that the restriction $\id_\Y \simeq r|_{\{1\}\times \Y} : \Y \to \Y$ factors through $\fU$, hence $\fU = \fY$.
\end{proof}

\begin{lem} \label{lem:quasi_compact}
Given a $\Theta$-deformation retract of an algebraic stack $\Y$ onto $\Z$, if $\pi : \Y \to \Z$ is locally finitely presented then it is quasi-compact.
\end{lem}
\begin{proof}
Let $S$ be a quasi-compact $B$-scheme and consider a map $\xi : S \to \Z$, then the stack $\Y' := \Y \times_{\pi,\Z,\xi} S$ is locally finitely presented over $S$. The definition of a $\Theta$-deformation retract is preserved under base change along a map to $\Z$, so $\Y'$ admits a $\Theta$-deformation retract onto $S$. We must show that $\Y'$ is quasi-compact. Because $S$ is quasi-compact one can find a smooth map from a quasi-compact scheme $Y \to \Y'$ such that $\sigma : S \to \Y'$ factors through the image of $Y \to \Y'$. \Cref{lem:full_opens} now implies that the open substack $\op{im}(Y \to \Y') \subset \Y'$ is all of $\Y'$, so $Y \to \Y'$ is surjective and hence $\Y'$ is quasi-compact.
\end{proof}

\begin{lem} \label{lem:retract_components}
Given a $\Theta$-deformation retract of an algebraic stack $\Y$ onto $\Z$, $\pi$ induces a bijection on connected components with inverse given by $\sigma$.
\end{lem}
\begin{proof}
Consider a map of stacks $\varphi : \Y \to \Y$ that is $\bA^1$-homotopic to the identity in the sense that there is a map $r : \bA^1 \times \Y \to \Y$ with $\{1\} \times \Y \to \Y$ isomorphic to the identity and $\{0\} \times \Y \to \Y$ isomorphic to $\varphi$. Then restricting $r$  to $\bA^1 \times \op{Spec}(k)$ for any $k$-point of $\Y$ shows that any point of $\Y$ lies on the same connected component as some point in the image of $\varphi$, so $\varphi : \pi_0(\Y) \to \pi_0(\Y)$ is surjective. More precisely, any $p \in |\Y|$ lies in the same connected component as $\varphi(p)$, which shows that if $p,q\in|\Y|$ are such that $\varphi(p)$ and $\varphi(q)$ lie on the same connected component, then $p$ and $q$ lie on the same connected component as well. This shows that $\varphi : \pi_0(\Y) \to \pi_0(\Y)$ is in fact bijective. The statement of the lemma now follows from this general fact applied to $\varphi = \sigma \circ \pi$ and from the fact that $\pi \circ \sigma \simeq \id$ and is thus also bijective on connected components.
\end{proof}

We now apply these lemmas in our situation. Consider the map $\Theta \times \Theta^n \to \Theta^n$ induced by the scalar multiplication map $\bA^1 \times \bA^n \to \bA^n$ taking $(t,v) \mapsto tv$, which is equivariant with respect to the group homomorphism $\Gm \times \Gm^n \to \Gm^n$ taking $(t,z_1,\ldots,z_n) \mapsto (tz_1,\ldots,tz_n)$. This defines a map
\begin{equation}\label{eqn:retract}
r : \Theta \times \filt[n]{\X} \to \filt[n]{\X}
\end{equation}
that takes any $S$-point, corresponding to a pair $(a : S \to \Theta, b : \Theta^n_S \to \X)$, to the $S$-point of $\filt[n]{\X}$ corresponding to the composition
\[
\Theta_S^n \xrightarrow{(a,\id)} \Theta \times \Theta_S^n \xrightarrow{\scriptstyle \txt{scalar\\multiplication}} \Theta_S^n \to \X.
\]

\begin{lem} \label{lem:retract}
The maps $\pi = \agr$, $\sigma$, and the map $r$ from \eqref{eqn:retract} define a $\Theta$-deformation retract of the stack $\Filt^n(\X)$ onto $\Grad^n(\X)$. Hence if $\X$ satisfies \ref{hyp2}, then $\agr$ is quasi-compact and induces a bijection on connected components whose inverse is given by $\sigma$.
\end{lem}
\begin{proof}
This is an immediate consequence of the definition of $r$ on $S$-points. If $a = 1$ is the constant map $S \to \Theta$, then the composition $\Theta_S^n \xrightarrow{(a,\id)} \Theta \times \Theta_S^n \to \Theta_S^n$ is isomorphic to the identity. If $a=0$ is the constant map, then this composition is isomorphic to the map $\Theta_S^n \to \Theta_S^n$ induced by the $\Gm^n$-equivariant map $\bA^n \mapsto \{0\} \subset \bA^n$. The fact that the map $r$ of \eqref{eqn:retract} commutes with the projection $\agr : \Filt^n(\X) \to \Grad^n(\X)$ follows from the fact that for any $a : S \to \Theta$ the composition
\[
(\{0\}/\Gm^n)_S \to \Theta^n_S \to \Theta \times \Theta_S^n \to \Theta_S^n
\]
is canonically isomorphic to the inclusion $(\{0\}/\Gm^n)_S \hookrightarrow \Theta_S^n$, and this isomorphism is natural in $S$.
\end{proof}

\subsubsection{More on connected components of $\Filt^n(\X)$}

Given a $k$-point of $\Filt^n(\X)$, corresponding to a map $f : \Theta^n_k \to \X$, one can consider the ``inertia" $k$-subgroup
\begin{equation} \label{eqn:inertia}
I_f := \ker( (\Gm^n)_k \to \Aut (f(0))).
\end{equation}
We use the notation $D_S(A)$ for the diagonalizable $S$-group scheme associated to a finitely generated abelian group $A$ as in \cite{conrad2014reductive}*{App.~B}. $I_f \subset (\Gm^n)_k$ is an fppf sub-group-scheme and is thus a diagonalizable group of the form $D_k(\bZ^n / N_f)$ for some sub-group $N_f \subset \bZ^n$ \cite{conrad2014reductive}*{Cor.~B.3.3, Proposition B.3.4}. We can use this observation to separate connected components of $\Filt^n(\X)$.

\begin{prop} \label{prop:components}
Let $\X$ be an algebraic stack satisfying \ref{hyp3}. Then the map $f \in \Filt^n(\X)(k) \mapsto N_f \subset \bZ^n$ is locally constant on $\Filt^n(\X)$.
\end{prop}
\begin{proof}
By construction the map $f \mapsto N_f$ factors through the projection $\agr : \Filt^n(\X) \to \Grad^n(\X)$, so it suffices to show that this function is locally constant on $\Grad^n(\X)$. By the same argument as in the proof of \Cref{prop:existence_mapping}, we can reduce to the case where $B = \Spec(R)$ is noetherian and affine and $\X$ is finite type over $B$.

It suffices to consider an integral $B$-scheme of finite type $S$ along with a family of graded objects $f : (\pt/\Gm^n)_S \to \X$. Consider the analogous subgroup
\[
I_f := \ker( (\Gm^n)_S \to \Aut(f_S)) \subset (\Gm^n)_S.
\]
As a kernel of a homomorphism to a separated group scheme, it is closed, and its restriction to each point of $S$ is the subgroup \eqref{eqn:inertia}. Let $1 \in S(K)$ denote the generic point and $0 \in S(k)$ any other point, and let $N_0,N_1 \subset \bZ^n$ be the corresponding subgroups such that $I_f|_{\op{Spec}(K)} = D_K(\bZ^n/N_1)$ and $I_f|_{\op{Spec}(k)} = D_k(\bZ^n/N_0)$. The proposition amounts to the claim that $N_0=N_1$.

\medskip
\noindent \textit{Proof that $N_0 \subset N_1$:}
\medskip

We can find an open subscheme $U \subset S$ for which $I_f|_U$ is flat and hence of multiplicative type \cite{conrad2014reductive}*{Cor.~B.3.3}, and in fact $I_f|_U \simeq D_U(\bZ^n/N_1) \subset D_U(\bZ^n) = (\Gm^n)_U$ by \cite{conrad2014reductive}*{Prop.~B.3.4}. Because $S$ is integral it is the scheme theoretic closure of $U$, and it follows that $D_S(\bZ^n/N_1)$ is the scheme theoretic closure of the open subscheme $D_U(\bZ^n/N_1) \subset D_S(\bZ^n/N_1)$ because $D_S(\bZ^n/N_1) \to S$ is flat and affine. As a consequence, $D_S(\bZ^n/N_1)$ is the scheme theoretic closure of $D_U(\bZ^n/N_1)$ in $(\Gm^n)_S$ and hence $D_S(\bZ^n/N_1) \subset I_f$. In particular $D_k(\bZ^n/N_1) \subset I_f|_{\op{Spec}(k)}$, which establishes that $N_0 \subset N_1$.

\medskip
\noindent \textit{Proof that $N_1 \subset N_0$:}
\medskip

We begin with some general observations about graded points of $\X$. Note that for $N \subset \bZ^n$, $D_k(N)$ is again a torus, and we have a map $p_N : (\pt/\Gm^n)_k \to BD_k(N)$ whose fiber is $BD_k(\bZ^n/N)$. \Cref{lem:graded_stack} implies that a graded point $f : (\pt/\Gm^n)_k \to \X$ factors through the map $p_N$ if and only if $D_k(\bZ^n/N) \subset \ker((\Gm^n)_k \to \Aut(f))$, and the factorization is unique in this case. Thus $N_f$ can be characterized as the largest subgroup $N \subset \bZ^n$ for which $f$ factors through $p_N$. Tannaka duality \cite{coherent_tannaka}*{Thm.~8.4} implies that this factorization exists if and only if the symmetric monoidal functor $f^\ast : \Coh(\X) \to \Coh((\pt/\Gm^n)_k)$ factors through the pullback functor $p_N^\ast : \Coh(BD_k(N)) \to \Coh((\pt/\Gm^n)_k)$.

Given a sheaf $F \in \Coh((\pt/\Gm^n)_k)$ we let $F_\chi$ denote the direct summand of $F$ on which $\Gm^n$ acts with the character $\chi \in \bZ^n$. The pullback functor $p_N^\ast$ is fully faithful, with essential image consisting of sheaves $F$ for which $F_\chi = 0$ for $\chi \notin N$. So, Tannaka duality implies that $f$ factors through $p$ if and only if
\[
\bigcup_{E \in \Coh(\X)} \left\{ \chi \in \bZ^n \text{ s.t. } (f^\ast E)_\chi \neq 0 \right\} \subset N.
\]
Under the characterization of $N_f$ above, we see that $N_f$ is the subgroup generated by those $\chi \in \bZ^n$ for which $\exists E \in \Coh(\X)$ with $(f^\ast E)_\chi \neq 0$.

Now to prove that $N_1 \subset N_0$, let $f_0 : (\pt/\Gm^n)_k \to \X$ and let $f_1 : (\pt/\Gm^n)_K \to \X$ be the graded points at a non-generic and generic point of $S$. Choose a finite generating set $\chi_1,\ldots,\chi_k \in N_1$ for which there exists $E_1,\ldots,E_k \in \Coh(\X)$ for which $(f_1^\ast E_i)_{\chi_i} \neq 0$. Nakayama's lemma applied to $(f^\ast E_i)_{\chi_i} \in \Coh(S)$ implies that $(f_0^\ast E_i)_{\chi_i} \neq 0$ for all $i$ as well, and in particular $\chi_i \in N_0$. This shows that $N_1 \subset N_0$.

\end{proof}

\subsubsection{The action of $\bN^\times$ on $\Filt(\X)$}

We define an action of $\bN^\times$ on $\Theta$ canonically commuting with the inclusion of the point ${1}$. For each $n \in \bN$, the morphism $(\bullet)^n : \Theta \to \Theta$ is defined by the map $\bA^1 \to \bA^1$ given by $z \mapsto z^n$, which is equivariant with respect to the group homomorphism $\Gm \to \Gm$ given by the same formula. Given a morphism $f : \Theta_S \to \X$, we let $f^n$ denote the precomposition of $f$ with $(\bullet)^n$, and we also use $(\bullet)^n$ to denote the precomposition morphism $\Filt(\X) \to \Filt(\X)$. This defines an action of the monoid $\bN^\times$ on $\Filt(\X)$ for which the morphism $\ev_1 : \Filt(\X) \to \X$ is canonically invariant.

\begin{rem}
It follows from the fact that $\ev_1$ is canonically $\bN^\times$-invariant that for any $\xi : T \to \X$ the $\bN^\times$-action induces an action on $\Flag(\xi)$.
\end{rem}

\begin{prop} \label{prop:N_action}
Let $\X$ be an algebraic stack satisfying \ref{hyp3}, and let $n > 0$ be an integer. Then the map $(\bullet)^n : \Filt(\X) \to \Filt(\X)$ is both an open and a closed immersion whose image consists of maps $f : \Theta_k \to \X$ for which the subgroup $I_{f}$ of \eqref{eqn:inertia} contains $(\mu_n)_k$.
\end{prop}

\begin{proof}
As in the proof of \Cref{prop:existence_mapping}, it suffices to prove the claim when $B = \Spec(R)$. Then we can write $\X$ as a filtered union of quasi-compact substacks and apply relative noetherian approximation to reduce to the case where $B=\Spec(R)$ is noetherian and affine and $\X$ is finite type over $B$. So we assume this for the remainder of the proof.

Let us denote the morphism $(\bullet)^n : \Theta \to \Theta$ as $p$. We first show that $p$ is a monomorphism. Under the Reese equivalence (\Cref{prop:quasicoh_theta}) for any scheme $S$, the pullback functor $p^\ast : \QCoh(\Theta_S) \to \QCoh(\Theta_S)$ maps a diagram $\cdots \to E_{w+1} \to E_w \to \cdots$ to the relabeled diagram $E'_\bullet$ with $E'_w = E_{\lfloor w/n \rfloor}$, so $p^\ast$ is fully faithful. Tannaka duality \cite{coherent_tannaka}*{Thm.~8.4(ii)} implies that when $S$ is locally noetherian a map $f : \Theta_S \to \X$ is uniquely determined by the corresponding symmetric monoidal functor $f^\ast : \QCoh(\X) \to \QCoh(\Theta_S)$. The fully faithfulness of $p^\ast$ thus implies the fully faithfulness of the map $\Map(\Theta_S,\X) \to \Map(\Theta_S,\X)$ taking $f \mapsto f \circ p$, so $(\bullet)^n : \Filt(\X) \to \Filt(\X)$ is a monomorphism.

Now consider a finite type $B$-scheme $S$ and a map $f : \Theta_S \to \X$. For any $p \in S(k)$, let $f_p$ denote the restriction $f|_{\Theta_k} : \Theta_k \to \X$. \Cref{prop:components} implies that the set of $p\in |S|$ for which $(\mu_n)_k \subset I_{f_p}$ is both open and closed. In particular there is a unique maximal open subscheme $U \subset S$ for which $f|_U : \Theta_U \to \X$ factors through $(\bullet)^n : \Theta_U \to \Theta_U$ by \Cref{lem:tannaka_factor} below, and the factorization of $f|_U$ is unique up to unique isomorphism by the previous paragraph. This implies that $(\bullet)^n : \Filt(\X) \to \Filt(\X)$ is an open immersion. However we have seen that the image, which is the set of points for which $(\mu_n)_k \subset I_{f_p}$, is closed as well, so $(\bullet)^n$ is an open immersion of locally Noetherian stacks whose image is closed, and is thus a closed immersion as well.
\end{proof}

We have used the following lemma in the previous proof.

\begin{lem} \label{lem:tannaka_factor}
Let $\X$ satisfy \ref{hyp2}, and consider a map $f : \Theta_S \to \X$, where $S$ is a Noetherian $B$-scheme. Then the following are equivalent:
\begin{enumerate}
\item $f$ factors through $(\bullet)^n : \Theta_S \to \Theta_S$,
\item $(\mu_n)_S \subset I_f := \ker((\Gm)_S \to \Aut(f|_{S \times \{0\}}))$,
\item for any $E \in \Coh(\X)$ and $\chi \in \bZ \setminus n\bZ$, we have $(f^\ast (E)|_{(\{0\}/\Gm)_S})_\chi = 0,$ and
\item for any $E \in \Coh(\X)$, $p \in S(k)$, and $\chi \notin n\bZ$, we have $$(f^\ast(E)|_{(\{0\} /\Gm)_{p}})_\chi = 0.$$
\end{enumerate}
In this case the factorization through $(\bullet)^n$ is unique up to unique isomorphism.
\end{lem}
\begin{proof}
First, note that $(3)$ and $(4)$ are equivalent by Nakayama's lemma. The implication $(1) \Rightarrow (2)$ follows from the fact that $(\mu_n)_S$ is the kernel of the map on automorphism groups induced by the restriction of $(\bullet)^n$ to $(\pt/\Gm)_S$. The implication $(2) \Rightarrow (3)$ is a consequence of the effect that pullback along $(\bullet)^n : (\pt/\Gm)_S \to (\pt/\Gm)_S$ has on sheaves $E \in \Coh((\pt/\Gm)_S)$: it simply scales the non-vanishing weights the complex by $n$.

Finally, the implication $(3) \Rightarrow (1)$ is a consequence of Tannaka duality, as has already been discussed in the proof of \Cref{prop:components} above: \cite{coherent_tannaka}*{Thm.~8.4} implies that (1) holds if and only if the symmetric monoidal functor $f^\ast : \Coh(\X_S) \to \Coh(\Theta_S)$ factors through the essential image of the fully faithful pullback functor $\Coh(\Theta_S)\to \Coh(\Theta_S)$ for the map $(\bullet)^n:\Theta_S \to \Theta_S$. This subcategory consists exactly of those sheaves $F$ for which the restriction $F|_{S \times \{0\}/\Gm}$ vanishes in any weight that is \emph{not} divisible by $n$.

\end{proof}

The following is an immediate corollary of \Cref{prop:components} and \Cref{prop:N_action}.

\begin{cor} \label{cor:free_component_action}
Let $\X$ satisfy \ref{hyp3}, and let $n > 0$. Then the morphism $(\bullet)^n : \Filt(\X) \to \Filt(\X)$ induces an isomorphism between connected components, acts injectively on the set of connected components, and only fixes those components for which the subgroup \eqref{eqn:inertia} is all of $\Gm$.
\end{cor}

\subsubsection{Change of base lemmas}

\begin{lem} \label{lem:DM_filtrations}
If $\X$ is a stack satisfying \ref{hyp2} that has quasi-finite relative inertia $I_{\X/B} \to \X$, then $\ev_1 : \Filt^n(\X) \to \X$ and the forgetful map $u : \Grad^n(\X) \to \X$ are equivalences whose inverse classifies the constant maps $\Theta^n \times \X \to \X$ and $(B\Gm^n) \times \X \to \X$ respectively. The map $\agr : \Filt^n(\X) \to \Grad^n(\X)$ is also an equivalence.
\end{lem}
\begin{proof}
Because $\X \to B$ is locally of finite presentation, it suffices to prove the claim for maps from $B$-schemes $S$ that are noetherian. Consider a map $f : \Theta^n_S \to \X_S$ relative to $S$. Then the induced group homomorphism $(\Gm^n)_S \to \underline{\Aut}_{\X/B}(f|_{\{0\}\times S})$ is trivial on geometric fibers, because $\X$ has quasi-finite inertia and $\Gm$ is geometrically connected. It follows that $f^\ast : \Coh(\X) \to \Coh(\Theta^n_S)$ factors through the full symmetric monoidal subcategory consisting of sheaves for which $\Gm^n$ acts trivially on $E|_{\{0\}\times S}$. This subcategory is the essential image of the fully faithful pullback functor along $\Theta^n_S \to S$, so $f$ factors uniquely through the projection $\Theta^n_S \to S$ by \cite{coherent_tannaka}*{Thm.~8.4}. Furthermore, the pullback along $\Theta^n_S \to S$ followed by restriction along $\{1\} \times S \to \Theta^n_S$ is equivalent to the identity functor on $\Coh(S)$, so the map $f$ is uniquely determined up to unique isomorphism by its restriction to $\{1\} \times S$. The arguments for the maps $u$ and $\agr$ are similar.
\end{proof}

\begin{rem}
Note the embedding $I_{\X/B} \hookrightarrow I_{\X} := I_{\X/\Spec(\bZ)}$, so if $\X$ has quasi-finite absolute inertia, for instance if it is a scheme, then the condition in the lemma, that $I_{\X/B} \to \X$ is quasi-finite, holds as well.
\end{rem}

\begin{cor} \label{cor:change_of_base_filt}
Let $\cB'$ be a stack satisfying \ref{hyp2} that has quasi-finite relative inertia $I_{\X/B} \to \X$, and let $\X$ be a stack satisfying \ref{hyp2} relative to $\cB'$. Then the canonical morphism $\iMap_{\cB'}(\Theta^n_{\cB'}, \X) \to \iMap_B(\Theta^n_B,\X)$ is an equivalence of stacks over $B$. The same is true with $\Theta^n$ replaced by $B\Gm^n$.
\end{cor}
\begin{proof}
Directly from the definition of the mapping stack, one can verify that the canonical map is an equivalence
\[
\iMap_{\cB'}(\Theta^n_{\cB'},\X) \cong \cB' \times_{\iMap_B(\Theta^n_B,\cB')} \iMap_B(\Theta^n_B,\X),
\]
where the morphism $\cB' \to \iMap_B(\Theta^n_B,\cB')$ is induced by the constant morphism $\Theta^n \times \cB' \to \cB'$. If $\cB'$ has quasi-finite inertia relative to $B$, then \Cref{lem:DM_filtrations} implies that $\cB' \to \iMap_B(\Theta^n_B,\cB')$ is an equivalence of stacks over $B$, and the claim follows.
\end{proof}

\begin{cor} \label{cor:representable_base_change}
Consider a cartesian diagram of stacks satisfying \ref{hyp2}
\[
\xymatrix{ \X' \ar[r] \ar[d] & \Y' \ar[d] \\ \X \ar[r] & \Y }
\]
in which $\Y$ and $\Y'$ have quasi-finite inertia relative to $B$. Then $\Filt^n(\X') \simeq \Filt^n(\X) \times_\X \X'$, where the fiber product is taken with respect to $\ev_1 : \Filt^n(\X) \to \X$, and $\Grad^n(\X') \simeq \Grad^n(\X) \times_\X \X'$, where the fiber product is take with respect to $u : \Grad^n(\X) \to \X$.
\end{cor}
\begin{proof}
For any $B$-scheme $S$, the functor $\Map(\Theta^n_S,-)$ commutes with homotopy limits, so the diagram
\[
\xymatrix{ \filt[n]{\X'} \ar[r] \ar[d] & \filt[n]{\Y'} \ar[d] \\ \filt[n]{\X} \ar[r] & \filt[n]{\Y} }
\]
is automatically cartesian. \Cref{lem:DM_filtrations} identifies $\filt[n]{\Y'} \cong \Y'$ and $\filt[n]{\Y} \cong \Y$ via the maps $\ev_1$. Therefore the right square and outermost square in the commutative diagram
\[
\xymatrix{ \Filt^n(\X') \ar[d] \ar[r]^{\ev_1} & \X' \ar[d] \ar[r] & \Y' \ar[d] \\ \Filt^n(\X) \ar[r]^{\ev_1} & \X \ar[r] & \Y },
\]
are cartesian, and it follows that the left square is cartesian as well. The argument for $\Grad^n(\X')$ is similar.
\end{proof}


\subsection{Graded and filtered points in global quotient stacks}

Here we compute the stacks of filtered and graded objects first for stacks of the form $X/\GL_N$ in \Cref{thm:describe_strata_quotient_GLN} and then for quotients by split reductive groups over a field in \Cref{thm:describe_strata_global_quotient}.

\subsubsection{Bia{\l}ynicki-Birula type theorems}

First we recall some facts about concentration under the action of $\Gm^n$.

\begin{prop} \label{prop:Hesselink_concentration}
Let $X \to B$ be a quasi-separated locally finitely presented map of algebraic spaces, and let $\Gm^n$ act on $X$ so that the structure map $X \to B$ is $\Gm^n$-invariant. Then the functors
\[
\begin{array}{c}
X^+(T) = \left\{ \Gm^n \text{-equivariant maps } \bA^n \times T \to X \text{ over } B\right\}, \text{ and} \\[3pt]
X^0(T) = \left\{\Gm^n \text{-equivariant maps } T \to X \text{ over } B\right\}
\end{array}
\]
are representable by algebraic spaces that are quasi-separated and locally of finite presentation over $B$. The forgetful map $X^0 \to X$ is a closed immersion, and the map $X^+ \to X^0$ induced by restriction to $0$ is quasi-compact.
\end{prop}

\begin{rem}
When $X$ is a scheme admitting a $\Gm$-equivariant affine open cover, then this is a special case of the main theorem of section 4 of \cite{He81} for which the ``center'' is $C = X$ and the ``speed'' is $m=1$. The general statement of the existence of $Y$ when $B$ is a field and $n=1$ is the main result of \cite{drinfeld2014theorem}. Also when $B$ is a field and $n=1$, a slightly more general version is proved in \cite{alper2015luna}*{Thm.~5.16}, where $X$ is allowed to be a Deligne-Mumford stack, and it is shown that the projection $X^+ \to X^0$ is affine.
\end{rem}

\begin{proof}
This follows from \Cref{prop:existence_mapping} and \Cref{prop:open_closed_filtrations}. Let $B \to \Filt^n((\pt/\Gm^n)_B)$ be the section classifying the projection map $\Theta^n \to \pt / \Gm^n$. We claim that we have a cartesian diagram of functors
\[
\xymatrix{ X^+ \ar[r] \ar[d] & \Filt^n(X/\Gm^n) \ar[d] \\ B \ar[r] & \Filt^n((\pt / \Gm^n)_B)}.
\]
Indeed, $S$ points of this fiber product correspond to sections of the pullback of the map $X/\Gm^n \to (\pt / \Gm^n)_B$ along the canonical projection map $\Theta^n \times S \to \pt/\Gm^n$, which is the same as $\Gm^n$-equivariant sections of the map $X \times \bA^n \to \bA^n$. The hypothesis that $X$ is quasi-separated implies that $X$ and $X/\Gm^n$ have quasi-affine diagonal \cite{stacks-project}*{\href{https://stacks.math.columbia.edu/tag/02LR}{Tag 02LR}}, so \Cref{prop:existence_mapping} implies that $\Filt^n(X/\Gm^n)$ is an algebraic stack locally of finite presentation over $B$ with quasi-affine diagonal, and \Cref{prop:open_closed_filtrations} implies that the map $\Filt^n(X/\Gm^n) \to \Filt^n(\pt/\Gm^n)$ is representable by an algebraic space, hence $X^+$ is representable by an algebraic space quasi-separated and locally of finite presentation over $B$. A similar discussion shows that $X^0$ is representable by an algebraic space quasi-separated and locally of finite presentation over $B$.

To show that the projection $\pi : X^+ \to X^0$ is quasi-compact, we claim that for any quasi-compact open subspace $U \subset X^0$, there is a quasi-compact open subspace $V \subset X^+$ such that $\pi^{-1}(U) \subset V$. Consider the map $X^0 \to X^+$ induced by the $\Gm^n$-equivariant projection $\bA^n_B \to B$. For any $p \in X^+(k)$, corresponding to an equivariant map $f : \bA^n_k \to X$, we can define a family of equivariant maps $f_t : \bA^n_k \to X$ by $f_t(z_1,\ldots,z_n) = f(tz_1,\ldots,tz_n)$. This defines a map $\bA^1_k \to X^+$ such that the image at $0$ is the image of $p$ under the composition $X^+ \to X^0 \to X^+$, and every other point lies in the same orbit for the natural action of $\Gm^n$ on $X^+$. It follows that for every $p \in \pi^{-1}(U)$, its orbit closure lies in the image of $U$ under the map $X^0 \to X^+$. Therefore if we let $U' \subset X^+$ be a quasi-compact open neighborhood of the image of $U$ under this map, and let $V$ be the $\Gm^n$-orbit of $U'$, we have $\pi^{-1}(U) \subset V$ as desired.

Finally, to show that $X^0 \to X$ is a closed immersion, we use \cite{ahr2}*{Thm.20.1} to construct a $\Gm^n$-equivariant \'etale over $X' \to X$, where $X'$ is a disjoint union of affine schemes. \Cref{cor:inertia_preserving} implies that $\Grad^n(X'/\Gm^n) \cong (X' / \Gm^n) \times_{X/\Gm^n} \Grad^n(X/\Gm^n)$, and analogously to the case of $\Filt^n$ above one has
\[
X^0/\Gm^n = \Grad^n(X/\Gm^n) \times_{\Grad^n(B\Gm^n)} B\Gm^n,
\]
and likewise for $(X')^0 / \Gm^n$. This shows that $(X')^0 \cong X^0 \times_X X'$, so it suffices to prove the claim for $X'$. But in the affine case, it is a direct computation that the fixed subscheme $\Spec(A)^0$ is represented by $\Spec(A/I)$, where $I \subset A$ is the ideal generated by all non-invariant homogeneous elements in the $\bZ^n$-graded algebra $A$.
\end{proof}

\begin{cor} \label{cor:contraction_functor}
In the context of \Cref{prop:Hesselink_concentration}, if $X$ is separated, then restriction of a map $\bA^n\times T \to X$ to  $\{1\}\times T \subset \bA^n \times T$ defines a locally finitely presented unramified monomorphism of algebraic spaces $X^+ \hookrightarrow X$ identifying $X^+(T)$ with
$$\left\{ f : T \to X \text{ s.t. }  \Gm^n \times T \xrightarrow{t \cdot f(x)} X \text{ extends to } \bA^n \times T \right\} \subset \Map_B(T,X).$$

\end{cor}
\begin{proof} Restriction to $\{1\}\times T$ identifies the set of equivariant maps $\Gm^n \times T \to X$ with $\op{Hom}(T,X)$. If the corresponding map extends to $\bA^n \times T$ it will be unique because $X$ is separated. Likewise the uniqueness of the extension of $\Gm^n \times \Gm^n \times T \to X$ to $\Gm^n \times \bA^n \times T \to X$ guarantees the $\Gm^n$ equivariance of the extension $\bA^n \times T \to X$. The map is locally finitely presented because both $X^+$ and $X$ are, and unramified by \cite{stacks-project}*{Tag 05VH}.
\end{proof}

\begin{rem}
$X^+$ typically has several connected components, and if $X$ admits an equivariant immersion $X \hookrightarrow \bP^n$ for some linear action of $\Gm$ on $\bP^n$, then $X^+$ is the disjoint union of the Bia{\l}ynicki-Birula strata for the action of $\Gm$ on $X$ and the map $j$ is a local immersion when restricted to each connected component. In general this need not be the case: Let $\Gm$ act on $\bP^1$ fixing two points $\{0\}$ and $\{\infty\}$. If $X$ is the nodal curve obtained by identifying these two points, then $\Gm$ acts on $X$ as well. In this case $X^+ = \bA^1$ and $j$ is the composition $\bA^1 \to \bP^1 \to X$. There is no neighborhood of $\{0\}$ in which $j$ is a local immersion.
\end{rem}

We observe the following strengthened version of the Bia{\l}ynicki-Birula theorem, which follows from the discussion above and the results of \Cref{S:deformation_theory}:
\begin{prop} \label{prop:stratum_regularity}
Let $X \to X'$ be a smooth $\Gm^n$-equivariant map of $B$-spaces that both satisfy the conditions of \Cref{prop:Hesselink_concentration}, and consider the map of algebraic spaces $X^+ \to (X')^+$ defined by composing an equivariant map $T \times \bA^n \to X$ with the map $X \to X'$. Then
\begin{enumerate}
\item At any point $p \in X^{0}(k)$, the fiber of the relative cotangent complex $\bL_{X/X',p}$ is canonically a representation of $(\Gm^n)_k$, and there are canonical isomorphisms
\[
\bL_{X^{0}/(X')^{0},p} \simeq (\bL_{X/X',p})^{\Gm^n} \quad \text{ and } \quad \bL_{X^+/(X')^+,p} \simeq \bL_{X/X',p}^{\rm{weight }\leq 0},
\]
where ``weight $\leq 0$'' denotes direct summands whose weight $(w_1,\ldots,w_n)$ with respect to $\Gm^n$ satisfies $w_i \leq 0$ for all $i$. In particular $X^{0} \to (X')^{0}$ and $X^+ \to (X')^+$ are smooth. \\
\item If furthermore $X \to B$ is a smooth schematic map, then the projection $X^+ \to X^{\Gm^n}$ is an \'etale locally trivial bundle of affine spaces with linear $\Gm^n$-action on the fibers.
\end{enumerate}
\end{prop}

\subsubsection{Stacks of the form $X / \GL_N$, where $X$ is an algebraic space}

We must first establish some notational book keeping. We regard $\Gm^N \subset \GL_N$ as the subgroup of diagonal matrices. For any sequence of integers $(w_1,\ldots,w_N)$ we define a one-parameter subgroup $\lambda_w : \Gm \to \GL_N$ given by $\lambda_w(t) = \op{diag}(t^{w_1},t^{w_2},\ldots,t^{w_N})$. More generally, homomorphisms $\Hom(\Gm^q,\Gm^N)$ correspond bijectively to $q$-tuples of one-parameter subgroups $(\lambda_1,\ldots,\lambda_q)$, or $N \times q$ matrices.

Given a $\psi \in \Hom(\Gm^n,\Gm^N)$, let $X^{\psi,0}$ denote the fixed locus for the $\Gm^n$-action induced via the $\GL_N$ action on $X$ and the homomorphism $\psi : \Gm^n \to \GL_N$. Let $X^{\psi,+}$ denote the algebraic space represented by the functor in \Cref{prop:Hesselink_concentration}, which we call the \emph{blade} corresponding to $\psi$.

\begin{rem}
Note that $X^{\psi,0}$ and $X^{\psi,+}$ typically have several connected components, and that the projection $X^{\psi,+} \to X^{\psi,0}$ is a bijection on connected components. The term ``blade" is sometimes used instead to refer to the connected components of $X^{\psi,+}$.
\end{rem}

Given $\psi \in \Hom(\Gm^n,\Gm^N)$, we let $P_\psi=(\GL_N)^{\psi,+}$ denote the blade for the action of $\GL_N$ on itself by conjugation. Concretely, $P_\psi \subset \GL_N$ is the closed subgroup of block matrices whose entries have nonnegative weight under the action of $\Gm$ on $N\times N$ matrices
\[
M \mapsto \psi(1,\ldots,1,t,1,\ldots,1) M \psi(1,\ldots,1,t,1,\ldots,1)^{-1}
\]
where the position of $t$ ranges over all $N$ possible positions. When $n=1$ this is a standard parabolic subgroup, and for $n>1$ it is an intersection of $n$ such parabolic subgroups. Likewise we define $L_\psi := (\GL_N)^{\psi,0}$ with respect to the conjugation action, i.e. the centralizer of $\psi$. It a closed subgroup of block diagonal matrices of a shape determined by $\psi$, and from the universal properties of the blade and fixed locus one has a canonical split surjective group homomorphisms $P_\psi \to L_\psi$.

Applying the blade construction to the group action map $G \times X \to X$, gives a map
\[
P_\psi \times X^{\psi,+} \simeq (G\times X)^{\psi,+} \to X^{\psi,+},
\]
which satisfies the axioms for a group action of $P_\psi$ on $X^{\psi,+}$. Finally, note that the symmetric group $S_N$ acts on $\Hom(\Gm^n,\Gm^N)$ by conjugation by permutation matrices, and for $w \in S_N \subset \GL_N$, $w \cdot X^{\psi,+} = X^{w\psi w^{-1},+}$ and $w P_{\psi} w^{-1} = P_{w \psi w^{-1}}$.

\begin{thm}\label{thm:describe_strata_quotient_GLN}
Let $X$ be an algebraic space with an action of $\GL_N$, and let $X \to B$ be a quasi-separated locally finitely presented $\GL_N$-invariant map to an algebraic space $B$. There are canonical isomorphisms
\begin{gather*}
\Filt^n(X/\GL_N) \simeq \bigsqcup_{\psi \in \Hom(\Gm^n,\Gm^N) / S_N} X^{\psi,+} / P_\psi \quad \text{and}\\
\Grad^n(X/\GL_N) \simeq \bigsqcup_{\psi \in \Hom(\Gm^n,\Gm^N) / S_N} X^{\psi,0} / L_\psi
\end{gather*}
where the notation $\psi \in \Hom(\Gm^n,\Gm^N) / S_N$ means that we choose a single representative for each $S_N$-orbit on the set of homomorphisms.

We have the following description of the universal maps \eqref{eqn:universal_maps} in this case: $\agr$ corresponds to the projection $X^{\psi,+} \to X^{\psi,0}$, which is equivariant with respect to the group homomorphism $P_\psi \to L_\psi$, and $\ev_1$ corresponds to the map $X^{\psi,+} \to X$, which is equivariant with respect to the inclusion of groups $P_\psi \subset \GL_N$.
\end{thm}

\begin{proof}

\medskip
\noindent \textit{The case where $X = \pt$:}
\medskip
The stack $\Filt^n(\pt/\GL_N)$ is the stack of equivariant vector bundles on the toric variety $\bA^n$, and our description agrees with that Payne's description \cite{Pa08}. The result states that for any $\psi \in \Hom(\Gm^n,\Gm^N)$, we associate the vector bundle $\cE_\psi = \cO^{\oplus N}_{\bA^n}$ with an equivariant structure in which $\Gm^n$ acts in the fiber directions via the homomorphism $\psi$. The Rees construction, \Cref{prop:quasicoh_theta}, identifies the category of equivariant vector bundles on $\bA^n$ with $\bZ^n$-weighted filtrations of the fiber at $1^n \in \bA^n$. Under this identification, the automorphisms of $\cE_\psi$ are precisely the automorphisms of the fiber that preserve this filtration, which is precisely the group $P_\psi$. Thus we have maps $\pt / P_\psi \to \Filt^n(\pt / \GL_N)$, and Payne's result says that each of these maps is an open and closed immersion. For completeness, we recall his argument at the end of this proof. The computation of $\Grad^n(\pt/\GL_N)$ is similar, but uses the identification of vector bundles on $\pt / \Gm^n$ with $\bZ^n$-graded vector spaces, rather than the Rees construction.

\medskip
\noindent \textit{The general case:}
\medskip

Consider the cartesian square
\[
\xymatrix{Y \ar[r] \ar[d]  & \cY \ar[r] \ar[d] & \Filt^n(\X) \ar[d] \\
\pt \ar[r] & \pt / P_\psi \ar[r]^-{\cE_\psi} & \Filt^n(\pt / \GL_N) }.
\]
We know from \Cref{prop:open_closed_filtrations} that $Y$ is actually representable by an algebraic space, and therefore $\cY = Y/P_\psi$ for some action of $P_\psi$ defined by the existence of this cartesian square. Note that if $\op{Frame}(\cE_\psi)$ denotes the frame bundle, then $\op{Frame}(\cE_\psi) \times_{\GL_N} X = \bA^n \times X$ with $\Gm^n$ acting simultaneously on the left via its standard action, and on the right via $\psi : \Gm^n \to \GL_N$. Unraveling the definitions, one can compute for any $B$-scheme $S$
\[
Y(S) = \left\{ \Gm^n\text{-invariant sections of the bundle } \bA^n \times X_S \to \bA^n \right\},
\]
which is equivalent to the functor represented by $X^{\psi,+}$ in \Cref{prop:Hesselink_concentration}. Under this identification $Y \simeq X^{\psi,+}$, the $P_\psi$ action on $Y$ agrees with the action on $X^{\psi,+}$ obtained from applying the blade construction to the $G$ action on $X$.

The argument that reduces the computation of $\Grad^n(X/\GL_N)$ to the computation of $\Grad^n(\pt/\GL_N)$ is identical, but with the $\Gm^n$-invariants construction replacing the blade construction.

\medskip
\noindent \textit{Recalling the computation in \cite{Pa08} of $\Filt^n(\pt/\GL_N)$.}
\medskip

Given a scheme $S$ and an equivariant locally free sheaf $\cV$ on $\Theta^n_S$, the restriction to $\{(0,\ldots,0)\}\times S / \Gm^n$ is a $\bZ^n$-graded locally free sheaf which splits into a direct sum of sub-bundles of constant weight with respect to $\Gm^n$. The ranks of each weight bundle are locally constant, and this data is encoded by a homomorphism $\psi : \Gm^n \to \GL_N$ up to conjugation, which we can assume comes from a homomorphism $\psi : \Gm^n \to \Gm^N$. Assume that $S$ is connected, so that $\psi$ is constant on $S$. Now for any point $s \in S$ one can choose an affine open neighborhood $s \in U \subset S$ and sections $v_1,\ldots,v_N \in \Gamma(\bA^n_U,\cV)$ that are eigenvectors for the $\Gm^n$ action and restrict to a basis in the fiber of $\cV$ over $s \times \{(0,\ldots,0)\}$. These sections must give a framing of the locally free sheaf in a $\Gm^n$-equivariant neighborhood of this point in $\bA^n_U$, and they define an isomorphism $\cE_\psi \simeq \cV|_{\Theta^n_{U'}}$ for some smaller open subset $s \in U' \subset U$. Because we get such an isomorphism in an open neighborhood of any point, we see that $\cV$ is the locally free sheaf on $\Theta^n_S$ associated to a principal $\Aut(\cE_\psi) \simeq P_\psi$-bundle over $S$.

\end{proof}

\subsubsection{Quotients by split algebraic groups over a field}

A result of Totaro \cite{totaro2004resolution} implies that any finite normal Noetherian stack with the resolution property is equivalent to a stack of the form $X/\GL_N$ for some quasi-affine scheme $X$. In particular this includes stacks of the form $X/G$, where $G$ is an algebraic group over a field $k$ with split maximal torus, and $X$ is a $G$-quasi-projective scheme. It will nevertheless be useful to have explicit descriptions of the stack of filtered and graded objects that more closely reflect the representation theory of $G$.

For any homomorphism $\psi : (\Gm^n)_k \to G$ we consider the blade $X^{\psi,+}$ for the action of $(\Gm^n)_k$ on $X$ via $\psi$ as in \Cref{prop:Hesselink_concentration}. Likewise we define $P_\psi$ to be the blade for the conjugation action of $(\Gm^n)_k$ on $G$ via $\psi$. If $G$ is reductive and $n = 1$, then $P_\psi \subset G$ is a standard parabolic, but $P_\psi$ need not be parabolic generally. As before, $P_\psi$ acts naturally on $X^{\psi,+}$ because the blade construction commutes with products. Similarly the closed subgroup $L_\psi \subset G$, the centralizer of $\psi$ in $G$, acts naturally on $X^{\psi,0}$.

We wish to define a map $X^{\psi,+} / P_\psi \to \Filt^n(X/G)$ by constructing a map $\Theta^n \times (X^{\psi,+} / P_\psi) \to X/G$. The latter map classifies a $(\Gm^n)_k \times P_\psi$-equivariant $G$-bundle over $\bA^n_k \times X^{\psi,+}$ along with a $G$-equivariant and $(\Gm)_k \times P_\psi$ invariant map to $X$. We use the trivial $G$-bundle $\bA^n_k \times X^{\psi,+} \times G$ equipped with a $(\Gm^n)_k \times P_\psi$-equivariant structure via the left action
$$(t,p) \cdot (z,x,g) = (tz,p\cdot x, \psi(tz) p \psi(z)^{-1} g)$$
This expression is only well defined when $z \neq 0$, but it extends to a regular morphism because $\lim_{z \to 0} \psi(z) p \psi(z)^{-1} = l$ exists. It is straightforward to check that this defines an action of $(\Gm^n)_k \times P_\psi$, that the action commutes with right multiplication by $G$, and that the map $\bA_k^n \times X^{\psi,+} \times G \to X$ defined by
$$(z,x,g) \mapsto g^{-1} \psi(z) \cdot x$$
is $(\Gm^n)_k \times P_\psi$-invariant.

It is simpler to construct a map $\pt / (\Gm^n)_k \times X^{\psi,0}/ L_\psi \to X/G$ and thus a map $X^{\psi,0}/ L_\psi \to \Grad^n(X/G)$. We simply use the inclusion of schemes $X^{\psi,0} \hookrightarrow X$, which is equivariant with respect to the group homomorphism $\Gm \times L_\psi \to G$ given by $(t,l) \mapsto \psi(t)l \in G$.

\begin{thm} \label{thm:describe_strata_global_quotient}
Let $\X = X/G$ be a quotient of a $k$-scheme $X$ by a smooth affine $k$-group $G$ with a split maximal torus $T \subset G$ and Weyl group $W$. The natural maps  $X^{\psi,+} / P_\psi \to \Filt^n(\X)$ and $X^{\psi,0} / L_\psi \to \Grad^n(\X)$ induce isomorphisms
\begin{gather*}
\Filt^n(\X) \simeq \bigsqcup_{\psi \in \Hom(\Gm^n,T) / W} X^{\psi,+} / P_\psi, \quad \text{and}\\
\Grad^n(\X) \simeq \bigsqcup_{\psi \in \Hom(\Gm^n,T)/W} X^{\psi,0} / L_\psi.
\end{gather*}
Furthermore, $\agr$ corresponds to the projection $X^{\psi,+} \to X^{\psi,0}$, which is equivariant with respect to the group homomorphism $P_\psi \to L_\psi$, and $\ev_1$ corresponds to the canonical map $X^{\psi,+} \to X$, which is equivariant with respect to the inclusion of groups $P_\psi \subset G$.
\end{thm}

The statement is essentially the same as \Cref{thm:describe_strata_quotient_GLN} and in principle can be reduced to it by choosing a linear embedding $G\hookrightarrow \GL_N$ and identifying $X/G \simeq \GL_N \times_{G} X / \GL_N$. We give a different, more direct proof in \Cref{appendix:filtrations_global_quotient} that does not make use of the Rees construction.

One application of \Cref{thm:describe_strata_global_quotient} is a concrete description of points of the flag scheme $\Flag(p)$ for $p \in X(k)$, when $X$ is separated. Its $k$-points are specified by the three pieces of data:
\begin{itemize}
\item a one parameter subgroup $\lambda \in \Hom(\Gm,T)/W$;
\item a point $q \in X$ such that $\lim_{t \to 0} \lambda(t) \cdot q$ exists; and
\item a $g \in G$ such that $g \cdot q = p$.
\end{itemize}
Where two sets of such data specify the same point of the fiber if and only if $(g',q') = (g h^{-1}, h \cdot q)$ for some $h \in P_\lambda$.

\begin{rem} \label{rem:points_in_fiber}
Alternatively, given such a datum we define the one parameter subgroup $\lambda'(t) := g \lambda(t) g^{-1}$, and $\lim_{t \to 0} \lambda'(t) \cdot p$ exists. The point in $\ev_1^{-1}(p)$ is uniquely determined by this data, thus we can specify a point in the fiber by one parameter subgroup $\lambda$, not necessarily in $T$, for which $\lim_{t \to 0} \lambda(t) \cdot p$ exists. Two one parameter subgroups specify the same point in the fiber if and only if $\lambda' = h \lambda h^{-1}$ for some $h \in P_\lambda$.
\end{rem}


\section{\texorpdfstring{$\Theta$}{Theta}-stratifications} \label{sect:stratifications}

In this section we introduce the notion of a (weak) $\Theta$-stratification of an algebraic stack $\X$, and we establish some general theorems for constructing such stratifications.

Our definition is motivated by the Harder-Narasimhan stratification of the stack $\X$ of vector bundles on a smooth curve (see \Cref{E:vector_bundles} and \Cref{sect:moduli_derived}). The unstable strata in $\X$ parameterize points $p \in \X(k)$ along with a Harder-Narasimhan (HN) filtration, which we regard as a canonical map $f: \Theta_k \to \X$ along with an isomorphism $f(1) \cong p$. More precisely, for each unstable point $p \in \X(k)$, there is a certain connected component of $\Filt(\X)$ such that $(\ev_1)^{-1}(p)$ consists of a single point in this component. This is the perspective we apply to more general algebraic stacks.

In general, a $\Theta$-stratification is determined uniquely by the subset of irreducible components $S \subset \op{Irred}(\Filt(\X))$ that containing some HN filtration, and an indexing map $\mu : S \to \Gamma$ to some totally ordered set $\Gamma$. The main result of this section, \Cref{thm:main_stratification}, gives necessary and sufficient conditions for such data to define a $\Theta$-stratification of $\X$. In later sections, we will establish stronger theorems along these lines (compare to \Cref{thm:main_improved}).

We also discuss conditions under which a $\Theta$-stratification of $\X$ induces a $\Theta$-stratification of a stack $\Y$ via a map $\Y \to \X$ (\Cref{defn:induced_stratum}), and we show that a $\Theta$-stratification of $\X$ always induces a $\Theta$-stratification of $\Grad(\X)$ (\Cref{prop:graded_strata}).

\subsection{Definition and first properties}

\begin{defn} \label{defn:theta_stratum}
Let $\X$ be a stack satisfying \ref{hyp2}. A \emph{$\Theta$-stratum} in $\X$ is a union of connected components $\S \subset \Filt(\X)$ such that the restriction $\ev_1 : \S \to \X$ is a closed immersion. We call $\S$ a \emph{weak $\Theta$-stratum} if $\ev_1$ is finite and radicial.
\end{defn}

\begin{defn} \label{defn:theta_stratification}
Under the same hypotheses, a \emph{(weak) $\Theta$-stratification} of $\X$ consists of:
\begin{enumerate}
\item a totally ordered set $\Gamma$ and a collection of open substacks $\X_{\leq c}$ for $c \in \Gamma$ such that $\X_{\leq c} \subset \X_{\leq c'}$ for $c<c'$ and $\X = \bigcup_c \X_{\leq c}$;
\item a (weak) $\Theta$-stratum in each $\X_{\leq c}$, $\ev_1 : \S_c \hookrightarrow \X_{\leq c}$, such that
\[
\X_{\leq c} \setminus \ev_1(\S_c) = \X_{<c} := \bigcup_{c'<c} \X_{\leq c'};\text{ and}
\]
\item for every point $x \in |\X|$, the set $\{c \in \Gamma | x \in |\X_{\leq c}|\}$ has a minimal element.
\end{enumerate}
We assume there is a minimal element $0 \in \Gamma$. We refer to the open substack $\X_{\leq 0}$ as the \emph{semistable locus} $\X^{\rm{ss}}$ and its complement in $|\X|$ as the \emph{unstable locus} $|\X|^{\rm{us}}$. We allow the situation $\X^{\rm{ss}} = \emptyset$.
\end{defn}

Note that (3) holds automatically if the index set $\Gamma$ is well-ordered.

\begin{rem}\label{rem:alternate_stratification_def}
Say that $\X$ is a stack satisfying \ref{hyp2} with a weak $\Theta$-stratification. We can regard the given substack $\S_c \subset \Filt(\X_{\leq c})$ as an open substack of $\Filt(\X)$ as well under the open immersion $\Filt(\X_{\leq c}) \subset \Filt(\X)$ induced by the open immersion $\X_{\leq c} \subset \X$ (see \Cref{prop:open_closed_filtrations}). Part (2) and (3) of \Cref{defn:theta_stratification} imply that
\[
|\X_{\leq c}| = |\X| \setminus \bigcup_{c'>c} \ev_1(\S_{c'}),
\]
so the stratification is completely determined by the open substack $\S := \bigcup_{c>0} \S_c \subset \Filt(\X)$ and the locally constant function $\mu : \S \to \Gamma$ that takes the value $c$ on $\S_c$. Therefore, an alternative definition of a (weak) $\Theta$-stratification is an open substack $\S \subset \Filt(\X)$ and a locally constant function $\mu : \S \to \Gamma$ such that for any $c \geq 0$ the subset $|\X_{\leq c}|$ above is open, and the corresponding open substacks $\cX_{\leq c} \subset \X$ and $\S_c = \mu^{-1}(c) \subset \Filt(\X)$ satisfy the conditions of \Cref{defn:theta_stratification}.
\end{rem}

\begin{lem}[HN filtrations] \label{lem:HN_filtrations}
Let $\X$ be a stack satisfying \ref{hyp2}, and let $\{\X_{\leq c}\}_{c \in \Gamma}$ be a weak $\Theta$-stratification of $\X$. Then for every unstable point $p \in \X(k)$, there is a unique $c \in \Gamma$ and point $f \in |\S_c|$ such that $p \in |\X_{\leq c}|$ and $\ev_1(f) = p \in |\X|$. Furthermore:
\begin{enumerate}
\item $f$ can be defined over a finite purely inseparable extension of $k$.
\item if $f$ is defined over $k$, and thus corresponds to a $k$-point of $\Filt(\X)$, then $\ev_1$ induces an isomorphism of underlying reduced algebraic groups $\Aut_{\Filt(\X)}(f)_{\rm red} \cong \Aut_\X(p)_{\rm red}$.
\end{enumerate}
\end{lem}

\begin{proof}
Note that the composition $|\S_c| \to |\X_{\leq c}| \to |\X|$ is a locally closed immersion, so we regard the former as a subset of the latter. Property (2) in \Cref{defn:theta_stratification} implies that $|\S_c|$ are disjoint for different $c$, and property (3) implies that every unstable point lies in $\S_{c^\ast}$, where $c^\ast = \min \{c \in \Gamma | x \in |\X_{\leq c}\}$. Thus we have existence and uniqueness of HN filtrations.

We therefore have some field extension $k'/k$ and a $k'$-point of $\Flag(p)$ representing the unique point of $|\Flag(p)|$ lying over the union of connected components $\S_c \subset \Filt(\X_{\leq c})$, one can take $k'$ to be a finite extension because $\Flag(p)$ is locally finite type over $\Spec(k)$ and thus any irreducible component of $\Flag(p)$ contains some finite type point. Furthermore we may assume that $k'/k$ is normal by replacing $k'$ with its normal closure. If $f : \Theta_{k'} \to \X$ is the HN filtration for $p \in \X(k)$, then $f$ is also the HN filtration for $f(1) \in \X(k')$. Uniqueness of the HN filtration implies that $f$ descends to a $k''$-point of $\Flag(p)$ for the purely inseparable extension $k \subset k'' := (k')^{\op{Gal}(k'/k)}$.

The claim (2) follows from the fact that the homomorphism of $k$-group schemes $\Aut_{\Filt(\X)}(f) \to \Aut_\X(p)$ is finite and radicial as a map of schemes, and that implies that it is an isomorphism on underlying reduced subgroups.
\end{proof}

\begin{defn}[HN filtration]
We refer to the filtration $f$ of \Cref{lem:HN_filtrations} as the \emph{Harder-Narasimhan (HN) filtration} of $p \in \X(k)$. Note that if $k$ is perfect, the HN filtration of $p \in \X(k)$ will be defined over $k$.
\end{defn}

Recall the map $\sigma : \Grad(\X) \to \Filt(\X)$ of \eqref{eqn:universal_maps}, and that a point of $\Filt(\X)$ is said to be \emph{split} if it lies in the image of $\sigma$.

\begin{defn} \label{defn:center}
Let $\S$ be a (weak) $\Theta$-stratum in $\X$. Then we define its \emph{center} $\Z^{\rm{ss}}$ to be the union of connected components $\sigma^{-1}(\S) \subset \Grad(\X)$.
\end{defn}

\Cref{lem:retract} guarantees that $\agr$ maps $\S$ to $\Z^{\rm{ss}}$ and hence the canonical morphisms in \eqref{eqn:universal_maps} induce canonical morphisms for any weak $\Theta$-stratum:
\[
\xymatrix{\Z^{\rm{ss}} \ar@/_/[r]_-{\sigma} & \S \ar[l]_-{\agr} \ar[r]^-{\ev_1} & \X}.
\]
In fact, the bijection of \Cref{lem:retract} implies that $\S = \agr^{-1}(\Z^{\rm{ss}}) \subset \Filt(\X)$, so a $\Theta$-stratum can be equivalently specified by the collection of connected components $\Z^{\rm{ss}} \subset \Grad(\X)$. In a $\Theta$-stratification $\X = \bigcup_c \X_{\leq c}$, each center $\Z_c^{\rm{ss}}$ is an open substack of $\Grad(\X)$ by \Cref{prop:open_closed_filtrations}, and the data of the $\Theta$-stratification is uniquely encoded by these open substacks.

\begin{lem} \label{lem:tangent_spaces}
A weak $\Theta$-stratification is a $\Theta$-stratification if and only if at every finite type point $f \in \S_c(k) \subset \Filt(\X)(k)$ that is split, either of the following equivalent conditions holds
\begin{enumerate}
\item the fiber of the relative cotangent complex $(\bL_{\S_{c}/\X})_f \in \APerf(\Spec(k))$ is $0$-connective, i.e. $H^0((\bL_{\S_{c}/\X})_f)=H^1((\bL_{\S_{c}/\X})_f) = 0$, or
\item the canonical map on Lie algebras $\op{Lie}(\Aut_{\S_c}(f)) \to \op{Lie}(\Aut_{\X} (f(1)))$ is surjective.
\end{enumerate}
\end{lem}
\begin{proof}

By replacing $\X$ with $\X_{\leq c}$ it suffices to prove the claim for a single $\Theta$ stratum $\S \hookrightarrow \X$. Note that $\bL_{\S/\X} \simeq \bL_{\S/\X}$, and the proper and radicial map $\ev_1 : \S \to \X$ is a closed immersion if and only if $\bL_{\S/\X}$ is acyclic in non-negative cohomological degree. This is an open condition on $\S$, and can be checked at finite type points, by Nakayama's lemma. Thus $\S$ is a $\Theta$-stratum if and only if (1) holds at all finite type points of $\S$, and the content of the claim is that it suffices to check only at \emph{split} points of $\S$.

The stratum $\S$ is identified with a union of connected components of $\Filt(\X)$, which corresponds under the bijection of \Cref{lem:retract_components} to the union of connected components $\Z^{\rm{ss}} \subset \Grad(\X)$. \Cref{lem:retract} implies that $\S$ admits a $\Theta$-deformation retract onto $\Z^{\rm{ss}}$. The split points of $\S$ are by definition the image of $\sigma : \Z^{\rm{ss}} \to \S$, so \Cref{lem:full_opens} implies that any open substack of $\S$ containing all split finite type points is all of $\S$. The first part of the claim follows.

\medskip
\noindent {\textit Proof that $(1) \Leftrightarrow (2)$:}
\medskip

We can identify the map of Lie algebras in (2) with the $k$-linear dual of the map
\[
H^1((\bL_\X)_f) \to H^1((\bL_{\S})_f),
\]
so we can reformulate the condition (2) as the property that this map is injective. We will consider the spectral stack $\X^{\rm{sp}}$ associated to $\X$, and we will denote by $\tilde{\S} \subset \Filt(\X^{\rm{sp}})$ the collection of connected components whose underlying classical stack is $\S$ under the equivalence of \Cref{lem:spectral_maps}. Note that $\tilde{\S}$ is a spectral $\Theta$-stratum in the sense that $\ev_1 : \tilde{\S} \to \X^{\rm{sp}}$ is a closed immersion, because a map of spectral algebraic stacks is a closed immersion if and only if the underlying map of classical stacks is a closed immersion. \Cref{lem:spectral_cotangent} shows that conditions (1) and (2) are equivalent to the analogous conditions formulated for the spectral $\Theta$-stratum $\tilde{\S} \to \X^{\rm{sp}}$, and we prove it in this context.

Consider the long exact sequence coming from the exact triangle
\[
\xymatrix{ H^0((\bL_{\tilde{\S}})_f) \ar[r]^-{\alpha} & H^0((\bL_{\tilde{\S}/\X^{\rm{sp}}})_f) \ar[r] & H^1((\bL_{\X^{\rm{sp}}})_f) \ar[r] & H^1((\bL_{\tilde{\S}})_f) \ar[r] & 0}.
\]
Note that because $f$ is split, there is a canonical non-trivial homomorphism $(\Gm)_k \to \Aut_{\tilde{S}}(f)$, and this sequence canonically descends to a sequence of representations of $(\Gm)_k$. By \Cref{lem:cotangent_complex_mapping} the complex $(\bL_{\tilde{\S}})_f$ has non-positive weights with respect to $(\Gm)_k$, and the complex $(\bL_{\tilde{\S}/\X^{\rm{sp}}})_f$ has strictly positive weights with respect to $(\Gm)_k$. It follows that $\alpha = 0$ automatically, so $H^0((\bL_{\tilde{\S}/\X^{\rm{sp}}})_f) = 0$ if and only if the map $H^1((\bL_{\X^{\rm{sp}}})_f) \to H^1((\bL_{\tilde{\S}})_f)$ is injective, which shows that (1) and (2) are equivalent.
\end{proof}

\begin{ex}
The relevant modular example for the failure of the map on tangent spaces to be injective is the failure of Behrend's conjecture for the moduli of $G$-bundles on a curve in finite characteristic \cite{heinloth2008bounds}. In that example, the moduli of $G$ bundles on a smooth projective curve $C$ is stratified by the type of the canonical parabolic reduction of a given unstable $G$-bundle. In finite characteristic there are examples where the map $H^1(C,\fp) \to H^1(C,\fg)$ is not injective, where $\fg$ is the adjoint bundle of a principle $G$-bundle and $\fp$ the adjoint bundle of its canonical parabolic reduction.
\end{ex}

\begin{cor} \label{cor:char_0_theta_strat}
Let $\X$ be a stack defined over a field of characteristic $0$ and satisfying \ref{hyp2}. Then any weak $\Theta$-stratification of $\X$ is a $\Theta$-stratification.
\end{cor}
\begin{proof}
We use the characterization (2) of \Cref{lem:tangent_spaces}. The map $\ev_1:\S_c \to \X_{\leq c}$ is by hypothesis, representable, proper, and radicial. It follows that for any $f \in \S_c(k)$, the canonical homomorphism of algebraic $k$-groups $\Aut_{\S_c}(f) \to \Aut_{\X}(f(1))$ is bijective on $k'$ points for any field extension $k'/k$. If $k$ has characteristic $0$, then this implies the map $\Aut_{\S_c}(f) \to \Aut_{\X}(f(1))$ is an isomorphism of group schemes and hence induces an isomorphism of Lie algebras.
\end{proof}


\subsection{First construction of \texorpdfstring{$\Theta$}{Theta}-stratifications}
\label{sect:general_HN_problem}

Let $\X$ be a stack satisfying \ref{hyp2} that has a weak $\Theta$-stratification. We have seen in \Cref{rem:alternate_stratification_def} that the weak $\Theta$-stratification is determined by an open substack $\S \subset \Filt(\X)$ and a locally constant function $\mu : \S \to \Gamma$. Our goal here is to replace $\S$ with data that is more set-theoretic.

For any $B$-stack $\Y$ satisfying \ref{hyp2}, let $\op{Irred}(\Y)$ denote the set of irreducible components of $|\Y|$. The open immersion induces inclusions of sets
\[
\op{Irred}(\S_c) \subset \op{Irred}(\Filt(\X_{\leq c})) \subset \op{Irred}(\Filt(\X)).
\]
The open and closed substack $\S_c \subset \Filt(\X_{\leq c})$ is completely determined by this subset of $\op{Irred}(\Filt(\X))$. Thus the weak $\Theta$-stratification of $\X$ is encoded by a collection of subsets of $\op{Irred}(\Filt(\X))$ labeled by $c \in \Gamma$. More formally, this shows that a weak $\Theta$-stratification is completely specified by the following data:
\begin{equation} \label{eqn:data}
\begin{array}{l}
\text{1) a set of irreducible components } S \subset \op{Irred}(\Filt(\X)), \text{ and}\\ \text{2) a map to a totally ordered set } \mu : S \to \Gamma.
\end{array}
\end{equation}

Given such data on an arbitrary stack, one can extend $\mu$ to a function $|\Filt(\X)| \to \Gamma$ by defining
\[
\mu(f) = \max \left( \{0\} \cup \{\mu(s) | f \text{ lies in irreducible component } s \in S\} \right)
\]
We may assume without loss of generality that suprema exist in $\Gamma$, and one can define a \emph{stability function} on $|\X|$ by
\begin{equation}\label{eqn:stability_function}
M^\mu(p) = \sup \left\{ \mu(f) \left| f \in |\Filt(\X)| \text{ s.t. } f(1) = p \right. \right\} \in \Gamma \cup \{\infty\}
\end{equation}
We regard $p$ as unstable if $M^{\mu}(p)>0$ and semistable otherwise.

\begin{defn} \label{defn:defined_stratification}
We say that the data \eqref{eqn:data} \emph{defines a (weak) $\Theta$-stratification} if the subsets
\begin{equation} \label{eqn:putative_stratification}
\begin{array}{c}
|\X|_{\leq c} := \{ p \in |\X| \text{ s.t. } M^\mu(p) \leq c \} \subset |\X| \\
|\Filt(\X)|_c := \{ f \in |\Filt(\X)| \text{ s.t. } f \text{ lies in } S \text{ and } \mu(f) = M^\mu(f) \}
\end{array}
\end{equation}
are open, and the corresponding open substacks $\X_{\leq c} \subset \X$ and $\S_c \subset \Filt(\X)$ are the data of a (weak) $\Theta$-stratification of $\X$ (\Cref{defn:theta_stratification}).
\end{defn}

If one starts with a weak $\Theta$-stratification and extracts the data \eqref{eqn:data} as above, then this data will define the original stratification in the sense of \Cref{defn:defined_stratification}. Of course, different data \eqref{eqn:data} can define the same weak $\Theta$-stratification, and not every datum \eqref{eqn:data} defines a $\Theta$-stratification in this way.

\begin{thm} \label{thm:main_stratification}
Let $\X$ be an algebraic stack satisfying \ref{hyp3} over a locally noetherian base stack $B$. Then data as in \eqref{eqn:data} define a weak $\Theta$-stratification of $\X$ if and only if the following conditions are satisfied:
\begin{enumerate}
\item \textbf{HN-property:} For all \emph{finite type} unstable points $p \in \X(k)$, $|\Flag(p)|$ contains a unique point $f$ lying over an irreducible component in $S$ with $\mu(f) = M^{\mu}(p)$. This is the \emph{Harder-Narasimhan (HN) filtration} of $p$.
\item \textbf{HN-specialization:} For any valuation ring $R$ with fraction field $K$ and residue field $k$ and any map $\xi : \Spec(R) \to \X$ whose generic point is unstable and a HN filtration $f_K \in \Flag(\xi)(K)$ of $\xi_K$, one has $$\mu(f_K) \leq M^\mu(\xi|_{\Spec(k)}),$$
and when equality holds there is a unique extension of $f_K$ to a filtration $f_R \in \Flag(\xi)(R)$.
\item \textbf{Open strata:} If $\Spec(R) \to \Filt(\X)$ is a map from a discrete valuation ring essentially of finite type over the base $B$ whose special point is an HN filtration, then its generic point is an HN filtration as well.
\item \textbf{Local finiteness:} For any map $\varphi : T \to \X$, with $T$ a finite type affine scheme, there is a finite subset of $S$ such that every unstable finite type point in $T$ has an HN filtration lying on one of these irreducible components.
\item \textbf{HN-consistency:} If $f : \Theta_k \to \X$ is a HN filtration for $f(1)$, then $M^\mu(f(0)) \leq \mu(f)$.
\end{enumerate}
Alternatively, one can replace conditions (2) and (4) with the following:
\begin{enumerate}[label=(S)]
\item \textbf{Simplified HN-specialization:} \label{princ:S} For any \emph{discrete} valuation ring $R$ \emph{essentially of finite type over the base $B$} with fraction field $K$ and residue field $k$, and for any map $\xi : \Spec(R) \to \X$ whose generic point is unstable and a HN filtration $f_K \in \Flag(\xi)(K)$ of $\xi_K$, one has $$\mu(f_K) \leq M^\mu(\xi|_{\Spec(k)}),$$
and when equality holds there is an extension of discrete valuation rings $R' \supset R$ with fraction field $K'$ such $f_K|_{K'}$ extends to a (not necessarily unique) filtration $f_{R'} \in \Flag(\xi)(R')$.
\item[(4')] \textbf{Strong HN-boundedness:} For any map $\xi : T \to \X$, with $T$ a finite type affine scheme, there is a quasi-compact subspace of $\Flag(\xi)$ that contains an HN filtration for every unstable finite type point of $T$.
\end{enumerate}
\end{thm}

The proof of this theorem is conceptually straightforward, but the theorem itself is a useful way of organizing the construction of $\Theta$-stratifications. In the rest of the paper we will restrict our attention to functions $\mu$ and stacks $\X$ for which many of these conditions hold \emph{automatically}. We can make some immediate simplifications before developing further theory:

\begin{simplification}[Locally constant $\mu$] \label{simp:locally_constant}
In all of the examples studied in this paper, the function $\mu$ will be \emph{locally constant}, in which case either ``HN-specialization" condition (2) or \ref{princ:S} immediately implies the ``open strata" condition (3), making the latter redundant. We say that $\mu$ is locally constant if it is induced by a pair
\[
S' \subset \pi_0(\Filt(\X)) \text{ and } \mu': S' \to \Gamma,
\]
in the sense that $S$ is the preimage of $S'$ under the canonical surjective map $\op{Irred}(\Filt(\X)) \to \pi_0(\Filt(\X))$, and $\mu$ is the restriction of $\mu'$ along this map.
\begin{proof}
The HN-specialization condition implies that the open strata condition only fails if for some map $\Spec(R) \to \Filt(\X)$ whose special point is a HN filtration, the image of $\Spec(K) \to \Filt(\X)$ does not lie on an irreducible component of $S$ or the value of $\mu$ is smaller at the generic point than at the special point of $\Spec(R)$. Neither of these things can happen if $\mu$ is locally constant.
\end{proof}
\end{simplification}

We say that a stack $\X$ satisfying \ref{hyp2} has \emph{quasi-compact flag spaces} if for any map $\xi : T \to \X$ from a finitely presented affine $B$-scheme and any connected component $\Y \subset \Filt(\X)$ the fiber product $\Y \times_{\ev_1,\X,\xi} T \subset \Flag(\xi)$ is quasi-compact. This property holds, for instance, for the stack of coherent sheaves on a projective scheme $X$ (see \Cref{ex:main2}), and we give a criterion that guarantees this property for a quasi-compact stack $\X$ in \Cref{prop:quasi_compact_flags}.

\begin{simplification}[Quasi-compact flag spaces]
If $\X$ has quasi-compact flag spaces, then the local finiteness condition (4) of \Cref{thm:main_stratification} is equivalent to the HN boundedness condition. It is also equivalent to requiring a finite set of \emph{connected components} of $\Filt(\X)$, rather than a finite set of irreducible components, such that every unstable point of $T$ has a HN filtration in $\Flag(\xi)$ lying over one of these connected components.

\end{simplification}

\begin{rem} \label{rem:HN-specialization}
Note that a stronger condition which immediately implies the ``HN-specialization" condition is the property that for any map $\xi : \Spec(R) \to \X$ with $M^\mu(\xi|_K) > -\infty$ and a HN filtration $f_K : \Theta_K \to \X$ of $\xi_K$, there is unique extension of $f_K$ to a filtration $f_R : \Spec(R) \to \Filt(\X)$ of $\xi$. This motivates the definition of $\Theta$-reductive stacks, \Cref{defn:reductive_stack} below, for which HN-specialization holds automatically for discrete valuation rings. In \Cref{defn:numerical_invariant} below we provide a method for constructing locally constant $\mu$, and we will show in \Cref{thm:theta_reductive_stratifications} that a simpler version of the strong HN-boundedness condition, \ref{princ:B2}, is equivalent to conditions (1)-(5) for such $\mu$ when $\X$ is $\Theta$-reductive.
\end{rem}

\begin{proof}[Proof of \Cref{thm:main_stratification}]
Let $T$ be a finite type affine scheme and let $\xi : T \to \X$ be a smooth map, and let $s_1,\ldots,s_n \in \op{Irred}(\Filt(\X))$ be the irreducible components capturing all of the HN filtrations of finite type points in $T$, whose existence is guaranteed by either condition $(4)$ or $(4')$. The fact that $\xi$ is smooth implies that $\Flag(\xi) \to \Filt(\X)$ is smooth, so there are finitely many irreducible components $Y_1,\ldots, Y_N \subset \Flag(\xi)$ in the preimage of some $s_i$, and we regard them as subspaces with their reduced structure. By relabeling and combining components with the same value of $\mu$, we will assume that each $Y_i$ is a finite union of irreducible components of $\Flag(\xi)$ such that $\mu = \mu_i$ on each irreducible component, and we will assume that $\mu_i < \mu_j$ for $i < j$.

Any point $p \in |T|$, not necessarily of finite type, must have $M^{\mu}(\xi(p)) \leq \mu_N$. Indeed if there were some point $f' \in |\Flag(\xi)|$ with $\mu(f')>\mu_N$, then the irreducible component containing $f'$ would have to contain finite type points as well with $\mu > \mu_N$, which would contradict the hypothesis that $Y_1,\ldots,Y_N$ contain HN filtrations for all finite type points of $\xi$. This implies that any point $f\in Y_N$ must be an HN filtration for $f(1)$.

The HN-specialization property $(2)$ now implies that for any valuation ring $R$ with fraction field $K$ and any map $\Spec(R) \to T$ with a lift $\Spec(K) \to Y_N$, there is a unique lift to a map $\Spec(R) \to \Flag(\xi)$. Because $Y_N$ is closed the lift is actually a map $\Spec(R) \to Y_N$. Thus $Y_N \to T$ satisfies the valuative criterion for properness, and \Cref{lem:quasi-compactness} below shows that $Y_N$ is quasi-compact and hence proper over $T$.

At this point we modify the argument if the strong HN-boundedness condition $(4')$ holds: we let $U \subset \Flag(\xi)$ be a quasi-compact open subspace containing an HN filtration for every unstable finite type point of $T$. Then $Y_N \subset U$, because every finite type point of $Y_N$ is a HN filtration. Hence $Y_N \to T$ is a finite type map of locally noetherian algebraic spaces, and it is injective on finite type points by the uniqueness of HN filtrations. The simplified HN-specialization property \ref{princ:S} corresponds to the variant of the valuative criterion for properness of \Cref{lem:valuative_variant} below, which implies that $Y_N \to T$ is a finite morphism.

Under either sets of hypotheses we have now shown that $S_N := \op{im} (Y_N \to T)$ is closed, and an inductive argument shows that $S_i := \bigcup_{j\geq i} \op{im}(Y_j \to T)$ is closed and the induced map
\[
\ev_1 : Y_i \setminus \ev_1^{-1}(S_{i+1}) \to T \setminus S_{i+1}
\]
is proper with image $S_i$. Note also that every point $p \in |T \setminus S_i|$, not necessarily of finite type, has $M^\mu(p) \leq \mu_{i-1} < \mu_i$, so every point with $M^\mu(p) = \mu_i$ has an HN filtration in $Y_i \setminus \ev_1^{-1}(S_{i+1})$. Applying this observation inductively shows that every unstable point of $T$ lies in the closed subscheme $S_1 \hookrightarrow T$, that all points of the quasi-compact locally closed subspace
\[
Y := \bigcup_i Y_i \setminus \ev_1^{-1}(S_{i+1}) \subset \Flag(\xi)
\]
are HN filtrations, and that $Y$ contains an HN filtration of every unstable point of $T$, not just the finite type points.

This shows that the set $|X_{\leq c}|$ is open as claimed, and thus defines an open substack $\X_{\leq c} \subset \X$. Now let $\S^{red}_c \subset \Filt(\X)$ denote the union of irreducible components corresponding to those $s_i$ for which $\mu(s_i) = c$, with reduced structure. \Cref{prop:open_closed_filtrations} implies that we have inclusions of open substacks $\Filt(\X_{\leq c}) = \agr^{-1}(\X_{\leq c}) \subset \ev_1^{-1}(\X_{\leq c}) \subset \Filt(\X)$, and the HN-consistency property $(5)$ implies that these inclusions become equalities after intersecting with $\S^{red}_c$. Thus we have
\[
\ev_1^{-1}(\X_{\leq c}) \cap \S^{red}_c = \Filt(\X_{\leq c}) \cap \S^{red}_c,
\]
which implies that the set of points of $|\Filt(\X)|$ corresponding to HN filtrations of points in $\X$ with $\mu=c$ is actually a subset of $|\Filt(\X_{\leq c})|$ that is closed in the subspace topology. Furthermore, this subset $|\ev_1^{-1}(\X_{\leq c}) \cap \S^{red}_c| \subset |\Filt(\X)|$ is a constructible subset of the locally finite type $B$-stack $\Filt(\X_{\leq c})$, so the open strata condition $(3)$ implies that this set is open as well. Thus there is a union of connected components $\S_c \subset \Filt(\X_{\leq c})$ realizing all HN filtrations of points in $\X$ for which $\mu=c$. The analysis of the previous paragraph shows shows that the map $\ev_1 : \S_c \to \X_{\leq c}$ is proper, and the uniqueness of the HN filtration implies that it is radicial, hence a weak $\Theta$-stratum.

\end{proof}

In the previous proof, we used two different variants of the valuative criterion for properness -- one that doesn't assume separatedness, and one that doesn't assume quasi-compactness.

\begin{lem}\label{lem:valuative_variant}
A finite type map of noetherian algebraic stacks $p : \X \to \Y$ is universally closed if for any discrete valuation ring $R$ essentially of finite type over $\Y$ and any diagram of solid arrows as below, there is an extension of discrete valuation rings $R' \supset R$ with fraction field $K' \supset K$ and a dotted arrow making the diagram commute:
\[
\xymatrix{
\Spec(K') \ar[r] \ar[d] & \Spec(K) \ar[r] \ar[d] & \X \ar[d]^p \\
\Spec(R') \ar@/^5pt/@{-->}[urr] \ar[r] & \Spec(R) \ar[r] & \Y } .
\]
Furthermore, if $p$ is representable by algebraic spaces and satisfies this condition, and for any finite type point $\Spec(k) \to \Y$ the base change $\X_k$ contains a single point, then $p$ is separated, hence finite.
\end{lem}
\begin{proof}
By \cite{HLPreygel}*{Lem.~2.4.6}, it suffices to check that for any finite type morphism $\Spec(A) \to \Y$, the base change $\X_A \to \Spec(A)$ is closed. By Chevalley's theorem, the image of the map $|\X_A| \to |\Spec(A)|$ is constructible, so it suffices to check that it is closed under specialization. Say $x$ lies in the image and specializes to $y$. By base changing first to the reduced closure $\overline{\{x\}}$ and then to the local ring of $y$, it suffices to prove that if $\X_A \to \Spec(A)$ is dominant and $A$ is a local noetherian domain, the closed point of $\Spec(A)$ lies in the image. The local ring $R$ at the exceptional divisor of the blowup of $\Spec(A)$ at the closed point is a DVR essentially of finite type over $A$ such that $\Spec(R) \to \Spec(A)$ is birational and preserves the respective closed points. The fact that $\X_A \to \Spec(A)$ is dominant implies that we may lift the morphism $\Spec(R) \to \Spec(A)$ at the generic point, after a suitable finite extension of DVR's, and the valuative criterion of the lemma guarantees that the special point lifts as well.

Now assume that for any finite type map from a field $\Spec(k) \to \Y$, the base change $|\X_k|$ contains a single point. Let us show that the diagonal $\X \to \X \times_\Y \X$ is universally closed. As this property is smooth local, we may again assume $\Y = \Spec(A)$ for some noetherian ring $A$ and that $\X=X$ is a finite type algebraic space over $\Spec(A)$. By part (3) of \cite{HLPreygel}*{Lem.~2.4.6} it suffices to show that the map $X \times \bA^n \to (X \times_{\Spec(A)} X) \times \bA^n$ is closed for all $n \geq 1$. After replacing $A$ with $A[x_1,\ldots,x_n]$ for various $n$, it thus suffices to show that $X \to X \times_{\Spec(A)} X$ is closed on topological spaces. Both $|X| \to |\Spec(A)|$ and $X\times_{\Spec(A)} X \to \Spec(A)$ are injective on finite type points, hence injective. We also know, because universally closed maps are stable under base change and composition, that both maps are closed. It follows that both $|X|$ and $|X \times_{\Spec(A)} X|$ are closed subspaces of $|\Spec(A)|$ and thus the surjection $|X| \to |X \times_{\Spec(A)}X|$ is a homeomorphism, hence closed. The same argument shows that the double diagonal $X \to X \times_{X \times_{\Spec(A)} X} X$ is a homeomorphism as well. This shows that $\X \to \X \times_\Y \X$ is separated and universally closed, hence $p$ is separated.

\end{proof}

\begin{lem}\label{lem:quasi-compactness}
Let $f : X \to Y$ be a map of algebraic spaces that is locally of finite presentation and satisfies the valuative criterion of properness (i.e., all valuation rings). If $X$ has finitely many irreducible components, then $f$ is quasi-compact, hence proper.
\end{lem}
\begin{proof}
If $X$ has finitely many irreducible components and $X' \to X$ is a smooth map, then $X'$ has finitely many irreducible components. This allows us to reduce to showing that if $Y = \Spec(A)$ is affine then $X$ is quasi-compact if it has finitely many irreducible components. It suffices to assume furthermore that $X$ is reduced and irreducible, and $A$ is an integral domain.

Let $K$ denote the function field of $X$, and let $RZ(K,A)$ denote the Riemann-Zariski space, whose points are valuation rings of $K$ that contain $A$. The valuative criterion for properness for the map $X \to \Spec(A)$ implies that for any valuation ring of $K$ containing $A$, there is a unique map $\Spec(R) \to X$ over $A$. Assigning every valuation ring to the image of its special point defines a map $RZ(K,A) \to |X|$ which is continuous: the preimage of an affine open $\Spec(B) \subset X$ is the subset $RZ(K,B) \subset RZ(K,A)$, which is open for the ``Zariski topology'' on $RZ(K,A)$. Furthermore, the map $RZ(K,A) \to |X|$ is surjective, so quasi-compactness of $|X|$ follows from the fact that $RZ(K,A)$ with its Zariski topology is quasi-compact \cite{matsumura_1987}*{Thm.~10.5}.
\end{proof}


\subsection{Induced stratifications}

Let $\pi : \Y \to \X$ be a map of algebraic stacks that satisfy \ref{hyp2}, and let $\Filt(\pi) : \Filt(\Y) \to \Filt(\X)$ denote the induced map.

\begin{defn} \label{defn:induced_stratum}
We say that a (weak) $\Theta$-stratum $\S \subset \Filt(\X) \to \X$ \emph{induces} a (weak) $\Theta$-stratum in $\Y$ if
\[
\Filt(\pi)^{-1} (\S) \subset \Filt(\Y) \xrightarrow{\ev_1} \Y
\]
is a (weak) $\Theta$-stratum whose set-theoretic image is $\pi^{-1}(\ev_1(\S)) \subset \Y$.

Likewise we say a (weak) $\Theta$-stratification $\{\X_{\leq c}\}_{c\in \Gamma}$ of $\X$ induces a (weak) $\Theta$-stratification of $\Y$ if each $\Theta$-stratum in $\X_{\leq c}$ induces a $\Theta$-stratum in $\Y_{\leq c} := \pi^{-1}(\X_{\leq c})$ for all $c \in \Gamma$.
\end{defn}

Note that in the previous definition if $\Y_{\leq c} := \pi^{-1}(\X_{\leq c})$ is an induced (weak) $\Theta$-stratification, then the stability function $M^\mu : |\Y| \to \Gamma$ defined as $M^\mu(p) = \min \{c \in \Gamma | p \in |\Y_{\leq c}|\}$ is the restriction of $M^\mu : |\X| \to \Gamma$ along the map $|\Y| \to |\X|$.

\begin{lem} \label{lem:induced_closed_open}
Let $\X$ be a stack satisfying \ref{hyp2} with a (weak) $\Theta$-stratification $\{ \X_{\leq c}\}_{c\in \Gamma}$, and let $\Y \to \X$ be either an open union of strata or a finite morphism. Then $\{\X_{\leq c}\}_{c \in \Gamma}$ induces a (weak) $\Theta$-stratification of $\Y$.
\end{lem}
\begin{proof}
The case of an open union of strata is a straightforward consequence of the definition, so we focus on the case of a closed immersion. It suffices to prove the claim for a single (weak) $\Theta$-stratum $\S \to \X$. In this case, the canonical map $\Filt(\Y) \to \Filt(\X) \times_\X \Y$ is a surjective closed immersion by \Cref{prop:affine_map}, and we have a diagram in which each square is Cartesian
\[
\xymatrix{ \Filt(\pi)^{-1}(\S) \ar[r] \ar[d] & \Filt(\Y) \ar@/^10pt/[rr]^{\ev_1} \ar[r] & \Filt(\Y) \times_{\Y} \X \ar[r] \ar[d]^{\Filt(\pi)} & \Y \ar[d]^\pi \\ \S \ar[rr] & & \Filt(\X) \ar[r]^{\ev_1} & \X }
\]
This diagram shows that if $\S \to \X$ is a closed immersion (respectively, a finite radicial morphism), then so is $\Filt(\pi)^{-1}(\S) \to \Y$, and the set theoretic image of the latter map is the preimage $\pi^{-1}(\ev_1(|\S|))$.
\end{proof}

\begin{lem} \label{lem:induced_base_change}
Consider a cartesian diagram of stacks satisfying \ref{hyp2}
\[
\xymatrix{\X' \ar[r]^{\pi'} \ar[d] & \X \ar[d] \\ \Y' \ar[r]^\pi & \Y }
\]
in which $\Y'$ and $\Y$ have quasi-finite inertia. Then any (weak) $\Theta$-stratification in $\X$ induces a (weak) $\Theta$-stratification in $\X'$. Conversely if $\Y' \to \Y$ is faithfully flat and $\S \subset \Filt(\X)$ is a union of connected components for which $\Filt(\pi')^{-1}(\S) \subset \Filt(\X')$ is a (weak) $\Theta$-stratum, then $\S$ is a (weak) $\Theta$-stratum.
\end{lem}

\begin{proof}
\Cref{cor:representable_base_change} implies that $\Filt(\X') \simeq \Filt(\X) \times_\X \X'$. As in the proof of \Cref{lem:induced_closed_open}, this implies that for any $\Theta$-stratum $\S \subset \Filt(\X)$ the substack $\Filt(\pi')^{-1}(\S) \subset \Filt(\X')$ is a $\Theta$-stratum as well. Conversely, if $\X' \to \X$ is faithfully flat then $\Filt(\pi') : \Filt(\X') \to \Filt(\X)$ is as well, because it is the base change of the map $\X' \to \X$. We know that $\Filt(\pi')^{-1}(\S) \simeq \S \times_{\X} \X'$, so faithfully flat descent implies that if $\ev_1 : \Filt(\pi')^{-1}(\S) \simeq \S \times_{\X} \X' \to \X'$ is a closed immersion (respectively finite and radicial), then so is $\S \to \X$.
\end{proof}

\begin{prop} \label{prop:graded_strata}
Let $\X$ be a stack satisfying \ref{hyp2} and consider the map
\[
u : \Grad(\X) \to \X
\]
induced by restriction along $\pt \to \pt/\Gm$. Any $\Theta$-stratification $\{\X_{\leq c}\}_{c \in \Gamma}$ of $\X$ induces a $\Theta$-stratification of $\Grad(\X)$.
\end{prop}

\begin{proof}
\Cref{cor:inertia_preserving} implies that $\Grad(\X_{\leq c}) = u^{-1}(\X_{\leq c})$ for all $c \in \Gamma$, so it suffices to focus on a single $\Theta$-stratum $\S_c \subset \Filt(\X_{\leq c}) \to \X_{\leq c}$, and we will drop the subscript $c$ from the notation. We first apply the functor $\Grad(-)$ to the sequence $\S \subset \Filt(\X) \xrightarrow{\ev_1} \X$ to obtain a commutative diagram
\[
\xymatrix{ \Grad(\S) \ar[r] \ar[d]^u & \Grad(\Filt(\X)) \ar[d]^{u_{\Filt(\X)}}\ar[r]^{\Grad(\ev_1)} & \Grad(\X) \ar[d]^u \\ \S \ar[r] & \Filt(\X) \ar[r] & \X }.
\]
Because $\S \subset \Filt(\X)$ is a union of connected components and $\S \to \X$ is a locally closed immersion, \Cref{cor:inertia_preserving} implies that both the outer square and the left square in this diagram are cartesian.

As an immediate consequence of their description as mapping stacks we can identify
\begin{equation} \label{eqn:flip_mapping_stacks}
\Filt(\Grad(\X)) \simeq \iMap(\Theta \times (\pt / \Gm), \X) \simeq \Grad(\Filt(\X)),
\end{equation}
and under this isomorphism we have a canonical identification between the map $u_{\Filt(\X)}$ in the commutative diagram above and the application of the functor $\Filt(-)$ to the map $u : \Grad(\X) \to \X$, in other words $u_{\Filt(\X)} \simeq \Filt(u_\X)$. Therefore, \eqref{eqn:flip_mapping_stacks} identifies the union of connected components $\Filt(u)^{-1}(\S)$ with $\Grad(\S)$, and the map $\ev_1 : \Filt(u)^{-1}(\S) \to \Grad(\X)$ is the base change of $\ev_1 : \S \to \X$ along $u : \Grad(\X) \to \X$, which establishes that the $\Theta$-stratum $\S \to \X$ induces a $\Theta$-stratum in $\Grad(\X)$.
\end{proof}

\begin{rem} \label{rem:equivariant_HN}
One concrete consequence is that if $g \in \Grad(\X)$ is a graded object whose underlying point $u(g) \in \X$ is unstable, then the HN filtration of $u(g)$ lifts uniquely to an HN filtration of $g$ in $\Grad(\X)$. This can be shown directly for any graded object for which $u(g)$ admits a unique HN filtration, even when $\mu$ does not define a $\Theta$-stratum. The general statement underlying this fact is that if $x \in \X(\bar{k})$ is a geometric point and $\S \subset \Filt(\X)$ a connected component such that $\ev_1^{-1}(x) \cap \S$ has a single $\bar{k}$ point $f$, then any homomorphism from a reduced $k$-group $G \to \Aut_\X(x)$ lifts uniquely to a homomorphism $G \to \Aut_{\Filt(\X)}(f)$.
\end{rem}

Let $\{\X_{\leq c}\}_{c \in \Gamma}$ be a $\Theta$-stratification of $\X$, and consider the center $\Z_c^{\rm{ss}}$ of the stratum $\S_c \hookrightarrow \X_{\leq c}$ (\Cref{defn:center}). Assume for simplicity that the stratification is induced by a locally constant $\mu$ as in \Cref{simp:locally_constant}, which will always be the case for the stratifications constructed in the remainder of this paper. For every $c \in \Gamma$, let $\Z_c \subset \Grad(\X)$ be a the union of all connected components meeting the open substack $\Z_c^{\rm{ss}} \subset \Grad(\X_{\leq c}) \subset \Grad(\X)$. The main consequence of \Cref{prop:graded_strata} is

\begin{cor} \label{cor:semistable_centers}
The center $\Z_c^{\rm{ss}}$ is the semistable locus of a $\Theta$-stratification $\{(\Z_c)_{\leq c'}\}_{c' \geq c}$ of $\Z_c$ induced by the forgetful map $u : \Z_c \to \X$, where $c'=c$ is the minimal index. In particular $(\Z_c)_{\leq c'}$ is the preimage of $\X_{\leq c'}$ under $u$, and $(\Z_c)_{\leq c'} = \emptyset$ for $c' < c$.
\end{cor}

\begin{proof}
The claim that $u : \Z_c \to \X$ induces a $\Theta$-stratification on $\Z_c$ follows from \Cref{prop:graded_strata} applied to $\Grad(\X) \to \X$ and \Cref{lem:induced_closed_open} applied to $\Z_c \to \Grad(\X)$.

Let $\bar{\S}_c \subset \Filt(\X)$ be the union of all connected components that meet the open substack $\S_c \subset \Filt(\X_{\leq c}) \subset \Filt(\X)$. Because we have assumed $\mu$ is locally constant, $\mu$ is defined and takes the value $c$ on all irreducible components of $\bar{\S}_c$. It follows that $M^\mu(p) \geq c$ for every point in $\ev_1(|\bar{\S}_c|) \subset \X$ and that $M^\mu(p) = c$ if and only if $p \in \ev_1(|\S_c|)$. In particular $\S_c$ is the preimage of $\X_{\leq c} \subset \X$ under the restriction of $\ev_1$ to $\bar{\S}_c \subset \Filt(\X)$.

Now making use of \Cref{lem:retract}, we see that $\Z_c = \sigma^{-1}(\bar{\S}_c) \subset \Grad(\X)$, and combining this with the fact that $\ev_1 \circ \sigma \simeq u$ shows that $\Z_c^{\rm{ss}} \subset \Z_c$ is the preimage of $\X_{\leq c} \subset \X$ under $u$. In order to identify $\Z_c^{\rm{ss}}$ with the ``semistable'' locus for the induced $\Theta$-stratification of $\Z_c$ in the sense of \Cref{defn:theta_stratification}, we must verify the preimage of $X_{\leq c'}$ is empty for any $c' < c$. This follows from the fact that $M^\mu(u(g)) = M^\mu(\ev_1(\sigma(g))) \geq c$ for any $g \in |\Z_c|$
\end{proof}


\section{Combinatorial structures in moduli theory} \label{sect:structures}

In this section, given a stack $\X$ and a $k$-point $p \in \X(k)$, we introduce the \emph{degeneration space} $\iDeg(\X,p)$ (\Cref{defn:degeneration_space}), which encodes all possible filtrations of $p$ up to the action of $\bN^\times$. Under fairly general hypotheses on $\X$ it is a compactly generated Hausdorff space by \Cref{prop:intersecting_simplices} and \Cref{prop:fans_cg}. $\iDeg(\X,p)$ is defined as the ``projective'' geometric realization of a formal fan $\Deg(\X,p)_\bullet$, a type of combinatorial object that we introduce in \Cref{defn:formal_fan}, which behaves somewhat like a semisimplicial set.

\begin{ex}
For the stack $\X = \bA^2 / \Gm^2$ and the point $p = (1,1)$, any pair of nonnegative integers, $(a,b)$, defines a group homomorphism $\Gm \to \Gm^2$ that extends to a map of quotient stacks $f:\Theta \to \bA^2/\Gm^2$ along with an isomorphism $f(1)\simeq (1,1)$, and the morphism $f$ corresponding to the pair $(na,nb)$ is the precomposition with the $n$-fold ramified cover $\Theta \to \Theta$. Thus non-degenerate filtrations of $(1,1) \in \bA^2/\Gm^2$ (up to ramified coverings) correspond to rational rays in the cone $(\bR_{\geq 0})^2$. We identify this, in turn, with the set of rational points in the unit interval $(\bR_{\geq 0}^2 \setminus 0) / \bR_{>0}^\times$, and $\iDeg(\X,p)$ is this interval. Note that $\iDeg(\X,p)$ parameterizes filtrations of $p$, but a path in $\iDeg(\X,p)$ does not correspond to an algebraic family of filtrations. In this case $\Flag(p)$ for $p = (1,1) \in \bA^2 / \Gm^2$ is an infinite disjoint union of points, one for every pair of nonnegative integers $(a,b)$.
\end{ex}

More generally, when $X$ is a toric variety and $\X = X/T$, $\Deg(\X,p)_\bullet$ encodes the classical fan in the cocharacter space of $T$ associated to $X$ (\Cref{lem:degeneration_fan_torus}), and when $\X = \pt /G$ for a split semisimple algebraic group $G$, $\iDeg(\X,p)$ is canonically homeomorphic to the spherical building of $G$ (\Cref{prop:spherical_building}). In addition, a proper representable morphism of stacks induces a homeomorphism of degeneration spaces, and an affine morphism of stacks induces a locally convex closed embedding of degeneration spaces (\Cref{prop:representable_maps}).

Our main technical result, \Cref{thm:perturbation}, identifies small perturbations of a given filtration $f : \Theta_k \to \X$ with small perturbations of the canonical filtration of the associated graded point $\agr(f)$. More precisely, $\agr(f)$ receives a canonical filtration $\canon$ as a point in $\Grad(\X)$, and the theorem identifies a neighborhood of $\canon \in \iDeg(\Grad(\X),\agr(f))$ with a neighborhood of $f \in \iDeg(\X,\ev_1(f))$. This concept plays an important role in \Cref{sect:construction}.

The second pair of objects we introduce is the component fan $\Comp(\X)_\bullet$ and its projective geometric realization $\iComp(\X)$, the component space (\Cref{defn:reduced_degeneration_space}). These parameterize all possible filtrations in $\X$ up to deformation equivalence, i.e., connected components of $\Filt(\X)$. The key fact is that $\iComp(\X)$ is compact when $\X$ is quasi-compact (\Cref{cor:bounded_comp_1}). We will also show how to construct continuous functions on $\iComp(\X)$ from cohomology classes in $H^{2n}(\X;\bR)$, which we will use to construct numerical invariants of filtrations in \Cref{sect:construction}. Finally, we use the component space to give a criterion in \Cref{prop:quasi_compact_flags} under which $\X$ has quasi-compact flag spaces.

\subsection{Formal fans and their geometric realizations}

The basic combinatorial objects that we study, formal fans, are analogous to semisimplicial sets.

\begin{defn} \label{defn:formal_fan}
We define a category of \emph{integral simplicial cones} $\Cones$ to have
\begin{itemize}
\item objects: positive integers $[n]$ with $n > 0$,
\item morphisms: a morphism $\phi : [k] \to [n]$ is an injective group homomorphism $\bZ^k \to \bZ^n$ that maps the standard basis of $\bZ^k$ to the cone spanned by the standard basis of $\bZ^n$.
\end{itemize}
We define the category of \emph{formal fans}
$$\op{Fan} := \op{Fun}(\Cones^{op}, \op{Set})$$
For $F \in \op{Fan}$ we use the abbreviated notation $F_n = F([n])$ and refer to this as the set of \emph{cones}. We refer to the elements of $F_1$ as \emph{rays}.
\end{defn}

Given a \gls{formal_fan} $F_\bullet$, in analogy with the theory of semisimplicial sets we define the \emph{face maps} $d_i : F_n \to F_{n-1}$ for $i = 1,\ldots,n$ to be the maps induced by the homomorphism $\bZ^{n-1} \to \bZ^n$ mapping $(x_1,\ldots,x_{n-1}) \mapsto (x_1,\ldots,x_{i-1},0,x_i,\ldots,x_n)$. We also introduce the \emph{vertex maps} $v_i : F_n \to F_1$ for $i=1,\ldots,n$ induced by the group homomorphism $\bZ \to \bZ^n$ that is the inclusion of the $i^{th}$ factor.

For any $F \in \op{Fan}$, we can define two notions of geometric realization. First form the comma category $(\Cones | F)$ whose objects are elements $\xi \in F_n$ and morphisms $\xi_1 \to \xi_2$ are given by morphisms $\phi : [n_1] \to [n_2]$ with $\phi^\ast \xi_2 = \xi_1$. There is a canonical functor $(\Cones|F) \to \op{Top}$ assigning $\xi \in F_n$ to the cone $(\bR_{\geq 0})^n$ spanned by the standard basis of $\bR^n$. Using this we can define the \emph{geometric realization} of $F$
$$|F| := \colim \limits_{(\Cones | F)} (\bR_{\geq 0})^n = \colim \limits_{[n] \in \Cones} F_n \times (\bR_{\geq 0})^n$$
This is entirely analogous to the geometric realization functor for semisimplicial sets.

The map $\bR^k_{\geq 0} \to \bR^n_{\geq 0}$ corresponding to a map $[k] \to [n]$ in $\Cones$ takes integral points $\bZ^k_{\geq 0} \subset \bR^k_{\geq 0}$ to $\bZ^n_{\geq 0}$. For any fan $F_\bullet$, we define the set of \emph{integral points} of $|F_\bullet|$ to be the image of the defining maps $\bZ^n_{\geq 0} \to |F_\bullet|$. Alternatively, the set of non-zero integral points of $|F_\bullet|$ is in bijection with the set of rays, where $\sigma \in F_1$ is identified with the image of $1 \in \bR_{\geq 0}$ under the map $\bR_{\geq 0} \to |F_\bullet|$ corresponding to $\sigma$.

Given a map $[k] \to [n]$ in $\Cones$, the corresponding linear map $\phi : \bR^k \to \bR^n$ is injective and equivariant with respect to scaling by $\bR^\times_{>0}$. $\phi$ descends to a an injective map $\simp{}^{k-1} \to \simp{}^{n-1}$, where $\simp{}^{n-1} = \left( (\bR_{\geq 0})^n - \{0\} \right) / (\bR_{>0})^\times$ is the standard $(n-1)$-simplex realized as the space of rays in $(\bR_{\geq 0})^n$. Thus for any $F \in \op{Fan}$ we have a functor $(\Cones | F) \to \op{Top}$ assigning $\xi \in F_n$ to $\simp{n-1}$. We define the \emph{projective realization} of $F$ to be
$$\bP(F) := \colim \limits_{\xi \in (\Cones|F)} \simp{\xi} = \colim \limits_{[n] \in \Cones} F_n \times \simp{}^{n-1}$$
where the notation $\simp{\xi}$ denotes a copy of the standard $n-1$ simplex formally indexed by $\xi \in F_n$.\footnote{When we wish to emphasize the dimension of the simplex, we will denote it $\simp{\xi}^{n-1}$.} We refer to the map $\xi : \simp{\xi} \to \bP(F_\bullet)$ as a \emph{rational simplex} of $\bP(F_\bullet)$ and when $n=1$ we call this a rational point of $\bP(F_\bullet)$.

\begin{ex}
For a representable fan, $h_{[n]}(\bullet) = \op{Hom}_{\Cones}(\bullet, [n])$, the category $(\Cones|h_{[n]})$ has a single terminal object, the identity map on $[n]$. Hence $|h_{[n]}| \simeq (\bR_{\geq 0})^n$ and $\bP(h_{[n]}) \simeq \simp{}^{n-1}$. For any other fan $F_\bullet$ and $\xi \in F_n$, the corresponding map on geometric realizations $\bP(h_{[n]}) \to \bP(F_\bullet)$ is the rational simplex $\xi : \simp{\xi} \to \bP(F_\bullet)$.
\end{ex}

\begin{ex}
One can define a fan $F_\bullet$ whose cones are the same as the cones of $h_{[2]}$, except that $F_1$ is the quotient of $h_{[2]}([1])$ by the equivalence relation that identifies the two boundary rays (corresponding the homomorphisms $\phi(e_0)=e_1$ and $\phi(e_0)=e_0$). Then $|F_\bullet|$ has the shape of a rolled paper cone, $\bP(F_\bullet)$ is a circle, and the rational simplex $\bP(h_{[2]}) \to \bP(F_\bullet)$ is the non-injective map from the closed unit interval to the circle obtained by identifying its endpoints.
\end{ex}

\begin{rem}
As an alternative description of $\bP(F_\bullet)$, we observe that injective maps $(\bR_{\geq 0})^k \to (\bR_{\geq 0})^n$ arising in the construction of $|F|$ are equivariant with respect to the action of $\bR_{>0}^\times$ on $(\bR_{\geq 0})^n$ by scalar multiplication, and thus $|F|$ has a canonical continuous $\bR_{>0}^\times$-action. Commuting colimits shows that
\[
\bP(F) \simeq (|F| - \{\ast\}) / \bR_{>0}^\times
\]
as topological spaces, where $\ast \in |F|$ is the cone point corresponding to the origin in each copy of $(\bR_{\geq 0})^n$.
\end{rem}

The terminology of formal fans and cones is motivated by the following construction, which establishes a relationship between formal fans and the classical notion of a fan in a vector space.

\begin{const} A subset $K \subset \bR^N$ that is invariant under multiplication by $\bR_{\geq 0}^\times$ is called a \emph{cone in $\bR^N$}. Given a set of cones $K_\alpha \subset \bR^N$, we define
\begin{equation} \label{eqn:tautological_R_set}
\ray[n] (\{K_\alpha\} ) := \left\{ \begin{array}{c} \text{injective homomorphisms } \phi : \bZ^n \to \bZ^N  \text{ s.t.} \\ \exists \alpha \text{ s.t. } \phi(e_i) \subset K_\alpha, \forall i \end{array} \right\}
\end{equation}
These sets naturally form an object $\ray (\{K_\alpha\}) \in \op{Fan}$.
\end{const}

\begin{rem}
We use the phrase \emph{classical fan} to denote a collection of rational polyhedral cones in $\bR^N$ such that a face of any cone is also in the collection, and the intersection of two cones is a face of each. If $K_\alpha \subset \bR^N$ are the cones of a classical fan $\Sigma$, it is possible to reconstruct the partially ordered set of cones in $\Sigma$ partially ordered by inclusion from the data of $\ray (\{K_\alpha\})$. Likewise one can recover the original $K_\alpha \subset \bR^n$ from the data of the inclusion $\ray(\{K_\alpha\}) \subset \ray(\{\bR^n\})$.
\end{rem}

\begin{lem} \label{lem:reconstruct_fan}
Let $K_\alpha \subset \bR^N$ be a finite collection of cones, each of which is a finite union of rational polyhedral cones. Then the canonical map $|\ray (\{K_\alpha\})| \to \bigcup K_\alpha$ is a homeomorphism. Furthermore, $\bP ( \ray (\{K_\alpha\}) ) \simeq S^{N-1} \cap \bigcup_{\alpha} K_\alpha$ via the evident quotient map $\bR^N -\{0\} \to S^{N-1}$.
\end{lem}
\begin{proof}
By considering all intersections of rational polyhedral cones and their faces appearing in the description of some $K_\alpha$, we may assume there is a classical fan $\Sigma = \{\sigma_i\}$ in $\bR^N$ such that each $\sigma_i$ is contained in some $K_\alpha$ and each $K_\alpha$ is the union of some collection of $\sigma_i$. By further refinement we may assume that $\Sigma$ is simplicial, and that the ray generators of each $\sigma_i$ form a basis for the lattice generated by $\sigma_i \cap \bZ^N$. Consider the formal fans $F_\bullet = \ray (\{K_\alpha\})$ and $F^\prime_\bullet = \ray (\{\sigma_i\})$. By hypothesis $F^\prime_\bullet$ is a subfunctor of $F_\bullet$. Hence we have a functor of comma categories $(\Cones | F^\prime) \to (\Cones | F)$ and thus a map of topological spaces $|F^\prime| \to |F|$ that commutes with the map to $\bR^N$.

We can reduce to the case when $\{K_\alpha\} = \{\sigma_i\}$. Indeed the map $|F^\prime| \to |F|$ is surjective because any point in $K_\alpha$ lies in some $\sigma_i$. If the composition $|F^\prime| \to |F| \to \bigcup_\alpha K_\alpha = \bigcup \sigma_i$ were a homeomorphism it would follow that $|F^\prime| \to |F|$ was injective as well, and one could use the inverse of $|F^\prime| \to \bigcup_\alpha K_\alpha$ to construct and inverse for $|F| \to \bigcup_\alpha K_\alpha$.

For a single simplicial cone $\sigma \subset \bR^N$ of dimension $n$ whose ray generators form a basis for the lattice generated by $\sigma \cap \bZ^N$, we have $\ray (\sigma) \subset \ray (\bR^N)$ is isomorphic to $h_{[n]}(\bullet)$. The terminal object of $(\Cones | h_{[n]})$ corresponds to the linear map $\bR^n \to \bR^N$ mapping the standard basis vectors to the ray generators of $\sigma$, and thus $| \ray (\sigma)| \to \sigma$ is a homeomorphism.

Let $\sigma^{\max}_1,\ldots,\sigma^{\max}_r \in \Sigma$ be the cones that are maximal with respect to inclusion and let $n_i$ be the dimension of each. Let $\sigma_{ij}^\prime := \sigma^{\max}_i \cap \sigma^{\max}_j \in \Sigma$, and let $n_{ij}$ denote its dimension. Then by construction
$$F := \ray (\{\sigma_i\}) \simeq \op{coeq} \left( \bigsqcup_{i,j} h_{n_{ij}} \rightrightarrows \bigsqcup_i h_{n_i} \right)$$
as functors $\Cones^{op} \to \op{Set}$. Our geometric realization functor commutes with colimits, so it follows that
$$|F| = \op{coeq} \left( \bigsqcup_{i,j} \sigma^\prime_{ij} \rightrightarrows \bigsqcup_i \sigma^{\max}_{ij} \right)$$
Which is homeomorphic to $\bigcup \sigma_i$ under the natural map $|F| \to \bR^N$. The final claim follows from the fact that $\bP(F) \simeq (|F| - \{0\}) / \bR_{>0}^\times$.
\end{proof}

\begin{ex} \label{ex:R_set_pathology}
Objects of $\op{Fan}$ describe a wider variety of structures than classical fans. For instance if $K_1$ and $K_2$ are two simplicial cones that intersect but do not meet along a common face, then $\ray (K_1,K_2)$ will not be equivalent to $\ray (\{\sigma_i\})$ for any classical fan $\Sigma = \{\sigma_i\}$.
\end{ex}

\begin{ex} \label{ex:fan_pathology}
The category $\op{Fan}$ is broad enough to include some pathological examples. For instance, if $K \subset \bR^3$ is the cone over a circle that is not contained in a linear subspace, then $|\ray (K)|$ consists of the rational rays of $K$ equipped with the discrete topology and is not homeomorphic to $K$.

Another example: if $K \subset \bR^2$ is a convex cone generated by two irrational rays, then $|\ray (K)|$ is the interior of that cone along with the origin. There are also examples of fans whose geometric realizations are not Hausdorff, such as multiple copies of the standard cone in $\bR^2$ glued to each other along the set of rational rays.
\end{ex}

\subsubsection{Some useful lemmas}

\begin{lem} \label{lem:equivalence_rel}
Let $F$ be a fan and let $x \in \bP(F)$ be a point lying in the image of two rational simplices $\simp{\sigma_i} \to \bP(F)$. Then there is a $\xi \in (\Cones|F)$ with maps $\xi \to \sigma_i$ such that $x$ lies in the image of $\simp{\xi} \to \bP(F)$.
\end{lem}

\begin{proof}
By the definition of $\bP(F)$ as a colimit, it is the quotient of the set $\bigsqcup_{\sigma \in (\cC|F)} \simp{\sigma}$ by the equivalence relation generated by $x \sim \phi(x)$ for the maps $\phi : \simp{\sigma_1} \to \simp{\sigma_2}$ induced by morphisms in $(\Cones|F)$. The lemma follows from an alternative description of this equivalence relation: points $x_1 \in \simp{\sigma_1}$ and $x_2 \in \simp{\sigma_2}$ are equivalent if there exists a diagram in $(\Cones|F)$ of the form $\sigma_1 \leftarrow \xi \to \sigma_2$ and a point $y \in \simp{\xi}$ mapping to $x_1$ and $x_2$. This relation is clearly symmetric and reflexive, so we must show transitivity.

Consider a diagram of the form $\sigma_1 \leftarrow \xi_1 \to \sigma_2 \leftarrow \xi_2 \to \sigma_3$ and points $y_i \in \simp{\xi_i}$ and $x_i \in \simp{\sigma_i}$ that are identified by these maps. The corresponding maps $\simp{\xi_1} \hookrightarrow \simp{\sigma_2}$ and $\simp{\xi_2} \hookrightarrow \simp{\sigma_2}$ are embeddings of simplices with rational vertices, and hence their intersection in $\simp{\sigma_2}$ is a polyhedron with rational vertices, which contains the point $y_1 = y_2 = x_2$. It follows that there is some simplex with rational vertices in $\simp{\sigma_2}$ containing $x_2$, which defines a $\xi' \in (\Cones|F)$ mapping to $\xi_1$ and $\xi_2$ and a point $y' \in \simp{\xi'}$ mapping to $y_1$ and $y_2$. Thus we have a diagram of the form $\sigma_1 \leftarrow \xi' \to \sigma_3$ exhibiting the equivalence $x_1 \sim x_3$.
\end{proof}

We can formulate the previous lemma a bit more generally.

\begin{lem} \label{lem:fiber_prod}
For any maps of fans $F,F' \to G$, the canonical map $\bP(F \times_G F') \to \bP(F) \times_{\bP(G)} \bP(F')$ is surjective.
\end{lem}

\begin{proof}
Let $S = \bP(F) \times_{\bP(G)} \bP(F')$. Any point in $S$ can be represented by a pair of rational simplices $\xi \in F_n$ and $\xi' \in F'_{n'}$ and two points $x \in \simp{\xi}$ and $x' \in \simp{\xi'}$ that map to the same point in $\bP(G)$. \Cref{lem:equivalence_rel} implies that we can find subcones of $\xi$ and $\xi'$ that map to the same cone in $G$ and whose corresponding rational simplex contains a point that maps to both $x$ and $x'$. This defines a point in $\bP(F \times_G F')$ mapping to our original point in $S$.
\end{proof}

\begin{lem} \label{lem:fans_surjectivity_base_change}
If $F \to G$ is a map of fans such that $\bP(F) \to \bP(G)$ is surjective, then for any fan $F'$ and map $F' \to G$, the map $\bP(F' \times_G F) \to \bP(F')$ is surjective.
\end{lem}
\begin{proof}This is an immediate corollary of the previous lemma and the fact that surjective maps of topological spaces are stable under base change.
\end{proof}

\begin{lem} \label{lem:fans_injective_surjective}
If $F_\bullet \to G_\bullet$ is an injective (respectively surjective) map of fans, then so is $\bP(F_\bullet) \to \bP(G_\bullet)$.
\end{lem}
\begin{proof}
\Cref{lem:equivalence_rel} implies that if two points in $x_0,x_1 \in \bP(F_\bullet)$, represented by cones $\sigma_1,\sigma_2$ in $F_\bullet$ with a choice of real ray in each, map to the same point in $\bP(G_\bullet)$, then there is a cone and real ray in $G_\bullet$ that is a subcone of the image of both $\sigma_1$ and $\sigma_2$. Because the map $F_\bullet \to G_\bullet$ is injective and its image is closed under taking subcones, this new cone must lie in $F_\bullet$ as well. The corresponding claim for surjectivity is immediate from the definition of $\bP(G_\bullet)$.
\end{proof}

\begin{cor} \label{cor:fans_image}
Let $f : F_\bullet \to G_\bullet$ be a map of fans and let $\bP(f) : \bP(F_\bullet) \to \bP(G_\bullet)$ be the induced map. Then
\begin{enumerate}
\item $\op{im}(\bP(f)) = \bP(\op{im}(f))$ where $\op{im}$ denotes the image of a map of fans or topological spaces, and
\item if $G'_\bullet \to G_\bullet$ is a map and $f' : F_\bullet \times_{G_\bullet} G'_\bullet \to G'_\bullet$ is the base change, then $\op{im}(\bP(f'))$ is the preimage of $\op{im}(\bP(f)) \subset \bP(G_\bullet)$ in $\bP(G'_\bullet)$.
\end{enumerate}
\end{cor}

\begin{proof}
The first claim is an immediate consequence of \Cref{lem:fans_injective_surjective} and the definition of the image in both categories as the unique factorization of a map as a surjection followed by an injection. The second then follows from this and \Cref{lem:fans_surjectivity_base_change}.
\end{proof}

Note that as a consequence of part (2) of this corollary, for any two subfans $F_\bullet, F'_\bullet \subset G_\bullet$, we have $\bP(F_\bullet) \cap \bP(F'_\bullet) = \bP(F_\bullet \cap F'_\bullet)$ as subsets of $\bP(G_\bullet)$.

\subsubsection{Boundedness of fans}

We shall sometimes restrict to a class of fans whose projective realizations exhibit fewer pathologies.

\begin{defn} \label{def:bounded_fan}
We say that a fan is \emph{bounded} if $\bP(F)$ is covered by finitely many rational simplices, and a map $F \to G$ is \emph{bounded} if for any cone $h_{[n]} \to G$, the fiber product $F \times_G h_{[n]}$ is bounded. We say that a fan $F$ is \emph{quasi-separated} if for any pair of cones $h_{[n_0]}, h_{[n_1]} \to F$, the fan $h_{[n_0]} \times_F h_{[n_1]}$ is bounded.
\end{defn}

\begin{rem}
Note that if $K \subset \bR^n$ is a rational polyhedral cone that is not simplicial, then $R_\bullet(\{K\})$ does not admit a surjection from a finite union of representable fans $h_{[n]}$, but $R_\bullet(\{K\})$ is still bounded in the above sense. On the other hand, a surjective map of fans $F_\bullet \to G_\bullet$ induces a surjection $\bP(F_\bullet) \to \bP(G_\bullet)$. Therefore $F_\bullet$ is bounded if it admits a surjection $\bigsqcup_{i=1}^N h_{[n_i]} \to F_\bullet$.
\end{rem}

\begin{lem}
A fan $G$ is quasi-separated if and only if for any bounded fan $F$, the map $F \to G$ is bounded.
\end{lem}
\begin{proof}
The ``if'' direction is immediate from the definitions, taking $F = h_{[n]}$ in the statement of the lemma. For the ``only if'' direction, let $h_{[m]} \to G$ be a cone, and let $\bigsqcup_i h_{[n_i]} \to F$ be a map that is surjective after applying $\bP(-)$. Then $\bigsqcup_i \bP(h_{[n_i]} \times_G h_{[m]}) \to \bP(F \times_G h_{[m]})$ is surjective by \Cref{lem:fans_surjectivity_base_change}, so the latter is bounded if the former is.
\end{proof}

\begin{rem}
The notion of boundedness is analogous to quasi-compactness in algebraic geometry, and the notion of a quasi-separated fan is analogous to that of a quasi-separated algebraic space, i.e., an algebraic space $X$ such that for any pair of maps $\Spec(A_i) \to X$ the fiber product $\Spec(A_1) \times_X \Spec(A_2)$ is a quasi-compact and quasi-separated algebraic space.

An algebraic space is quasi-separated if and only if the diagonal $X \to X \times X$ is quasi-compact. For fans, however, the condition that $F \to F \times F$ is bounded is weaker than $F$ being quasi-separated. The equivalence fails because $h_{[n]} \times h_{[m]}$ is not bounded for $m\neq n$. This is related to the failure of $\bP(-)$ to commute with products. Presumably this deficiency can be addressed by reformulating the theory using a category $\Cones$ that allows degenerate maps between cones, just as the theory of simplicial sets addresses similar deficiencies in the theory of semisimplicial sets.
\end{rem}

\begin{lem} \label{lem:fans_closed_map}
Let $\psi : F \to G$ be a bounded map of fans. Then $\psi : \bP(F) \to \bP(G)$ is a closed map whose fibers are finite sets. The image of this map restricted to any rational simplex of $\bP(G)$ is a finite union of closed rational sub-simplices.
\end{lem}

\begin{proof}
Given a rational simplex $\xi : \simp{\xi}^n \to \bP(G)$, one can find a surjection from a finite union of rational simplices
\[
\bigsqcup_i \simp{\sigma_i} \twoheadrightarrow \bP(F \times_G h_{[n+1]}) \twoheadrightarrow \bP(F) \times_{\bP(G)} \simp{\xi}^n,
\]
by \Cref{lem:fiber_prod}. Let $p_1$ and $p_2$ denote the projection from $\bigsqcup_i \simp{\sigma_i}$ to $\bP(F)$ and $\simp{\xi}^n$ respectively. Then for any closed set $S \subset \bP(F)$,
\[
\xi^{-1}(\psi(S)) = p_2(p_1^{-1}(S)).
\]
All three claims follow from the fact that $p_2$ is a closed map with finite fibers.
\end{proof}

\begin{lem} \label{lem:fans_bounded_surjectity}
If $\psi : F_\bullet \to G_\bullet$ is a bounded map of fans, then $\bP(\psi) : \bP(F_\bullet) \to \bP(G_\bullet)$ is injective (resp. surjective) if and only if it is injective (resp. surjective) on rational points.
\end{lem}
\begin{proof}
Assume that $\bP(\psi)$ is surjective on rational points. \Cref{lem:fans_closed_map} implies that the image of $\bP(\psi)$ restricted to each rational simplex $\simp{\xi} \to \bP(G_\bullet)$ is a finite union of rational subsimplices. If this finite union contains all rational points of $\simp{\xi}$, then it must cover $\simp{\xi}$, so $\bP(\psi)$ is surjective.
\end{proof}

\begin{cor} \label{cor:fans_closed_embedding}
If $F \subset G$ is a bounded injective map of fans, then $\bP(F) \to \bP(G)$ is a closed embedding, and its preimage under any rational simplex $\simp{\xi} \to \bP(G)$ is a finite union of rational sub-simplices.
\end{cor}
\begin{proof}
Combine \Cref{lem:fans_closed_map} and \Cref{lem:fans_injective_surjective}.
\end{proof}

\begin{prop} \label{prop:fans_cg}
Let $F$ be a quasi-separated fan. Then any rational simplex $\simp{\xi} \to \bP(F)$ is a closed map, and $\bP(F)$ is a Hausdorff space whose topology is compactly generated by the images of rational simplices.
\end{prop}

\begin{proof}
The fact that an rational simplex is a closed map is \Cref{lem:fans_closed_map}. Once we know this, we see from the definition of the colimit topology that a subset of $\bP(F)$ is closed if and only if its intersection with the compact subset $\xi(\simp{\xi})$ is closed for all rational simplices, so the topology on $\bP(F)$ is compactly generated.

We show that $\bP(F)$ is Hausdorff by showing that the diagonal $\bP(F) \hookrightarrow \bP(F) \times \bP(F)$ is a closed subset. Writing $F = \bigcup_\alpha F_\alpha$ as a union of bounded subfans, we have that $\bP(F) = \bigcup \bP(F_\alpha)$, because $\bP(-)$ commutes with colimits. This colimit is filtered, so $\bP(F) \times \bP(F) = \bigcup_{\alpha, \beta} \bP(F_\alpha) \times \bP(F_\beta)$ with the colimit topology. It therefore suffices to show that for each $\alpha$ and $\beta$ the preimage of the diagonal under the map
\[
\bP(F_\alpha) \times \bP(F_\beta) \to \bP(F) \times \bP(F)
\]
is a closed subset.

The fans $F_\alpha$ are bounded by hypotheses and quasi-separated because they are sub-fans of a quasi-separated fan, so $\bP(F_\alpha)$ admit closed surjections from a finite disjoint union of rational simplices. This allows one to reduce the claim to showing that for any two rational simplices $\xi_i : \simp{\xi_i} \to \bP(F)$, the preimage of the diagonal under the map $\simp{\xi_1} \times \simp{\xi_2} \to \bP(F) \times \bP(F)$ is a closed subset. We can identify this preimage with $\simp{\xi_1} \times_{\bP(F)} \simp{\xi_2}$, which by \Cref{lem:fiber_prod} is covered by finitely many rational simplices because $h_{[n_1]} \times_F h_{[n_2]}$ is bounded. Hence the preimage of the diagonal in $\simp{\xi_1} \times \simp{\xi_2}$ is closed.
\end{proof}

In fact we can be even more precise:

\begin{lem} \label{lem:fans_quasi_sep_relation}
Let $G$ be a quasi-separated fan. Then for any finite collection of cones $\xi_i \in G_{n_i}$, the image of $\bigsqcup \simp{\xi_i} \to \bP(G)$ is the quotient of $\bigsqcup \simp{\xi_i}$ by a closed equivalence relation
\[
R \subset (\bigsqcup_i \simp{\xi_i}) \times (\bigsqcup_i \simp{\xi_i}),
\]
which is a union of finitely many closed subsets $\simp{}^k \hookrightarrow (\bigsqcup_i \simp{\xi_i}) \times (\bigsqcup_i \simp{\xi_i})$ obtained from a pair of rational sub-simplices $\simp{}^k \hookrightarrow \simp{\xi_i}$ and $\simp{}^k \hookrightarrow \simp{\xi_j}$.
\end{lem}
\begin{proof}
Writing $\bigsqcup_i \simp{\xi_i} = \bP(\bigsqcup h_{[n_i]})$, $R = \bP(\bigsqcup h_{[n_i]}) \times_{\bP(G)} \bP(\bigsqcup h_{[n_i]})$ is covered by the image of
$$\bigsqcup_{i,j} \bP(h_{[n_i]} \times_F h_{[n_j]}) \to \bP(\bigsqcup h_{[n_i]}) \times \bP(\bigsqcup h_{[n_i]}),$$
by \Cref{lem:fiber_prod}, which in turn is covered by finitely many rational simplices by hypotheses. For each rational simplex $\simp{\sigma}^k \to \bP(h_{[n_i]} \times_F h_{[n_j]})$, the projection onto each factor $\bP(\bigsqcup h_{[n_i]})$ is a closed embedding, and hence so is the map from $\simp{\sigma}^k$ to the product.
\end{proof}

\begin{cor} \label{cor:fans_rational_injectivity}
Let $\psi : F_\bullet \to G_\bullet$ be a map of fans, with $F_\bullet$ quasi-separated. Then $\bP(\psi) : \bP(F) \to \bP(G)$ is injective if and only if it is injective on rational points.
\end{cor}
\begin{proof}
Assume that $\bP(\psi)$ is injective on rational points. Consider two points $x,x' \in \bP(F)$ that map to the same point in $\bP(G)$. Using \Cref{lem:equivalence_rel} one can find two cones $\xi,\xi' \in F_n$ that map to the same cone in $G_n$ along with a point $y \in \simp{}^{n-1}$ such that $x$ and $x'$ are the image of $y$ under the two rational simplices $\xi,\xi' : \simp{}^{n-1} \to \bP(F)$. By \Cref{lem:fans_quasi_sep_relation} the image of the map $\simp{\xi}^{n-1} \sqcup \simp{\xi'}^{n-1} \to \bP(F)$ is the quotient of these two simplices by a closed equivalence relation $R \subset (\simp{\xi}^{n-1} \sqcup \simp{\xi'}^{n-1}) \times (\simp{\xi}^{n-1} \sqcup \simp{\xi'}^{n-1})$. Because $\bP(\psi)$ is injective on rational points, the relation $R$ identifies corresponding rational points on $\simp{\xi}^{n-1}$ and $\simp{\xi'}^{n-1}$, but because it is closed it must identify non-rational points as well. Hence $x = x' \in \bP(F)$.
\end{proof}

\subsubsection{Realizable sets}

\begin{defn} \label{defn:realizable_subset}
Let $F_\bullet$ be a fan. We say that a subset of $\bP(F_\bullet)$ is \emph{realizable} if it is of the form $\bP(G_\bullet) \subset \bP(F_\bullet)$ for some sub-fan $G_\bullet \subset F_\bullet$.
\end{defn}

By part (1) of \Cref{cor:fans_image}, a subset of $\bP(F_\bullet)$ is realizable if and only if it is a union of the images of some collection of rational simplices $\simp{\xi} \to \bP(F)$. In addition, a subset $U \subset \bP(F_\bullet)$ is realizable if and only if its preimage under any rational simplex $\simp{\xi} \to \bP(F_\bullet)$ is realizable.

\begin{ex}
Any open subset of $\bP(F_\bullet)$ is realizable, because any open subset of a simplex $\simp{\xi}$ is a union of sub-simplices.
\end{ex}

Note that for a realizable set $U \subset \bP(F_\bullet)$, there is a \emph{canonical} subfan $G_\bullet \subset F_\bullet$ for which $\bP(G_\bullet) = U$ -- $G_\bullet$ consists of all cones $\xi$ for which $\simp{\xi} \to \bP(F_\bullet)$ factors through $U$. The canonical fan realizing $U$ is also the largest fan realizing $U$.

\begin{lem}
The set of realizable subsets of $\bP(F_\bullet)$ is closed under finite intersections and arbitrary unions.
\end{lem}
\begin{proof}
This follows from the fact that $\bP(-)$ commutes with arbitrary colimits and finite intersections of subfans (see \Cref{cor:fans_image}).
\end{proof}

\subsubsection{Convexity and fans}

Let $U \subset \bR^n_{\geq 0}$ be a convex subset. For any map $[m] \to [n]$, the preimage in $\bR_{\geq 0}^m$ of $U$ under the corresponding linear map $\bR_{\geq 0}^m \to \bR_{\geq 0}^n$ is also convex. For any $F_\bullet \in \op{Fan}$, the space $|F_\bullet|$ consists of copies of $\bR_{\geq 0}^n$ glued along linear maps.

\begin{defn} \label{defn:convex_subset}
For $F_\bullet \in \op{Fan}$, a subset $U \subset |F_\bullet|$ is said to be \emph{locally convex} if for any cone $\sigma \in F_n$, the preimage of $U$ under the corresponding map $\bR_{\geq 0}^n \to |F_\bullet|$ is convex. Likewise a subset $U \subset \bP(F_\bullet)$ is said to be \emph{locally convex} if its preimage under the projection $|F_\bullet| - \{\ast\} \to \bP(F_\bullet)$ is convex. We say that $U \subset \bP(F_\bullet)$ is \emph{convex} if it is locally convex, realizable, and any two rational points are connected by a rational $1$-simplex. \footnote{The conventional notion of a convex subset of the standard simplex $\simp{}^n$ does not include the realizability condition. Realizability is necessary to avoid certain pathologies that are introduced by only considering only line segments with rational endpoints. For instance, consider two cones $\sigma_1,\sigma_2 \subset \bR^2$ that are disjoint away from the origin, and let $F_\bullet = \ray{\{\sigma_1,\sigma_2\}}$. Let $U \subset \bP(F_\bullet)$ be the union of a non-rational point in $\simp{\sigma_1}$ and a non-rational point in $\simp{\sigma_2}$. This $U$ satisfies the definition of convexity without the realizability condition, but it does not have the usual properties of a convex set.}
\end{defn}

Many of the formal properties of convex sets are local and thus generalize to locally convex subsets of $|F_\bullet|$ and $\bP(F_\bullet)$. For instance, the whole space is a locally convex subset, and the intersection of any collection of locally convex subsets is locally convex, so for every $S \subset \bP(F_\bullet)$ or $S \subset |F_\bullet|$ there is a unique smallest locally convex subset containing $S$, which we call the \emph{locally convex hull} of $S$.

\begin{lem} \label{lem:convexity}
A subset $U \subset \bP(F_\bullet)$ is locally convex if and only if for any rational simplex $\simp{}^{n} \to \bP(F_\bullet)$ and points $x_1,x_2 \in \simp{} \to \bP(F_\bullet)$ mapping to $U$, every point in the line segment joining $x_1$ and $x_2$ in $\simp{}^n$ also maps to $U$.\footnote{For the purposes of defining the line segment joining two points, we regard $\simp{}^n$ as embedded in $\bR^{n+1}_{\geq 0}$ as the points $(p_0,\ldots,p_n)$ where $p_0+\cdots+p_n=1$. This notion of line segment does not depend on the choice of affine linear embedding $\simp{}^n \hookrightarrow \bR^m$.}
\end{lem}
\begin{proof}
From the definition and the fact that colimits commute with quotients reduces the lemma to the following claim: if $\pi : \bR^{n+1}_{\geq 0} - \{0\} \to \simp{}^n$ is the quotient projection and $\iota : \simp{}^n \hookrightarrow \bR^{n+1}_{\geq 0}$ is the standard embedding, then for any subset $U \subset \simp{}^n$ the subset $\iota(U) \subset \bR^{n+1}$ is convex if and only if $\pi^{-1}(U) \subset \bR^{n+1}$ is convex. This follows from the fact that the image of a line segment under $\pi$ is a line segment in $\simp{}^n$.
\end{proof}

\begin{rem}
If $\{K_\alpha\}$ is a collection of cones in $\bR^n$ and $U \subset |R_\bullet(\{K_\alpha\})| \simeq \bigcup_{\alpha} K_\alpha$ is convex in the usual sense for subsets of $\bR^n$, then it is locally convex in the sense of \Cref{defn:convex_subset}. The converse is not true, though. For example if $\bigcup K_\alpha = \bR^n$, then the open subset $\bR^n - \{0\} \subset \bR^n$ is locally convex.
\end{rem}

Recall that a function $\Phi$ defined on a convex subset $U \subset \bR^n_{\geq 0}$ is concave if $\Phi(t x + (1-t) y) \geq t \Phi(x) + (1-t) \Phi(y)$ for distinct points $x,y \in U$ and $t \in (0,1)$, and it is strictly concave if strict inequality holds. This notion is preserved by restriction along a linear map $\bR^m_{\geq 0} \to \bR^n_{\geq 0}$, so we may define a function $\Phi(x)$ on a locally convex subset $U \subset |F_\bullet|$ to be \emph{(strictly) concave} if for every cone $\sigma \in F_n$ the restriction of $\Phi(x)$ along the canonical map $\bR^n_{\geq 0} \to |F_\bullet|$ is a (strictly) concave function on the preimage of $U$.

\begin{warning}
The line segment joining two points in $\simp{}^n$ does not have a canonical parameterization. Consider two vectors $v_0,v_1 \in \bR^{n+1}_{\geq 0} -\{0\}$ with distinct image under the quotient map $\pi : \bR^{n+1}_{\geq 0} -\{0\} \to \simp{}^n$. Define $v_t = t v_1 + (1-t) v_0$, so that $\{\pi(v_t) | t \in [0,1]\}$ is the line segment joining the two points in $\simp{}^n$. If we instead let $v'_i = \lambda_i v_i$ for $i=0,1$ and define $v'_t = t v'_1 + (1-t)v'_0$, then
\[
\pi(v'_t) = \pi \left( \frac{t\alpha_1}{t\alpha_1+(1-t)\alpha_0} v_1 +  \frac{(1-t)\alpha_0}{t\alpha_1+(1-t)\alpha_0} v_0 \right) = \pi(v_\tau),
\]
where $\tau = t\alpha_1 /(t\alpha_1+(1-t)\alpha_0)$. Both $\pi(v'_t)$ and $\pi(v_\tau)$ are parameterizations of the same line segment in $\simp{}^n$, but the change of coordinate from $\tau$ to $t$ is not affine-linear and does not preserve concave functions, so \emph{there is no notion of a concave function on $\bP(F_\bullet)$}. Informally, while $|F_\bullet|$ has a canonical affine structure, the space $\bP(F_\bullet)$ only has a canonical convex structure.
\end{warning}

\begin{defn} \label{defn:quasi-concave_function}
Let $U \subset \bP(F_\bullet)$ be a locally convex subset, and let $\Gamma$ be a totally ordered set. We say that a function $\Phi : U \to \Gamma$ is \emph{quasi-concave} if for any rational simplex $\simp{}^n \to \bP(F_\bullet)$ and any three distinct points $x_0,x_1,x_2 \in \simp{}^n$ mapping to $U$ such that $x_1$ lies in the interior of the line segment joining $x_0$ and $x_2$, we have
\[
\Phi(x_1) \geq \min(\Phi(x_0),\Phi(x_2)),
\]
with $\Phi(x_0)=\Phi(x_2)$ whenever equality holds.\footnote{Some definitions elsewhere of quasi-concavity do not require $\Phi(x_0)=\Phi(x_2)$ when equality holds, but this will be useful for our applications.} We say $\Phi$ is \emph{strictly quasi-concave} if this inequality holds strictly for all such points.
\end{defn}

One nice property of quasi-concave functions is that if $\Phi$ is quasi-concave and $m : \Gamma \to \Gamma$ is strictly monotone increasing, the $m\circ \Phi$ is quasi-concave.

\begin{lem} \label{lem:max_quasi_concave}
Let $U \subset \bP(F_\bullet)$ be a convex subspace, let $\Gamma$ be a totally ordered set, and let $\Phi : U \to \Gamma$ be a quasi-concave function that is lower semi-continuous, i.e., $\{ x \in U | \Phi(x) \leq c\}$ is closed for any $c \in \Gamma$. Then any rational point that is a local maximum for $\Phi$ is a global maximum of $\Phi$, and there can be at most $1$ such rational point if $\Phi$ is strictly quasi-concave.
\end{lem}
\begin{proof}
Let $x \in U$ be a rational point that is a local maximum. If $y \in U$ is another rational point, then we can find a rational $1$-simplex in $U$ that connects $x$ and $y$. The restriction of $\Phi$ to this $1$-simplex is a quasi-concave function with a local maximum at the endpoint corresponding to $x$, and it follows from the definition that $x$ is the global maximum along this $1$-simplex, and it is the unique global maximum of $\Phi$ is strictly quasi-concave, so $\Phi(x) \geq \Phi(y)$ with strict inequality if $\Phi$ is strictly quasi-concave. Now if $y \in U$ is any point, the fact that $U$ is realizable implies that we can find a sequence of rational points converging to $y$, and then the fact that $\Phi$ is lower semi-continuous implies that $\Phi(x) \geq \Phi(y)$ as well. Hence $x$ is a global maximum.

\end{proof}

We note that without further hypotheses, $\Phi$ can achieve a maximum at more than one non-rational point.

\begin{ex}
Let $F_\bullet$ be the quotient of $h_{[2]} \sqcup h_{[2]}$ by the equivalence relation identifying corresponding $1$-cones in each summand. Then $\bP(F_\bullet)$ is two copies of $\simp{}^1$ glued along corresponding pairs of rational points, and $U = \bP(F_\bullet)$ is convex. Choose a strictly quasi-concave function on $\simp{}^1$ that achieves its maximum at a non-rational point. Then this function defines a strictly quasi-concave function $\Phi : U \to \bR$ which has two global maxima.
\end{ex}

More generally, we have the following variant of the extreme value theorem:

\begin{lem} \label{lem:continuous_maximizer}
Let $F_\bullet$ be a fan, let $U \subset \bP(F_\bullet)$ be a compact subspace, and let $\Gamma$ be a totally ordered set. If $\Phi : U \to \Gamma$ is an upper semi-continuous function, i.e., $\{x \in U | \Phi(x)< c\}$ is open for all $c \in \Gamma$, then: 
\begin{enumerate}
\item $\Phi$ obtains a global maximum at some point $x \in U$; and
\item If furthermore $F_\bullet$ is quasi-separated, $U \subset \bP(F_\bullet)$ is convex, and $\Phi$ is continuous and strictly quasi-concave, then $x$ is also the unique local maximum of $\Phi$.
\end{enumerate}
\end{lem}
\begin{proof}

It suffices to embed $\Gamma \subset \bar{\Gamma}$ into a totally ordered set such that every subset of $\bar{\Gamma}$ has a supremum, and then to prove the claims for the composition $U \to \bar{\Gamma}$. We will therefore assume that suprema exist in $\Gamma$.

\medskip
\noindent{\textit{Proof of (1):}}
Because $U$ is compact and covered by the open subsets $\Phi^{-1}((-\infty,c))$, it must be equal to one of these subsets. It follows that $\Phi$ is bounded above, so $M := \sup(\Phi(U)) \in \Gamma$ exists.

Suppose that $M \notin \Phi(U)$. If $y \in U$ had $\Phi(y) \geq c$ for all $c<M$, then $\Phi(y)$ would be an upper bound for $\Phi(U)$ that is $<M$. Thus the open subsets $\Phi^{-1}((-\infty,c))$ for $c<M$ cover $U$. Compactness of $U$ then implies that $U = \Phi^{-1}((-\infty,c))$ for some $c<M$, which contradicts the definition of $M$ as the supremum of $\Phi(U)$. Therefore $M \in \Phi(U)$, i.e., the supremum is achieved at some point in $U$.

\medskip
\noindent{\textit{Proof of (2):}}
Let $x \in U$ be a \emph{local} maximum of $\Phi$ with value $m := \Phi(x)$.

First we show that $S := \Phi^{-1}([m,\infty))$ is connected. If $m' := \sup(\{c \in \Gamma | c<m\}) < m$, then $S = \Phi^{-1}((m',\infty))$ is an open subset that is convex because $\Phi$ is quasi-concave, and this implies that it is connected. If $m'=m$, then we have
\[
S = \bigcap_{c < m} \overline{\Phi^{-1}((c,\infty))},
\]
where the closure is taken in $U$. The sets $\overline{\Phi^{-1}((c,\infty))}$ have convex, non-empty interiors, so they are connected. It follows from the fact that $U$ is Hausdorff by \Cref{prop:fans_cg} and compact by hypothesis that $S$ is connected \cite{engelking1989general}*{Thm.~6.1.18}.

For any rational simplex $\simp{\xi} \to U$ and point $\tilde{x}$ mapping to $x$, the strictly quasi-concave function $\Phi|_{\simp{\xi}}$ obtains a local maximum at $\tilde{x}$ and hence $\tilde{x}$ is the unique global maximum of $\Phi|_{\simp{\xi}}$. It follows that if $\xi^{-1}(S) \subset \simp{\xi}$ contains a point in the preimage of $x$, it cannot contain any other point. This implies that the singleton $\{x\} \subset S$ is a closed subset whose complement is also closed. Because $S$ is connected, this implies $S = \{x\}$. Hence $x$ is the unique global maximum.
\end{proof}


\subsection{The degeneration space \texorpdfstring{$\iDeg(\X,p)$}{Deg(X,p)}}

For any field $k$, we introduce the $2$-category of pointed $k$-stacks, $\stack{k}$ as follows: objects are $k$-stacks $\X$ whose inertia group is representable and separated over $\X$ along with a fixed $k$-point $p : \Spec k \to \X$ over $\Spec(k)$. A 1-morphism in this category is a 1-morphism of $k$-stacks $\psi : \X \to \X^\prime$ along with an isomorphism $p^\prime \simeq \psi \circ p$, subject to the constraint that $\psi_\ast : \Aut(q) \to \Aut(\psi(q))$ has finite kernel for all $q \in \X(k)$.\footnote{When $\X$ and $\Y$ are algebraic stacks, we interpret this to mean the kernel of the homomorphism of group schemes $\Aut(q) \to \Aut(\psi(q))$ is a finite $k$-group scheme. More generally for any map $\psi : \X \to \Y$ whose inertia $I_\psi := \X \times_{\X \times_{\Y} \X} \X$ is representable over $\X$, one should require that the fiber of $I_\psi$ over $q \in \X(k)$ is a finite $k$-group scheme. If one would like to consider this definition for non-algebraic stacks $\X$ and $\Y$, then one should instead require $\ker(\Aut(q) \to \Aut(f(q)))$ to be finite after base change to an arbitrary extension $k'/k$.} A $2$-isomorphism is a $2$-isomorphism $\psi \to \tilde{\psi}$ that is compatible with the identification of marked points.

We let the point $1^n$ denote the canonical $k$-point in $\Theta^n_k$ coming from $(1,\ldots,1) \in \bA^n_k$, and regard $\Theta_k^n$ as an object of $\stack{k}$.
\begin{lem} \label{lem:fan_object_in_stacks}
The assignment $[n] \mapsto (\Theta_k^n,1^n)$ extends to a functor $\Cones \to \stack{k}$.
\end{lem}
\begin{proof}
A morphism $\phi : [k] \to [n]$ in $\Cones$ is represented by a matrix of nonnegative integers $\phi_{ij}$ for $i=1,\ldots,n$ and $j = 1,\ldots,k$. One has a map of stacks $ \phi_\ast : \Theta^k \to \Theta^n$ defined by the map $\bA^k \to \bA^n$
$$(z_1,\ldots,z_k) \mapsto (z_1^{\phi_{11}}\cdots z_k^{\phi_{1k}},\ldots,z_1^{\phi_{n1}}\cdots z_k^{\phi_{nk}}),$$
which is intertwined by the group homomorphism $\Gm^k \to \Gm^n$ defined by the same formula. It is clear that these $1$-morphisms are compatible with composition in $\Cones$, and that this construction gives canonical identifications $\phi_\ast (1^k) \simeq 1^n \in \Theta^n(k)$, hence these are morphisms in $\stack{k}$.
\end{proof}

We call a $\bZ^n$-weighted filtration of a $k$-point $f : \Theta_k^n \to \X$ \emph{non-degenerate} if the resulting homomorphism $(\Gm)_k^n \to \Aut(f(0))$ of group sheaves over $\Spec(k)$ has finite kernel. A non-degenerate filtration defines a map $(\Theta_k^n,1^n) \to (\X,f(1^n))$ in $\stack{k}$. Morphisms between pointed $k$-stacks form a groupoid in general, but \Cref{rem:theta_stack_fibers} shows that $\Map_{\stack{k}}((\Theta_k^n,1^n),(\X,p))$ is equivalent to a set whenever the inertia of $\X$ is representable by separated algebraic spaces, and we regard it as such.

\begin{defn} \label{defn:degeneration_space}
Let $\X$ be a stack with representable separated inertia over a fixed base stack $B$. Any $k$-point $p: \Spec(k) \to \X$ induces a $k$-point of $B$ and a canonical pointed $k$-stack $\X_p := \Spec(k) \times_{B} \X$. We define the \emph{\gls{degeneration_fan}}
\begin{equation} \label{eqn:degeneration_space}
\Deg (\X,p)_n := \Map_{\stack{k}}((\Theta_k^n,1^n), (\X_p,p)),
\end{equation}
which defines a functor $\Cones^{op} \to \op{Set}$ via precomposition with the morphisms of \Cref{lem:fan_object_in_stacks}: a morphism $\phi : [k] \to [n]$ of $\Cones$ defines a map $\Deg(\fX,p)_n \to \Deg(\fX,p)_k$ by $f \mapsto f \circ \phi_\ast$. We define the \emph{\gls{degeneration_space}} $\iDeg(\X,p) := \bP(\Deg(\X,p)_\bullet)$. We will sometimes abbreviate the notation to $\Deg(p)_\bullet$ and $\iDeg(p)$ respectively, as the stack $\X$ is implicitly specified by the map $p : \Spec(k) \to \X$.
\end{defn}

Concretely, elements of the set $\Deg(\X,p)_n$ consist of non-degenerate filtrations $f : \Theta_k^n \to \X$ along with an isomorphism $f(1^n) \simeq p$. We will give a complete description of the degeneration space of a point in a quotient stack in \Cref{sect:degeneration_space_quotient}.

\begin{rem}
When $\X$ is an algebraic stack satisfying \ref{hyp3}, \Cref{prop:representable_flags} implies that \eqref{eqn:degeneration_space} is the set of $k$-points of the algebraic space $\Flag^n(p)$ whose underlying filtration is non-degenerate (See also \eqref{eqn:base_changing_mapping}), so this would have been an equivalent definition of $\Deg(\X,p)_n$. As an alternative, one could consider the fan whose $n$-cones consist of \emph{all} finite type points of $|\Flag^n(p)|$ and use this as a basis for the theory of stability. The resulting fan, however, would be equivalent to $\Deg(\X,\bar{p})$ where $\bar{p} : \Spec(\bar{k}) \to \Spec(k) \to \X$ is a geometric point associated to $p$, so \Cref{defn:degeneration_space} is more general.
\end{rem}

\begin{rem}
When $\X$ is not an algebraic stack, we can still discuss the groupoid of $k$-points of
$$\filt[n]{\X}_p := \filt[n]{\X} \times_{\ev_1,\X,p} \Spec(k),$$
but this groupoid need not be equivalent to a set. For instance, if $\X$ is the stack that associates every $k$-scheme $T$ to the groupoid of coherent sheaves on $T$ and isomorphisms between them, the homotopy fiber of $\ev_1 : \filt{\X} \to \X$ has non-trivial automorphism groups. These automorphisms arise as automorphisms of coherent sheaves on $\Theta$ that are supported on the origin in $\bA^1$.
\end{rem}

\subsubsection{The degeneration fan and toric geometry}

The simplest degeneration space to compute is the following:

\begin{ex} \label{ex:degenration_fan_torus}
If $G=T$ is a torus, then $\Deg(\pt/T,\pt)_n$ is the set of all injective homomorphisms $\bZ^n \to \bZ^r$ where $r = \op{rank}T$. This fan is equivalent to $\ray(\bR^r)$ where $\bR^r \subset \bR^r$ is thought of as a single cone. Because this cone admits a simplicial subdivision, \Cref{lem:reconstruct_fan} implies that $|\Deg(\pt/T,\pt)| \simeq \bR^r$ and $\bP(\Deg(\pt/T,\pt)) \simeq S^{r - 1}$.
\end{ex}

Now let us compute the degeneration space for the action of a split torus $T$ on a finite type separated $k$-scheme $X$. Let $p \in X(k)$, and define $T^\prime = T / \Aut(p)$. Define $Y \subset X$ to be the closure of $T \cdot p$ and $\tilde{Y} \to Y$ its normalization. The subgroup $\Aut(p)$ acts trivially on $Y$ and $\tilde{Y}$. $\tilde{Y}$ is a toric variety for the torus $T^\prime$ and thus defines and is defined by a classical fan consisting of rational polyhedral cones $\sigma_i \subset N^\prime_\bR$, where $N^\prime$ is the cocharacter lattice of $T^\prime$. 

\begin{lem} \label{lem:degeneration_fan_torus}
Let $\pi : N_\bR \to N^\prime_\bR$ be the linear map induced by the surjection of cocharacter lattices $N \to N^\prime$ corresponding to the quotient homomorphism $T \to T^\prime$. Then the cones $\pi^{-1} \sigma_i \subset N_{\bR}$ define a classical fan, and the canonical map
$$\Deg(X/T,p)_\bullet \to \Deg(\pt/T,\pt)_\bullet \simeq \ray (N_{\bR})$$
identifies $\Deg(X/T,p)_\bullet$ with the sub fan $\ray \{\pi^{-1} (\sigma_i)\} \subset \ray (N_{\bR})$.
\end{lem}

\begin{proof}
\Cref{cor:embedding_flag_space} implies that because $\tilde{Y} \to X$ is finite, the canonical composition map induces an isomorphism
\[
\Deg(\tilde{Y}/T,p)_\bullet \xrightarrow{\simeq} \Deg(X/T,p)_\bullet,
\]
where on the left $p$ denotes the unique lift of $p$ to the normalization $\tilde{Y}$. So it suffices to prove the claim when $X = \tilde{Y}$ is normal and $p$ has a dense orbit.

\Cref{thm:describe_strata_global_quotient}, along with the fact that for any $\psi \in \Hom(\Gm^n,T)$ the map $X^{\psi,+} \to X$ is a monomorphism (\Cref{cor:contraction_functor}) implies that $p$ has at most one preimage in each component $X^{\psi,+} / T$ of $\Filt^n(X/T)$. In particular $\Deg(X/T,p)_n$ is precisely the set of non-degenerate $\psi \in \Hom(\Gm^n,T)$ for which there is an equivariant map $\bA_k^n \to X$ mapping $1^n \mapsto p$. Such an equivariant map exists if and only if the same is true when we regard $X$ as a $T'$ scheme and consider the composition $\psi' : (\Gm^n)_k \to T \to T'$.

Equivariant maps between toric varieties preserving a marked point in the open orbit are determined by maps of lattices such that the image of any cone in the first lattice is contained in some cone of the second \cite{Fu93}. Applying this to the toric variety $\bA^n_k$ under the torus $\Gm^n$ and to $X$ under the torus $T$, pointed non-degenerate equivariant maps from $\bA^n_k$ to $X$ correspond exactly to non-degenerate homomorphisms $\phi : \bZ^n \to N$ such that the composition $\bR^n \to N'_\bR \to N'_\bR$ maps $\bR_{\geq 0}^n$ to some cone $\sigma_i \subset N'_\bR$. This is exactly the sub-fan $\ray \{\pi^{-1}(\sigma_i)\} \subset \ray (N_\bR)$.
\end{proof}

\begin{ex}
Let $X$ be an affine toric variety defined by a rational polyhedral cone $\sigma \subset \bR^n$ and let $p \in X$ be generic. Then $\Deg(X/T,p)_\bullet \simeq \ray (\sigma)$ as defined in \eqref{eqn:tautological_R_set}. For instance, $\Deg(\bA^n/\Gm^n,1^n)_\bullet \simeq h_{[n]}$ is represented by the object $[n] \in \Cones$.
\end{ex}

The previous lemma and \Cref{lem:reconstruct_fan} implies that $|\Deg(X/T,p)_\bullet| \simeq \bigcup_i \sigma_i \subset N_\bR$ is the support of the fan of the toric variety $\tilde{Y}$, and $\iDeg(X/T,p)$ can be identified with the intersection of this set and the unit sphere in $N_\bR$ with respect to some norm.

\subsubsection{Properties of the degeneration space}

The construction of $\Deg(-)_\bullet$ is functorial in three ways. First, for any map $\psi : \X \to \Y$ with quasi-finite inertia and any $p \in \X(k)$, composition of a filtration $f : \Theta_k^n \to \X$ with $\psi$ defines a functorial map of fans
\[
\psi_\ast : \Deg(\X,p)_\bullet \to \Deg(\Y,\psi(p))_\bullet.
\]
Second, for any pair $p,p' \in \X(k)$ and an isomorphism of $k$-points $p \xrightarrow{\sim} p'$, one has a functorial map
\[
\Deg(\X,p)_\bullet \to \Deg(\X,p')_\bullet
\]
that maps a filtration $f : \Theta_k^n \to \X$ along with an isomorphism $f(1) \simeq p$ to the same filtration with the composed isomorphism $f(1) \simeq p \simeq p'$. In particular, the group $\Aut(p)$ acts on $\Deg(\X,p)_\bullet$.

Finally, one also has canonical base change maps. Given a field extension $k'/k$ and $p \in \X(k)$, extending the map $\Theta_k^n \to \X_p$ to $k'$ defines a map of fans $\Deg(\X,p)_\bullet \to \Deg(\X,p')_\bullet$, where $p' \in \X(k')$ is the image of $p \in \X(k)$. Faithfully flat descent implies

\begin{lem}
$\Deg(\X,p)_\bullet \to \Deg(\X,p')_\bullet$ is injective. If $k'/k$ is Galois, then $\Deg(\X,p)_\bullet$ is the sub-fan of cones in $\Deg(\X,p')_\bullet$ that are fixed points for the natural $\op{Gal}(k'/k)$-action on $\Deg(\X,p')_\bullet$.
\end{lem}

Consider a pair of maps of stacks $\X \to \Y$ and $\X' \to \Y$, and let $\Z := \X \times_\Y \X'$. Fix a point $z \in \Z(k)$ corresponding to points $p \in \X(k), p'\in \X'(k),$ and $q\in \Y(k)$ along with isomorphisms $p \simeq q \simeq p'$ in $\Y(k)$. Assume that the resulting maps of pointed $k$-stacks $\X_p \to \Y_q$ and $\X'_{p'} \to \Y_q$ lie in $\stack{k}$. Then we have
\begin{equation} \label{eqn:fiber_prod}
\Deg(\Z,z)_\bullet = \Deg(\X,p)_\bullet \times_{\Deg(\Y,q)_\bullet} \Deg(\X',p')_\bullet
\end{equation}
from the definition of $\Deg(\Z,z)_\bullet$. The first consequence of this observation is the following:

\begin{prop} \label{prop:intersecting_simplices}
If $\X$ is an algebraic stack satisfying \ref{hyp3}, then $\Deg(\X,p)_\bullet$ is quasi-separated.
\end{prop}

We first prove some initial lemmas. The first observation is that maps $\Theta \to \X$ are the same as maps from the formal completion of $\Theta$ at $\{0\}$. We shall denote the $r^{th}$ infinitesimal neighborhood of $\{0\}$ in $\Theta$ by
\[
\cQ_r := \Spec(\bZ[t]/(t^{r+1})) / \Gm
\]
where $t$ has weight $-1$ with respect to $\Gm$. Note the canonical closed immersions $B\Gm^n = \cQ^n_0 \hookrightarrow \cQ^n_1 \hookrightarrow \cQ^n_2 \hookrightarrow \cdots \hookrightarrow \Theta^n$.

\begin{lem} \label{lem:theta_maps}
For any stack $\X$ satisfying \ref{hyp2} and any integer $n \geq 1$, the restriction map
\[
\Map(\Theta^n,\X) \to \op{holim} \limits_r \Map(\cQ^n_r, \X)
\]
is an equivalence of stacks on the category of finitely presented $B$-schemes.
\end{lem}

\begin{proof}
By \Cref{lem:base_change_mapping}, we may assume that the base $B$ is affine. Note that for any $\Spec(R)$ over $B$ and any open substack $\cU \subset \X$, a map $\Theta^n_R \to \X$ factors through $\cU$ if and only if the restriction to $(\cQ^n_r)_R$ does for any $r$, so we may replace $\X$ with a quasi-compact open substack. Then, we may use relative noetherian approximation to replace $B$ with a noetherian affine scheme and assume $\X$ is of finite presentation over $B$.

For any finitely presented affine $B$-scheme $\Spec(R)$, coherent Tannaka duality \cite{coherent_tannaka}*{Thm.~8.4.ii.(a,b)} implies that
\[
\Map(\Theta^n_R,\X) \cong \Hom_{c,\otimes,\cong}(\QCoh(\X),\QCoh(\Theta^n_R)),
\]
where the latter denotes the category of colimit preserving symmetric monoidal functors and natural equivalences between them. Both categories are compactly generated, which implies that the restriction map is also an equivalence
\[
\Map(\Theta^n_R,\X) \cong \Hom_{c,\otimes,\cong}(\Coh(\X),\Coh(\Theta^n_R)).
\]
The same argument implies $\Map((\cQ^n_r)_R,\X) \cong \Hom_{c,\otimes,\cong}(\Coh(\X),\Coh((\cQ^n_r)_R))$ for any $r \geq 0$. The claim now follows from the fact that the restriction functor $\Coh(\Theta^n_R) \cong \lim_{r} \Coh((\cQ^n_r)_R)$ is an equivalence of symmetric monoidal abelian categories, by \cite{ahr2}*{Prop.~5.1}.
\end{proof}

\begin{lem} \label{lem:almost_toric}
Let $Z$ be a qc.qs. algebraic space of finite type over a field $k$ with an action of $T=(\Gm^n)_k$. Then for any $z \in (Z/T)(k)$, $\Deg(Z/T,z)_\bullet$ is bounded.
\end{lem}
\begin{proof}

We can replace $Z$ with the orbit closure of $z$, and therefore assume $T \cdot z$ is dense. Because $Z/T$ has finitely many points, \cite{ahr2}*{Thm.~20.1} implies that there is a $T$-equivariant Nisnevich cover $W \to Z$ with $W$ affine.

The map $W\to Z$ induces an isomorphism on the reduced identity component of every stabilizer group because it is \'etale. It follows that any map $(B\Gm^n)_k \to Z/T$ admits a lift to $W/T$. Now imagine one is given a morphism $f : \Theta_k^n \to Z/T$ and a lift $\cQ_0^n = (B\Gm^n)_k \to W/T$ of the restriction of $f$ to $\cQ_0^n$. Because $W/T \to Z/T$ is \'etale and the closed immersions $\cQ_r^n \hookrightarrow \cQ^n_{r+1}$ are nilpotent, this lift extends uniquely to a lift $\cQ_r^n \to W/T$ for every $r>0$. \Cref{lem:theta_maps} then implies that this compatible systems of lifts extends uniquely to a lift $\Theta^n_k \to W/T$ of $f$.

Let $w_1,\ldots,w_r \in W(k)$ be the points in the fiber over $z \in Z(k)$. Then the discussion of the previous paragraph shows that
\[
\bigsqcup_{i=1}^r \Deg(W/T,w_i)_\bullet \to \Deg(Z/T,z)_\bullet
\]
is surjective. The left hand side is bounded by \Cref{lem:degeneration_fan_torus}, which implies that $\Deg(Z/T,z)_\bullet$ is bounded.
\end{proof}

\begin{proof}[Proof of \Cref{prop:intersecting_simplices}]
Consider the maps $f_i : \Theta^{n_i+1} \to \X$ with isomorphisms $f_i(1) \simeq p$ corresponding to two cones of $\Deg(\X,p)_\bullet$. Form the fiber product and consider following the diagram, where every square is Cartesian
$$\xymatrix{Y \ar[r] \ar[d]  & Y_0 \ar[d] \ar[r] & \bA^{n_0+1} \ar[r] \ar[d] & \pt \ar[d] \\ Y_1 \ar[r] \ar[d] & \Y \ar[d] \ar[r] & \Theta^{n_0+1} \ar[d] \ar[r] & \pt/\Gm^{n_0+1} \\ \bA^{n_1+1} \ar[r] & \Theta^{n_1+1} \ar[r] & \X & }$$
In other words, there is $\Gm^{n_0+n_1+2}$-torsor over $\Y$ whose total space is $Y$, which is an algebraic space of finite presentation over $k$ because $\X \to B$ is quasi-separated. Hence $\Y$ is a global quotient of an algebraic $k$-space of finite presentation by a torus. The claim now follows from \Cref{lem:almost_toric} and \eqref{eqn:fiber_prod}, which show that $\bP(h_{[n_1]} \times_{\Deg(\X,p)_\bullet} h_{[n_2]}) = \bP(\Y,q)$ is covered by finitely many simplices.

\end{proof}

A second consequence of \eqref{eqn:fiber_prod} is the following:
\begin{prop} \label{prop:representable_maps}
Let $\pi : \X \to \Y$ be a morphism of $B$-stacks that is representable by algebraic spaces, and let $\Y$ (and hence $\X$) have representable separated inertia. Let $p \in \X(k)$, let $q = \pi(p)$, and let $\pi_\ast: \Deg(\X,p)_\bullet \to \Deg(\Y,q)_\bullet$ be the induced map of fans.
\begin{enumerate}
\item If $\pi$ is finite, then $\pi_\ast$ is an isomorphism;
\item If $\pi$ is separated, then $\pi_\ast$ is injective.
\end{enumerate}
If we assume further that $\pi$ is finite type and quasi-separated, then $\pi_\ast$ is bounded, and:
\begin{enumerate}[resume]
\item If $\pi$ is separated, then $\iDeg(\X,p) \to \iDeg(\Y,q)$ is a closed embedding whose restriction to each rational simplex of $\iDeg(\Y,q)$ is a finite union of rational sub-simplices;
\item If $\pi$ is proper, then $\iDeg(\X,p) \to \iDeg(\Y,q)$ is a homeomorphism; and
\item If $\pi$ is affine, then the restriction of the closed subspace $\iDeg(\X,p) \to \iDeg(\Y,q)$ to any rational simplex of $\iDeg(\Y,q)$ is a convex rational polytope.
\end{enumerate}
\end{prop}

\begin{proof}
The fact that $\pi_\ast : \Deg(\X,p)_\bullet \to \Deg(\Y,q)_\bullet$ is an isomorphism if $\pi$ is finite follows immediately from \Cref{prop:affine_map}, which implies that $\Filt^n(\X) \to \Filt^n(\Y) \times_\Y \X$ is a surjective closed immersion and thus universally bijective, and the fact that a filtration $f : \Theta^n_k \to \X$ is non-degenerate if and only if $\pi \circ f$ is non-degenerate.

Consider a rational simplex $\xi : \simp{\xi} \to \iDeg(\Y,q)$, corresponding to a non-degenerate pointed map $f : \Theta^n_k \to \Y$. Let $(\Z,z)$ denote the fiber product $\Theta^n \times_\Y \X$ along with its canonical $k$-point. From \eqref{eqn:fiber_prod} we have an equivalence
\[
\Deg(\Z,z)_\bullet \simeq h_{[n]} \times_{\Deg(\Y,q)_\bullet} \Deg(\X,p)_\bullet.
\]
Furthermore the map $\Z \to \Theta_k^n$ is representable, so $\Z \simeq Z/(\Gm^n)_k$ for an algebraic space $Z$. The set of cones in $\Deg(\X,x)_n$ mapping to the given cone of $\Deg(\Y,y)_n$ is in bijection with the set of equivariant sections of the equivariant map $Z \to \bA^n_k$. If $\pi$ is separated, then $Z$ is separated, so a section is uniquely determined by its restriction to the open orbit. This implies that $\Deg(\X,x)_\bullet \to \Deg(\Y,y)_\bullet$ is injective

Now assume further that $\pi$ is finite type and quasi-separated. Then $Z$ is finite type and quasi-separated, so \Cref{lem:almost_toric} implies that $\Deg(\Z,z)_\bullet$ is bounded. Therefore, the map $\Deg(\X,p)_\bullet \to \Deg(\Y,q)_\bullet$ is bounded. \Cref{lem:fans_closed_map} implies that the map $\iDeg(\X,x) \to \iDeg(\Y,y)$ is closed and its image restricted to any rational simplex of $\iDeg(\Y,y)$ is a finite union of rational sub-simplices. In the case where $f$ is separated, $\iDeg(\X,x) \to \iDeg(\Y,y)$ is then a closed embedding by \Cref{cor:fans_closed_embedding}.

If $\pi$ is proper and we consider the case of $n=1$, the map $Z \to \bA^1_k$ is proper, so the point $z \in Z(k)$ over $1 \in \bA^1_k$ extends uniquely to an equivariant section over $\bA^1_k$ by the valuative criterion. This implies that $\iDeg(\X,x) \to \iDeg(\Y,y)$ is also surjective on rational points, hence a homeomorphism by \Cref{lem:fans_bounded_surjectity}.

Finally, if $\pi$ is affine, then $Z$ is an affine. By part (1) of the proposition, we can replace $Z$ with the normalization of the orbit closure of $z$ without affecting $\Deg(Z/T,z)_\bullet$, so we may assume $Z$ is a normal affine toric variety with a toric map to $\bA^n_k$. It follows from our computation in \Cref{lem:degeneration_fan_torus} that the sub-fan $\Deg(\Z,z)_\bullet \hookrightarrow \Deg(\Theta_k^n,1^n)_\bullet$ has the form $\ray(\sigma) \subset \ray(\bR^n_{\geq0})$ for some rational polyhedral cone $\sigma \subset \bR^n_{\geq0}$ corresponding to $Z$. Thus $\iDeg(\Z,z) \subset \simp{}^{n-1}$ is a rational polytope.
\end{proof}

\begin{rem}
In the proof of \Cref{prop:representable_maps}, we explicitly described the preimage of the closed subspace $\iDeg(\X,x) \hookrightarrow \iDeg(\Y,y)$ under an rational simplex $\simp{\xi} \to \iDeg(\Y,y)$ when $\X \to \Y$ is representable by separated finite type schemes. If the rational simplex is represented by a map $\Theta_k^n \to \Y$, then we consider the fiber product $Z := \bA^n_k \times_\Y \X$, a separated finite type scheme with a $(\Gm^n)_k$-action. There is a canonical orbit in $Z$ lying above $1^n \in \bA^n_k$, and we can let $Z'$ be the normalization of the closure of this orbit. Then the restriction of $\iDeg(\X,x)$ to $\simp{\xi}$ is the projectivization of the support of the fan of the toric variety $Z'$ inside $\bR^n_{\geq 0}$.
\end{rem}

\subsubsection{Separated flag spaces}

When an algebraic stack has separated flag spaces, which is always the case when $\X$ has affine diagonal by \Cref{prop:representable_flags}, we can say more about the degeneration space. 

\begin{lem} \label{lem:unique_fillings}
Let $\X$ be a stack satisfying \ref{hyp3} such that $\ev_1 : \Filt(\X) \to \X$ is separated, and let $p \in \X(k)$. The map $(v_0,\ldots,v_n) : \Flag^{n+1}(p)(k) \to \prod_{i=0}^n \Flag(p)(k)$ is injective, so any cone in $\Deg(\X,p)_\bullet$ is uniquely determined by its vertices.
\end{lem}
\begin{proof}
This follows by induction from the claim that the map
\[
(v_0,d_0) : \Flag^{n+1}(p)(k) \to \Flag(p)(k) \times \Flag^n(p)(k)
\]
is injective. Using the equivalence $\Filt^{n+1}(\X) \simeq \Filt^n(\Filt(\X))$, we can identify $\Flag^{n+1}(p)$ with the set of points $f \in \Filt(\X)(k)$ with an isomorphism $f(1) \simeq p$ and a pointed map $\Theta_k^n \to \Filt(\X)$, where the latter is regarded as pointed stack with point $f$. The latter data is equivalent to a pointed map $\sigma : \Theta_k^n \to \X$ along with a section of the projection
\[
\Theta^n_k \times_{\X} \Filt(\X) \to \Theta^n_k
\]
taking the point $1^n \in \Theta_k^n$ to the point $(1^n,f) \in \Theta^n_k \times_{\X} \Filt(\X)$. The point $1^n$ is open and dense in $\Theta_k^n$, so if the map $\ev_1$ is separated, any section is uniquely determined by its restriction to a dense open subset. It follows that the existence of a section is a condition, and not data, so the pointed map $\sigma : \Theta^{n+1}_k \to \X$ is uniquely determined by $f = v_0(\sigma)$ and $d_0(\sigma) \in \Flag^{n}(p)(k)$.
\end{proof}

\begin{lem} \label{lem:closed_embedding}
Let $\X$ be a stack satisfying \ref{hyp3} such that $\ev_1 : \Filt(\X) \to \X$ is separated. Then for any point $p \in \X(k)$ any rational simplex $\simp{\xi} \to \iDeg(\X,p)$ is a closed embedding.
\end{lem}
\begin{proof}
It suffices by \Cref{cor:fans_rational_injectivity} and \Cref{lem:fans_closed_map} to show that $\simp{\xi}^n \to \iDeg(\X,p)$ is injective on rational points. Assume that two points $x_0,x_1 \in \simp{\xi}$ map to the same rational point in $\iDeg(\X,p)$. Let $f_0,f_1 \in \bZ_{>0}^{n+1}$ represent the points $x_0,x_1$, and correspond to non-degenerate pointed maps $\Theta_k \to \Theta_k^{n+1}$. $f_0$ and $f_1$ together define a unique non-degenerate pointed map $\Theta_k^2 \to \Theta_k^{n+1}$. Composing this with the map $\xi : \Theta_k^{n+1} \to \X$ gives a (possibly degenerate) point $\gamma \in \Flag^2(p)(k)$ such that $v_0(\gamma) = f_0$ and $v_1(\gamma)=f_1$. By hypothesis after replacing $f_i$ with a positive multiple we can arrange that the two compositions
\[
\Theta_k \xrightarrow{f_i} \Theta_k^{n+1} \xrightarrow{\xi} \X
\]
are the \emph{same} pointed map. It follows from the uniqueness of fillings, \Cref{lem:unique_fillings}, that $\gamma$ must correspond to a $2$-dimensional filtration that factors through the projection $\Theta_k^2 \to \Theta_k$ corresponding to the homomorphism $\phi : \bZ^2 \to \bZ$ defined by $\phi(a,b) = a+b$. As a result, every $1$-dimensional sub-cone of $\gamma$ corresponds to the same rational point in $\iDeg(\X,p)$ as well. In other words the map $\simp{\xi} \to \iDeg(\X,p)$ collapses all rational points in the $1$-simplex connecting $x_0,x_1 \in \simp{\xi}$ to a single point in $\iDeg(\X,p)$. This contradicts the finiteness of the fibers of the map $\simp{\xi} \to \iDeg(\X,p)$, established in \Cref{lem:fans_closed_map}.
\end{proof}


\subsection{\texorpdfstring{$\iDeg(X/G,p)$}{Deg(X/G,p)} is a generalized spherical building} \label{sect:degeneration_space_quotient}
We have already seen in \Cref{lem:degeneration_fan_torus} that for a toric variety $X$, the degeneration fan $\Deg(X/T,\one)_\bullet$ is essentially equivalent to a more classical construction, the toric fan associated to $X$. In this section we provide a more classical description of $\Deg(X/G,p)_\bullet$ when $G$ is not abelian.

First we address the case where $\X = \pt / G$ for a semisimple group $G$. Recall that the spherical building $\Delta(G)$ is a simplicial complex whose vertices are the maximal parabolic subgroups of $G$, where a set of vertices spans a simplex if and only if the corresponding maximal parabolics contain a common parabolic subgroup of $G$.
\begin{prop} \label{prop:spherical_building}
Let $\X = \pt / G$ where $G$ is a split semisimple group, and let $p \in \X(k)$ correspond to the trivial $G$-bundle. Then $\iDeg(\X, p)$ is homeomorphic to the spherical building $\Delta(G)$.
\end{prop}

\begin{proof}
\Cref{thm:describe_strata_global_quotient} implies that $\Flag^n(p) \simeq \bigsqcup_{\psi} G / P_\psi$ where $\psi$ ranges over $\Hom((\Gm^n)_k,T) / W$, and $\Deg(\X,p)_n$ is the set of non-degenerate $k$-points of this flag space. Thus to every cone $\xi \in \Deg(\X,p)_n$ we associate a $\psi \in \Hom((\Gm^n)_k,T)$ up to conjugation by $W$ and a subgroup conjugate to $P_\psi \subset G$. $P_\psi$ need not be parabolic if $n>1$. We denote $F = \Deg(\X,p)_\bullet$, and $F^{par} \subset F$ the sub-fan consisting of all cones such that the associated $P_\psi$ is parabolic.

Fix a split maximal torus $T \subset G$ with cocharacter lattice $N$, and let $\cW$ denote the classical fan in $N_\bR$ defined by the Weyl chambers. We also choose a dominant Weyl chamber and an ordering on the set of minimal generators $\{v_1,\ldots,v_r\} \in N$ for the rays of the dominant Weyl chamber. Every subset of $\{v_1,\ldots,v_r\}$ of cardinality $n$ defines a homomorphism $\psi : \Gm^n \to G$ such that $P_\psi$ is parabolic. Let $\cS \subset (\Cones|F^{par})$ be the full subcategory of cones $\xi \in \Deg(\fX,p)_n$ whose associated conjugacy class of homomorphism $\psi : \Gm^n \to G$ is one of these. The inclusions $\cS \subset (\Cones|F^{par}) \subset (\Cones|F)$ induce canonical maps
$$\colim \limits_{\xi \in \cS} \simp{\xi} \xrightarrow{(a)} \colim \limits_{\xi \in (\Cones|F^{par})} \simp{\xi} \xrightarrow{(b)} \colim \limits_{\xi \in (\Cones|F)} \simp{\xi}.$$

The statement of the proposition follows from the claim that $(a)$ and $(b)$ are homeomorphisms, because we may identify the first colimit with the spherical building $\Delta(G)$. Indeed, $k$-cones in $\cS$ correspond bijectively to proper parabolic subgroups of $G$ that are contained in $k$ distinct maximal parabolics, i.e., $k-1$-simplices in $\Delta(G)$, and morphisms in $\cS$ correspond to containment of parabolics, which corresponds to the containment of faces in $\Delta(G)$.

\medskip
\noindent \textit{Step 1:  $F^{par} \subset F$ is a bounded inclusion of fans that is surjective on $1$-cones, so the map $(b)$ is a homeomorphism by \Cref{lem:fans_bounded_surjectity} and \Cref{cor:fans_closed_embedding}.}
\medskip

Surjectivity on $1$-cones follows from the fact that for any $\lambda : \Gm \to G$, the subgroup $P_\lambda$ is parabolic. To show that the map $F^{par} \to F$ is bounded, consider a $\xi : h_{[n]} \to F$, corresponding under \Cref{thm:describe_strata_global_quotient} to a $\psi \in \Hom((\Gm^n)_k,T)$ and a $k$-point of $G/P_\psi$. For any sub-cone $h_{[m]} \subset h_{[n]}$, corresponding to a homomorphism $\bZ^m \hookrightarrow \bZ^n$ with nonnegative matrix coefficients, the composition $h_{[m]} \to h_{[n]} \to F$ lies in $F^{par}$ if and only if for the resulting homomorphism
\[
\psi' : (\Gm^m)_k \to (\Gm^n)_k \xrightarrow{\psi} T
\]
the subgroup $P_\psi$ is parabolic. This happens if and only if the image of $\bR_{\geq 0}^n$ under the homomorphism on cocharacter lattices $\bZ^m \to N$ is contained in some cone of the fan $\cW$. In other words
\[
h_{[n]} \times_F F^{par} = \ray{\cW'} \subset \ray(\bR^n_{\geq 0}) = h_{[n]},
\]
where $\cW'$ is the classical fan in $\bR^n$ obtained by taking the preimage of $\cW$ under the homomorphism $\bR^n \to N_\bR$ induced by $\psi$ and intersecting with the nonnegative cone $\bR^n_{\geq 0}$.

\medskip
\noindent \textit{Step 2: $\cS \subset (\Cones|F^{par})$ is cofinal, so the map $(a)$ is a homeomorphism.}
\medskip

It suffices to show that the inclusion $\cS \subset (\Cones|F^{par})$ admits a left adjoint. Let $\xi \in F_n^{par}$, and let $\psi : \Gm^n \to T$ represent the conjugacy class of homomorphism associated to $\xi$. Because $P_\psi$ is parabolic, the cone in $N_\bR$ associated to $\psi$ must be in $\ray{\cW}$. Up to conjugation by an element of the Weyl group we may assume that $\psi$ maps $(\bR_{\geq 0})^n$ to the dominant Weyl chamber of $\cW$.

Let $\sigma \in \cW$ be the smallest cone containing the image of $(\bR_{\geq 0})^n$. $P_\sigma = P_\psi$, and because cones in the dominant Weyl chamber classify conjugacy classes of parabolics, $\sigma$ is thus uniquely determined by the conjugacy class of $\psi$. As the ray generators of $\sigma$ form a basis for the subgroup generated by $\sigma \cap N$, it follows that there is a unique $n$-cone $\xi^\prime \in \cS$ and a morphism $\phi : [n] \to [n]$ such that $\xi = \phi^\ast \xi^\prime$. A similar argument shows that any morphism in $(\Cones|F^{par}_\bullet)$ from $\xi$ to an element of $\cS$ factors uniquely through this $\xi^\prime$, and thus that the assignment $\xi \mapsto \xi'$ is a left adjoint for the inclusion $\cS \subset (\Cones|F^{par})$.
\end{proof}

\Cref{prop:spherical_building} justifies the following
\begin{defn}
For an algebraic group $G$ over a field $k$, the \emph{generalized spherical building} of $G$ is the degeneration space $\Delta(G) := \iDeg(\pt/G,\pt)$.
\end{defn}
For non-semisimple $G$, this space does not canonically have the structure of a simplicial complex. For instance, when $G$ is a torus we have seen in \Cref{ex:degenration_fan_torus} that $\Delta(G)$ is a sphere. Nevertheless for a split reductive group $G$ of rank $r$, $\Delta(G)$ has a canonical cover by rational simplices after fixing a maximal torus $T \subset G$. For any Borel subgroup $B \subset G$ and rational basis $\{v_1,\ldots,v_p\}$ of $N_\bR^W$, where $N$ is the cocharacter lattice of $T$ and $W$ is the Weyl group, we shall construct a rational simplex $\simp{B,v_1,\ldots,v_p}^{r-1} \subset \Delta(G)$ :

The Borel subgroup $B$ selects a dominant Weyl chamber in $N'_\bR := N_\bR / N_\bR^W$. There is a unique and hence canonical $W$-equivariant splitting $N_\bR = N_\bR^W \oplus N'_\bR$, so we can using generators for the extremal rays of the dominant Weyl chamber in $N'_\bR$ we can extend the vectors above to a rational basis $\{v_1,\ldots,v_r\}$ of $N_\bR$. Clearing denominators so that $v_i$ are integral, they define a homomorphism $\psi : \Gm^r \to G$ for which $P_\psi$ is conjugate to $B$. Hence $B$ defines a $k$ point of $G/P_\psi$ and thus an $r$-cone $\xi \in \Deg(\pt/G,\pt)_r$. The basis $v_1,\ldots,v_r$ is uniquely defined from the original data $(B,v_1,\ldots,v_p)$ up to positive rescaling of the $v_i$, which implies that the resulting rational simplex
\[
\simp{B,v_1,\ldots,v_p} \hookrightarrow \iDeg(\pt/G,\pt)
\]
is determined uniquely by this data. These simplices are top dimensional, and they cover $\Delta(G)$. When $G$ is semisimple, $p=0$ and these are the simplices arising in the usual description of the spherical building.

\begin{rem} \label{rem:gerbes}
One can deduce from the above discussion that if $G' := G/Z(G)$, then $\Delta(G)$ is a sphere bundle over $\Delta(G')$. It seems likely that in general if $\pi : \X \to \Y$ is a gerbe for the group $\Gm^n$, then one can canonically identify $\iDeg(\X,p)$ with a fiber bundle over $\iDeg(\Y,\pi(p))$ with fiber $S^n$. Note that the map $\pi$ is not quasi-finite, and therefore does not induce a map of fans $\Deg(\X,p)_\bullet \to \Deg(\Y,\pi(p))_\bullet$.
\end{rem}

We can now complete the description of $\iDeg(X/G,p)$ when $X$ is a separated finite type $k$-scheme and $p \in X(k)$. We regard it as a closed subspace $\iDeg(X/G,p) \subset \iDeg(\pt/G,\pt)$ using \Cref{prop:representable_maps}. It suffices to describe the intersection with any rational simplex. \Cref{thm:describe_strata_global_quotient} shows that a pointed map $\xi : \Theta^n_k \to \pt/G$ corresponds to a choice of $\psi : \Hom(\Gm^n,T)/W$ and a point $g \in G/P_\psi$. One can identify the fiber product
$$\xymatrix{(\bA^n \times X) / \Gm^n \ar[r] \ar[d] & X/G \ar[d] \\ \Theta^n \ar[r]^{\xi} & \pt/G},$$
where the action of $\Gm^n$ on $\bA^n \times X$ is the diagonal action induced by the standard action of $\Gm^n$ on $\bA^n$ and the action on $X$ via $g \psi g^{-1}$. The point $p \in X(k)$ defines a point $(1^n,p) \in \bA^n \times X$. Let $\tilde{X}$ be the normalization of the closure of $\Gm^n \cdot (1^n,p)$ in $\bA^n \times X$. Then \Cref{prop:representable_maps} and the discussion leading up to it imply that the intersection of $\iDeg(X/G,p) \subset \Delta(G)$ with the simplex $\simp{\xi}^{n-1} \subset \Delta(G)$ is the degeneration space $\iDeg(\tilde{X}/\Gm^n,(1^n,p))$. By \Cref{lem:degeneration_fan_torus} this is the support of the fan of the toric variety $\tilde{X}$ modulo the scaling action of $\bR^\times_{>0}$, where the toric map $\tilde{X} \to \bA^n$ is used to identify the support of this fan with a union of subcones of $(\bR_{\geq 0})^n$. In particular $\iDeg(\tilde{X}/\Gm^n,(1^n,p))$ is a sub-polyhedron of $\simp{\xi}^{n-1}$.

\begin{ex} \label{ex:affine_quotient}
Let $V$ be an affine scheme and $G$ be a split reductive $k$-group. For any $k$-point in $V$, the morphism $V/G \to BG$ identifies $\iDeg(V/G,p)$ with a closed subspace of the spherical building $\Delta(G)$ whose intersection with any rational simplex $\simp{\xi} \subset \Delta(G)$ is a rational polytope $P \subset \simp{\xi}$. We call such a subset of $\Delta(G)$, whose intersection with each $\simp{\xi}$ is a rational polytope, a rational polytope in $\Delta(G)$. This suggests a non-abelian analog of toric geometry where one encodes a normal affine $G$-variety with dense open orbit by a rational polytope in $\Delta(G)$, analogous to the rational polyhedral cone that encodes a normal affine toric variety.
\end{ex}


\subsection{A result on extension of filtrations}

In this section we study the relationship between $n+1$-dimensional filtrations $\sigma$ of a point $p\in \X(k)$ with $v_0(\sigma) = f$, and $n$-dimensional filtrations of $\agr(f)$ as a point in $\Grad(\X)$. Our main result, \Cref{prop:extension}, is essentially the first half of the main theorem on perturbation of filtrations in the next section \Cref{thm:perturbation}.

Given a graded object $g \in \Grad(\X)(k)$ and a filtration $\nu \in \Deg(\Grad(\X),g)_1$, we consider the corresponding map
\[
\nu : (\pt/\Gm)_k \times \Theta_k \simeq \bA^1_k / (\Gm^2)_k \to \X.
\]
We regard the coordinate $t$ on $\bA^1_k$ as having weight $(0,-1)$ with respect to the torus $(\Gm^2)_k$, and consider the weight decomposition of the representation of $(\Gm^2)_k$
\begin{equation} \label{eqn:cotangent_weights}
H^0((\pt/\Gm^2)_k, \nu^\ast \bL_\X|_{\{0\}}) \oplus H^1((\pt/\Gm^2)_k, \nu^\ast \bL_\X|_{\{0\}}) = \bigoplus_{a_0,a_1 \in \bZ} W_{a_0,a_1}.
\end{equation}
We say that $(a_0,a_1)$ is a \emph{cotangent weight} for $\nu \in \Deg(\Grad(\X),g)_1$ if $W_{a_0,a_1} \neq 0$ in the decomposition above.

\begin{defn} \label{defn:near_canon}
Let $\X$ be an algebraic stack, and let $g \in \Grad(\X)(k)$. Then we define the sub-fan
\[
\Deg(\Grad(\X),g)^\canon_\bullet \subset \Deg(\Grad(\X),g)_\bullet
\]
to consist of cones $\sigma$ such that if $\nu \in \Deg(\Grad(\X),g)_1$ is a sub-cone of $\sigma$, then any cotangent weight $(a_0,a_1)$ of $\nu$ with $a_0<0$ has $a_1 \leq 0$. We denote the projective realization $\iDeg(\Grad(\X),g)^\canon$.
\end{defn}

\medskip

For a general map of stacks $\phi : \X \to \X'$ and point $p \in \X(k)$, the composition map $\Flag^n(p) \to \Flag^n(\phi(p))$ does not preserve non-degenerate filtrations.\footnote{We have shown that this is the case for maps $\phi$ with quasi-finite inertia.} Nevertheless we have

\begin{lem} \label{lem:ev_0_fans}
For any stack $\X$ with separated inertia, the map of stacks $\agr : \Filt(\X) \to \Grad(\X)$ preserves non-degenerate $n$-dimensional filtrations and thus induces a map of fans
\[
\agr : \Deg(\Filt(\X),f)_\bullet \to \Deg(\Grad(\X),\agr(f))_\bullet
\]
for any $f \in \Filt(\X)(k)$.
\end{lem}
\begin{proof}
This is just a matter of unraveling definitions. An $n$-dimensional filtration of $f$ is a map $\sigma : \Theta^{n+1}_k = \bA_k^{n+1} / (\Gm^{n+1})_k \to \X$ with an isomorphism $v_0(\sigma) \simeq f$. Denote the coordinates on $\bA_k^{n+1}$ by $(t_0,\ldots,t_n)$. Then $\agr(\sigma)$ is the restriction of $\sigma$ to the closed substack $\{t_0=0\} / (\Gm^{n+1})_k \simeq (\pt/\Gm)_k \times \Theta_k^n \subset \Theta_k^{n+1}$. Non-degeneracy of either filtration $\sigma$ or $\agr(\sigma)$ amounts to the same condition: that the homomorphism $(\Gm^{n+1})_k \to \Aut(\sigma(0^{n+1}))$ induced by $\sigma$ at $0^{n+1} \in \Theta_k^{n+1}$ has finite kernel when restricted to the subtorus $\{1\} \times (\Gm^n)_k$.
\end{proof}

We can now state our main extension result:

\begin{prop} \label{prop:extension}
Let $\X$ be an algebraic stack satisfying \ref{hyp3}, and let $f \in \Filt(\X)(k)$ be a non-degenerate filtration. Every cone in $\Deg(\Grad(\X),\agr(f))^\canon_\bullet$ lifts uniquely to $\Deg(\Filt(\X),f)_\bullet$, so there is a unique map of fans $\op{ext}$ making the following diagram commute
\[
\xymatrix{ & \Deg(\Filt(\X),f)_\bullet \ar[d]_{\agr} \\ \Deg(\Grad(\X),\agr(f))^\canon_\bullet \ar@{-->}[ur]^{\op{ext}} \ar@{^{(}->}[r] & \Deg(\Grad(\X),\agr(f))_\bullet }.
\]
\end{prop}

The proof of this proposition will appear at the end of this section, and is a straightforward application of a slightly more technical result, which we now formulate. Say we are given the following data:
\begin{enumerate}
\item a filtration $f : \Theta_k \to \X$,
\item an $\bZ^n$-weighted filtration $f_0 : \Theta_k^n \to \Grad(\X)$, and
\item an isomorphism $f_0(1) \simeq \agr(f) \in \Grad(\X)$.
\end{enumerate}
We can reinterpret this data in terms of certain substacks of $\Theta_k^{n+1}$. Consider the following subschemes of $\bA_k^{n+1}$, with coordinates $t_0,\ldots,t_n$:
\[
Y_0 := \{ t_0 = 0\} \quad \text{and} \quad U := \{t_1 \cdots t_n \neq 0\}.
\]
Note that $U / \Gm^{n+1} \simeq \Theta_k$ with coordinate $t_0$, and that a map
\[
f_0 : Y_0 / \Gm^{n+1} \simeq (\pt/\Gm)_{t_0} \times \Theta_k^n \to \X
\]
classifies an $n$-dimensional filtration $\Theta_k^n \to \Grad(\X)$ of the graded object classified by the restriction of $f_0$ to a map $\{(0,1,\ldots,1)\} / (\Gm)_{t_0} \to \X$. Note also that the inclusion of this point induces an equivalence of stacks $\{(0,1,\ldots,1)\} / (\Gm)_{t_0} \simeq U \cap Y_0 / \Gm^{n+1}$. Thus the data above is equivalent to specifying:
\begin{enumerate}[label=(\arabic*')]
\item a map $f : U / \Gm^{n+1} \to \X$,
\item a map $f_0 : Y_0 / \Gm^{n+1} \to \X$,
\item an isomorphism between $f$ and $f_0$ after restricting to $U \cap Y_0 / \Gm^{n+1} \simeq \{(0,1,\ldots,1)\}/(\Gm)_{t_0}$.
\end{enumerate}
We can now formulate the unusual gluing question of when the data (1')-(3') results from restricting a $\bZ^{n+1}$-weighted filtration $f' : \Theta^{n+1} \to \X$ to $Y_0/\Gm^{n+1}$ and $U/\Gm^{n+1}$. In other words, we want to know then there is a dotted arrow that makes the following diagram commute:
\begin{equation} \label{eqn:extension}
\xymatrix{ & Y_0/\Gm^{n+1} \ar@/^/[drrr]^{f_0} \ar[dr] & & & \\
Y_0 \cap U / \Gm^{n+1} \ar[ur] \ar[dr] & & \bA_k^{n+1} / \Gm^{n+1} \ar@{-->}[rr]^{f'} & & \X \\
& U / \Gm^{n+1} \ar[ur] \ar@/_/[urrr]_{f} & & & }
\end{equation}
If such an extension $f'$ exists, then we can restrict $f'$ to $\Theta_k^n \simeq \{t_0 \neq 0\} / \Gm^{n+1} \subset \Theta_k^{n+1}$ to obtain a new $\bZ^n$-weighted filtration of $f(0)$. This is the perturbation of the filtration $f$ discussed above in the case of coherent sheaves.

\begin{prop} \label{prop:extension_technical}
Let $\X$ be a stack satisfying \ref{hyp2}, and assume we are given data (1')-(3'). Consider the space $H^i((f_0^\ast \bL_\X)|_{\{(0,\ldots,0)\}})$ as a representation of $\Gm^{n+1}$, and let $H^i((f_0^\ast \bL_\X)|_{\{(0,\ldots,0)\}})_{a_0,\ldots,a_n}$ denote the summand of weight $(a_0,\ldots,a_n)$ with respect to $\Gm^{n+1}$. If
\[
H^i((f_0^\ast \bL_\X)|_{\{(0,\ldots,0)\}})_{a_0,\ldots,a_n} = 0
\]
whenever $i=0,1$, $a_0<0$, and $a_j>0$ for some $j$, then there exists an extension $f'$ making the diagram \eqref{eqn:extension} commute, and it is unique up to unique isomorphism.
\end{prop}

First we establish the following:

\begin{lem} \label{lem:extension_basic}
Let $\X$ be a stack satisfying \ref{hyp2}, and consider the open subscheme
\[
V = \{(t_0,\ldots,t_n) | t_n \neq 0\} \subset \bA_k^{n+1}.
\]
Fix maps $f : V/\Gm^{n+1} \to \X$, $f_0 : Y_0 /\Gm^{n+1} \to \X$, and an isomorphism between $f$ and $f_0$ after restricting to $V \cap Y_0 / \Gm^{n+1}$. Assume that for all $m<0,l>0$ and for $i=0,1$ we have
\[
\left[ H^i((f_0^\ast \bL_\X)|_{\{(0,\ldots,0)\}}) \right]_{\begin{subarray}{c}t_0\text{-weight } m, \\ t_n\text{-weight }l \\ \text{ summand}\end{subarray}} = 0.
\]
Then $f,f_0$, and the given isomorphism are induced by restriction from a map $f' : \bA_k^{n+1}/\Gm^{n+1} \to \X$, which is unique up to unique isomorphism.
\end{lem}

\begin{rem} The extension problem in this lemma is identical to that of \eqref{eqn:extension}, but with $V$ instead of $U$. The gluing data in this case is more complicated, though, because $Y_0 \cap V / \Gm^{n+1} \simeq (\pt/\Gm)_{t_0} \times (\Theta^{n-1})_{t_1,\ldots,t_{n-1}}$.
\end{rem}

\begin{proof}
Let $Y_m = \Spec (k[t_0,\ldots,t_n] / (t_0^{m+1}))$ be the $m^{th}$ infinitesimal neighborhood of $Y_0 = \{0\} \times \bA^n \hookrightarrow \bA^{n+1}$. Then $Y_m / \Gm^{n+1} \simeq \cQ_m \times \Theta_k^n$, where $\cQ_m = \Spec(k[t]/(t^{m+1}))/\Gm$. Combining \Cref{lem:theta_maps} with the canonical isomorphism $\filt[n]{\filt{\X}} \simeq \filt[n+1]{\X}$, we conclude that the restriction map
\[
\Map(\Theta_k^{n+1},\X) \to \op{holim}_m \Map(Y_m / \Gm^{n+1},\X)
\]
is an equivalence of groupoids. We likewise define $V_m = V \cap Y_m$. The same reasoning as above implies that the canonical map
\[
\Map(V / \Gm^{n+1},\X) \to \op{holim}_m \Map(V_m / \Gm^{n+1}, \X)
\]
is an equivalence of groupoids.

We can now rephrase the goal: we are given a map $f_0 : Y_0 / \Gm^{n+1} \to \X$ and an extension of the restriction of $f_0$ to $V_0/\Gm^{n+1} \subset Y_0 / \Gm^{n+1}$ to a pro-system of maps $V_m / \Gm^{n+1} \to \X$, and we wish to extend this pro-system by filling in the dotted arrows in the diagram
\[
\xymatrix@R=7pt{V_0 \ar@{^{(}->}[r] \ar@{}|-{\cap}[d] & V_1 \ar@{^{(}->}[r] \ar@{}|-{\cap}[d] & V_2 \ar@{^{(}->}[r] \ar@{}|-{\cap}[d] & V_3 \ar@{^{(}->}[r] \ar@{}|-{\cap}[d] & \cdots \\
Y_0 \ar@{^{(}->}[r]^{i_0} \ar[dd]_{f_0} & Y_1 \ar@{^{(}->}[r]^{i_1} \ar@{-->}[ddl]_{f_1} & Y_2 \ar@{^{(}->}[r]^{i_2} \ar@{-->}[ddll]_{f_2} & Y_3 \ar@{^{(}->}[r]^{i_3} \ar@{-->}[ddlll]^{f_3} & \cdots \\
 & & & & \\
\X & & & & }
\]
The maps $f_m$ can be constructed inductively, so it suffices to show that the dotted arrow has a unique filling making the following diagram commute:
\begin{equation} \label{eqn:def_problem}
\xymatrix{
    V_{m-1} / \Gm^{n+1} \ar@{^{(}->}[r] \ar@{}|-{\cap}[d] & V_m/\Gm^{n+1} \ar@{}|-{\cap}[d] \ar@/^/[ddr] & \\
    Y_{m-1} / \Gm^{n+1} \ar@/_/[drr]_{f_{m-1}} \ar@{^{(}->}[r]^{i_{m-1}} & Y_m / \Gm^{n+1} \ar@{-->}[dr]_{f_m} & \\
    & & \X
}
\end{equation}

This now amounts to a deformation theory problem. Each inclusion $i_{m-1} : Y_{m-1} / \Gm^{n+1} \to Y_{m} / \Gm^{n+1}$ is a square-zero extension of algebraic stacks by the coherent sheaf $i_\ast \cO_{Y_0} \twist{m,0,\ldots,0}$, where we abuse notation slightly to denote the composed inclusion
\[
i : Y_0 / \Gm^{n+1} \to Y_{m-1}/\Gm^{n+1}.
\]
Given an extension of $f_0$ to a map $f_{m-1} : Y_{m-1} / \Gm^{n+1} \to \X$, the map can be extended to $Y_m/\Gm^{n+1}$ if and only if the composition
\[
f_{m-1}^\ast \bL_\X \to \bL_{Y_{m-1}/\Gm^{n+1}} \to i_\ast (\cO_{Y_0} \twist{m,0,\ldots,0}[1])
\]
vanishes as an element of $H^0 (R\Hom(f_{m-1}^\ast \bL_\X, i_\ast (\cO_{Y_0} \twist{m,0,\ldots,0}[1])))$. Furthermore, if this map vanishes, then the set of extensions to $Y_{m}/\Gm^{n+1}$ up to isomorphism is a torsor for the group $$H^0 (R\Hom(f_{m-1}^\ast \bL_\X, i_\ast (\cO_{Y_0} \twist{m,0,\ldots,0}))).$$

The same analysis applies to the square zero extension $V_{m-1} / \Gm^{n+1} \hookrightarrow V_m / \Gm^{n+1}$, so the extension $f_m$ in \eqref{eqn:def_problem} exists and is unique as long as the restriction map
\begin{multline*}
R\Hom_{Y_{m-1}/\Gm^{n+1}}(f^\ast_{m-1} \bL_\X, i_\ast (\cO_{Y_0} \twist{m,0,\ldots,0})) \to \\
R\Hom_{V_{m-1}/\Gm^{n+1}}(f^\ast_{m-1} \bL_\X|_{V_{m-1}},i_\ast (\cO_{V_0} \twist{m,0,\ldots,0}))
\end{multline*}
is injective in cohomological degree $1$ and an equivalence in degree $0$, for every $m \geq 1$. Using the adjunction and the isomorphism $f_0 \simeq f_{m-1} \circ i$, we can rewrite the $R\Hom$ complexes above only in terms of $f_0$ and $Y_0/\Gm^{n+1}$, identifying the map above with the restriction map
\begin{multline*}
R\Hom_{Y_0/\Gm^{n+1}}(f^\ast_0 \bL_\X, \cO_{Y_0} \twist{m,0,\ldots,0} ) \to \\
R\Hom_{V_0/\Gm^{n+1}}(f^\ast_0 \bL_\X|_{V_0},\cO_{V_0} \twist{m,0,\ldots,0}).
\end{multline*}

We apply the long exact sequence in local cohomology for the closed substack $\cS_0 := \{t_n = 0\} / \Gm^{n+1} \subset Y_0/\Gm^{n+1}$, whose complement it $V_0/\Gm^{n+1}$. This sequence implies that for the map above to be injective in degree $1$ and bijective in degree $0$, it suffices to show that for all $m>0$
\begin{equation} \label{eqn:local_cohomology}
H^i R\Hom_{Y_0 / \Gm^{n+1}}(f_0^\ast \bL_\X, R\Gamma_{\cS_0} \cO_{Y_0} \otimes \cO_{Y_0} \twist{m,0,\ldots,0})) = 0 \text{ for } i=0,1,
\end{equation}
where $R\Gamma_{\cS_0} \cO_{Y_0} \simeq k[t_1,\ldots,t_{n-1},t_n^{\pm}] / k[t_1,\ldots,t_n] [-1]$ is the derived subsheaf with supports on $\{t_n=0\}$. Note that this is supported in cohomological degree $1$ and $f_0^\ast \bL_\X$ is supported in cohomological degree $\leq 1$, so the $R\Hom$ complex in \eqref{eqn:local_cohomology} is supported in cohomological degree $\geq 0$.

Now observe that $R\Gamma_{\cS_0} \cO_{Y_0}$ has a filtration whose associated graded is isomorphic to $\bigoplus_{l>0} \cO_{\cS_0} \twist{0,\ldots,0,-l}[-1]$. A spectral sequence argument shows that for the vanishing of \eqref{eqn:local_cohomology} to hold for all $m>0$, it suffices to show that
\[
H^i \left(\bigoplus_{l,m>0} R\Hom_{\cS_0}((f_0^\ast \bL_\X)|_{\cS_0} [1], \cO_{\cS_0} \twist{m,0,\ldots,0,-l}) \right) = 0 \text{ for } i=0,1.
\]
Now observe that $\cS_0 = \Spec(k[t_1,\ldots,t_{n-1}]) / \Gm^{n+1}$ where the first and last copies of $\Gm$ act trivially. This implies that the complex $f_0^\ast \bL_\X|_{\S_0}$ decomposes into a direct sum of weight complexes under the action of $(\Gm)_{t_0} \times (\Gm)_{t_n}$. For the homological vanishing condition above, it is sufficient to check that the $t_0$-weight $-m$ and $t_n$-weight $l$ summand of $f_0^\ast \bL_\X|_{\S_0}$ is supported in cohomological degree $<0$. By Nakayama's lemma, and using the fact that every closed substack of $\cS_0$ meets the origin, it suffices to check that the $t_0$-weight $-m$ and $t_n$-weight $l$ summand of the homology space $H^i (f_0^\ast \bL_\X|_{\{0,\ldots,0\}})$ is zero for $i = 0,1$ and all $l,m>0$.
\end{proof}

\begin{proof}[Proof of \Cref{prop:extension_technical}]
We use \Cref{lem:extension_basic} to argue by induction on $n$. The base case is $n=1$, where the statement of \Cref{prop:extension_technical} is the same as \Cref{lem:extension_basic}. Now let $n>1$ and $V = \{t_n \neq 0\}$ and $U := \{t_1\cdots t_n \neq 0\}$ as defined above. We expand the diagram \eqref{eqn:extension} to the following commutative diagram
\[
\xymatrix{
Y_0 \cap U / \Gm^{n+1} \ar@{}[dr]|*+[o][F-]{A} \ar@{}|-{\subset}[r] \ar@{^{(}->}[d] & Y_0 \cap V / \Gm^{n+1} \ar@{}|-{\subset}[r] \ar@{^{(}->}[d] \ar@{}[dr]|*+[o][F-]{B} & Y_0 / \Gm^{n+1} \ar@/^/[ddr]^{f_0} \ar@{^{(}->}[d] & \\
U / \Gm^{n+1} \ar@{}|-{\subset}[r] \ar@/_/[drrr]_f & V / \Gm^{n+1} \ar@{}|-{\subset}[r] \ar@{-->}[drr]^{f'} & \bA_k^{n+1} / \Gm^{n+1} \ar@{-->}[dr]^{f''} & \\
& & & \X
}
\]
Note that $(\Gm)_{t_n}$ acts freely on $V$, so we have a natural equivalence
\[
\bA_k^n / \Gm^n \simeq \{t_n = 1\} / \Gm^n \simeq V / \Gm^{n+1}.
\]
This equivalence identifies the square $(A)$ above with the square in \eqref{eqn:extension}. Semi-continuity of fiber homology of coherent complexes implies that
\[
\bigoplus_{a_n \in \bZ} H^i((f_0^\ast \bL_\X)|_{\{(0,\ldots,0)\}})_{a_0,\ldots,a_n} = 0 \Rightarrow H^i((f_0^\ast \bL_\X)|_{\{(0,\ldots,0,1)\}})_{a_0,\ldots,a_{n-1}} = 0.
\]
We use this to verify the weight conditions needed to apply the inductive hypothesis to the square $(A)$, which implies the existence and uniqueness of the extension $f'$. Once we have $f'$, we are in the situation of \Cref{lem:extension_basic}, and the appropriate weight hypotheses are satisfied, so we can conclude the existence and uniqueness of $f''$.
\end{proof}

\subsubsection{Proof of \Cref{prop:extension}}

We have seen in the proof of \Cref{lem:ev_0_fans} that for $\sigma \in \Flag^n(f)(k)$, the non-degeneracy of $\sigma$ is equivalent to the non-degeneracy of $\agr(\sigma)$. So the proof of \Cref{prop:extension} amounts to showing that for any $f_0 \in \Deg(\agr(f))^\canon_\bullet$, there is a unique filtration $\op{ext}(f_0) \in \Flag^{n+1}(p)(k)$ such that $v_0(\op{ext}(f_0)) = f$ and $\agr(\op{ext}(f_0)) = f_0$. This follows immediately from \Cref{prop:extension_technical} and the following:

\begin{lem}\label{lem:extension_rephrase}
For any $f_0 \in \Deg(\Grad(\X),\agr(f))^\canon_n$, the resulting map $f_0 : (\pt/\Gm) \times \Theta^n_k = \{t_0=0\} / (\Gm^{n+1})_k \to \X$ satisfies the weight conditions of \Cref{prop:extension_technical}.
\end{lem}

\begin{proof}
We shall prove the contrapositive. Consider the weight decomposition of the representation of $(\Gm^{n+1})_k$
\begin{multline} \label{eqn:cotangent_weight_decomposition}
H^0((\pt/\Gm^{n+1})_k, f_0^\ast \bL_\X|_{\{(0,\ldots,0)\}}) \oplus H^1((\pt/\Gm^{n+1})_k, f_0^\ast \bL_\X|_{\{(0,\ldots,0)\}}) = \\ \bigoplus_{a_0,\ldots,a_n \in \bZ} W_{a_0,\ldots,a_n}.
\end{multline}
The condition in \Cref{prop:extension_technical} is that if $(a_0,\ldots,a_n)$ is a weight for which $W_{a_0,\ldots,a_n} \neq 0$ and $a_0 < 0$, then $a_i \leq 0$ for all $i =1,\ldots,n$.

Let us compute the cotangent weights \eqref{eqn:cotangent_weights} corresponding to a rational ray in the \emph{interior} of the cone $f_0$. Such a ray is determined by a vector with non-zero entries $(r_1,\ldots,r_n) \in \bZ^n_{> 0}$, corresponding to a pointed map $\Theta_k \to \Theta^n_k$ that maps $0 \mapsto 0^n$. If we let $\tilde{W}_{a_0,a_1}$ denote the weight spaces in the decomposition \eqref{eqn:cotangent_weights} corresponding to the composition $\Theta_k \to \Theta^n_k \to \Grad(\X)$, then
\[
\tilde{W}_{b_0,b_1} = \bigoplus_{r_1 a_1 + \cdots r_n a_n = b_1} W_{b_0,a_1,\ldots,a_n}.
\]
If $W_{a_0,\ldots,a_n} \neq 0$ with $a_0 < 0$ and $a_i>0$ for some $i$, then we can choose $r_i \gg r_j$ for $j \neq i$, and we will have $\tilde{W}_{a_0,b_1} \neq 0$ with $a_0<0$ and $b_1 = r_1 a_1 + \cdots + r_n a_n > 0$. So if $f_0$ does not satisfy the weight condition of \Cref{prop:extension_technical}, we can find a $1$-dimensional sub-cone of $f_0$ that does not lie in $\Deg(\agr(f))_1^\canon$.
\end{proof}


\subsection{Local structure of \texorpdfstring{$\iDeg(\X,p)$}{Deg(X,p)}}
\label{sect:local_structure}

Let $0 \subset \cdots \subset \cE_{i+1} \subset \cE_i \subset \cdots \subset \cE$ be a $\bZ$-weighted filtered coherent sheaf on some scheme, and let $\cG := \bigoplus_w \gr_w(\cE_\bullet)$ denote its associated graded coherent sheaf. Let $F_\bullet \cG$ be a filtration of $\cG$ in the category of graded coherent sheaves, which simply means a filtration of each weight summand $\cG_w := \gr_w(\cE_\bullet) \subset \cG$. We make the following assumption on $F$: if $F_i \cG_w \neq 0$, then $F_i \cG_v = \cG_v$ for all $v>w$. Then we may induce a new filtration on $\cE$ by defining:
\[
F'_i \cE := \sum_{w \text{ s.t. } F_i \cG_w \neq 0} \{ \text{preimage of } F_i \cG_w \subset \cE_w / \cE_{w+1} \text{ in } \cE_w \} \subset \cE.
\]
This construction admits a slight generalization in which a $\bZ^n$-weighted filtration of $\cG$ induces a $\bZ^n$-weighted filtration of $\cE$. We regard the filtration $0 \cdots \subset F'_{i+1} \cE \subset F'_i \cE \subset \cdots \subset \cE$ as a ``perturbation" of the original filtration of $\cE$ by the filtration $F$ of the graded coherent sheaf $\gr_\bullet(\cE_\bullet)$. We will justify this terminology in the next section. In this section we describe this construction intrinsically and extend it to filtrations in an arbitrary algebraic stack $\X$.

Observe that the associated graded coherent sheaf $\cG$ above acquires a canonical filtration, whose weight-$i$ subsheaf is
\[
F_i \cG := \bigoplus_{j \geq i} \op{gr}_j (\cE_\bullet) \subset \cG.
\]
The condition placed on the filtration $F_\bullet \cG$ above can be phrased as saying that it is ``close" to the canonical filtration in the degeneration space of $\cG$. In fact, we will identify a neighborhood of the filtration of the original filtration in the degeneration space of $\cE$ with a neighborhood of the canonical filtration in the degeneration space of $\cG$ as a graded coherent sheaf.

To generalize the canonical filtration of a graded coherent sheaf above, consider a stack $\X$ and a point $g \in \Grad(\X)(k)$. Then $g$ has a canonical filtration classified by the composition
\begin{equation} \label{eqn:canonical_comp}
\xymatrix{
\Theta_k \times (\pt/\Gm)_k \ar[r] & (\pt/\Gm)_k \ar[r]^g & \X
},
\end{equation}
where the first map is determined by the group homomorphism $(\Gm)_{t} \times (\Gm)_z \to (\Gm)_z$ given by $(t,z) \mapsto tz$ and the equivariant projection $\bA^1_t \to \pt$. The restriction of this map to $\{1\} \times (\pt/\Gm)_k$ is canonically isomorphic to $g$. We denote the resulting point $\canon \in \Flag(g)(k)$. If $g$ corresponds to a non-degenerate map $(\pt/\Gm)_k \to \X$, then $\canon$ will be a non-degenerate filtration of $g$ in $\Grad(\X)$, so we may regard it as an element of $\Deg(\Grad(\X),g)_1$. By a slight abuse of notation use the same notation for the corresponding rational point $\canon \in \iDeg(\Grad(\X),g)$.

\begin{rem}
We will refer to the fan $\Deg(\Grad(\X),g)^\canon_\bullet$ as the fan of simplices that are \emph{near} $\canon$, justifying the notation of \Cref{defn:near_canon}. Note that all of the cotangent weights $(a_0,a_1)$ of $\canon \in \Deg(\Grad(\X),g)_1$ have $a_1=a_0$, so $\canon \in \Deg(\Grad(\X),g)^\canon_1$. Note also that if $f \in \Deg(\Grad(\X),g)_1$ is antipodal to $\canon$, then all of its cotangent weights $(a_0,a_1)$ have $a_0 = -a_1$, so $\iDeg(\Grad(\X),g)^\canon$ contains no points that are antipodal to $\canon$.
\end{rem}

Informally stated, for a non-degenerate filtration $f$, our main results identify small perturbations of the canonical filtration of $\agr(f)$ with small perturbations of $f$ itself. The comparison uses the map of fans defined as composition $\bT := (\ev_1)_\ast \circ \op{ext}$,
\begin{equation} \label{eqn:define_T}
\xymatrix{\Deg(\Grad(\X),\agr(f))^\canon_\bullet \ar[rr]^{\op{ext}}_{\substack{see \\ \Cref{prop:extension}}} \ar@/^20pt/[rrr]^{\bT} & & \Deg(\Filt(\X),f)_\bullet \ar[r]^{(\ev_1)_\ast} & \Deg(\X,f(1))_\bullet}.
\end{equation}
Concretely, given an $\sigma \in \Deg(\agr(f))^\canon_n$ corresponding to a map $\sigma : (\pt/\Gm)_k \times \Theta_k^n \to \X$, the cone $\bT(\sigma)$ is constructed by filling in the vertical arrows in the following commutative diagram from left to right:
\[
\xymatrix@C=40pt{(\pt/\Gm)_k \times \Theta_k^n \ar@{^{(}->}[r]^-{t_0 = 0} \ar[dr]_{\sigma} & \Theta_k^{n+1} \ar@{-->}[d]_{\op{ext}(\sigma)} & \Theta_k^n \ar@{_{(}->}[l]_-{t_0 = 1} \ar@{-->}[dl]^{\bT(\sigma)} \\ & \X & }
\]
Using this one can see that $\bT(\canon) = f \in \Deg(\X,f(1))_1$. For simplicity we first state our results under the hypothesis that $\X$ has separated flag spaces, which holds when $\X$ has affine diagonal (\Cref{prop:representable_flags}).

\begin{figure}
\includegraphics[height=6cm]{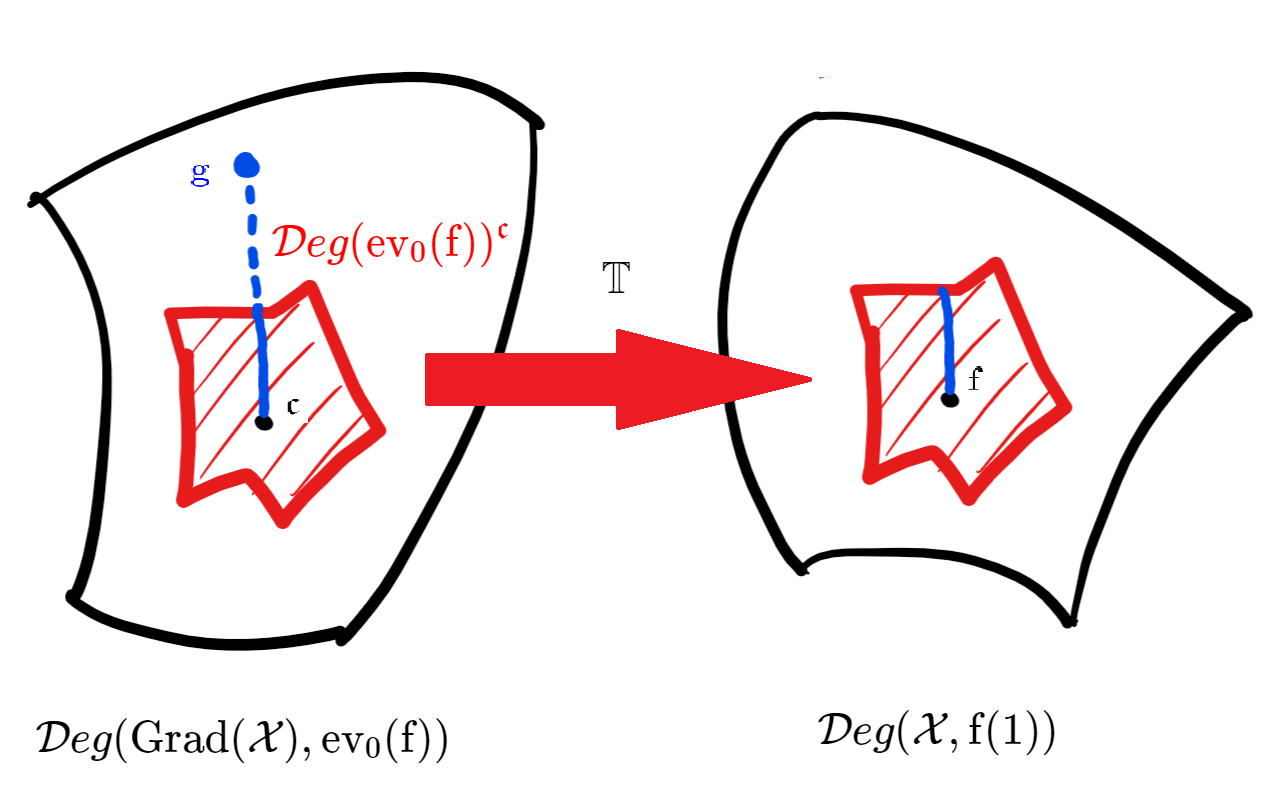}

\caption{Visual summary of \Cref{thm:perturbation}, \Cref{prop:graded_degenerations}, and \Cref{thm:perturbation_light}: For any non-degenerate filtration in $\X$, the map $\bT$ identifies a neighborhood of the canonical filtration of $\agr(f)$ with a neighborhood of the filtration $f$ in $\iDeg(\X,f(1))$. Furthermore, for any point $g \in \iDeg(\Grad(\X,\agr(f)))$ that is not antipodal to $\canon$, there is a canonical rational 1-simplex connecting $g$ with $\canon$ that passes through this neighborhood.}
\end{figure}

\begin{prop} \label{prop:graded_degenerations}
Let $\X$ be an algebraic stack satisfying \ref{hyp3} that has separated flag spaces, and let $g \in \Grad(\X)(k)$. Then
\begin{enumerate}
  \item the inclusion $\Deg(\agr(f))_\bullet^\canon \subset \Deg(\agr(f))_\bullet$ is bounded,
  \item $\canon$ lies in the interior of the closed subset $\iDeg(\agr(f))^\canon \subset \iDeg(\agr(f))$, and
  \item for any rational point $x \in \iDeg(\agr(f))$ that is not antipodal to $\canon$, there is a unique rational $1$-simplex $$\simp{\sigma}^1 \to \iDeg(\agr(f))$$ with $v_0(\sigma) = \canon$ and $v_1(\sigma) = x$. The set of points in $\simp{\sigma}^1$ mapping to $\iDeg(\agr(f))^\canon$ is a closed subinterval with non-empty interior containing the vertex $v_0$.
\end{enumerate}
\end{prop}

\begin{thm}[Perturbation of filtrations] \label{thm:perturbation}
Let $\X$ be an algebraic stack satisfying \ref{hyp3} that has separated flag spaces, and let $f \in \Filt(\X)(k)$ be a non-degenerate filtration. The map of fans $\bT : \Deg(\agr(f))_\bullet^\canon \to \Deg(f(1))_\bullet$ of \eqref{eqn:define_T} is injective and maps $\canon$ to $f$. The subset $\bT(\iDeg(\agr(f))^\canon) \subset \iDeg(f(1))$ contains an open neighborhood of $f \in \iDeg(f(1))$. Furthermore
\begin{enumerate}
\item $\bT$ is ``bounded near $f$'' in the sense that or any rational simplex $\simp{\xi} \to \iDeg(f(1))$ mapping $x \in \simp{\xi}$ to $f$, the preimage of $\bT(\iDeg(\agr(f))^\canon)$ in $\simp{\xi}$ is a finite union of rational simplices; and
\item If $f \in \Filt(\X)(k)$ is contained in a closed substack $\Y \subset \Filt(\X)$ such that $\ev_1 : \Y \to \X$ is quasi-compact, then $\bT$ is a bounded inclusion, hence $\bT(\iDeg(\agr(f))^\canon) \subset \iDeg(f(1))$ is closed (\Cref{cor:fans_closed_embedding}).
\end{enumerate}
\end{thm}

Note that the stronger hypothesis on $\X$ in the second part of \Cref{thm:perturbation} is satisfied when $\X$ has quasi-compact flag spaces (See \Cref{defn:quasi-compact_flags} below). \Cref{prop:graded_degenerations} and \Cref{thm:perturbation} are synopses of the lemmas that we prove in the remainder of this section. The lemmas also lead to a more light-weight version of the statement that does not require $\ev_1$ to be separated. We state it here for reference, as we will use it below:
\begin{thm} \label{thm:perturbation_light}
Let $\X$ be an algebraic stack satisfying \ref{hyp3}, and let $f \in \Filt(\X)(k)$ be a non-degenerate filtration.
\begin{enumerate}
\item For any rational point in $\iDeg(\agr(f))$ that is not antipodal to $\canon$, there is a canonical rational $1$-simplex $\simp{\sigma}^1 \to \iDeg(\agr(f))$ with $v_0(\sigma) = \canon$ and $v_1(\sigma) = x$, and the preimage of $\iDeg(\agr(f))^\canon \subset \iDeg(\agr(f))$ contains an open neighborhood of $v_0$.
\item Conversely, for any rational $1$-simplex $\simp{\sigma}^1 \to \iDeg(f(1))$ with $v_0(\sigma) = f$, there is a rational $1$-simplex $\simp{\sigma'}^1 \to \iDeg(\agr(f))^\canon$ such that $v_0(\sigma') = \canon$ and $\bT(\sigma')$ is a subcone of $\sigma$ containing $v_0$.
\end{enumerate}
\end{thm}

\begin{rem}
Note that \Cref{thm:perturbation} essentially identifies a closed neighborhood of $f \in \Deg(\X,p)$ with non-empty interior such that every point in this neighborhood is connected to $f$ by a canonical $1$-simplex. If $\X$ admits a positive definite class in $H^4(\X;\bR)$, one can canonically parameterize this $1$-simplex by arc length using the spherical metric on $\Deg(\X,p)$. It should be possible to use this to show that this neighborhood of $f$ is contractible, hence $\Deg(\X,p)$ is a locally contractible space.
\end{rem}

\subsubsection{Proof of \Cref{prop:graded_degenerations}: Lemmas on the degeneration fan of a graded point}

\begin{lem} \label{lem:near_canon_bounded}
Let $\X$ be an algebraic stack satisfying \ref{hyp3}, and let $g \in \Grad(\X)(k)$. The inclusion $\Deg(g)^{\canon}_\bullet \subset \Deg(g)_\bullet$ is bounded, so $\iDeg(g)^{\canon} \to \iDeg(g)$ is a closed subspace by \Cref{cor:fans_closed_embedding}.
\end{lem}
\begin{proof}
Given a cone $\xi : h_{[n]} \to \Deg(\agr(f))_n$, we must show that the preimage of $\Deg(g)_\bullet^\canon \subset \Deg(g)_\bullet$ is a bounded sub-fan of $h_{[n]}$. Because cones in $\Deg(g)_\bullet^\canon$ are characterized by their $1$-dimensional sub-cones, it suffices to show that the preimage of $\iDeg(g)^\canon \to \iDeg(g)$ under a rational simplex $\simp{\xi} \to \iDeg(g)$ is a finite union of rational simplices.

Say that $\simp{\xi} \to \iDeg(g)$ corresponds to a map $f : (\pt/\Gm)_k \times \Theta_k^n \to \X$, we consider the representation $\bigoplus_{a_0,\ldots,a_n} W_{a_0,\ldots,a_n}$ of $(\Gm^{n+1})_k$ described in \eqref{eqn:cotangent_weight_decomposition}. Points in $\simp{\xi}$ correspond to $(r_1,\ldots,r_n) \in \bR_{\geq 0}^n - \{0\}$ up to positive scale. One can check that the set of points in the interior of $\simp{\xi}$ that map to $\iDeg(g)^\canon$ are the points $(r_1,\ldots,r_n)$ with all $r_i>0$ for which
\begin{equation} \label{eqn:near_canon_constraints}
r_1 a_1 + \cdots + r_n a_n \leq 0 \quad \forall (a_1,\ldots,a_n) \text{ s.t. } a_0< 0 \text{ for some } W_{a_0,\ldots,a_n} \neq 0.
\end{equation}
Note that these constraints are a finite list of linear inequalities. Rational points in the boundary of $\simp{\xi}$ correspond to maps $\Theta_k \to \Theta_k^n$ that do not map $0$ to $0^n$. However, semicontinuity for the fiber cohomology of the complex $f^\ast \bL_{\X}$ implies that any point in the boundary of $\simp{\xi}$ satisfying the above inequalities still maps to $\iDeg(g)^\canon$.

We have produced a rational polytope in $\simp{\xi}$ mapping to $\iDeg(g)^\canon$ and containing every point in the interior of $\simp{\xi}$ that maps to $\iDeg(g)^\canon$. Applying this argument to the boundary of $\simp{\xi}$ inductively shows that the set of points in $\simp{\xi}$ mapping to $\iDeg(g)^\canon$ is a finite union of rational polytopes, and can thus be covered by finitely many simplices.
\end{proof}

\begin{lem} \label{lem:limits_qsep_fans}
Let $F_\bullet$ be a quasi-separated fan, and let $\{x_i\}_{i=0}^\infty$ be a sequence of points in $\bP(F_\bullet)$ converging to $x \in \bP(F_\bullet)$. Then there is a rational simplex $\simp{\sigma} \to \bP(F_\bullet)$ and a subsequence $\{x_{i_j}\}_{j=1}^\infty$ that lifts to a sequence $\{\tilde{x}_{i_j}\} \subset \simp{\sigma}$ that converges to a lift $\tilde{x} \in \simp{\sigma}$ of $x$.
\end{lem}
\begin{proof}
There must be a rational simplex $\simp{\sigma} \to \bP(F_\bullet)$ whose image contains infinitely many of the $x_i$. If $x = x_i$ for infinitely many $i$, then we chose any rational simplex containing $x$. If not, then passing to a subsequence we may assume that $x \neq x_i$ for any $i$. Now there must be a rational simplex whose image contains infinitely many of the $x_i$, or else the subset $\{x_i\} \subset \bP(F_\bullet)$ would be closed, implying that any limit point appears in the sequence.

Fix such a rational simplex $\simp{\sigma} \to \bP(F_\bullet)$. Select the subsequence of points that lift to $\simp{\sigma}$, and choose lifts. Because $\simp{\xi}$ is compact we can pass to a smaller subsequence that converges to some point $\tilde{x} \in \simp{\xi}$. Because $\bP(F_\bullet)$ is Hausdorff by \Cref{prop:fans_cg}, the image of $\tilde{x}$ under the map $\simp{\xi} \to \bP(F_\bullet)$ must be $x$.
\end{proof}

\begin{lem} \label{lem:graded_deg_interior}
If $\X$ is an algebraic stack satisfying \ref{hyp3}, then $\canon$ lies in the interior of the closed subset $\iDeg(g)^{\canon} \subset \iDeg(g)$.
\end{lem}
\begin{proof}
To show that $\canon$ lies in the interior of $\iDeg(g)^\canon$, it suffices to show that no net of points in the complement of $\iDeg(g)^\canon$ converges to $\canon$. For this, it suffices to show that any sequence in $\iDeg(g)$ converging to $\canon$ must contain points in $\iDeg(g)^\canon$. By \Cref{prop:intersecting_simplices} and \Cref{lem:limits_qsep_fans} we may assume that this sequence lifts to a sequence $\{x_i\}_{i=1}^\infty \subset \simp{\sigma}$ for some rational simplex $\simp{\xi} \to \iDeg(g)$ converging to a point that maps to $\canon$. By subdividing $\simp{\xi}$ we may assume that $v_0(\xi) = \lim_{i \to \infty} x_i$ maps to $\canon$. By definition of the canonical filtration $\canon$, the fact that $v_0(\xi) = \canon \in \iDeg(g)$ implies that in the linear inequalities \eqref{eqn:near_canon_constraints} defining the preimage of $\iDeg(g)^{\canon} \subset \iDeg(g)$ in $\simp{\xi}$, the integer $a_1$ is always some positive multiple of $a_0$. It follows that the constraints are satisfied in some open neighborhood of the point corresponding to $(r_1,\ldots,r_n) = (1,0,\ldots,0)$. Thus the sequence $\{x_i\}$, which converges to $(1,0,\ldots,0)$, must map to $\iDeg(g)^\canon$ for $i$ sufficiently large.
\end{proof}

Finally, we show that both the degeneration space $\iDeg(g)$ and $\iDeg(g)^\canon$ are ``star shaped" with respect to the canonical point $\canon$.

\begin{lem} \label{lem:star_shaped}
If $\X \to B$ is a morphism of algebraic stacks with separated relative inertia $I_{\X/B} \to \X$, then for any rational point $x \in \iDeg(g)$ that is not antipodal to $\canon$, there is a canonical rational $1$-simplex $\simp{\sigma}^1 \to \iDeg(g)$ with $v_0(\sigma) = \canon$ and $v_1(\sigma) = x$. The set of points in $\simp{\sigma}^1$ mapping to $\iDeg(g)^\canon$ is a closed subinterval with non-empty interior containing the vertex $v_0$. If furthermore $\X \to B$ satisfies \ref{hyp2} and has separated flag spaces, then $\sigma$ is unique.
\end{lem}
\begin{proof}
The point $x$ is represented by some map $\nu : (\pt/\Gm)_k \times \Theta_k \to \X$. We can compose this with the map $(\pt/\Gm)_k \times \Theta_k^2 = (\Theta_k \times (\pt/\Gm)_k) \times \Theta_k \to (\pt/\Gm)_k \times \Theta_k$ that is the identity on the first factor and the canonical map \eqref{eqn:canonical_comp} on the second factor. It follows from \Cref{lem:contract_degenerate} below that if the resulting map $\Theta^2_k \to \Grad(\X)$ is non-degenerate, then either $v_0(\sigma),v_1(\sigma) \in \Deg(\Grad(\X),g)_1$ are positive multiples of one another, or they are antipodal, so under the hypothesis of the lemma $\sigma$ is non-degenerate. If $\X$ satisfies \ref{hyp2} and has separated flag spaces, then $\sigma$ is uniquely determined by $v_0$ and $v_1$ by \Cref{lem:unique_fillings}.

Finally, we know from the proof of \Cref{lem:near_canon_bounded}, that the condition for a point in the interior of $\simp{\sigma}$ to map to $\iDeg(g)^\canon$ is given by a finite list of linear inequalities that hold strictly at $v_0$ and thus hold in a neighborhood of $v_0$. These inequalities define a closed interval with non-empty interior containing $v_0$. To complete the proof we must argue that if $x \in \iDeg(g)^\canon$, then all of $\sigma$ must map to $\iDeg(g)^\canon$ as well. Say $x$ is represented by $\nu$ as above, and consider a point in the interior of $\simp{\sigma}^1$, represented by a nonvanishing pair $(r_0,r_1) \in \bZ^2_{>0}$ which defines a map
\[
\nu' : \Theta_k \times (\pt/\Gm)_k \xrightarrow{(r_0,r_1)} \Theta_k^2 \times (\pt/\Gm)_k \xrightarrow{\sigma} \X.
\]
We can directly compare the two representations $\nu^\ast \bL_\X|_{(0,0)}$ and $(\nu')^\ast \bL_\X|_{(0,0)}$ of $(\Gm)^2_k$: the latter is the pullback of the former along the homomorphism $\Gm^2 \to \Gm^2$ given by
\[
(z_0,z_1) \mapsto (z_0 z_1^{r_0}, z_1^{r_1}).
\]
It follows that $(a_0,a_1)$ is a cotangent weight of $\nu$ if and only if $(a_0, a_0 r_0 + a_1 r_1)$ is a cotangent weight of $\nu'$. In particular if $\nu \in \Deg(g)^\canon_1$, then any cotangent weight of $\nu'$ with $a_0<0$ is of the form $(a_0, a_0 r_0 + a_1 r_1)$ with $a_1 \leq 0$, so $\nu' \in \iDeg(g)^\canon_1$ as well.
\end{proof}

The following lemma was used in the previous proof.

\begin{lem} \label{lem:contract_degenerate}
Let $\X$ be an algebraic stack, and let $f : \Theta_k^2 \to \X$ be a morphism over a $k$-point of $B$ such that $\ker ((\Gm)_k^2 \to \Aut_\X f(0))$ is positive dimensional. Then one of the following holds:
\begin{enumerate}
\item the two pointed maps $\Theta_k \to \X$ induced by restricting $f$ to $(\bA_k^1 \setminus 0) \times \bA_k^1 / \Gm^2 \cong \Theta_k$ and to $\bA_k^1 \times (\bA_k^1 \setminus \{0\}) / \Gm^2 \cong \Theta_k$ are antipodal, or
\item $f$ factors uniquely through a morphism $\pi : \Theta_k^2 \to \Theta_k$ given in coordinates by $(z_1,z_2) \mapsto z_1^a z_2^b$ for some $a,b \geq 0$.
\end{enumerate}

\end{lem}

\begin{proof}
The claim only depends on $\X \times_B \Spec(k)$, so we may assume $B= \Spec(k)$. As $G = \ker((\Gm^2)_k \to \Aut_\X f(0))$ is a positive dimensional subgroup, it contains a subtorus $G_0$ of rank $1$. $G_0$ has the form $\{(t^a,t^b) | t \in \Gm\}$ for some pair of coprime integers $a,b$ with $a\geq 0$. We prove the claim in three cases:

\medskip
\noindent \textit{Case 1:} $a b<0$
\medskip

In this case $b<0$. Consider the morphism $\pi : \bA_k^2 \to \bA^1_k$ taking $(z_1,z_2) \mapsto z_1^{-b}z_2^a$. This identifies $\bA^1_k$ with the good quotient $\bA^2_k /\!/ G_0 = \Spec(k[z_1,z_2]^{G_0})$. Note that $\pi$ is equivariant for the group homomorphism $\Gm^2 \to \Gm$ given by the same formula and thus induces a map $\pi : \Theta_k^2 \to \Theta_k$ that is a relative good moduli space morphism. The universal property for good moduli space morphisms \cite{ahr2}*{Thm.~17.2} implies that $f$ factors uniquely through $\pi$, because $G_0 = \ker((\Gm^2)_k \to \Aut_{\Theta_k}(\pi(0)))$ is contained in $\ker((\Gm^2)_k \to \Aut_\X(f(0)))$ by construction.

\medskip
\noindent \textit{Case 2:} $a b > 0$
\medskip

In this case $b>0$, and the good quotient $\bA^2_k /\!/ G_0 \cong \Spec(k)$. So the morphism $\pi : \Theta_k^2 \to B(\bG_m)_k$ induced by the projection $\bA^2_k \to \Spec(k)$ and the homomorphism $\phi : \Gm^2 \to \Gm$ given by $(z_1,z_2) \mapsto z_1^{-b} z_2^a$ is a relative good moduli space morphism. Note that $\phi$ maps the two factors of $\Gm^2$ to opposite cocharacters of $\Gm$ (up to positive multiple). Again the universal property \cite{ahr2}*{Thm.~17.2} implies that $f$ factors uniquely through $\pi$, hence the two $k$-points of $\Filt(\X)$ determined by $f|_{\bA^2-\{0\} / \Gm^2}$ are antipodal.

\medskip
\noindent \textit{Case 3:} $a b=0$
\medskip

We may assume, by permuting coordinates, that $a = 0$. Then let $\pi : \Theta_k^2 \to \Theta_k$ be the projection onto the first factor, given in coordinates by $(z_1,z_2) \mapsto z_1$. As above, $\pi$ is a good moduli space morphism, and the universal property implies that $f$ factors uniquely through $\pi$.

\end{proof}

\subsubsection{Proof of \Cref{thm:perturbation}}

\begin{lem} \label{lem:ext_bounded}
If $\X$ is a stack satisfying \ref{hyp3}, then the map
\[
\op{ext} :\Deg(\agr(f))^\canon_\bullet \to \Deg(f)_\bullet
\]
is bounded, and $\op{ext}(\canon)$ lies in the interior of the closed subset $\iDeg(\agr(f))^\canon \subset \iDeg(f)$.
\end{lem}
\begin{proof}
Consider a cone $\xi : h_{[n]} \to \Deg(f)_\bullet$, then \Cref{prop:extension} implies that a subcone $\xi' : h_{[k]} \to h_{[n]} \to \Deg(f)_\bullet$ lies in the image of $\op{ext}$ if and only if $\agr(\xi') \in \Deg(\agr(f))^\canon_n$. Hence the preimage of $\Deg(\agr(f))_\bullet^\canon \subset \Deg(f)_\bullet$ is bounded by \Cref{lem:near_canon_bounded}. Similarly, image of $\op{ext}$ is precisely the preimage of $\Deg(\agr(f))_\bullet^\canon \subset \Deg(f)_\bullet$, so $\op{ext}(\canon)$ lies in the interior of the image of $\iDeg(\agr(f))^\canon$ by \Cref{lem:graded_deg_interior}.
\end{proof}

\begin{lem} \label{lem:canonical_lift}
For any $B$-stack $\X$ with representable and separated relative inertia, any $p \in \X(k)$, and any cone $\sigma \in \Deg(\X,p)_n$ with $v_0(\sigma) = f \in \Deg(\X,p)_1$, there is a canonical cone $\sigma' \in \Deg(\Filt(\X),f)_n$ with $(\ev_1)_\ast(\sigma') = \sigma$.
\end{lem}
\begin{proof}
Consider the product map $p:\Theta_k \times \Theta_k \to \Theta_k$ given in coordinates by $(t_0,t_1) \mapsto t_0 t_1$. The composition
\[
\Theta_k^{n+1} \xrightarrow{p \times \id_{\Theta_k^{n-1}}} \Theta_k^n \xrightarrow{\sigma} \X
\]
can be regarded as a non-degenerate $n$-dimensional filtration of $f \in \Filt(\X)(k)$ lifting $\sigma \in \Deg(\X,f(1))_n$.
\end{proof}

The following should be interpreted as saying that the map of fans $\bT$ is ``bounded near $f$''

\begin{lem} \label{lem:T_almost_bounded}
Let $\X$ be a stack satisfying \ref{hyp3} with separated flag spaces. Then $\bT$ is injective, and $f$ is an interior point of the subset $\bT(\iDeg(\agr(f))^\canon) \subset \iDeg(f(1))$. For any rational simplex $\xi : \simp{\xi} \to \iDeg(f(1))$ mapping $x \in \simp{\xi}$ to $f$, the preimage $\xi^{-1}(\bT(\iDeg(\agr(f))^\canon)) \subset \simp{\xi}$ is a finite union of sub-simplices.
\end{lem}
\begin{proof}
We have seen in \Cref{prop:representable_maps} that representable separated maps of stacks induce injective maps on degeneration fans, so $(\ev_1)_\ast : \Deg(\Filt(\X),f)_\bullet \to \Deg(\X,f(1))_\bullet$ is injective. $\op{ext}$ is injective because $\agr \circ \op{ext}$ is injective, and hence $\bT = (\ev_1)_\ast \circ \op{ext}$ is injective. First we show that for any rational simplex $\xi : \simp{\xi} \to \iDeg(f(1))$ mapping $x \in \simp{\xi}$ to $f$, the preimage $\xi^{-1}(\bT(\iDeg(\agr(f))^\canon)) \subset \simp{\xi}$ is a finite union of sub-simplices that contains an open neighborhood of $x$.

By barycentric subdivision of $\simp{\xi}$ centered at $x$, one can reduce to proving the claim when $x$ is the $0^{th}$ vertex of $\simp{\xi}$, so we may assume that for our $\xi \in \Deg(\X,f(1))_n$ we have $v_0(\xi) = f^m$ for some $m>0$. \Cref{prop:N_action} implies that the $m^{th}$-power map $\Filt(\X) \to \Filt(\X)$ is an isomorphism on connected components, so if $\Theta_k^{n-1} \to \Filt(\X)$ is a map such that the image of the generic point admits an $m^{th}$ root, then the entire map admits and $m^{th}$ root. The resulting filtration $\Theta_k^n \to \X$ represents the same simplex in $\iDeg(\X,f(1))$ set theoretically. We may therefore assume that $v_0(\xi)=f$.

In this case $\xi$ lifts to $\Deg(f)_\bullet$ by \Cref{lem:canonical_lift}. Injectivity of $\ev_1$ implies that for any $\xi : h_{[n]} \to \Deg(f(1))_\bullet$ that lies in the image of $(\ev_1)_\ast : \Deg(f)_\bullet \to \Deg(f(1))_\bullet$, we have
\[
h_{[n]} \times_{\Deg(\X,f(1))_\bullet} \Deg(\agr(f))_\bullet = h_{[n]} \times_{\Deg(\Filt(\X),f)_\bullet} \Deg(\agr(f))_\bullet.
\]
\Cref{lem:ext_bounded} therefore implies that $\xi^{-1}(\iDeg(\agr(f))^\canon) \subset \simp{\xi}$ is a finite union of rational sub-simplices that contains $x$ in its interior.

If $f \in \bT(\iDeg(\agr(f))^\canon)$ were not an interior point, then one could use \Cref{lem:limits_qsep_fans} to find a sequence of $x_i$ in some rational simplex $\simp{\xi} \to \iDeg(f(1))$ that converge to a point $x \in \simp{\xi}$ such that $x$ maps to $f$, but no $x_i$ maps to $\bT(\iDeg(\agr(f))^\canon)$. This contradicts the fact that $\xi^{-1}(\bT(\iDeg(\agr(f))^\canon))$ contains an open neighborhood of $x$. Therefore the subset $\bT(\iDeg(\agr(f))^\canon)$ must contain an open neighborhood of $f$.
\end{proof}

Under stronger hypotheses still, we can prove that $\bT$ is bounded.

\begin{lem}
If $\X$ satisfies \ref{hyp3} and has separated flag spaces, and if $f$ is contained in a closed substack $\Y \subset \Filt(\X)$ for which $\ev_1 : \Y \to \X$ is quasi-compact, then the injective map $\bT$ is bounded, and $f$ lies in the interior of the closed subset $\iDeg(\agr(f))^\canon \subset \iDeg(\X,f(1))$.
\end{lem}
\begin{proof}
We know that the map $\ev_1 : \Deg(\Filt(\X),f)_\bullet \to \Deg(\X,f(1))_\bullet$ is bounded by \Cref{prop:representable_maps}, and $\op{ext}$ is bounded by \Cref{lem:ext_bounded}.
\end{proof}


\subsection{The component space \texorpdfstring{$\iComp(\X)$}{Comp(X)}}

We introduce a topological space that is much smaller than $|\Deg(\X,p)_\bullet|$ and plays a key role in our construction of numerical invariants below.

\begin{defn} \label{defn:reduced_degeneration_space}
Let $\X$ be a stack satisfying \ref{hyp2}, and let $\xi : S \to \X$ be a map from a scheme. We say that a connected component of $\Filt^n(\X)$ is \emph{non-degenerate} if it contains a non-degenerate point (note that if $I_{\X/B} \to \X$ is separated, then \Cref{prop:components} implies every point in the component is non-degenerate). We define a formal fan $\Comp(\xi)_\bullet$ and the \emph{\gls{component_fan}} $\Comp(\X)_\bullet$ by
\begin{gather*}
\Comp(\X)_n := \left\{\text{non-degenerate } \alpha \in \pi_0 \Filt^n(\X) \right\}, \text{and} \\
\Comp(\xi)_n := \left\{ \alpha \in \pi_0 (\Flag^n(\xi)) \text{ whose image in } \pi_0 \Filt^n(\X) \text{ is non-degenerate} \right\}.
\end{gather*}
We define the \emph{\gls{component_space}} of $\X$ (respectively $\xi$) to be $\iComp(\X):= \bP(\Comp(\X)_\bullet)$ (respectively $\iComp(\xi) = \bP(\Comp(\xi)_\bullet)$).
\end{defn}

Note that the rational points of $\iComp(\X)$ correspond exactly to the set $\pi_0 (\Filt(\X)) / \bN^\times$.

\begin{rem}
For an arbitrary stack $\Y$, one can define the set of connected components $\pi_0(\Y)$ as the set of points $|\Y|$ modulo the smallest equivalence relation identifying any two points in the image of a map $S \to \Y$ where $S$ is a connected scheme. Using this one can define $\Comp(\X)_\bullet$ and $\Comp(\xi)_\bullet$ even for stacks over $B$ that are not algebraic.
\end{rem}

For any scheme $S$ and $S$-point $\xi : S \to \X$, the map $\Flag^n(\xi) \to \Filt^n(\X)$ defines a canonical map of fans $\Comp(\xi)_\bullet \to \Comp(\X)_\bullet$. Furthermore for $p \in \X(k)$ we can define a map of fans $\Deg(\X,p)_\bullet \to \Comp(p)_\bullet$ by assigning a $k$-point of $\Flag^n(p)$ to its connected component in $\pi_0(\Flag^n(p))$.\footnote{Note that the group of $k$-points of $\Aut(p)$ also acts on $\Deg(\X,p)_\bullet$ and if $G$ is connected, then the morphism $\Deg(\X,p)_\bullet \to \Comp(p)_\bullet$ factors through $\Deg(\X,p)_\bullet / \Aut(p)$. One can define a fan that is still coarser than $\Comp(\xi)_\bullet$ whose set of $n$-cones is $\op{im}(\Comp(\X,\xi)_n \to \pi_0 \Filt^n(\X))$, but $\Comp(\xi)_\bullet$ will suffice for our purposes.}

\subsubsection{The component space of a quotient stack}

In this subsection $G$ will denote either a split reductive group over a field $k$ or $G=\GL_N$ over $\bZ$, $T \subset G$ will denote a split maximal torus, and $X$ will denote a $G$-quasi-projective scheme over a base scheme $B$. We shall compute the component fan $\Comp(X/G)_\bullet$ by considering the map $X/T \to X/G$. Using the description in \Cref{thm:describe_strata_global_quotient} of the stack of filtered objects in $X/G$ and $X/T$, we can identify the map $\Filt^n(X/T) \to \Filt^n(X/G)$ with the canonical surjective map
\[
\bigsqcup_{\psi \in \Hom(\Gm^n,T)} X^{\psi,+} / T \to \bigsqcup_{\psi \in \Hom(\Gm^n,T) / W} X^{\psi,+} / P_\psi
\]
The Weyl group $W$ acts naturally on $\Filt^n(X/T)$ by the canonical identification $w \cdot X^{\psi,+} \simeq X^{w \psi w^{-1},+}$ for $w \in W$, and the map $\Filt^n(X/T) \to \Filt^n(X/G)$ is invariant with respect to this action.

\begin{lem} \label{lem:component_fan_quotient}
If $G$ is geometrically connected, then the $W$-invariant map $\Comp(X/T)_\bullet \to \Comp(X/G)_\bullet$ induces an isomorphism
\[
\Comp(X/T)_\bullet / W \simeq \Comp(X/G)_\bullet.
\]
\end{lem}
\begin{proof}
Observe that $X^{\psi,+} / T \to X^{\psi,+} / P_\psi$ induces a bijection on connected components. This follows from the fact that $P_\psi$ is connected, so the set of connected components of both stacks corresponds bijectively to the set of connected components of $X^{\psi,+}$. The components of $\Filt^n(X/G)$ and $\Filt^n(X/T)$ that are non-degenerate correspond to components of $X^{\psi,+}$ for $\psi$ that have a finite kernel, so the equivalence $\pi_0(\Filt^n(X/G)) \simeq \pi_0(\Filt^n(X/T)) / W$ preserves the non-degenerate components and gives an equivalence of component fans.
\end{proof}

Because projective geometric realization of fans commutes with colimits, and the quotient is a colimit, it follows that $$\iComp(X/G) \simeq \iComp(X/T) / W$$ as well. Consider the set $\{(T_i \subset T, Z_i \subset X^{T_i})\}_{i \in I}$ of all sub-tori that arise as the reduced identity component of the stabilizer of some point in $X$ along with a choice of connected component $Z_i$ of the fixed locus $X^{T_i}$. The index set $I$ is finite, as can be shown by reducing to the case of a linear action of $T$ on projective space. We therefore have
\begin{lem} \label{lem:bounded_comp_quotient}
The fans $\Comp(X/T)_\bullet$ and $\Comp(X/G)_\bullet$ are bounded.
\end{lem}
\begin{proof}
By \Cref{lem:component_fan_quotient} it suffices to prove this for $X/T$. For each $i$, choose a finite type point in $Z_i(k_i)$. Consider the map of stacks
\[
\X' := \bigsqcup_i \Spec(k_i) / T_i \to \bigsqcup_i Z_i / T_i \to X/T.
\]
\Cref{lem:graded_stack} implies that $\Grad^n(\X') \to \Grad^n(X/T)$ is surjective for all $n$, and \Cref{lem:retract} implies that $\Grad^n(-)$ and $\Filt^n(-)$ have the same connected components. It follows that $\Comp(\X')_\bullet \to \Comp(X/T)_\bullet$ is surjective. The claim now follows from the observation that $\Comp(\Spec(k_i)/T_i)_\bullet \simeq \Deg(\Spec(k_i)/T_i,\pt)_\bullet$ is bounded (see \Cref{lem:degeneration_fan_torus}).
\end{proof}

We can give an even more concise description of $\iComp(X/T)$. Partially order the set $\{(T_i,Z_i)\}_{i\in I}$ above by the rule
\[
(T_i,Z_i) \prec (T_j, Z_j) \text{ if } T_i \subset T_j \text{ and } Z_j \subset Z_i.
\]
Let $N = \Hom(\Gm,T)$ be the cocharacter lattice. For each $i$ we consider the subgroup $N_i := \Hom(\Gm,T_i) \subset N$. Note that with our indexing convention the $N_i$ are not distinct as subgroups of $N$, but we regard them as distinct abstract groups. For any $i \prec j$ we have a canonical embedding $N_i \subset N_j$. Then $|\Comp(X/T)_\bullet|$ is a union of the vector spaces $(N_i)_\bR$ along the resulting inclusions $(N_i)_\bR \subset (N_j)_\bR$. The proof of \Cref{lem:component_fan_quotient} actually implies the following
\begin{lem} \label{lem:describe_component_fan}
$|\Comp(X/T)_\bullet| \simeq \colim \limits_{i \in I} (N_i)_\bR$ and $\iComp(X/T) \simeq \colim \limits_{i \in I} ((N_i)_\bR - \{0\}) / \bR_{>0}^\times$.
\end{lem}
Thus $|\Comp(X/T)_\bullet|$ is a finite union of real vector spaces along linear embeddings, and $|\Comp(X/G)_\bullet|$ is the quotient of this space by an action of $W$ that permutes indices and acts linearly on each vector space.

\begin{rem}
Note that the inclusions of vector spaces come from an inclusion of lattices $N_i \hookrightarrow N_j$. The image of these lattices in $\iComp(X/T) = \bigcup (N_i)_\bR$ is precisely the set of integral points, and the action of $W$ preserves integral points. The set of non-zero integral points of $\iComp(X/G)$, which is in bijection with the set of non-degenerate connected components of $\Filt(X/G)$, is the union of the image of the maps $N_i\setminus \{0\} \to \iComp(X/G)$.
\end{rem}

\subsubsection{The component space of a quasi-compact stack}

\begin{lem} \label{lem:stratified_comp}
Let $\X = \bigcup \X_i$ be a set-theoretic disjoint union of locally closed substacks. Then the map of fans $\bigsqcup_i \Comp(\X_i)_\bullet \to \Comp(\X)_\bullet$ is surjective.
\end{lem}
\begin{proof}
We know from \Cref{lem:retract} that $\Grad^n(\X) \to \Filt^n(\X)$ induces a bijection on connected components, so it suffices to show that for any $n$, $\bigsqcup_i \Grad^n(\X_i) \to \Grad^n(\X)$ induces a surjection on connected components. This follows from \Cref{cor:closed_open_graded}, which implies that this map is universally bijective.
\end{proof}

\begin{cor} \label{cor:bounded_comp_1}
If $\X$ is quasi-compact satisfying \ref{hyp2}, then the fan $\Comp(\X)_\bullet$ is bounded, and hence $\iComp(\X)$ is compact.
\end{cor}
\begin{proof}
Because $\X$ is quasi-compact, we can find a smooth map $\Spec(R) \to B$ such that the base change $\X_R \to \X$ is surjective. By \Cref{cor:base_change_grad} the maps $\Grad^n_R(\X_R) \to \Grad^n_B(\X)$ are surjective and hence surjective on $\pi_0$. It follows that the map $\Comp(\X_R)_\bullet \to \Comp(\X)_\bullet$, where the former is computed relative to $\Spec(R)$ and the latter relative to $B$, is surjective. $\X_R$ is quasi-compact because $B$ is quasi-separated, so we may therefore assume that $B = \Spec(R)$. We can write any noetherian stack with affine stabilizers as a set-theoretic union of locally closed substacks that are quotient stacks \cite{kresch1999cycle}*{Prop.~3.5.9}\cite{hall2014coherent}*{Prop.~8.2}. The component fan of a quotient stack over the base $\Spec(R)$ is bounded, by the computations above, so it follows from \Cref{lem:stratified_comp} that $\Comp(\X)_\bullet$ is bounded.
\end{proof}

\begin{rem}
After developing some more sophisticated methods, we will show in \Cref{lem:qc_bounded_relative_comp} that for any quasi-compact stack $\X$ satisfying \ref{hyp2} and any map $\xi : S \to \X$ from a quasi-compact space $S$, the fan $\Comp(\X,\xi)_\bullet$ is bounded. 
\end{rem}


\subsection{Cohomology and functions on the component space}

We now discuss a method of constructing continuous functions on $\iDeg(\X,p)$ and $|\Deg(\X,p)_\bullet|$ from cohomology classes on the stack $\X$ with values in some fixed coefficient ring $A \subset \bR$. In order for our framework and results to be as flexible as possible -- for instance to work over fields other than $\bC$ -- we will axiomatize the properties we need. These properties hold for many common cohomology theories. We use $H^\ast$ to denote any contravariant functor from some subcategory of the homotopy category of $B$-stacks,\footnote{By homotopy category we mean the category whose objects are stacks and whose maps are $1$-morphisms of stacks up to $2$-isomorphism.} that must at least contain the stacks $\Theta^n_S$ for any finitely presented $B$-scheme $S$, to graded $A$-modules that satisfies the following axioms.
\begin{enumerate}
\item For fields $k$ of finite type over $B$, there is a canonical isomorphism $H^\ast(\Theta_k^n) \simeq A[u_1,\ldots,u_n]$ with $u_i$ in cohomological degree $2$. We regard the ring $A[u_1,\ldots,u_n]$ as polynomial $A$-valued functions on $(\bR_{\geq 0})^n$.
\item For any $\phi : [m] \to [n]$ in $\Cones$, the restriction homomorphism $H^\ast(\Theta_k^n) \to H^\ast(\Theta_k^m)$ induced by the morphism $\Theta_k^m \to \Theta_k^n$ agrees with the restriction of polynomial functions along the corresponding inclusion $(\bR_{\geq 0})^m \subset (\bR_{\geq 0})^n$.
\item For any integral affine finitely presented $B$-scheme $S$, the composed homomorphism $H^\ast(\Theta_S^n) \to H^\ast(\Theta^n_{k(s)}) \simeq A[u_1,\ldots,u_n]$ is independent of the finite type point $s \in S$.
\end{enumerate}

We regard the choice of cohomology theory $H^\ast(-)$ satisfying (1)-(3) as fixed throughout this paper. We sometimes use the phrase \emph{rational cohomology classes} to refer to classes in a cohomology theory with coefficient ring $A=\bQ$.

\begin{ex}
When $B$ is locally of finite type over $\bC$, we may discuss the topological stack underlying the analytification of any $\X$ satisfying \ref{hyp2}. This topological stack is defined by taking a presentation of $\X$ by a groupoid in schemes and then taking the analytification, which is a groupoid in topological spaces. The cohomology is then defined as the cohomology of the classifying space of this topological stack \cite{No12}. For global quotient stacks $\X = X/G$ this agrees with the equivariant cohomology $H_{G^{an}}^\ast (X^{an};A)$, which agrees with $H_K^\ast (X^{an};A)$ when $K \subset G$ is a maximal compact subgroup weakly homotopy equivalent to $G$. The properties (1)-(3) are well known in this setting.
\end{ex}

\begin{rem}
The stack $\Theta^n = \bA^n / \Gm^n$ deformation retracts onto $\pt / \Gm^n$, so standard computations in equivariant cohomology show that $H^\ast(\bA^n/\Gm^n) \simeq A[u_1,\ldots,u_n]$, where $u_i,\ldots,u_n \in H^2$. Note that while $H^\ast (\pt / \Gm^n) \simeq A[u_1,\ldots,u_n]$ as well, the generators in $H^2$ are only canonical up to the action of $GL_n(\bZ)$ via automorphisms of $\pt / \Gm^n$. Automorphisms of $\Theta^n$ correspond to $M \in \GL_n(\bZ)$ such that $M$ and $M^{-1}$ both fix $(\bR_{\geq 0})^n$, which implies that $M$ is a permutation matrix. It follows that the generators of $H^2(\Theta^n)$ are canonical up to permutation. We encode this distinction in (1) by regarding $A[u_1,\ldots,u_n]$ as functions on $(\bR_{\geq 0})^n$ rather than $\bR^n$.
\end{rem}

\begin{ex}
When $B$ is locally finite type over some other field, one can take $H^\ast$ to be operational equivariant Chow cohomology \cites{edidin1998equivariant, kresch1999cycle}. For properties (1) and (2) see \cite{Pa06}. Applying \cite{gonzales2015equivariant}*{Thm.~3.5} to the special $(\Gm)_k^n$-linear scheme $\bA_k^n$ gives a version of the K\"unneth formula $H^\ast(\Theta_S^n) \simeq H^\ast(S) \otimes H^\ast(\Theta_k^n)$ for any finite type $B$-scheme $S$. This combined with the fact that $H^\ast(\Spec(k)) = A$ for any field $k$ of finite type over $B$ implies (3).
\end{ex}

\begin{lem} \label{lem:cohomology_functions}
Let $\X$ be a stack and let $\eta \in H^{2l}(\X)$. Then there is a unique continuous function $\hat{\eta} : |\Comp(\X)_\bullet| \to \bR$ defined by the property that for any $\xi \in \Comp(\X)_n$ the restriction of $\hat{\eta}$ along the map $\xi : \bR^n_{\geq 0} \to |\Comp(\X)_\bullet|$ is the homogeneous degree $l$ function $f^\ast (\eta) \in H^\ast(\Theta^n_k) \simeq A[u_1,\ldots,u_n]$, where $f : \Theta_k^n \to \X$ is a finite type point of $\Filt^n(\X)$ lying on the connected component corresponding to $\xi$.
\end{lem}
\begin{proof}
Property (3) guarantees that $f^\ast (\eta) \in A[u_1,\ldots,u_n]$ only depends on the connected component on which $f \in \Filt^n(\X)$ lies. Consequently $\eta \in H^{2l}(\X)$ defines a function $\pi_0 \Filt^n(\X) \to A[u_1,\ldots,u_n]_{\rm{deg}-l}$, which we regard as a space of real valued degree $l$ polynomials on $\bR^n_{\geq 0}$. The geometric realization is a colimit, so a continuous function $|F| \to \bR$ is defined by a family of continuous functions $(\bR_{\geq 0})^n \to \bR$ for each $\xi \in F_n$ that is compatible with the continuous maps $(\bR_{\geq 0})^k \to (\bR_{\geq 0})^n$ for each morphism in $(\Cones | F)$. Properties (1) and (2) give exactly such a family of continuous functions on $(\bR_{\geq 0})^n$.
\end{proof}

For a stack $\X$ with representable separated inertia and $p \in \X(k)$, we will overload notation by using $\hat{\eta}$ to denote the function on $|\Deg(\X,p)|$ obtained by restricting $\hat{\eta}$ along the canonical map $|\Deg(\X,p)| \to |\Comp(\X)|$. A filtration $f : \Theta_k \to \X$ defines an integral point in both $|\Deg(\X,f(1))|$ and $|\Comp(\X)|$, the image of $1 \in \bR_{\geq 0}$ under the rational ray corresponding to $f$,
\[
\bR_{\geq 0} \to |\Deg(\X,f(1))| \to |\Comp(\X)|.
\]
We denote this point by $f \in |\Deg(\X,f(1))|$ as well. Given a class $\eta \in H^{2l}(\X)$, we can evaluate the function $\hat{\eta}$ at this point,
\begin{equation} \label{eq:cohomology_function}
\hat{\eta}(f) = \frac{1}{q^l} f^\ast(\eta) \in A,
\end{equation}
where $q^l$ is the canonical generator of the cohomology $H^\ast(\Theta_k;A) \simeq A[\![q]\!]$.

\begin{rem}
Operational Chow cohomology should satisfy property (3) in mixed characteristic as well, but the literature on operational Chow cohomology of stacks in mixed characteristic is not as developed. In practice, however, we will always be considering Chern classes of locally free sheaves (or complexes) on $\X$. For these classes, one can check property (3) directly by relating the Chern classes of a locally free sheaf on $\Theta_k^n$ to the weights of the sheaf restricted to $(\pt/\Gm^n)_k$, so the conclusion of \Cref{lem:cohomology_functions} still applies.
\end{rem}

\subsubsection{Positive definite classes}

One useful notion we will make frequent use of is the following
\begin{defn} \label{defn:positive_definite}
We say that a class $b \in H^4(\X)$ is \emph{positive definite} if for all non-degenerate maps $g : \pt/\Gm \to \X$, the pullback $\gamma^\ast(b) \in H^4(\pt/\Gm) \simeq A \cdot u^2 \subset \bR[u]$ is a positive multiple of $u^2$.
\end{defn}
Note that the generator of $H^2(\pt/\Gm)$ is only canonical up to sign, but $u^2$ is a canonical generator for $H^4(\pt/\Gm)$, so the notion of sign is well-defined. We leave the proof of the following to the reader:

\begin{lem}
The following are equivalent:
\begin{enumerate}
  \item $b \in H^4(\X)$ is positive definite,
  \item the function $\hat{b} : |\Comp(\X)_\bullet| \to \bR$ is positive away from the cone point, and
  \item for any $p \in \X(k)$, the function $\hat{b} : |\Deg(\X,p)_\bullet| \to \bR$ is positive away from the cone point.
\end{enumerate}
\end{lem}

Note that $b$ is positive definite if and only if for every non-degenerate map $g : (\pt/\Gm)^n \to \X$ the resulting quadratic form $g^\ast(b) \in \op{Sym}^2(\bR^n)$ is positive definite. In fact, the previous lemma shows that if the fan $\Comp(\X)_\bullet$ is bounded, then there is a finite set of non-degenerate morphisms
$$g_i : (pt / \Gm)^{n_i} \to \X, \text{ for } i=1,\ldots,N$$
with finite kernels such that $b \in H^4(\X)$ is positive definite if and only if $g_i^\ast(b) \in H^4((\pt / \Gm)^{n_i}) \simeq \op{Sym}^2 (\bR^{n_i})$ is a positive definite bilinear form for all $i=1,\ldots,N$. In this case the set of positive definite classes (if non-empty) is the interior of a convex cone of full dimension in $H^4(X/G)$, and small perturbations of a positive definite class remain positive definite.

\begin{rem}
Let $b \in H^4(\X)$ be positive definite. In this case given $\gamma \in \Deg(\X,p)_2$ we can define the length of the $\gamma$ as
$$\op{length} (\gamma) = \arccos \left( \frac{b(\phi^\ast_{1,1} \gamma) - b(\phi^\ast_{1,0} \gamma) - b(\phi^\ast_{0,1} \gamma)}{2 \sqrt{b(\phi^\ast_{0,1} \gamma) b(\phi^\ast_{1,0} \gamma)}}\right),$$
where $\phi^\ast_{m,n} : \Deg(\X,p)_2 \to \Deg(\X,p)_1$ for $m,n\geq 0$ is the map induced by the morphism in $\Cones$ corresponding to the homomorphism $\bZ \to \bZ^2$ taking $\phi_{m,n} : 1 \mapsto (m,n)$. One can check that the formula for $\op{length}(\gamma)$ is invariant under any homomorphism $\bZ^2 \to \bZ^2$ preserving the two positive coordinate rays in $\bR^2$, and thus $\op{length}$ is a well-defined function of a rational $1$-simplex in $\iDeg(\X,p)$. We can then define a \emph{spherical metric} on $\iDeg(\X,p)$ by the formula
$$d(f,g) := \inf \left\{ \sum_i \op{length}(\gamma^{(i)}) \right\}$$
Where the infimum is taken over all piecewise linear paths, meaning sequences $\gamma^{(0)},\ldots,\gamma^{(n)}$ of rational $1$-simplices in $\iDeg(\X,p)$ such that $\phi_{1,0}^\ast \gamma^{(0)} = f$, $\phi_{0,1}^\ast \gamma^{(n)} = g$ and $\phi_{0,1}^\ast \gamma^{(i)} = \phi_{1,0}^\ast \gamma^{(i+1)}$. This is a generalization of the spherical metric discussed in \cite{MFK94}*{Sect.~2}
\end{rem}


\subsection{A criterion for quasi-compactness of flag spaces}

\begin{defn} \label{defn:quasi-compact_flags}
If $\X$ is a stack satisfying \ref{hyp2}, then we say $\X$ \emph{has quasi-compact flag spaces} if for any connected component $\Y \subset \Filt(\X)$, the map $\ev_1 : \Y \to \X$ is quasi-compact.
\end{defn}

\begin{prop} \label{prop:quasi_compact_flags}
Let $\X$ be a quasi-compact algebraic stack satisfying \ref{hyp2}. Let $\Phi : |\Comp(\X)_\bullet| \to \bR$ be a function such that for any $c \in \bR$ the preimage of $\Phi^{-1}((-\infty,c])$ in any rational cone $\bR^n_{\geq 0} \to |\Comp(\X)_\bullet|$ is bounded. Then $\ev_1 : \Filt^n(\X) \to \X$ is quasi-compact on each connected component of $\Filt^n(\X)$.
\end{prop}

\begin{proof}
Because $\X$ is quasi-separated, it suffices to show that the connected components of $\Filt^n(\X)$ itself are quasi-compact. The image of $\X$ is contained in the image of some smooth map $\Spec(R) \to B$, which is quasi-compact because $B$ is quasi-separated. Then $\Filt_B(\X)_R  \cong \Filt_R(\X_R)$ by \Cref{lem:base_change_mapping}, and $\Filt_B(\X)_R \to \Filt_B(\X)$ is surjective and quasi-compact, so it suffices to replace $B$ with $\Spec(R)$ and $\X$ with $\X_R$.

The projection $\agr : \Filt^n(\X) \to \Grad^n(\X)$ is quasi-compact and a bijection on connected components by \Cref{lem:retract}, so it suffices to show that the connected components of $\Grad^n(\X)$ are quasi-compact.

We first address the case $n=1$. The function $\Phi : |\Comp(\X)| \to \bR$ defines a locally constant function on $\Filt(\X)$ and hence on $\Grad(\X)$ as well. The value on each connected component of $\Filt(\X)$ is the value of $\Phi$ at the image of $1 \in \bR_{\geq 0}$ under the corresponding rational ray $\bR_{\geq 0} \to |\Comp(\X)|$ (See the formula \eqref{eq:cohomology_function}). By convention we let $\Phi = -\infty$ on the connected components of $\Filt(\X)$ classifying trivial filtrations (they are isomorphic to $\X$). We claim that the substack of $\Grad(\X)$ on which $\Phi \leq c$ is quasi-compact.

Write $\X$ as a set theoretic union $\X = \bigcup \X_i$ of locally closed substacks of the form $X_i/\GL_{n_i}$ for a $\GL_{n_i}$-quasi-projective scheme. Then $\Grad(\X) = \bigcup_i \Grad(\X_i)$ is a set-theoretic union of locally closed substacks as well. Consider the restriction of $\Phi$ along the map $|\Comp(\X_i)| \to |\Comp(\X)|$. In \Cref{lem:describe_component_fan} we described $\iComp(X_i/\GL_{n_i})$ as a finite union of real vector spaces associated to various subgroups of the cocharacter group of $\Gm^{n_i}$, modulo the action of the Weyl group $S_{n_i}$. The hypothesis that $\Phi^{-1}((-\infty,c])$ is bounded in any rational cone of $|\Comp(\X)|$ implies that there are only finitely many integral points in $|\Comp(\X_i)|$ for which $\Phi \leq c$. We know from \Cref{thm:describe_strata_quotient_GLN} that $\Grad(\X_i)$ has quasi-compact connected components, so it follows that the substack of $\Grad(\X_i)$ on which $\Phi \leq c$ is quasi-compact. This holds for all $i$, so the substack of $\Grad(\X)$ on which $\Phi \leq c$ is quasi-compact as well.

The argument for $n>1$ is very similar, so we omit it. The only difference is that one should consider the set of connected components of $\Filt^n(\X)$ such that each of the vertex maps $v_i : \Filt^n(\X) \to \Filt(\X)$ maps to the locus where $\Phi \leq c$ and show that this subset of $\Filt^n(\X)$ is quasi-compact.
\end{proof}

\begin{ex}
When $\X$ is quasi-compact and admits a positive definite class $b \in H^4(\X;\bR)$, the function $\hat{b} : |\Comp(\X)| \to \bR$ of \Cref{lem:cohomology_functions} satisfies the hypothesis of \Cref{prop:quasi_compact_flags}. More generally, \Cref{prop:quasi_compact_flags} implies that any stack $\X$ satisfying \ref{hyp2} that is quasi-compact and admits a norm on graded points in the sense of \Cref{defn:norm_graded_points} below has quasi-compact flag spaces.
\end{ex}


\section{Construction of \texorpdfstring{$\Theta$}{Theta}-stratifications} \label{sect:construction}

In this section we strengthen the main existence result, \Cref{thm:main_stratification} in several ways, leading ultimately to relatively simple necessary and sufficient conditions for the existence of a $\Theta$-stratification in \Cref{thm:main_improved}.

Our $\Theta$-stratifications will be defined by a numerical invariant, which consists (\Cref{defn:numerical_invariant}) of a realizable subset $\cU \subset \iComp(\X)$ and a continuous function $\mu : \cU \to \bR$. For any $p \in \X(k)$, $\mu$ then induces a function on the preimage of $\cU$ under the projection $\iDeg(\X,p) \to \iComp(\X,p)$. We say that $p$ is unstable if there is a point $f \in \iDeg(\X,p)$ with $\mu(f)>0$, and the HN problem is to show that if $p$ is unstable, then $\mu$ achieves a maximum at a unique rational point of $\iDeg(\X,p)$, which we call the HN filtration.

The rough idea behind our approach to the HN problem is that a maximum for $\mu$ exists because $\iComp(\X,p)$ is compact, and it is unique because the locus of maximizers of $\mu$ on $\iDeg(\X,p)$ is convex and $\mu$ is locally strictly quasi-concave on $\iDeg(\X,p)$.

\subsection{Numerical invariants and the Harder-Narasimhan problem}

In \Cref{thm:main_stratification} we encoded the data of a $\Theta$-stratification of a stack $\X$ satisfying \ref{hyp3} as a subset $S \subset \op{Irred}(\Filt(\X))$ and a function $\mu : S \to \Gamma$ for some totally ordered set $\Gamma$. In the remainder of the paper we shall specialize to a more restrictive situation in which we can use the structures of \Cref{sect:structures} to analyze the Harder-Narasimhan problem. Recall the notion of the component space $\iComp(\X)$ (\Cref{defn:reduced_degeneration_space}).

\begin{defn} \label{defn:numerical_invariant}
Let $\X$ be a stack, let $\Gamma$ be a totally ordered set such that every bounded above subset admits a supremum, and fix a distinguished element $0 \in \Gamma$. A \emph{\gls{numerical_invariant}} on $\X$ with values in $\Gamma$ consists of a realizable subset $\cU \subset \iComp(\X)$ (see \Cref{defn:realizable_subset}) and a function $\mu: \cU \to \Gamma$. Given a numerical invariant $\mu$ we define the \emph{stability function} $M^\mu : |\X| \to \Gamma \cup \{\infty\}$ as
$$M^\mu(p) = \sup \left\{ \mu(f) \left| f \in \cU \text{ with } f(1) = p \in |\X| \right. \right\},$$
where we are abusing notation by using $f \in \cU$ to denote both a filtration $f$ and the point in $\iComp(\X)$ corresponding to the connected component of $\Filt(\X)$ that contains $f$, and by letting $\mu(f)$ denote the corresponding value of $\mu$. We say that $p \in |\X|$ is \emph{unstable} if $M^\mu(p)>0$ and \emph{semistable} otherwise.
\end{defn}

\begin{rem} \label{rem:alternate_numerical}
In the special case that $\cU = \iComp(\X)$, which is the situation in most examples, one can identify \Cref{defn:numerical_invariant} with the more explicit definition of a numerical invariant given in the introduction, \Cref{defn:numerical_invariant_simple}, i.e., $\mu$ consists of a family of scale-invariant functions $\mu_\gamma : \bR^n \setminus 0 \to \Gamma$ for any $p \in \X(k)$ and homomorphism $\gamma : (\bG_m)_k \to \Aut_\X(p)$ with finite kernel, and these functions are locally constant in algebraic families and compatible with field extension and restriction along homomorphisms of tori.
\end{rem}
\begin{proof}
Any $p \in \X(k)$ and homomorphism $\gamma: (\Gm^n)_k \to \Aut_\X(p)$ with finite kernel defines a non-degenerate point of $\Grad^n(\X)$, and hence a non-degenerate connected component $\alpha \in \pi_0(\Grad^n(\X))$. \Cref{lem:retract} gives a bijection $\pi_0(\Filt^n(\X)) \simeq \pi_0(\Grad^n(\X))$, so we can identify cones of $\Comp(\X)_\bullet$ with non-degenerate connected components of $\Grad^n(\X)$. The group $(\bZ/2\bZ)^n$ acts by coordinate-wise inversion on $B\Gm^n$ and hence acts on $\Grad^n(\X)$. If we regard $\alpha \in \pi_0(\Grad^n(\X))$ as a cone in $\Comp(\X)_n$, then the $(\bZ/2\bZ)^n$-orbit of $\alpha$ gives a collection of cones that naturally glue together to form the orthants of a continuous $\bR_{>0}^\times$-equivariant map $\bR^n \to |\Comp(\X)_\bullet|$. This defines a continuous map $(\bR^n \setminus 0) / \bR_{>0}^\times \to \iComp(\X)$, and $\mu_\gamma$ is the restriction of $\mu : \iComp(\X) \to \Gamma$ along the composition $\bR^n \setminus 0 \to \iComp(\X)$. It is not hard to see that this induces an equivalence between the data of \Cref{defn:numerical_invariant_simple} and that of a numerical invariant with $\cU=\iComp(\X)$.
\end{proof}

Let $\mu : \cU \subset \iComp(\cX) \to \Gamma$ be a numerical invariant. The set of rational points of $\cU$ on which $\mu>0$ corresponds to a subset of $\pi_0 (\Filt(\X)) / \bN^\times$, or equivalently an $\bN^\times$-invariant union of connected components, which we denote by $\Y \subset \Filt(\X)$. We can choose a subset $S \subset \op{Irred}(\Y)$ that is a complete set of $\bN^\times$-orbit representatives of irreducible components, and consider the induced function $\mu : S \to \Gamma_{\geq 0} = \{c \in \Gamma | c \geq 0\}$. This is an instance of the data \eqref{eqn:data}.

\begin{defn}
We say that $\mu$ \emph{defines} a (weak) $\Theta$-stratification of $\X$ if the subset $S \subset \op{Irred}(\Filt(\X))$ and $\mu : S \to \Gamma_{\geq 0}$ constructed above define a (weak) $\Theta$-stratification in the sense of \Cref{defn:defined_stratification}.
\end{defn}

Fix a numerical invariant $\mu : \cU \subset \iComp(\X) \to \Gamma$. For any $p \in \X(k)$, the canonical continuous map $\iDeg(p) \to \iComp(\X)$ allows us to consider the preimage $\cU_p \subset \iDeg(p)$ of $\cU$ and restrict $\mu$ to a function on $\cU_p$, which we also denote by $\mu$. Note that this is compatible with the map $\iDeg(p) \to \iDeg(p')$ if $k'/k$ is a field extension and $p' \in \X(k')$ the $k'$-point induced by $p$. This allows us to reformulate the Harder-Narasimhan problem as follows:
\begin{problem}[Harder-Narasimhan] \label{problem:hn}
Given an unstable $p \in \X(k)$, is there a unique rational point $f \in \cU_p \subset \iDeg(\X,p)$ that maximizes $\mu(f)$?
\end{problem}
If such an $f$ exists, it corresponds to a filtration $f : \Theta_k \to \X$ of $p$ up to composition with the ramified covering maps $(\bullet)^n : \Theta_k \to \Theta_k$. If $f$ is still a maximizer for $\mu$ after base change to the algebraic closure $\bar{k}$, it will be a Harder-Narasimhan filtration in the sense of part (1) of \Cref{thm:main_stratification}.

\begin{defn}[HN filtration] \label{defn:HN_filtration}
An HN filtration of a point $p \in \X(k)$ is a non-degenerate filtration $f : \Theta_k \to \X_k$ along with an isomorphism $f(1) \cong p \in \X_k(k)$ such that if $\bar{p} \in \X(\bar{k})$ is the base change of $p$ to an algebraic closure of $k$, then $f_{\bar{k}} \in \iDeg(\X;\bar{p})$ lies in $\cU_{\bar{p}}$ and maximizes $\mu$.
\end{defn}

\subsubsection{Properties of numerical invariants}

\begin{defn} \label{defn:locally_quasi-concave_invariant}
Let $\X$ be a stack over a base stack $B$, and let $\mu : \cU \subset \iComp(\X)$ be a numerical invariant. For any field $k$ over $B$, any $p \in \X(k)$, and any rational simplex $\simp{\sigma}^n \to \cU_p \subset \iDeg(\X,p)$, let us denote the restriction of $\mu$ by $\mu_\sigma : \simp{\sigma}^n \to \Gamma$. We say that $\mu$ is:
\begin{itemize}
\item \emph{locally quasi-concave} if for any $p$ and $\simp{\sigma}^n$ as above, $\mu_{\sigma}^{-1}([0,\infty)) \subset \simp{\sigma}^n$ is convex, and for any distinct points $x_0,x_1,x_2 \in \mu_\sigma^{-1}([0,\infty))$ with $x_1$ on the line segment joining $x_0$ and $x_2$, one has
\begin{equation} \label{eq:quasi-concave_numerical}
\mu_\sigma(x_1) \geq \min(\mu_\sigma(x_0),\mu_\sigma(x_2)),
\end{equation}
with $\mu_\sigma(x_0)=\mu_\sigma(x_2)$ whenever equality holds;\\
\item \emph{quasi-concave} if it is locally quasi-concave and for any $p \in \X(k)$ as above $\cU_p^{\mu>0}:=\{x \in \cU_p | \mu(x)>0\} \subset \iDeg(\X,p)$ is convex in the sense of \Cref{defn:convex_subset}; and\\
\item \emph{(locally) strictly quasi-concave} if it is (locally) quasi-concave, and the restriction of $\mu$ to $\cU_p^{\mu>0}$ is strictly quasi-concave in the sense of \Cref{defn:quasi-concave_function}, i.e., the inequality \eqref{eq:quasi-concave_numerical} holds strictly whenever $x_0,x_2 \in \mu_\sigma^{-1}((0,\infty))$.
\end{itemize}
\end{defn}

The condition that $\mu$ is strictly quasi-concave is very strong, as the following Corollary shows. Nevertheless we will see many examples in which this condition holds in \Cref{sect:theta_reductive}.

\begin{cor}[Uniqueness of HN filtrations] \label{cor:hn_uniqueness}
Let $\X$ be a stack and let $\mu : \cU \subset \iComp(\X) \to \Gamma$ be a strictly quasi-concave numerical invariant, then an HN filtration of any unstable point $p \in \X(k)$ is unique up to the action of $\bN^\times$, if it exists.
\end{cor}
\begin{proof}
By definition the subset $\cU_p^{\mu > 0} \subset \iDeg(\X,p)$ is convex and the function $\mu$ is strictly quasi-concave, so at most one rational point in $\cU_p^{\mu>0}$ can maximize $\mu$ by \Cref{lem:max_quasi_concave}.
\end{proof}

In order to formulate the standard hypotheses for a numerical invariant, we need the following:
\begin{defn} \label{defn:antipodal_points}
We will say that $f_1,f_2 \in \Filt(\X)(k)$ are \emph{antipodal} if there is a 2-commutative diagram of the form
$$\xymatrix@R=5pt{\Theta_k \ar[dr] \ar@/^/[drr]^{f_1} & & \\ & (\pt/\Gm)_k \ar[r] & \X \\ \Theta_k \ar[ur] \ar@/_/[urr]_{f_2} & & }$$
such that the cocharacters $(\Gm)_k \to (\Gm)_k$ induced by the two maps $\Theta_k \to (\pt/\Gm)_k$ have opposite sign. Note in particular that if $f_1$ and $f_2$ are antipodal, then they are both split filtrations. We say that two points of $\iComp(\X)$ are antipodal if the corresponding connected components of $\Filt(\X)$ contain a pair of antipodal points.
\end{defn}

\begin{ex}
When $\X = B\SL_2$, the degeneration space $\iDeg(\X,\ast)$ is an infinite disjoint union of points, one for each Borel subgroup of $\SL_2$, i.e. one point for each line in $\bA^2_k$ regarded as a two step filtration $L \subset k^2$. For any two distinct lines, there is a unique (up to multiples) one parameter subgroup of $\SL_2$ acting with positive weight on $L_1$ and negative weight on $L_2$. This one parameter subgroup corresponds to a map $\pt/\Gm \to \X$, and the maps $\Theta \to \X$ corresponding to the two Borel subgroups factor through this $\pt/\Gm$. Thus any two distinct points of $\iDeg(B\SL_2,\ast)$ are antipodal.
\end{ex}

\begin{defn}[\gls{standard_numerical_invariant}] \label{defn:standard_invariant}
We say that a numerical invariant $\mu : \cU \subset \iComp(\X) \to \Gamma$ is \emph{standard} if it is locally strictly quasi-concave and $\cU^{\mu>0}$ does not contain a pair of antipodal points.
\end{defn}

Note that any strictly quasi-concave numerical invariant satisfies the property that $\cU^{\mu>0}$ does not contain a pair of antipodal points, so it is automatically standard.

Given a point $p \in \X(k)$ for an algebraically closed field $k$, the existence of a maximizer of $\mu$ on $\cU_p \subset \iDeg(p)$ is weaker than the existence of a HN filtration for $p$, because the maximum need not occur at a rational point. We will therefore often restrict our focus to numerical invariants that satisfy the following condition:
\begin{enumerate}[label=(R)]
\item \label{princ:R} For any rational simplex $\simp{} \to \cU \subset \iComp(\X)$, if the restriction of $\mu$ to $\simp{}$ has a point with $\mu>0$, then $\mu$ obtains a maximum at a rational point of $\simp{}$.
\end{enumerate}
\begin{lem} \label{lem:princ_R}
Let $\mu : \cU \subset \iComp(\X) \to \Gamma$ be a numerical invariant satisfying \ref{princ:R}, let $p \in \X(k)$, and let $\bar{p} \in \X(\bar{k})$ be the induced point over the algebraic closure. Then $p$ has a HN filtration if and only if $\mu : \cU_{\bar{p}} \subset \iDeg(\bar{p}) \to \Gamma$ obtains a maximum.
\end{lem}
\begin{proof}
Because the flag space $\Flag(p)$ is locally finite type over $k$ and is compatible with base change, $p$ has a HN filtration if and only if $\bar{p}$ does. Condition \ref{princ:R} means that $\mu$ obtains a maximum on $\cU_{\bar{p}}$ if and only if it obtains a maximum at a rational point.
\end{proof}

\subsubsection{Numerical invariants from cohomology classes}

\begin{defn}[Norm] \label{defn:norm_graded_points}
A \emph{(semi-)norm on graded points} of a stack $\X$ assigns a (semi-)norm $\lVert \bullet \rVert_\gamma : \bR^n \to \bR$ for any $p \in \X(k)$ and homomorphism $\gamma : (\bG_m^n)_k \to \Aut_\X(p)$ with finite kernel. This data is required to satisfy the compatibility conditions (1), (2) and (3) of \Cref{defn:numerical_invariant_simple}.
\end{defn}

As in \Cref{rem:alternate_numerical}, this data is equivalent to the data of a function $|\Comp(\X)_\bullet| \to \bR$, which we also denote $\lVert \bullet \rVert$, that restricts to a (semi-)norm along any of the maps $\bR^n_\gamma \to |\Comp(\X)_\bullet|$ corresponding to a homomorphism $\gamma : (\bG_m^n)_k \to \Aut_\X(p)$ with finite kernel. More explicitly, for any filtration $f : \Theta_k \to \X$, one defines $\lVert f \rVert := \lVert 1 \rVert_{\agr(f) : (B\bG_m)_k \to \X}$.

\begin{ex} \label{ex:norm_from_H4}
Let $b \in H^4(\X;\bR)$ be positive definite (\Cref{defn:positive_definite}). Then the function $\sqrt{\hat{b}} : |\Comp(\X)_\bullet| \to \bR$ is a norm on graded points, where $\hat{b}$ is the function assigned to $b$ in \Cref{lem:cohomology_functions}.
\end{ex}

\begin{defn}[Associated numerical invariant] \label{defn:induced_invariant}
Let $\ell \in H^2(\X;\bR)$, and let $\lVert \bullet \rVert : |\Comp(\X)_\bullet| \to \bR$ be a semi-norm on graded points. Then the \emph{\gls{numerical_invariant} associated to $\ell$ and $\lVert \bullet \rVert$} is the pair
\[
\cU := \left\{x \in \iComp(\X) \left| \lVert \tilde{x} \rVert > 0  \right. \right\} \quad \text{and} \quad \mu(x) = \frac{\hat{\ell}(\tilde{x}) }{\lVert \tilde{x} \rVert} \in \bR,
\]
where $\tilde{x} \in |\Comp(\X)_\bullet|$ is some lift of $x \in \iComp(\X)$, and $\hat{\ell}$ is defined in \Cref{lem:cohomology_functions}.
\end{defn}

Note that both $\cU$ and $\mu$ are well-defined because both $\hat{\ell}$ and $\lVert \bullet \rVert$ are homogeneous of weight $1$ under the scaling action of $\bR_{>0}^\times$ on $|\Comp(\X)_\bullet|$. $\cU$ is open and hence realizable, so \Cref{defn:induced_invariant} indeed defines a numerical invariant. Note also that $\mu$ is \emph{continuous}. We will almost always be interested in the case when $\lVert \bullet \rVert$ is a norm and hence $\cU = \iComp(\X)$.

\begin{lem} \label{lem:induced_invariant}
The numerical invariant of \Cref{defn:induced_invariant} is locally quasi-concave, and $\cU^{\mu>0}$ does not contain a pair of antipodal points. It is standard (\Cref{defn:standard_invariant}) if $\lVert \bullet \rVert$ is locally a norm in which the unit ball is strictly convex, such as the norm on graded points associated to a positive definite $b \in H^4(\X;\bR)$.
\end{lem}
\begin{proof}
First note that $\cU^{\mu>0}$ does not contain a pair of antipodal points, because $\mu(f) = -\mu(f')$ for any two antipodal filtrations.

Consider the restriction $\mu_\sigma$ of $\mu$ along the map $\bR^n_{\geq 0} \to |\Deg(\X,p)_\bullet| \to |\Comp(\X)_\bullet|$ defined for a $p \in \X(k)$ and $\sigma \in \Deg(p)_n$. The function $\mu_\sigma$ is defined wherever $\lVert x \rVert_\sigma>0$, and it has the form $\ell_\sigma \cdot x / \lVert x \rVert_\sigma$ for some $\ell_\sigma \in \bR^n$. To show local quasi-concavity of $\mu$, we must show that $\mu_\sigma$ is quasi-concave along the line segment joining any two points $y_0,y_1 \in \bR_{\geq 0}^n$ that lie on distinct rays and satisfy $\ell_\sigma \cdot y_i \geq 0$ and $\lVert y_i\rVert_\sigma > 0$ for $i=0,1$. By the scale-invariance of $\mu$, it suffices to show this after rescaling the points by a positive number so that $\ell_\sigma \cdot y_0 = \ell_\sigma \cdot y_1 = c \geq 0$. It follows from the concavity of $\ell_\sigma \cdot (-)$ that $\ell_\sigma \cdot y=c$ for all points on the line segment joining $y_0$ and $y_1$. Thus on this line segment, we have $\mu(y) = c / \lVert y \rVert_\sigma$. The triangle inequality implies that for $\lambda \in (0,1)$,
\[
\lVert \lambda y_0 + (1-\lambda) y_1 \rVert_\sigma \leq \lambda \lVert y_0 \rVert_\sigma + (1-\lambda) \lVert y_1 \rVert_\sigma \leq \max(\lVert y_0 \rVert_\sigma, \lVert y_1 \rVert_\sigma),
\]
with $\lVert y_0 \rVert_\sigma = \lVert y_1 \rVert_\sigma$ whenever equality holds, which verifies the quasi-concavity of $\mu_\sigma$ on the segment joining $y_0$ and $y_1$. If the norm $\lVert \bullet \rVert_\sigma$ is strictly convex, then the first inequality above is always strict, hence $\mu_\sigma$ is strictly quasi-concave.

\end{proof}

\begin{lem} \label{lem:induced_invariants_R}
If $\ell \in H^2(\X;\bQ)$, and $\lVert \bullet \rVert : |\Comp(\X)_\bullet| \to \bR$ is a norm on graded points whose restriction to any cone is either 1) a piecewise rational linear function, or 2) a rational quadratic norm (such as the norm associated to a positive definite $b \in H^4(\X;\bQ)$), then the induced numerical invariant from \Cref{defn:induced_invariant} satisfies condition \ref{princ:R}.
\end{lem}
\begin{proof}
Let $\simp{\sigma} \to \cU$ be a rational simplex such that $\mu>0$ at some point in $\simp{\sigma}$. Restricted to $\simp{\sigma} = (\bR_{\geq 0}^n \setminus 0) / \bR_{>0}^\times$, $\mu$ has the form $\ell \cdot x / \lVert x \rVert$ where $\ell \in \bQ^n$ and $\lVert x \rVert$. We know that $\mu(x)$ achieves a maximum on $\bR^n_{\geq 0}$ which is $>0$, so the affine hyperplane $\{\ell \cdot x = 1\}$ meets the cone $\bR^n_{\geq 0}$. Maximizing $\ell \cdot x / \lVert x \rVert$ on $\bR^n_{\geq 0}$ is thus equivalent to minimizing $\lVert x \rVert$ on the rational polyhedron $\bR^n_{\geq 0} \cap \{\ell \cdot x = 1\}$.

We prove the claim in the case where $\lVert \bullet \rVert = \sqrt{b(x)}$, where $b$ is a rational quadratic form such that $b(x)>0$ for nonzero $x \in \bR_{\geq 0}^n$. The proof in the case where $\lVert \bullet \rVert$ is a rational piecewise linear function is similar.

The Kuhn-Tucker condition\footnote{The Kuhn-Tucker condition is analogous to Lagrange multiplier equations, but for optimization problems with inequality constraints.} for this convex minimization problem states that $x \in \bR^n$ is a global minimum for this constrained optimization problem if and only if it satisfies the constraints and there is a Lagrange multiplier $\lambda \in \bR$ and a Kuhn-Tucker multiplier $\tau \in \bR^n_{\geq 0}$ such that
\[
0 = \nabla b(x) - \tau + \lambda \ell \text{ and } \tau \cdot x = 0.
\]
These two equality constraints on $(x,\tau,\lambda)$ combined with the equality constraint $\ell \cdot x =0$ define a rational linear subspace of $\bR^n \times \bR^n \times \bR$. Thus the set of points $(x,\tau,\lambda)$ representing a solution of the Kuhn-Tucker condition is the intersection if a rational linear subspace with the rational polyhedral cone $\{x \geq 0, \tau \geq 0\}$, and such a set always contains a rational point if it is nonempty.
\end{proof}

\begin{ex}[GIT semistability] \label{ex:git_numerical_invariant}
Let $X$ be a projective over affine $k$-variety with a reductive group action. Let $\cL$ be a $G$-linearized ample invertible sheaf on $X$, and let $|\bullet|$ be a Weil-group-invariant positive definite bilinear form on the cocharacter lattice of $G$. \Cref{thm:describe_strata_global_quotient} identifies a filtration $f$ of a point $p \in X(k)$ with a one parameter subgroup $\lambda : \Gm \to G$ such that $q := \lim_{t \to 0} \lambda(t) \cdot p$ exists, taken up to conjugation by an element of $P_\lambda(k)$. One can choose a class $l \in H^2(\X;\bQ)$ and a positive definite $b \in H^4(\X;\bQ)$ such that
\begin{equation} \label{eqn:original_GIT_numerical_invariant}
\mu(f) = \frac{f^\ast l}{\sqrt{f^\ast b}} = \mu(p,\lambda) = \frac{-1}{|\lambda|} \op{weight}_\lambda \cL|_q \in \bR.
\end{equation}
is the normalized Hilbert-Mumford numerical invariant \cite{DH98}. Thus $\Theta$-stability agrees with GIT stability, and the Harder-Narasimhan filtration of an unstable point corresponds Kempf's canonical destabilizing one parameter subgroup \cite{Ke78}.
\end{ex}

\begin{proof}
We define $l = c_1(\cL) \in H^2(\X;\bQ)$, so that $$f^\ast l = c_1(f^\ast \cL) = - \op{weight} ( (f^\ast \cL)_{\{0\}}) \cdot q,$$ which follows from the fact\footnote{Every invertible sheaf on $\Theta$ is of the form $\cO_\Theta(n)$, which corresponds to the free $k[t]$ module with generator in degree $-n$. Note that the isomorphism $\op{Pic}(\Theta) \simeq \bZ$ is canonical, because $\Gamma(\Theta,\cO_\Theta(n)) = 0$ for $n>0$ whereas $\Gamma(\Theta,\cO_\Theta(n)) \simeq k$ for $n\leq 0$. This holds in contrast to $\op{Pic}(\pt / \Gm)$. The stack $\pt / \Gm$ has an automorphisms exchanging $\cO(1)$ and $\cO(-1)$.} that $f^\ast \cL \simeq \cO_{\Theta}(w)$ where $w = - \op{weight} (f^\ast \cL)_{\{0\}}$. For the denominator, $|\bullet|$ can be interpreted as a class in $H^4(\pt / G ; \bC)$ under the identification $H^\ast(\pt / G; \bC) \simeq (\op{Sym}(\fg^\dual))^G$, and we let $b$ be the image of this class under the map $H^4(\pt / G) \to H^4(X/G)$. For a morphism $f : \Theta \to X/G$, $f^\ast b$ is the pullback of the class in $H^4(\pt / G)$ under the composition $\Theta \to X/G \to \pt / G$. We therefore have $f^\ast b = |\lambda|^2 q^2 \in H^4(\Theta)$.
\end{proof}

\subsubsection{Induced numerical invariants and stratification of \texorpdfstring{$\Grad(\X)$}{Grad(X)}}

\begin{defn} \label{defn:restricted_invariant}
Let $\pi : \X \to \Y$ be a morphism of stacks with quasi-finite relative inertia, such as a representable morphism, and let $\mu : \cU \subset \iComp(\X) \to \Gamma$ be a numerical invariant. The \emph{restriction of $\mu$} along $\pi$ is the function
\[
\mu \circ \pi_\ast : (\pi_\ast)^{-1}(\cU) \to \Gamma,
\]
where $\pi_\ast : \iComp(\Y) \to \iComp(\X)$ is the associated map on component spaces.
\end{defn}

We are particularly interested below in the case where $\Y = \Grad(\X)$ and $u :\Grad(\X) \to \X$ is the canonical forgetful map.

\begin{lem}
Let $\pi : \Y \to \X$ be a map of stacks that satisfy \ref{hyp2}, and let $\mu : \cU \subset \iComp(\X) \to \Gamma$ be a numerical invariant that defines a (weak) $\Theta$-stratification. If this (weak) $\Theta$-stratification on $\X$ induces a (weak) $\Theta$-stratification on $\Y$ in the sense of \Cref{defn:induced_stratum}, then this (weak) $\Theta$-stratification is also defined by the induced numerical invariant $\mu \circ \pi_\ast$ on $\Y$.
\end{lem}

\begin{proof}
Recall from \eqref{eqn:putative_stratification} that a weak $\Theta$-stratification (with locally constant $\mu$) can be uniquely recovered from a subset $S \subset \pi_0(\Filt(\X))$ and a map $\mu : S \to \Gamma$ by defining $|\X_{\leq c}| = \{p \in |\X| | M^\mu(p) \leq c\}$ and $|\S_c| = \{f \in |\Filt(\X)| \text{ with }\mu(f) = M^\mu(f) \text{ and that lie on a component in } S\}$. In our case $S$ is a set of orbit representatives for the action of $\bN^\times$ on
\[
\pi_0(\Filt(\X))_{\cU} := \{ \alpha \in \pi_0(\Filt(\X)) | \alpha \in \cU \subset \iComp(\X)\},
\]
and $\mu$ is given by the numerical invariant $\mu : \pi_0(\Filt(\X))_{\cU} \to \Gamma$. If $\pi : \Y \to \X$ induces a stratification of $\Y$, then one can likewise encode this stratification by data $\mu' : S' \subset \pi_0(\Filt(\Y)) \to \Gamma$. In this case $S'$ is the preimage of $S$ under $\Filt(\pi) : \Filt(\Y) \to \Filt(\X)$, and $\mu'$ is the composition $S' \to S \to \Gamma$. This is by definition the data of the induced numerical invariant of $\Y$.
\end{proof}

As a consequence we have:

\begin{cor}
If $\X$ is a stack satisfying \ref{hyp2} and a numerical invariant $\mu$ defines a (weak) $\Theta$-stratification of $\X$, then the restriction of $\mu$ to $\Grad(\X)$ along the map $u : \Grad(\X) \to \X$ defines a (weak) $\Theta$-stratification, and this is the same as the induced (weak) $\Theta$-stratification of \Cref{prop:graded_strata}.
\end{cor}

We will see in the next section that if $p \in \Grad(\X)$ is a graded point such that the corresponding split filtration $\sigma(p)$ is destabilizing for the underlying point $u(p) \in \X$, then the entire connected component containing $p$ is unstable for the induced $\Theta$-stratification of $\Grad(\X)$. Therefore, many connected components of $\Grad(\X)$ contain no semistable points.


\subsection{The Recognition Theorem for HN filtrations}

In this section we generalize the original definition of the Harder-Narasimhan filtration of a vector bundle on a curve \cite{HN} as the unique filtration such that the associated graded bundles are semistable, with strictly increasing slopes. We will characterize an HN filtration $f:\Theta_k \to \X$ in terms of a notion of semistability for the associated graded point $\agr(f) \in \Grad(\X)(k)$.

Under the identification $\pi_0(\Grad(\X)) = \pi_0(\Filt(\X))$, a numerical invariant $\mu$ on a stack $\X$ defines a locally constant function on $\Grad(\X)$. For clarity let us call this
\[
\mu_{\rm can} : \Grad(\X) \to \Gamma.
\]
The value of $\mu_{\rm can}$ at $p \in \Grad(\X)(k)$ agrees with the value of the restriction of $\mu$ along the forgetful map $u : \Grad(\X) \to \X$ (see \Cref{defn:restricted_invariant}) on the canonical filtration of $p$, which is the same as the value of $\mu$ on the split filtration $\sigma(p) : \Theta_k \to \X$ of the underlying point $u(p) \in \X(k)$. So
\[
M_{\Grad(\X)}^{\mu \circ u_\ast}(p) \geq \mu_{\rm can}(p),
\]
and any component of $\Grad(\X)$ on which $\mu_{\rm can}>0$ is unstable with respect to the restricted numerical invariant. We therefore use an alternate notion of semistability on $\Grad(\X)$.

\begin{defn}[Graded-semistability] \label{defn:graded_semistable}
We say that a point $p \in \Grad(\X)$ is \emph{graded-semistable} if $M_{\Grad(\X)}^{\mu \circ u_\ast}(p) = \mu_{\rm can}(p)$, i.e., the canonical filtration of $p$ maximizes the value of the numerical invariant among all filtrations of $p$.
\end{defn}

If $\mu$ takes values in a totally ordered abelian group $\Gamma$, then we can define a new numerical invariant on $\Grad(\X)$ by the formula
\[
\mu'(f) = \mu(u_\ast(f)) - \mu_{\rm can}(f(1))
\]
for any filtration $f : \Theta_k \to \Grad(\X)$. With this modified numerical invariant, a point $p \in \Grad(\X)(k)$ is graded-semistable with respect to $\mu$ if and only if it is semistable with respect to $\mu'$. Combining this with \Cref{cor:semistable_centers}, we see that if $\cZ_c \subset \Grad(\X)$ is the union of the connected components that meet the center $\cZ_c^{\rm{ss}}$ of the $\Theta$-stratum $\cS_c \hookrightarrow \X$, then the center $\cZ_c^{\rm{ss}}$ is precisely the semistable locus for the numerical invariant induced by $\mu'$ on $\cZ_c$.

\begin{rem}
This notion of graded-semistability generalizes the inductive description of the KN stratification in geometric invariant theory given in \cite{Ki84}, where it is shown that the centers of the strata induced by a reductive group acting on a projective variety are themselves semistable loci for a reductive subgroup acting on a closed subvariety.
\end{rem}

\begin{thm}[Recognition theorem, preliminary version] \label{thm:recognition}
Let $\X$ be a stack satisfying \ref{hyp3}, let $\mu : \cU \subset \iComp(\X) \to \Gamma$ be a locally quasi-concave \gls{numerical_invariant} for which $\cU^{\mu>0}$ does not contain a pair of antipodal points, and let $f : \Theta_k \to \X$ be a filtration in $\X$ with $\mu(f) > 0$. If $f$ is an HN filtration, then $\agr(f) \in \Grad(\X)(k)$ is graded-semistable. Conversely, if $\agr(f)$ is graded-semistable, then $f$ is an HN filtration of $f(1)$ if either of the following hold:
\begin{enumerate}
\item $\mu$ is strictly quasi-concave; or
\item $\mu$ defines a weak $\Theta$-stratification of $\X$.
\end{enumerate}
\end{thm}

We will need the following observation:

\begin{lem} \label{lem:embedded_numerical_invariant}
Let $\X$ be a stack satisfying \ref{hyp3}, and let $f \in \Filt(\X)(k)$. Then the two compositions agree
\[
\xymatrix{\Deg(\agr(f))_\bullet^\canon \ar@/^/[r]^\bT \ar@/_/[r]_{u_\ast} & \Deg(f(1))_\bullet \ar[r] & \Comp(\X)_\bullet},
\]
where $u_\ast$ is the map induced by the forgetful map $u : \Grad(\X) \to \X$, and $\bT$ is the transfer map of \Cref{thm:perturbation}. \end{lem}
\begin{proof}
Considering the concrete description of the map $\bT$ of \eqref{eqn:define_T} following its definition, for any cone $\xi \in \Deg(\agr(f))^\canon_n$ we have a map $\op{ext}(\xi) : \Theta_k \times \Theta_k^n \to \X$ whose restriction to $\{t_0=0\}$ is a map $(\pt/\Gm)_k \times \Theta^n_k \to \X$ representing $\xi$, and the restriction to $\{t_0=1\}$ represents $\bT(\xi)$. We regard the restriction of $\op{ext}(\xi)$ to $\bA^1_k \times \Theta^n_k \to \X$ as a map $\bA^1_k \to \Filt^n(\X)$ that maps $\{0\}$ to the filtration underlying the cone $u_\ast(\xi)$ and maps $\{1\}$ to the the filtration underlying the cone $\bT(\xi)$. Hence these two filtrations lie on the same connected component of $\Filt^n(\X)$, which implies that the two maps take $\xi$ to the same cone in $\Comp(\X)_\bullet$.
\end{proof}

\begin{proof}[Proof of \Cref{thm:recognition}]
Let $\mu'$ denote the restriction of $\mu$ to $\Grad(\X)$, and let $f'$ be a filtration of $\agr(f) \in \Grad(\X)(k)$ for which $\mu'$ is defined and $\mu'(f')>\mu'(\canon)$. The filtration $f'$ is not antipodal to $\canon$, because it lies over $\cU^{\mu>0}$.  Let $\simp{\sigma}^1 \to \iDeg(\Grad(\X),\agr(f))$ be the canonical rational $1$-simplex with $v_0(\sigma) = \canon$ and $v_1(\sigma) = f'$ from \Cref{thm:perturbation_light}. By hypothesis $\mu'(f') > \mu'(\canon) = \mu(f)$ and $\mu'$ is locally quasi-concave, so the maximum of $\mu'$ restricted to $\simp{\sigma}$ occurs either at $f'$ or at some point in the interior of $\simp{\sigma}$, and $\mu'$ must be strictly increasing in a neighborhood of $\canon \in \simp{\sigma}$.

By \Cref{thm:perturbation_light} we may choose a small subinterval of $\simp{\sigma}$ containing $v_0$ and lying in $\iDeg(\agr(f))^\canon$ on which $\mu'$ is strictly increasing. The map $\bT$ identifies this interval with a rational $1$-simplex in $\iDeg(\X,f(1))$ with $v_0 = f$, and \Cref{lem:embedded_numerical_invariant} implies that under this identification the restriction of $\mu' : \iDeg(\agr(f)) \to \bR$ corresponds to restriction of $\mu: \iDeg(f(1))$. Hence we have found a rational $1$-simplex in $\Deg(\X,f(1))$ for which $v_0 = f$ and along which $\mu$ is strictly increasing, which shows that $f$ is not a HN filtration. The contrapositive: if $f$ is a HN filtration, then $\agr(f)$ is graded-semistable.

For the converse, first assume that $\mu$ is strictly quasi-concave. Replace $k$ with a field extension if necessary and let $f' : \Theta_{k} \to \X$ be a filtration of $f(1)$ with $\mu(f')>\mu(f)$. By \Cref{defn:locally_quasi-concave_invariant}, the realizable subset $\cU_p^{\mu>0} \subset \iDeg(\X,p)$ is convex, and so $f,f' \in \cU_p^{\mu>0}$ are connected by a rational $1$-simplex $\simp{\sigma} \to \cU_p^{\mu>0} \subset \iDeg(\X,f(1))$ with $v_0(\sigma) = f$ and $v_1(\sigma)=f'$. By the same argument as in the previous paragraph, $\mu$ must be strictly increasing in a neighborhood of $v_0$. Again by \Cref{thm:perturbation_light} we can identify a small subinterval of $\simp{\sigma}$ containing $v_0$ with a rational $1$-simplex in $\iDeg(\agr(f))^\canon$ with $v_0 = \canon$. The function $\mu$ is strictly increasing on this interval, so $\agr(f)$ is not graded-semistable.

Finally, assume that $\mu$ defines a weak $\Theta$-stratification of $\X$ and let $f'$ be the HN-filtration of $f(0)$ defined by $\mu$. If $f$ is not an HN filtration of $f(0)$, then $\mu(f')>\mu(f)$. Furthermore, $f'$ is $\Gm$-equivariant by part (2) of \Cref{lem:HN_filtrations} (see also \Cref{rem:equivariant_HN}), so it lifts to a filtration in $\Grad(\X)$ with $\mu(f')>\mu(\canon)$. This shows that $\gr(f)$ is not graded-semistable.

\end{proof}


\subsection{Specialization condition for the uniqueness of HN filtrations}

\begin{prop} \label{prop:specialization_implies_uniqueness}
Let $\X$ be a stack satisfying \ref{hyp3}, and let $\mu : \cU \subset \iComp(\X) \to \Gamma$ be a \gls{standard_numerical_invariant} that satisfies the Simplified HN-Specialization property \ref{princ:S} of \Cref{thm:main_stratification}. Then any $HN$-filtration of an unstable point is unique up to the action of $\bN^\times$.
\end{prop}

\begin{proof}
Let $p \in \X(k)$, and consider two $k$-points $f_1,f_2$ of $\Flag(p)$ such that $f_1^a \neq f_2^b$ for any $a,b>0$, i.e., they represent distinct rational points in $\iDeg(\X,p)$, and such that $f_1$ and $f_2$ lie in $\cU_p$ and are both HN filtrations for $p$. First, assume there is a rational $1$-simplex $\Delta_g \to \iDeg(\X,p)$ whose endpoints map to $f_1$ and $f_2$. By hypothesis $\cU_p^{\mu>0} \subset \iDeg(\X,p)$ is locally convex, so the image of $\Delta_g \to \iDeg(\X,p)$ must lie in $\cU^{>0}_p$, and also by hypothesis the restriction of $\mu$ to $\Delta_g$ is strictly quasi-concave, so it can achieve its maximum at exactly one point (see \Cref{lem:max_quasi_concave}). This contradicts the hypothesis that $f_1$ and $f_2$ were distinct points of $\iDeg(\X,p)$. Therefore the claim of the proposition follows from the existence of such a $\Delta_g$, which we show in the remainder of this proof.

Because $\X$ is locally finitely presented over $B$, both $f_1$, $f_2$, and the isomorphism $f_1(1) \simeq f_2(1)$ are induced via base change from a field that is essentially of finite type over $B$. Therefore we may assume that $k$ is essentially finite type over $B$.

Let $U = \bA_k^1 -\{0\}$. We consider $f_1$ and $f_2$ as morphisms $U \times \bA_k^1 / (\Gm)_k^2 \to \X$ and $\bA_k^1 \times U / (\Gm)_k^2 \to \X$ respectively with a fixed isomorphism of their restrictions to $U\times U / (\Gm)_k^2 \simeq \pt$, so we can glue them to define
$$f_1 \cup f_2 : (\bA_k^2 - \{0\}) / (\Gm)_k^2 \to \X$$
And the data of the morphism $f_1 \cup f_2$ is equivalent to the data of the pair $f_1,f_2$. We will first show that if $f_1$ and $f_2$ are both HN filtrations, then after passing to a finite extension of $k$ and replacing $f_1$ with a positive multiple, $f_1 \cup f_2$ extends uniquely to a morphism $\bA_k^2 / (\Gm)_k^2 \to \X$.\footnote{We can rephrase the construction of $f_1 \cup f_2$ above more formally and more generally as the observation that restriction defines an equivalence of stacks
\[
\inner{\op{Map}} \left( (\bA^2 - \{0\}) / \Gm^2 , \X \right) \xrightarrow{\simeq} \Filt(\X) \times_{\X} \Filt(\X),
\]
where the fiber product is taken with respect to $\ev_1$.}

\medskip
\noindent{\textit{Extending the morphism $f_1 \cup f_2$:}}
\medskip

Finding a morphism $g : \Theta^2_{k} \to \X$ along with an isomorphism of the restriction of $g$ to $(\bA^2_{k} - \{0\}) / \Gm^2$ with $f_1 \cup f_2$ is equivalent to finding a lift in the diagram
\[
\xymatrix{ (\bA^1_{k} - \{0\}) / \Gm \simeq \Spec k \ar[r] \ar[d] & \filt{\X} \ar[d]^{\ev_1} \\ \bA^1_{k} / \Gm \ar@{-->}[ur]^g \ar[r]^{f_1} & \X }
\]
where the left vertical morphism is the inclusion of the point $\{1\}$, and the top horizontal morphism classifies the morphism $f_2 : \Theta_{k} \to \X$.

By \Cref{prop:representable_flags}, we can write the fiber product with respect to $\ev_1$ as
\[X/\Gm := \bA^1_k / \Gm \times_{\X} \filt{\X},\]
where $X$ is some $\Gm$-equivariant algebraic space that is locally finite type and quasi-separated over $\bA^1_k$. The top vertical arrow in the lifting diagram defines a $\Gm$-invariant section $s^\circ$ of $X \to \bA^1_k$ over $\bA^1 - \{0\}$, and the existence of a dotted arrow $g$ is equivalent to the extension of this to a $\Gm$-invariant section $s : \bA^1_k \to X$.

If we restrict $f_1$ to the completion $\Spec(k[\![t]\!]) \to \bA^1_k / \Gm$, the condition \ref{princ:S} implies that the given lift $\Spec(k(\!(t)\!)) \to \filt{\X}$ of the map to $\X$, which classifies a HN filtration, extends to a lift $\Spec(k[\![t]\!]) \to \filt{\X}$ after passing to a finite extension and adjoining a root of $t$. This implies that after passing to an extension of $k$ and a ramified covering $\Theta_k \to \Theta_k$, we can assume the section $s^\circ$ extends to a section over the completion of $\bA^1_k$ at the origin, and hence by the Beauville-Laszlo theorem to a section $s : \bA^1_k \to X$.

If $X$ were separated, $s$ would automatically be $\Gm$-equivariant, but we do not know this in general. Instead, we choose a quasi-compact $\Gm$-equivariant open subspace $U$ containing the image of $s$, and let $Y$ denote the reduced closure of the image of $s^\circ$ in $U$. Note that the non-equivariant map $s : \bA^1_k \to X$ still factors through $Y$, and $Y$ is flat over $\bA^1_k$, hence quasi-finite. So, we have a reduced $(\Gm)_k$-space $Y$ with a $(\Gm)_k$ equivariant quasi-finite map $Y \to \bA^1_k$ and an equivariant section $s^\circ$ over $\bA^1_k \setminus 0$ that extends non-equivariantly to a section over $\bA^1_k$, and we wish to show that $s^\circ$ extends equivariantly.

Using Sumihiro's theorem for algebraic spaces \cite{ahr2}*{Thm.~20.1}, we may choose a surjective $(\Gm)_k$-equivariant \'etale map $\Spec(A) \to Y$. We claim that there is a unique scheme theoretic disjoint union $\Spec(A) = \Spec(A') \sqcup \Spec(B)$ such that $\Spec(A')$ is finite over $\bA^1_k$ and the image of $\Spec(B) \to \bA^1_k$ does not contain $0$. It suffices prove that for each connected component, if the image in $\bA^1_k$ contains $0$, then it is finite over $\bA^1_k$. Also, it suffices to prove this separately for each irreducible component, so we may assume $A$ is an integral domain. The case where $\Spec(A) \to \bA^1_k$ lies over the origin is immediate, so we can assume that $\Spec(A) \to \bA^1_k$ is flat and surjective. Because $A$ is integral $A^{\Gm}=k'$ for some finite extension $k'$ of $k$, and there is a unique closed orbit in $\Spec(A)$. By the hypothesis that $0$ is in the image of $\Spec(A) \to \bA^1_k$, this closed orbit must lie over $0$, and thus it consists of a single fixed point with ideal $\mathfrak{m} \subset A$. If $\mathfrak{m}$ contained homogeneous elements of positive and negative weight, then their product would lie in $A^{\Gm} \cong k'$ and therefore must vanish (or else it would be a unit). It follows that $\mathfrak{m}$ is generated by homogeneous elements of negative weight only, and thus $A$ is non-positively graded. A graded version of Nakayama's lemma implies that because $A$ is non-positively graded and $A \otimes_{k[t]} k$ is finite dimensional, $A$ is finite as a $k[t]$-module.

Now consider the section $s : \bA^1_k \to Y$. After passing to a finite extension of $k$ and a ramified covering $\Theta_k \to \Theta_k$, this section can be lifted to $\Spec(A)$ over the formal completion $\Spec(k[\![t]\!])$. This section must lie in the piece $\Spec(A')$ of the decomposition $\Spec(A) = \Spec(A') \bigsqcup \Spec(B)$ that is finite over $\bA^1_k$. This means that the $k$-point in $Y/\Gm$ corresponding to the equivariant section $s^\circ$ lies in the image of the \'etale map $\Spec(A')/\Gm \to Y/\bG_m$, i.e., after a finite extension of $k$ the section $s^\circ$ lifts to an \emph{equivariant} section $\bA^1_k \setminus 0 \to \Spec(A')$. Because $\Spec(A')$ is finite over $\bA^1_k$, this extends uniquely to an equivariant section $\bA^1_k \to \Spec(A')$, and thus the section $s^\circ$ extends to an equivariant section $s : \bA^1_k \to Y$.

\medskip
\noindent{\textit{Completing the proof of \Cref{prop:specialization_implies_uniqueness}}}
\medskip

We have now constructed an extension of $f_1 \cup f_2$ to a morphism $g : \Theta_k^2 \to \X$. To complete the proof, we show that if $f_1$ and $f_2$ are non-degenerate and do not have common multiple under the action of $\bN^\times$, then $g$ is also non-degenerate. Because $f_1$ and $f_2$ are non-degenerate, the only point at which the maps of stabilizer groups can have a positive dimensional kernel is the origin, and \Cref{lem:contract_degenerate} implies that this can not happen if $f_1$ and $f_2$ are distinct, unless they are antipodal. Finally $f_1$ and $f_2$ are not antipodal because no two antipodal filtrations both have $\mu>0$ for a standard numerical invariant, by \Cref{defn:standard_invariant}. This shows that as points of $\iDeg(\X,p)$, $f_1$ and $f_2$ are the endpoints of a rational $1$-simplex $\simp{g}$.
\end{proof}


\subsection{Boundedness condition for the existence of HN filtrations}
\label{sect:boundedness}

Our first method for establishing the existence of HN filtrations is fairly unsophisticated. Nevertheless it can clarify the question. For any $p \in \X(k)$ and $\cU \subset \iComp(\X)$, let $\cU_{\iComp(p)} \subset \iComp(\X,p)$ denote the preimage of $\cU$ under the canonical map $\iComp(\X,p) \to \iComp(\X)$.

\begin{lem} \label{lem:principle_B}
Let $\mu : \cU \subset \iComp(\X) \to \Gamma$ be a numerical invariant on a stack $\X$, let $p \in \X(k)$ be an unstable point, and let $\bar{p} \in \X(\bar{k})$ be the extension of $p$ to an algebraic closure $k \subset \bar{k}$. If $\mu$ is upper semi-continuous, then $\mu$ obtains a maximum on $\iDeg(\bar{p})$ if and only if the following boundedness condition is satisfied:
\begin{enumerate}[label=(B${}_p$)]
\item \label{princ:B1} $\exists$ a bounded sub-fan $F_\bullet \subset \Comp(\X,p)_\bullet$ such that $\bP(F_\bullet) \cap \cU_{\iComp(p)} \subset \bP(F_\bullet)$ is closed and $\forall x \in \cU_{\iComp(p)}$ with $\mu(x)>0$, there is another point $y \in \bP(F_\bullet) \cap \cU_{\iComp(p)}$ with $\mu(y) \geq \mu(x)$.
\end{enumerate}
Furthermore, if $\mu$ satisfies \ref{princ:R} and either i) $\mu$ is upper semi-continuous, or ii) $\cU = \bP(F_\bullet)$ for some bounded embedding of fans $F_\bullet \subset \Comp(\X)_\bullet$ (\Cref{def:bounded_fan}), then \ref{princ:B1} is equivalent to the existence of an HN filtration of $p$.
\end{lem}
\begin{proof}
Note that $\iDeg(\X,\bar{p}) \to \iComp(\X,p)$ is surjective, and therefore so is $\cU_p \to \cU_{\iComp(p)}$. It follows that $\mu$ obtains a maximum on the former if and only if it obtains a maximum on the latter. The ``only if'' direction of the claim follows immediately from choosing $F_\bullet$ to be the subfan generated by a rational simplex in $\cU_p \subset \iComp(\X,p)$ containing a maximizer of $\mu$. The ``if'' direction follows from the fact that $\bP(F_\bullet) \cap \cU_{\iComp(p)}$ can be covered by a finite disjoint union of closed subsets of standard $n$-simplices of various dimensions, because $F_\bullet$ is bounded, so \Cref{lem:continuous_maximizer} implies the existence of a maximizer if $\mu$ is upper semi-continuous. One can also deduce the existence of a maximizer directly from the condition \ref{princ:R} if $\bP(F_\bullet) \cap \cU_{\iComp(p)}$ is bounded, i.e., covered by finitely many rational simplices. The fact that under condition \ref{princ:R} the existence of a maximizer of $\mu$ on $\iDeg(\bar{p})$ is equivalent to the existence of an HN filtration is \Cref{lem:princ_R}.
\end{proof}

In the context of \Cref{thm:main_stratification}, however, one needs to show additionally that only finitely many HN types appear in a family over a quasi-compact scheme. Recall that in the statement of \Cref{thm:main_stratification}, this is formalized by the condition:
\begin{enumerate}[label=(4')]
\item For any map $\xi : T \to \X$, with $T$ an affine scheme of finite type over $B$, there is a quasi-compact substack\footnote{In the statement of \Cref{thm:main_stratification}, $\Flag(\xi)$ is an algebraic space because it is assumed that $\X$ has separated inertia relative to $B$ (\Cref{prop:representable_flags}), but we will not require that in this section.} $Y \subset \Flag(\xi)$ that contains an HN filtration for every unstable finite type point of $T$.
\end{enumerate}
Our main result in this section is to establish a more flexible alternative version of this ``HN-boundedness'' criterion, \ref{princ:B2}, that applies to numerical invariants.

\begin{prop}\label{prop:principle_B2}
Let $\X$ be a stack satisfying \ref{hyp2} over a locally noetherian base stack $B$. Let $\mu : \cU \subset \iComp(\X) \to \Gamma$ be a \gls{numerical_invariant} satisfying \ref{princ:R}. Assume that either:
\begin{enumerate}[label=\alph*)]
\item $\cU \subset \iComp(\X)$ is closed and $\mu$ is upper semi-continuous, i.e., $\{x \in \cU | \mu(x)<c\}$ is open for all $c \in \Gamma$; or
\item $\cU = \bP(F_\bullet)$ for some bounded inclusion of fans $F_\bullet \subset \Comp(\X)_\bullet$.
\end{enumerate}
Then condition (4') above is equivalent to the following condition:
\begin{enumerate}[label=(B)]
\item \label{princ:B2} \textbf{HN-boundedness:} For any map from a finite type affine scheme $\xi : T \to \X$, $\exists$ a quasi-compact substack $\X' \subset \X$ such that $\forall$ finite type points $p \in T(k)$ and $f \in \Flag(p)$ with $\mu(f)>0$, there is another filtration $f' \in \Flag(p)$ with $\mu(f') \geq \mu(f)$ and $f'(0) \in \X'$.
\end{enumerate}
Furthermore: i) these conditions are equivalent to the stronger form of (4') that asserts that $Y \subset \Flag(\xi)$ contains an HN filtration for every unstable point of $T$, not just finite type points; and ii) \ref{princ:B2} implies that the function $M^\mu : |\X| \to \Gamma$ of \Cref{defn:numerical_invariant} is constructible.
\end{prop}

\begin{rem}
Although condition (4') and \ref{princ:B2} have a similar appearance, \ref{princ:B2} is a significant improvement. For instance, condition \ref{princ:B2} holds automatically if $\X$ is quasi-compact, whereas (4') does not. Also, condition (4') presupposes the existence of HN filtrations for finite type points, whereas \ref{princ:B2} does not.
\end{rem}
We shall prove \Cref{prop:principle_B2} after establishing several lemmas.

\begin{lem} \label{lem:surj_filt}
Let $\X \to \Y$ be any smooth map of stacks that satisfy \ref{hyp2}. If $\Grad(\X) \to \Grad(\Y)$ is surjective, then so is $\Filt(\X) \to \Filt(\Y)$.
\end{lem}

\begin{proof}
\Cref{lem:theta_maps} implies that given a point of $f \in |\Filt(\Y)|$ and a lift of $\agr(f) \in \Grad(\Y)$ to $\Grad(\X)$, one can lift $f$ to $\Filt(\Y)$ by solving an infinite sequence of infinitesimal lifts along the map $\X \to \Y$. Therefore the fact that $\X \to \Y$ is formally smooth implies that $\Filt(\X) \to \Filt(\Y)$ is surjective.
\end{proof}

For the next lemma, we will need the following local structure theorem generalizing the main results of \cites{alper2015luna, ahr2}, whose proof will appear in \cite{slice3}:

\begin{thm}[Special case of \cite{slice3}*{Thm.~5.1}] \label{thm:slice2}
Let $\X$ be a quasi-separated algebraic with affine stabilizers, $\X_0 \hookrightarrow \X$ a closed substack, $Y_0$ an affine scheme with $\Gm^n$-action, and $f_0 : Y_0 / \Gm^n \to \X_0$ a smooth (resp. \'etale) morphism. Then there is an affine $\GL_{s}$-scheme $Y$ for some $s\geq 0$ and a smooth (resp. \'etale) morphism $f : Y/\GL_{s} \to \X$ such that $f|_{\X_0} \cong f_0$.
\end{thm}

\begin{lem} \label{lem:theta_surj}
Let $\X$ be a quasi-compact stack satisfying \ref{hyp2}. Then there is an affine scheme $X$ with an action of $\Gm^n$ for some $n\geq 0$ along with a smooth representable surjective map $X/\Gm^n \to \X$ such that the induced map $\Filt^r(X/\Gm^n) \to \Filt^r(\X)$ is smooth, representable, and surjective for all $r$.
\end{lem}
\begin{proof}
By \Cref{lem:surj_filt} it suffices to find a smooth surjective map $X/\Gm^n \to \X$ such that $\Grad^r(X/\Gm^n) \to \Grad^r(\X)$ is surjective. Because $\X$ is quasi-compact we may assume without loss of generality that the base $B$ is quasi-compact, and choose a smooth surjective map $\Spec(R) \to B$. \Cref{cor:base_change_grad} implies that $\Grad^r_B(\X_R) \to \Grad^r_B(\X)$ is smooth and surjective. Furthermore for any stack $\Y$ over $\Spec(R)$, \Cref{cor:grad_over_affine} implies the equivalence of mapping stacks $\Grad^r_B(\Y) \simeq \Grad^r_{\Spec(R)}(\Y)$. We may therefore replace $\X$ with $\X_R$ and prove the claim under the assumption that $B=\Spec(R)$.

For any stack of the form $\X = X/\GL_n$ for a finitely presented $\GL_n$-space $X$, the map $X/\Gm^n \to X/\GL_n$ suffices, by the explicit computation of $\Grad(\X)$ relative to the scheme $\Spec(R)$ in \Cref{thm:describe_strata_quotient_GLN}. In general, we may assume that $\X$ is reduced, because $\Grad^r(\X^{\rm{red}}) \to \Grad^r(\X)$ is a surjective closed immersion, by \Cref{cor:closed_open_graded}. Let us stratify $\X$ by locally closed substacks $\X_i$ of the form $X_i/\GL_{n_i}$ where $X_i$ is quasi-affine (See \cite{kresch1999cycle}*{Prop.~3.5.9} when the base is a field and \cite{hall2014coherent}*{Prop.~8.2} in general). Then \Cref{cor:inertia_preserving} implies that for all $i$, $\Grad^r(\X_i) \simeq \Grad^r(\X) \times_{\X} \X_i$. So we have a stratification by locally closed substacks
\[
\Grad^r(\X) = \bigcup_i \Grad^r(\X_i).
\]
We let $Y_i$ be a disjoint union of open affine $\Gm^{n_i}$-equivariant subschemes that cover $X_i$. Note that $\Grad^r(Y_i/\Gm^{n_i}) \to \Grad^r(X_i/\Gm^{n_i})$ is surjective. Each $Y_i$ is finitely presented over $\Spec(R)$ by hypothesis, so we may apply \Cref{thm:slice2} to the smooth map $Y_i/\Gm^{n_i} \to X_i/\GL_{n_i} \simeq \X_i$. The result is an affine $\GL_{s_i}$-scheme $U_i$ and a smooth map $U_i / \GL_{s_i} \to \X$ such that $(U_i / \GL_{s_i})\times_{\X} \X_i \cong Y_i/\Gm^{n_i}$ over $\X$. Again \Cref{thm:describe_strata_quotient_GLN} implies that $\Grad^r(U_i / \Gm^{s_i}) \to \Grad^r(U_i / \GL_{s_i})$ is surjective, so it follows from \Cref{cor:inertia_preserving} that the composition
\[
\Grad^r(U/\Gm^{s_i})\times_{\X} {\X_i} \to \Grad^r(Y_i/\Gm^{n_i}) \to \Grad^r(X_i/\Gm^{n_i}) \to \Grad^r(\X_i)
\]
is surjective. Choosing $n$ larger than $s_i$ for all $i$, we let $X = \bigsqcup_i U_i \times \Gm^{n-s_i}$, and the resulting map $X/\Gm^n \to \X$ satisfies the claim.
\end{proof}

It will be convenient to further reduce matters to a quotient stack over a scheme:

\begin{lem} \label{lem:theta_surj_2}
Let $\X$ be a quasi-compact stack satisfying \ref{hyp2}. For any smooth map $\Spec(R) \to B$ whose image contains the image of $\X \to B$, one can construct a stack $X/\Gm^n \to \Spec(R)$, where $X$ is an affine $R$-scheme of finite presentation, and a smooth representable surjective map $X/\Gm^n \to \X$ over $B$ such that the induced map $\Filt_R^r(X/\Gm^n) \to \Filt_B^r(\X)$ is smooth, representable, and surjective for all $r>0$.
\end{lem}
\begin{proof}
Consider the smooth representable map $X/\Gm^n \to \X$ over $B$ from \Cref{lem:theta_surj}. If we let $X'/\Gm^n = (X/\Gm^n)\times_B \Spec(R)$, then $\Filt_R(X'/\Gm^n) \simeq \Spec(R) \times_B \Filt_B(X/\Gm^n)$ by \Cref{lem:base_change_mapping}, so $\Filt_R^r(X'/\Gm^n) \to \Filt_B(X/\Gm^n)$ is smooth and surjective as well.
\end{proof}

\begin{lem} \label{lem:qc_hn}
Let $\X$ be a quasi-compact stack satisfying \ref{hyp2}. For any map from a quasi-compact space $\xi : S \to \X$ and any union of connected components $Y \subset \Flag(\xi)$, there is a quasi-compact substack $Y' \subset Y$ such that
\[
\ev_1(Y') = \ev_1(Y) \subset |S|.
\]
\end{lem}
\begin{proof}
Let $X/\Gm^n \to \X$ be a smooth surjective representable map such that $\Filt_R(X/\Gm^n) \to \Filt_B(\X)$ is surjective, as constructed in \Cref{lem:theta_surj_2}. Let $\xi' : S' \to X/\Gm^n$ be the base change of $\xi$ along the map $\pi : X/\Gm^n \to \X$. Because $\pi$ is representable, we know that $S'$ is a space. We can therefore consider the commutative (non-cartesian) diagram
\[
\xymatrix{
\Flag_R(\xi') \ar[r] \ar[d]^{\ev_1} & \Flag_B(\xi) \ar[d]^{\ev_1} \\ S' \ar[r]^\pi & S
}
\]
where the horizontal maps are smooth and surjective. If we denote the preimage of $Y \subset \Flag_B(\xi)$ as $Y' \subset \Flag_R(\xi')$, then $\pi(\ev_1(Y')) = \ev_1(Y)$.

It therefore suffices to prove the claim for a stack of the form $X/\Gm^n$ over an affine scheme $\Spec(R)$. The claim is now a consequence of the explicit description of the stack of filtrations $\Filt_R(X/\Gm^n)$ given in \Cref{thm:describe_strata_quotient_GLN}. The connected components of $\Filt(X/\Gm^n)$ are indexed by a one parameter subgroup $\lambda : \Gm \to \Gm^n$ and a connected component of the fixed locus $Z \subset X^{\lambda,0}$. The connected components of $\Filt(X/\Gm^n)$ are $X^{\lambda,+}_{Z} / \Gm^n$, where $X^{\lambda,+}_{Z}$ is the preimage of $Z$ under the projection $X^{\lambda,+} \to X^{\lambda,0}$. Note that $\ev_1$, corresponding to the canonical map $X^{\lambda,+} \to X$, is a locally closed immersion on each connected component, and despite $\Filt(X/\Gm^n)$ having infinitely many connected components, there are only finitely many locally closed substacks of $X/\Gm^n$ arising in this way.

Now any map $\xi : S \to X / \Gm^n$ can be presented as $S'/\Gm^n \to X/\Gm^n$ for some $\Gm^n$-equivariant map $S' \to X$. The connected components of $\Flag(\xi)$ are of the form $Y' / \Gm^n$, where $Y'$ is a connected component of $S' \times_X X^{\lambda,+}_{Z}$. There are finitely many locally closed subschemes of $S$ arising in this way, so given an infinite collection of connected components of $\Flag(\xi)$, finitely many suffice to cover their image in $S$.
\end{proof}

\begin{cor} \label{cor:constructible_filtrations}
If $\X$ is a quasi-compact stack satisfying \ref{hyp2}, then for any map from a quasi-compact quasi-separated algebraic space $\xi : S \to \X$ and any union of connected components $Y \subset \Flag(\xi)$, $\ev_1(|Y|) \subset |S|$ is constructible.
\end{cor}
\begin{proof}
This is an immediate consequence of \Cref{lem:qc_hn} and the fact that $\ev_1$ is locally of finite type and quasi-separated.
\end{proof}

\begin{lem} \label{lem:qc_bounded_relative_comp}
Let $\X$ be a quasi-compact stack satisfying \ref{hyp2}. Then for any map $\xi : S \to \X$ from a quasi-compact quasi-separated algebraic space $S$, the \gls{component_fan} $\Comp(\xi)_\bullet$ is bounded.
\end{lem}
\begin{proof}
From the base change properties of mapping stacks, \Cref{lem:base_change_mapping}, we may replace $\X$ with $\X \times_B S$ and therefore assume that $B = S$ and the map $\xi : S \to \X$ is a section of the structure map $\X \to S$. Let $\pi : X/\Gm^n \to \X$ be a smooth surjective representable map such that $\Filt^r(X/\Gm^n) \to \Filt^r(\X)$ is surjective, as constructed in \Cref{lem:theta_surj}. Let $\xi' : S' \to X/\Gm^n$ be the base change of $\xi$ along the $\pi$. Because $\pi$ is representable, we know that $S'$ is a space, and as in the proof of \Cref{lem:qc_hn} we have a smooth surjection $\Flag^r(\xi') \to \Flag^r(\xi)$ for all $r$. It follows that $\Comp(\xi')_\bullet \to \Comp(\xi)_\bullet$ is surjective, so we may assume that $\X = X/\Gm^n$ over a base space $S$ and that $\xi$ is a section of the map $X/\Gm^n \to S$.

We use the explicit description from \Cref{thm:describe_strata_quotient_GLN}:
\[
\Flag^r(\xi) = \bigsqcup_{\psi \in \Hom(\Gm^r,\Gm^n)} \xi^{-1}(X^{\psi,+}).
\]
Note that each connected component of $\Flag^r(\xi)$ is a locally closed subspace of $S$ under the projection $\Flag^r(\xi) \to S$. Also, because finitely many locally closed subspaces arise as blades $X^{\psi,+} \hookrightarrow X$ associated to some homomorphism $\psi : \Gm^r \to \Gm^n$, only finitely many locally closed subspaces of $S$ arise as connected components of $\Flag^r(\xi)$. It follows that there is a finite set of geometric points $p_i \in S(\bar{k}_i), i=1,\ldots,N$ such that
\[
\bigsqcup_i \Deg(X/\Gm^n,\xi(p_i))_\bullet \to \Comp(X/\Gm^n,\xi)_\bullet
\]
is surjective. By \Cref{lem:almost_toric} the left hand side is bounded, so the right hand side is bounded as well.
\end{proof}

\begin{proof}[Proof of \Cref{prop:principle_B2}]

It is immediate that (4') implies \ref{princ:B2}, so we will prove the other direction. We must show that given a map $\xi : T \to \X$ from a finite type $B$-scheme $T$, there is a quasi-compact substack $Y \subset \Flag(\xi)$ that contains an HN filtration for every unstable point of $T$.

We first prove the claim when $\X$ is quasi-compact, in which case the condition \ref{princ:B2} holds tautologically. In this case \Cref{lem:qc_bounded_relative_comp} implies that \ref{princ:B1} holds for any point $p \in |\X|$, not just finite type points, by letting $F_\bullet = \Comp(\X,p)_\bullet$. Hence by \Cref{lem:principle_B} every unstable point has an HN filtration.

We shall use Noetherian induction to construct $Y \subset \Flag(\xi)$. Consider the class of closed subsets of $T$
\[
\cC = \left\{Z \subset |T| \left| \begin{array}{c} \exists \text{ quasi-compact open subset } Y \subset |\Flag(\xi)| \\ \text{containing an HN filtration for every point in } Z \end{array} \right. \right\}.
\]
Our goal is to show that $|T| \in \cC$. For the inductive argument it suffices to consider an irreducible $Z \in \cC$, because $\cC$ is closed under union of subsets. So we consider a closed irreducible subset $Z \subset |T|$, and assume that $Z' \in \cC$ for every proper closed subset $Z' \subsetneq Z$.

Let $U \subset |\Flag(\xi)|$ denote the union of all connected components that contain an HN filtration for some point in $Z$. By \Cref{cor:constructible_filtrations} then $\ev_1(U) \cap Z \subset |T|$, the set of unstable points in $Z$, is constructible. If its closure is a proper subset of $Z$, we can apply the inductive hypothesis to conclude that $Z \in \cC$.

Otherwise, we may assume that $\ev_1(U) \cap Z$ is dense in $Z$. Let $V \subset |\Flag(\xi)|$ be a quasi-compact connected open subset containing an HN filtration of the generic point of $Z$, which we know exists because $\X$ is quasi-compact. Let $U^{\mu > \mu(V)} \subset U$ be the union of all connected components on which $\mu > \mu(V)$. We know from \Cref{cor:constructible_filtrations} that the image $\ev_1(U^{\mu > \mu(V)}) \cap Z$ is constructible, and it does not contain the generic point because $M^\mu=\mu(V)$ at the generic point of $Z$. It follows that
\[
Z' := \overline{Z \setminus \ev_1(V)} \cup \overline{\ev_1(U^{\mu > \mu(V)}) \cap Z} \subset Z
\]
is a closed proper subset. The inductive hypothesis provides a quasi-compact open subspace $Y \subset |\Flag(\xi)|$ containing an HN filtration for each point of $Z'$, and by construction $V$ contains an HN filtration for every point of $Z \setminus Z'$. Thus $Y \cup V$ is a quasi-compact open subset containing an HN filtration for every point of $Z$, and $Z \in \cC$.

\medskip
\noindent \textit{Completing the proof for non-quasi-compact $\X$:}
\medskip

Given $\xi : T \to \X$ as above, condition \ref{princ:B2} implies that we can find a quasi-compact open substack $\X' \subset \X$ such that: 1) the map $\xi$ factors through $\X'$; and 2) for the purposes of finding HN filtrations of unstable finite type points, it suffices to consider filtrations in the open substack $\Filt(\X') \subset \Filt(\X)$. We may therefore apply the proposition for the quasi-compact stack $\X'$ to find a quasi-compact substack $Y \subset \Flag_{\X'}(\xi) \subset \Flag_{\X}(\xi)$ that contains an HN filtration for every unstable finite type point of $T$. This already establishes condition (4'), but we claim further that $Y$ contains an HN filtration for all unstable points of $T$. Note that this immediately implies the constructibility of $M^\mu|_T$, because $Y$ is quasi-compact.

Suppose $p \in |T|$ is a non-finite type point whose image in $\X$ is unstable, and let $S \subset \Flag_{\X}(\xi)$ be a quasi-compact connected locally closed substack containing a destabilizing filtration of $p$. $Y$ contains an HN filtration for every unstable finite type point of $T$, and every point in the image of $\ev_1 : S \to T$ is unstable, so the fiber product $Y \times_T S \to S$ is a finite type morphism that is surjective on finite type points. It follows that $Y \times_T S \to S$ is surjective, and in particular $p$ lies in the image of $\ev_1 : |Y| \to |T|$ as well.

We claim that if we let $c \in \Gamma$ be the largest value of $\mu$ obtained among all connected components of $Y$ that meet $\ev_1^{-1}(p)$, and then choose a point in $f \in \ev_1^{-1}(p)$ with $\mu(f)=c$, then $f$ will be an HN filtration for $p$. To see this, assume $f' \in |\Flag_{\X}(p)|$ were a filtration with $\mu(f')>c$. We could then let $S \subset \Flag_\X(\xi)$ be a locally closed subspace containing $f'$, and let $Y^{\mu > c} \subset Y$ be the open and closed substack on which $\mu>c$. As in the previous paragraph, $S \times_T Y^{\mu>c} \to S$ would be surjective, hence there would be a point in $Y^{\mu>c}$ mapping to $p$, contradicting the maximality of $c$.

\end{proof}


\subsection{Main theorem on the existence of \texorpdfstring{$\Theta$}{Theta}-stratifications}

We now reformulate \Cref{thm:main_stratification} by incorporating the simplifications that we have established above.

\begin{thm} \label{thm:main_improved}
Let $\X$ be a stack satisfying \ref{hyp3} over a locally noetherian base stack $B$, and let $\mu : \cU \subset \iComp(\X)\to \Gamma$ be a \gls{standard_numerical_invariant} satisfying \ref{princ:R}. Assume that either:
\begin{enumerate}[label=\alph*)]
\item $\cU \subset \iComp(\X)$ is closed and $\mu$ is upper semi-continuous, i.e., $\{x \in \cU | \mu(x)<c\}$ is open for all $c \in \Gamma$; or
\item $\cU = \bP(F_\bullet)$ for some bounded inclusion of fans $F_\bullet \subset \Comp(\X)_\bullet$.
\end{enumerate}
Then $\mu$ defines a weak $\Theta$-stratification if and only if it satisfies the simplified HN-specialization property \ref{princ:S} of \Cref{thm:main_stratification} and the boundedness condition \ref{princ:B2} of \Cref{prop:principle_B2}.
\end{thm}
\begin{proof}
We verify the hypotheses of \Cref{thm:main_stratification}. The open strata property holds because $\mu$ is locally constant (see \Cref{simp:locally_constant}). The HN-boundedness condition combined with the existence of HN filtrations for finite type points is equivalent to the condition \ref{princ:B2} by \Cref{prop:principle_B2}. The uniqueness of the HN-filtration up to the action of $\bN^\times$ is \Cref{prop:specialization_implies_uniqueness}.

All that is left to verify is the HN-consistency condition, namely that if $f : \Theta_k \to \X$ is an HN filtration, then $M^\mu(f(0)) \leq \mu(f)$. Let $g : \Theta_k \to \X$ be the HN filtration of $f(0)$, whose existence and uniqueness (up to $\bN^\times$-action) we have already established. The uniqueness of $g$ implies that it is isolated in its fiber for the map $\ev_1 : \Filt(\X) \to \X$, hence it is fixed by the action of $\Gm$ on that fiber coming from the canonical $\Gm$-action on $f(0)$. After replacing $k$ with an extension, we may thus lift $g$ to a filtration $\tilde{g} : \Theta_k \to \Grad(\X)$ with $\tilde{g}(1) \cong \agr(f)$ and $u \circ \tilde{g} \cong g$. The Recognition Theorem, \Cref{thm:recognition}, implies that $\agr(f)$ is graded-semistable, so by \Cref{defn:graded_semistable} we have
\[
\mu(f) = M_{\Grad(\X)}^{\mu\circ u_\ast}(\tilde{g}) = \mu(g).
\]
This establishes the claim, because $\mu(g) = M^{\mu}(f(0))$ by hypothesis.
\end{proof}

In practice, the only hypotheses of \Cref{thm:main_improved} that are challenging to check are \ref{princ:S} and \ref{princ:B2}, because the others are all local, in the sense that they only depend on the functional form of the numerical invariant restricted to a rational simplex $\simp{} \to \iComp(\X)$.

\begin{ex} \label{ex:main_hypotheses_associated_invariant}
Any numerical invariant on a stack $\X$ associated in \Cref{defn:induced_invariant} to a class $\ell \in H^2(\X;\bQ)$ and a positive definite class $b \in H^4(\X;\bQ)$ is standard by \Cref{lem:induced_invariant} and satisfies \ref{princ:R} by \Cref{lem:induced_invariants_R}. Also, $\cU = \iComp(\X)$ for such a numerical invariant, so the hypothesis (b) holds tautologically.
\end{ex}


\section{Beyond geometric invariant theory}
\label{sect:beyond_git}

\Cref{thm:main_improved} shows that, under some mild hypotheses, a numerical invariant $\mu$ on a stack $\X$ defines a weak $\Theta$-stratification of $\X$ if and only if two conditions hold: \ref{princ:S} and \ref{princ:B2}. In this section we investigate two conditions that automatically imply \ref{princ:S} and are easier to check in many examples: the first is the condition that the stack $\X$ is $\Theta$-reductive  (\Cref{defn:reductive_stack}), and the second is the condition that the numerical invariant $\mu$ is strictly $\Theta$-monotone (\Cref{defn:theta_monotone}).

This leads to an alternative, and arguably more user-friendly, criterion for $\Theta$-stratifications in \Cref{T:monotone_stratifications}. As an application, we show in \Cref{prop:strat_proper_over_reductive} the existence of weak $\Theta$-stratifications defined by certain polynomial valued numerical invariants that are positive relative to a proper map $\X \to \Y$.

After the notion of $\Theta$-reductivity was introduced in the first version of this paper, it was discovered in \cite{AHLH} that this condition is one of a short list of necessary and sufficient conditions for a stack to have a good moduli space. We use this to show here that under suitable hypotheses, strict $\Theta$-monotonicity and a related notion of strict {\textsf{S}}-monotonicity imply the existence of a separated good moduli space for $\X^{\rm ss}$. As another interesting point of interaction between the theory of $\Theta$-stratifications and good moduli spaces, we formulate an intrinsic version of the main results regarding ``variation of GIT quotients'' in \Cref{thm:main_GIT}

\subsection{\texorpdfstring{$\Theta$}{Theta}-reductive stacks}
\label{sect:theta_reductive}


\begin{defn} \label{defn:reductive_stack}
If $\cC$ denotes a class of valuation rings over $B$ (e.g. discrete valuation rings, complete DVR's, DVR's essentially of finite type over $B$,...) we say that a stack $\X$ over $B$ is \emph{$\Theta$-reductive} with respect to $\cC$ if the morphism $\ev_1 : \Filt_B(\X) \to \X$ satisfies the valuative criterion for properness for all valuation rings in $\cC$. We say simply that $\X$ is $\Theta$-reductive if $\cC$ is all valuation rings.
\end{defn}

By the valuative criterion for properness in \Cref{defn:reductive_stack}, we mean the condition that for any valuation ring $R$ with function field $K$, any commutative diagram of the form
\[
\xymatrix{\Spec(K) \ar[d] \ar[r] & \Filt(\X) \ar[d]^{\ev_1} \\ \Spec(R) \ar[r] \ar@{-->}[ur] & \X},
\]
has a unique dotted arrow making the diagram commute. This is the ``usual" valuative criterion for maps of algebraic spaces. In words, it states that for any family over $\Spec(R)$, any filtration of the generic point extends uniquely to a filtration of the family.

Another way to state this property: For any valuation ring $R$ over $B$ in $\cC$ and any map $\Theta_R - \{(0,0)\} \to \X$ relative to $B$ there exists a unique extension to a map $\Theta_R \to \X$ relative to $B$, where $(0,0) \in \Theta_R$ denotes the closed point and $\Theta_R$ is regarded as a $B$-stack via the composition $\Theta_R \to \Spec(R) \to B$. In other words, every commutative diagram of the form
\[
\xymatrix{\Theta_R \setminus \{(0,0)\} \ar[r] \ar[d] & \Theta_R \ar[r] \ar@{-->}[ld] & \Spec(R) \ar[d] \\
\X \ar[rr] & & B }
\]
has a unique dotted arrow making the diagram commute.

\begin{rem}
Note that if $\X$ is an algebraic stack satisfying \ref{hyp3} over a locally noetherian base $B$, and if for any $\xi : \Spec(R) \to \X$, $\Flag(\xi)$ has quasi-compact irreducible components, then the weakest form of the valuative criterion, where $\cC$ is the class of DVR's essentially of finite type over $B$, already implies that the irreducible components of $\Flag(\xi)$, which are quasi-separated algebraic spaces locally of finite type over $\Spec(R)$ (\Cref{prop:representable_flags}), are proper \cite{stacks-project}*{Tag 0ARI}. Hence all other variants of the valuative criterion hold as well.
\end{rem}

\begin{rem} \label{rem:reductive_locality}
\Cref{defn:reductive_stack} is local in the sense that a stack $\X$ satisfying \ref{hyp2} is $\Theta$-reductive over $B$ if and only if its base change along every map $\Spec(R) \to B$ is $\Theta$-reductive (either relative to $R$ or in the absolute sense). This follows from the canonical identifications $\Spec(R) \times_B \Filt_B(\X) \cong \Filt_{\Spec(R)}(\X_R) \cong \Filt_{\Spec(\bZ)}(\X_R)$ (see \Cref{cor:change_of_base_filt} and \Cref{cor:representable_base_change}).
\end{rem}

\begin{warning} \label{warn:theta_reductive}
When $B$ is an algebraic space, then any map $\Theta_R \to B$ factors uniquely through a map $\Spec(R) \to B$, because $\Theta_R \to \Spec(R)$ is a good moduli space. It follows that \Cref{defn:reductive_stack} agrees with the relative notion of $\Theta$-reductivity of a morphism introduced in \cite{AHLH}*{Defn.~3.10}, applied to the morphism $\X \to B$. More generally, when $B$ is an algebraic stack with positive-dimensional automorphism groups, \Cref{defn:reductive_stack} is weaker than \cite{AHLH}*{Defn.~3.10}. The latter is not local over $B$ in the sense of \Cref{rem:reductive_locality}.
\end{warning}

\begin{ex}
The stack $\X = \pt/G$ is $\Theta$-reductive for any reductive group $G$ over a field $k$. Indeed, the formation of $\Filt(\pt/G)$ commutes with base change to an algebraic closure $\bar{k} / k$, so we may assume $G$ is split. Applying \Cref{thm:describe_strata_global_quotient}, we see that the fiber of $\ev_1 : \Filt(\pt/G) \to \pt/G$ over the point is an infinite disjoint union of generalized flag varieties $G/P_\lambda$, which are proper.
\end{ex}

\begin{ex} \label{ex:bad_stratum}
Let $V = \Spec (k[x,y,z])$ be a linear representation of $\Gm$ where $x,y,z$ have weights $-1,0,1$ respectively, and let $X = V - \{0\}$ and $\X = X / (\Gm)_k$. The fixed locus is the punctured line $Z = \{x=z=0\} \cap X$, and the connected component of $\Filt(\X)$ corresponding to the cocharacter $\lambda(t)=t$ is the quotient $S/(\Gm)_k$ where
\[
S = \{(x,y,z) | z = 0 \text{ and } y \neq 0\}
\]
$S \subset X$ is not closed. Its closure contains the points where $x \neq 0$ and $y = 0$. These points would have been attracted by $\lambda$ to the missing point $\{0\} \in V$. It follows that $S / (\Gm)_k \to X/(\Gm)_k$ is not proper, and hence $X/(\Gm)_k$ is not $\Theta$-reductive.
\end{ex}

\begin{lem} \label{lem:affine_over_reductive}
Let $f : \X \to \Y$ be a representable affine morphism of stacks. If $\Y$ is $\Theta$-reductive with respect to some class of valuation rings $\cC$, then so is $\X$.
\end{lem}

\begin{proof}
\Cref{prop:affine_map} states that the canonical map $\Filt(\X) \to \Filt(\Y) \times_\Y \X$ is a closed immersion, and if $\Y$ is $\Theta$-reductive, then $\Filt(\Y) \times_\Y \X \to \X$ satisfies the valuative criterion for properness with respect to valuation rings in $\cC$. It follows that the composition $\Filt(\X) \to \X$ satisfies the valuative criterion for properness with respect to valuation rings in $\cC$.
\end{proof}

\begin{ex}
\Cref{lem:affine_over_reductive} implies that any stack of the form $V/G$, where $G$ is a reductive group over a field $k$ acting on an affine $k$-scheme $V$, is $\Theta$-reductive.
\end{ex}

For further examples of $\Theta$-reductive stacks, see \Cref{thm:monotone_moduli_spaces} and \Cref{prop:theta_reductive_S_complete} below.

The main advantage of a $\Theta$-reductive stack is that the specialization condition \ref{princ:S} holds for any numerical invariant, because for any family over a DVR, any filtration over the generic point extends to a filtration of the family. One consequence is that HN filtrations, if they exist, are unique (\Cref{prop:specialization_implies_uniqueness}). In addition, our main existence theorem for $\Theta$-stratifications becomes:
\begin{thm} \label{thm:theta_reductive_stratifications}
Let $\X$ be a $\Theta$-reductive stack and $\mu$ a numerical invariant satisfying the hypotheses of \Cref{thm:main_improved}. Then $\mu$ defines a weak $\Theta$-stratification of $\X$ if and only if \ref{princ:B2} holds.
\end{thm}

\begin{ex}
Recall from \Cref{ex:main_hypotheses_associated_invariant} that for the numerical invariant associated to a class $\ell \in H^2(\X;\bQ)$ and a positive definite $b \in H^4(\X;\bQ)$ in \Cref{defn:induced_invariant}, all of the hypotheses of \Cref{thm:theta_reductive_stratifications} hold automatically. Also, the condition \ref{princ:B2} holds tautologically if $\X$ is quasi-compact. So the theorem implies that if $\X$ satisfies \ref{hyp3} over a locally noetherian base $B$, and $\X$ is $\Theta$-reductive and quasi-compact, then any such numerical invariant defines a $\Theta$-stratification of $\X$.
\end{ex}

Another consequence of $\Theta$-reductivity is the following:

\begin{lem} \label{lem:theta_reductive_convex}
Let $\X$ be a $\Theta$-reductive stack satisfying \ref{hyp3}, and let $p \in \X(k)$. Then any two distinct rational points in $\iDeg(p)$ are either antipodal or connected by a unique rational 1-simplex in $\iDeg(p)$. In particular any standard numerical invariant $\mu : \cU \subset \iComp(\X) \to \Gamma$ on $\X$ is strictly quasi-concave, i.e., for any $p \in \X(k)$ the subset $\cU_p^{\mu>0}$ is convex.
\end{lem}
\begin{proof}
The method of proof of \Cref{prop:specialization_implies_uniqueness} applies verbatim to show that for any two filtrations $f_1,f_2$ in $\iDeg(\X,p)$, the map $f_1 \cup f_2 : \bA^2_k \setminus 0 / \Gm^2 \to \X$ extends uniquely to a map $\Theta^2_k \to \X$, and the latter is non-degenerate if $f_1$ and $f_2$ are not antipodal or equivalent up to the action of $\bN^\times$. This exhibits a rational simplex in $\iDeg(\X,p)$ whose endpoints are $f_1$ and $f_2$. If $\mu$ is a standard numerical invariant, then this rational simplex has to lie in $\cU_p^{\mu>0}$.
\end{proof}


\subsection{\texorpdfstring{$\Theta$}{Theta}-monotone numerical invariants}
\label{sect:proper_maps}

In this section we formulate a replacement for the specialization condition \ref{princ:S} that is easier to check in practice. Recall that a reduced rational curve $C$ over a field $k$ with non-trivial $\Gm$-action has two fixed points, $\{0\} := \lim_{t\to 0} t\cdot x$ and $\{\infty\} := \lim_{t\to \infty} t \cdot x$ for a general point $x$.

\begin{defn} [$\Theta$-monotone numerical invariant] \label{defn:theta_monotone}
Let $\mu : \cU \subset \iComp(\X) \to \Gamma$ be a numerical invariant on a stack $\X$ over $B$. We say that $\mu$ is \emph{$\Theta$-monotone} if for any discrete valuation ring $R$ over $B$ with fraction field $K$ and any map $f : \Theta_R \setminus \{(0,0)\} \to \X$ over $B$ such that $f_K := f|_{\Theta_K}$ is non-degenerate, $f_K \in \cU$, and $\mu(f_K) \geq 0$, there is a commutative diagram of the form
\[
\xymatrix{
& \cW \ar[r]_{\tilde{f}} \ar[d]^p & \X \\
\Theta_R \setminus \{(0,0)\} \ar@{(=}[r] \ar[ur] \ar@/^30pt/[urr]^f & \Theta_R & 
}
\]
such that:
\begin{enumerate}
\item  $\cW$ is a reduced and irreducible algebraic stack, the morphism $\tilde{f}$ has quasi-finite relative inertia, the morphism $p$ is proper and relatively representable by Deligne-Mumford stacks, and $p$ is an isomorphism over $\Theta_R \setminus \{(0,0)\}$. (In particular, the lift $\Theta_R \setminus \{(0,0)\} \to \cW$ is uniquely determined.)
\item For any non-degenerate graded point $(B\Gm)_k \to \cW$, the filtration obtained by composition $\Theta_k \to (B\Gm)_k \to \cW \to \X$ lies in $\cU \subset \iComp(\X)$.
\item For any commutative diagram of the form
\[
\xymatrix{\bP^1_k / \Gm \ar[r]^\phi \ar[d] & \cW \ar[d] \\ (B\Gm)_k \ar[r] & \Theta_R },
\]
where the lower horizontal morphism is a \emph{positive} multiple of the canonical graded point $\{(0,0)\} / \Gm \hookrightarrow \Theta_R$ and the induced map $\bP^1_k / \Gm \to (B\Gm)_k \times_{\Theta_R} \cW$ is finite, one has an inequality for the induced graded points of $\X$,
\begin{equation}\label{E:monotonicity}
\mu(\tilde{f} \circ \phi_0 : \{0\}/\Gm \to \X) \leq \mu(\tilde{f} \circ \phi_\infty : \{\infty\}/\bG_m \to \X).
\end{equation}
\end{enumerate}
We say that $\mu$ is \emph{strictly $\Theta$-monotone} if $\cW$ and $\tilde{f}$ can always be chosen so that the inequality \eqref{E:monotonicity} is strict.
\end{defn}

\begin{ex} \label{E:reductive_implies_monotone}
Any numerical invariant on a $\Theta$-reductive stack $\X$ is strictly $\Theta$-monotone, because any morphism $(\bA^1_R \setminus (0,0)) /\Gm \to \X$ over $B$ extends uniquely to a morphism $\bA^1_R/\Gm \to \X$ over $B$. If we let $\cW = \Theta_R$, then conditions (2) and (3) hold vacuously.
\end{ex}

Our main result is the following strict generalization of \Cref{thm:theta_reductive_stratifications}:

\begin{thm} \label{T:monotone_stratifications}
Let $\X$ be an algebraic stack satisfying \ref{hyp3} over a locally noetherian base stack $B$, and let $\mu$ be a numerical invariant on $\X$ satisfying the hypotheses of \Cref{thm:main_improved}. If $\mu$ is strictly $\Theta$-monotone, then it satisfies condition \ref{princ:S}. In particular, it defines a weak $\Theta$-stratification if and only if it satisfies condition \ref{princ:B2} of \Cref{prop:principle_B2}.
\end{thm}

\subsubsection{Some lemmas supporting \Cref{defn:theta_monotone}.}

In order to explain and eventually apply the notion of $\Theta$-monotonicity, we will establish some preliminary results.

\begin{lem} \label{lem:surface_structure}
Let $X$ be a regular noetherian Nagata algebraic space with $\dim(X)=2$, let $T$ be a split torus acting on $X$, and let $\Sigma \to X$ be a $T$-equivariant birational proper $T$-equivariant map from a reduced algebraic space $\Sigma$.
\begin{enumerate}
\item There is a finite set of closed $T$-fixed points $\{p_1,\ldots,p_N\} \subset |X|$ such that $\Sigma \to X$ is an isomorphism over $X\setminus \{p_1,\ldots,p_N\}$.\\
\item For each $p \in \{p_1,\ldots,p_N\}$ the fiber $\Sigma_{p} = C_1 \cup \cdots \cup C_{n}$ over $p$ is a connected union of rational curves over $k(p)$.\\
\item For any cocharacter $\lambda : \Gm \to T$ that acts on the tangent space $T_p X$ with weights $a \leq 0 < b$, $\lambda$ acts non-trivially on each component $C_i$, and one can order the $\lambda$-fixed locus $\Sigma_p^{\Gm} = \{x_0,\ldots,x_n\}$ such that for any non-fixed $x \in C_i$, $\lim_{t\to 0} \lambda(t) \cdot x = x_{i-1}$ and $\lim_{t\to \infty} \lambda(t) \cdot x = x_i$.\\
\item If $D \hookrightarrow X$ is a $T$-equivariant divisor that is regular at $p$ and such that $T_p D \subset T_p X$ is the eigenspace of $T_p X$ of weight $a$ (respectively $b$), then the strict transform $D' \hookrightarrow \Sigma$ of $D$ contains $x_0$ (respectively $x_n$).

\end{enumerate}
\end{lem}

\begin{proof}
The construction of \cite{stacks-project}*{Tag 0BHM}, which can be carried out $T$-equivariantly, provides a finite set of closed points (which must be $T$-invariant) $\{p_1,\ldots,p_N\} \subset |X|$ satisfying (1) and a $T$ equivariant map $\Sigma' \to \Sigma$ relative to $X$. $\Sigma'$ is constructed via a sequence of blowups at closed $T$-invariant points lying over some $p_j$, starting with $X$. The irreducibility of $\Sigma'$ implies that $\Sigma' \to \Sigma$ is birational and surjective. The claims for $\Sigma$ can now be reduced to the equivalent claims for $\Sigma'$. For $\Sigma'$ they follow from an inductive argument based on the following observation:

If $X$ is a regular $2$-dimensional noetherian algebraic space with a $\bG_m$-action and $x \in X^{\bG_m}$ is a fixed point such that the tangent space $T_x X$ has eigenspaces with two $\bG_m$-weights $a \leq 0 < b$, then the blowup $\tilde{X} \to X$ at $X$ has an exceptional fiber which is $\bP^1$ with a non-trivial $\bG_m$ action. If $p$ is a general point in the exceptional divisor, then at the fixed point $0 = \lim_{t \to 0} t \cdot p$, the tangent space $T_{0} \tilde{X}$ has weights $a \leq 0 < b-a$, and at the fixed point $\infty = \lim_{t\to \infty} t \cdot p$, the tangent space $T_\infty \tilde{X}$ has weights $a-b < 0 < b$.
\end{proof}

We also consider the generalization of \Cref{lem:surface_structure} in which $\Sigma$ is a proper Deligne-Mumford stack over $X$. One formalizes this as a proper birational morphism of algebraic stacks $\cW \to X/T$ with finite relative inertia, where $\cW$ is regarded as the quotient of $\Sigma := \cW \times_{X/T} X$ by $T$.
\begin{lem} \label{lem:relative_coarse_space}
Let $X$ and $T$ be as in \Cref{lem:surface_structure}, let $\cW$ be a reduced algebraic stack, and let $\cW \to X/T$ be a proper birational morphism with finite relative inertia. Then there is a reduced algebraic space $\Sigma$, a $T$-equivariant proper birational morphism $\Sigma \to X$, and a morphism $\cW \to \Sigma/T$ over $X/T$ that is universal for maps from $\cW$ to algebraic stacks that are representable by algebraic spaces relative to $X/T$.
\end{lem}
\begin{proof}
Consider the level-wise base change of $\cW \to X/T$ to the action groupoid $T \times X \rightrightarrows X$. This gives a groupoid object in algebraic stacks $\cW_1 \rightrightarrows \cW_0$ along with a cartesian map of groupoid stacks $(\cW_1 \rightrightarrows \cW_0) \to (T \times X \rightrightarrows X)$, i.e., both the source and target map $\cW_1 \to \cW_0$ are cartesian over the corresponding map $T \times X \to X$. Let $W_i$ be the coarse moduli space for $\cW_i$, for $i=1,2$. By the properties of functoriality and stability under base change for coarse moduli spaces, one gets a groupoid algebraic space $W_1 \rightrightarrows W_0$ along with a cartesian map of groupoid algebraic spaces $(W_1 \rightrightarrows W_0) \to (T \times X \rightrightarrows X)$. This data corresponds to a canonically defined $T$-action on $\Sigma := W_0$, such that the projection $\Sigma \to X$ is $T$-equivariant. The resulting map to the quotient stack $\cW \to [W_0/W_1] \cong \Sigma / T$ is the universal morphism -- we leave it to the reader to verify that it has the properties stated in the lemma.
\end{proof}

Finally, we have

\begin{lem} \label{lem:stacky_surface_structure}
Let $X$ and $T$ be as in \Cref{lem:surface_structure}, and let $\cW \to X/T$ be a proper birational morphism of algebraic stacks with finite relative inertia. Then there is a finite set of closed $T$-fixed points $\{p_1,\ldots,p_N\} \subset |X|$ such that $\cW \to X/T$ is an isomorphism over $( X\setminus \{p_1,\ldots,p_N\} ) / T$. Furthermore, if $p \in \{p_1,\dots,p_N\}$ and $\lambda : \Gm \to T$ is a cocharacter that acts with weights $a \leq 0 < b$ on the tangent space $T_p X$, then there is a commutative diagram of the form
\[
\xymatrix{
C / \Gm \ar[r]^-{\phi} \ar[d] & \cW \ar[d] \\
(B\Gm)_k \ar[r]^\gamma & X/T},
\]
such that:
\begin{enumerate}
\item $C = C_1 \cup \cdots C_n$ is a connected union of curves, where each $C_i$ is isomorphic to $\bP^1_k$ with a non-trivial $\Gm$ action. The fixed points can be ordered
\[
(C_1 \cup \cdots \cup C_n)^{\Gm} = \{x_0,\ldots, x_n\}
\]
such that for a non-fixed $x \in C_i$, $\lim_{t\to 0} t \cdot x = x_{i-1}$ and $\lim_{t\to \infty} t \cdot x = x_i$;\\
\item The graded point $\gamma \in |\Grad(X/T)|$ is a positive multiple of the graded point $\{p\}/\Gm \to X/T$ associated to $\lambda$;\\
\item The arrow $\phi$ is induces a representable, finite surjection over $(B\Gm)_k$
\[
C / \Gm \to (B\Gm)_k \times_{X/T} \cW.
\]
In particular, every graded point of $\cW$ lying over $\gamma$ lifts to a graded point of $C/\Gm$ after replacing it with some positive multiple. \\
\item If $D \hookrightarrow X$ is a $T$-equivariant divisor that is regular at $p$ and such that $T_p D \subset T_p X$ is the eigenspace of $T_p X$ of weight $a$ (respectively $b$), then the strict transform $D' \hookrightarrow \cW$ of $D$ contains the image of $x_0$ (respectively $x_n$).
\end{enumerate}
\end{lem}
\begin{proof}
Consider the map $\cW \to \Sigma / T$ of \Cref{lem:relative_coarse_space}. By \Cref{lem:surface_structure}, the fiber $\Sigma_p$ is a union of (possibly singular) rational curves, each with non-trivial $\Gm$-action. We can replace each component of $\Sigma_p$ with its normalization and glue at the fixed points to obtain a finite $\Gm$-equivariant map $C = C_1 \cup \cdots \cup C_n \to \Sigma_p$ over the residue field $k(p)$. \Cref{lem:surface_structure} implies that $C$ has the properties listed in (1), after possibly extending the field so that each $C_i$ contains a $k$-rational point.

If $(\Gm)_k$ acts non-trivially on $\bP^1_k$ for some field $k$ with fixed points at $0$ and $\infty$, then one has a ramified covering map $\bP^1_{k} \to \bP^1_k$ given by $[z_0:z_1] \mapsto [z_0^m:z_1^m]$ on homogeneous coordinates. This map is equivariant with respect to the group homomorphism $(\bullet)^m : \Gm \to \Gm$. Let $\psi_m : C/(\Gm)_k \to C/(\Gm)_k$ be the morphism of stacks obtained by performing this operation on each component of $C$ at the same time.

We would like to find an arrow $C/ (\Gm)_{k} \to \cW$ lifting the map $C/(\Gm)_k \to \Sigma_p / \Gm$. This might not be possible, but after replacing $k$ with a finite extension we will construct a dotted arrow making the following diagram commute for some sufficiently divisible $m \geq 0$:
\[
\xymatrix{ & & \cW \ar[d] \\
C_k / (\Gm)_k \ar[r]^{\psi_m} \ar@{-->}@/^/[urr]^-{\phi} & C / (\Gm)_{k} \ar[r] & \Sigma/T }.
\]
The stack $(C \setminus \{x_0,\ldots,x_n\}) / \Gm$ is a disjoint union of copies of $B\mu_i$ for various $i$, so for sufficiently divisible $m$, $\psi_m$ induces a trivial map on isotropy groups on $C \setminus \{x_0,\ldots,x_n\}$. It follows that after passing to a finite extension of $k$, a lift exists over $(C \setminus \{x_0,\ldots,x_n\} ) / \Gm$ for sufficiently divisible $m$. The morphism $\cW \to \Sigma/T$ is proper, so for each component $\bP^1_k / (\Gm)_k \subset C_k / (\Gm)_k$, a lift $\bP^1_k / (\Gm)_k \to \cW$ over $\Sigma/T$ extending the given lift $(\bP^1_k \setminus \{0,\infty\}) / (\Gm)_k \to \cW$ is unique if it exists, and it must exist after passing to a suitable ramified extension of the local ring at $0$ and $\infty$. These ramified extensions are realized by a further extension of ground field and precomposition with the morphism $\psi_m$, so the dotted arrow above exists for $m$ sufficiently divisible.

From the construction, it is clear that $\phi$ induces a proper surjective map $C \to \cW \times_{X/T} \Spec(k)$ that does not contract any components and is thus finite. The claim (4) follows from the corresponding claim for $\Sigma/T \to X/T$ in \Cref{lem:surface_structure}
\end{proof}

The (strict) monotonicity condition (3) of \Cref{defn:theta_monotone} is equivalent to saying that for the chain of rational curves $C$ of \Cref{lem:stacky_surface_structure}, with fixed points $\{x_0,\ldots,x_n\}$, the value of the numerical invariant on the corresponding graded points $\mu(\{x_i\}/(\Gm)_k \to \X)$ is (strictly) monotone increasing in $i$.

\subsubsection{Proof of \Cref{T:monotone_stratifications}}

We show that any $\Theta$-monotone numerical invariant satisfies the HN-specialization condition \ref{princ:S}, at which point \Cref{T:monotone_stratifications} follows immediately from \Cref{thm:main_improved}. Let $R$ be a DVR over $B$ with fraction field $K$ and residue field $\kappa$, let $\xi : \Spec(R) \to \X$ be a morphism, and let $f_K : \Theta_K \to \X$ be a filtration over $B$ along with an isomorphism $f(1) \simeq \xi_K$.

We can glue this data to define a morphism $\xi \cup_{\xi_K} f_K : \Theta_R \setminus \{(0,0)\} \to \X$ over $B$. Let $\cW \to \Theta_R$ be a proper morphism with finite relative inertia that is an isomorphism over $\Theta_R \setminus \{(0,0)\}$, and let $\tilde{f} : \cW \to \X$ be the extension of $\xi \cup_{\xi_K} f_K$ to $\cW$ provided by \Cref{defn:theta_monotone}. Let us apply \Cref{lem:stacky_surface_structure} to the morphism $\cW \to \Theta_R$ to obtain a map $\phi : C/\Gm \to \cW$ over a graded point $\gamma : (B\Gm)_k \to \Theta_R$, where $C = C_1 \cup \cdots \cup C_n$ is a connected chain of copies of $\bP^1_k$ for some field $k$ each with non-trivial $\Gm$-action, and $\gamma$ is a positive multiple of the canonical graded point at the origin $\{(0,0)\}/\Gm \hookrightarrow \Theta_R$.

By part (1) of \Cref{lem:stacky_surface_structure}, the fixed locus $C^{\Gm} = \{x_0,\ldots,x_n\}$ is ordered such that for any non-fixed $x \in C_i$, $\lim_{t \to 0} t \cdot x = x_{i-1}$ and $\lim_{t \to \infty} t \cdot x = x_i$. Part (3) of \Cref{defn:theta_monotone} implies that the value of the numerical invariant
\[
\mu(\tilde{f} \circ \phi : \{x_i\}/(\Gm)_k \to \X) \in \Gamma
\]
is strictly monotone increasing in $i$.

The fact that $\cW \to \Theta_R$ is proper, implies that after passing to a field extension and composing with a ramified covering $(\bullet)^m : \Theta_{\kappa'} \to \Theta_{\kappa'}$, the resulting filtration $\Theta_{\kappa'} \to \Theta_R$ can be lifted to a filtration $\Theta_{\kappa'} \to \cW$. Parts (3) and (4) of the lemma implies that the associated graded point $(B\Gm)_{\kappa'} \to \cW$ agrees with the graded point $\{x_n\}/(\Gm)_k \to \cW$ after passing to a positive multiple of both. Thus if we denote the composition $f' : \Theta_{\kappa'} \to \cW \to \X$, which is a filtration of $f'(1) =\xi_\kappa \in |\X|$, then
\[
\mu \left(\agr(f') : (B\Gm)_{\kappa'} \to \X \right) = \mu \left(\tilde{f} \circ \phi : \{x_n\}/(\Gm)_k \to \X \right).
\]
Likewise, after passing to a finite extension of DVR's $R \subset R'$, one can lift the divisor $\Spec(R) \times (\{0\}/\Gm) \to \Theta_R$ to a map $\Spec(R') \times (\{0\}/\Gm) \to \cW$, and the graded point at the closed point of $\Spec(R')$, $\Spec(\kappa') \times (\{0\}/\Gm) \to \cW$ agrees with the graded point $\{x_0\}/(\Gm)_k \to \cW$ after replacing both with a positive multiple. Because the numerical invariant is stable under field extension and raising a graded point to a positive power, and it is locally constant for algebraic families of graded points, one has
\[
\mu \left(\agr(f_K) : (B\Gm)_K \to \X\right) = \mu\left(\tilde{f} \circ \phi : \{x_0\} / (\Gm)_k \to \X \right).
\]
Because $\mu(\{x_i\}/\Gm \to \X)$ is \emph{strictly} increasing in $i$, it follows that $\mu(f_K) \leq \mu(f') \leq M^{\mu}(\xi_\kappa)$, and if equality holds then $x_0 = x_n$. To complete the proof, we will show that in the latter case, if $\cW' \to \cW$ is the normalization, then $\cW' \to \Theta_R$ is an isomorphism, so the restriction of $\tilde{f}$ to $\cW'$ exhibits an extension of $\xi \cup_{\xi_K} f_K$ to all of $\Theta_R$.

Consider the base change $W := \bA^1_R \times_{\Theta_R} \cW' \to \bA^1_R$. It suffices to show that this map is an isomorphism. $W$ is a stack with finite inertia, so the map $W \to \bA^1_R$ factors through a coarse moduli space $W \to \Sigma \to \bA^1_R$, but the map $\Sigma \to \bA^1_R$ must be an isomorphism because it is finite and birational and $\bA^1_R$ is normal. We know that after base change to an \'etale neighborhood $U \to \bA^1_R$ of $(0,0)$, the stack $W$ is a quotient of an affine scheme by a finite flat group scheme, and in particular there is a locally free sheaf $E$ on $W_U$ on which the isotropy group of the zero fiber $W_{(0,0)}$ acts faithfully.

Because $W$ is normal and dimension $2$, $E$ is uniquely determined up to isomorphism by its restriction to the open substack $W \setminus W_{(0,0)} \cong \bA^1_R \setminus \{(0,0)\}$. On the other hand there is a unique vector bundle $E'$ on $\bA^1_R$ such that $E'|_{\bA^1_R \setminus \{(0,0)\}} \cong E$. This implies that $E$ is the pullback of $E'$ along the map $W \to \bA^1_R$. Because the isotropy of the zero fiber $W_{(0,0)}$ acts faithfully on $E$, this implies that it is in fact trivial, i.e., $W$ is an algebraic space. Thus the coarse moduli space morphism $W \to \bA^1_R$ is an isomorphism.

\qed

\subsection{\texorpdfstring{$\Theta$}{Theta}-monotonicity and proper morphisms}

Consider a morphism $\pi : \X \to \Y$ of stacks that is relatively representable by proper Deligne-Mumford stacks. Consider a commutative diagram of the form
\begin{equation} \label{eqn:positive}
\xymatrix{ \bP_k^1 / \Gm \ar[d] \ar[r]^\phi & \X \ar[d] \\ (B\Gm)_k \ar[r] & \Y },
\end{equation}
where $\Gm$ acts non-trivially on $\bP_k^1$ and the induced map $\bP_k^1/\Gm \to \X \times_{\Y} (B\Gm)_k$ is finite. We define the point $0 \in \bP_k^1$ to be the limit $\lim_{t\to 0} t\cdot x$ for a general point $x \in \bP_k^1$, and we let $\infty \in \bP_k^1$ be the second fixed point.

We regard $\bQ[\epsilon]$ as an ordered vector space by regarding $\epsilon$ as a formal infinitesimal positive number, i.e., $a_0 + a_1 \epsilon + \cdots \geq 0$ if $a_n \geq 0$ for the lowest $n$ for which $a_n \neq 0$. For any $\ell \in H^2(\X;\bQ[\epsilon])$ we can restrict $\phi^\ast(\ell)$ to $H^2(\{0\}/\Gm;\bQ)[\epsilon] \simeq \bQ[\epsilon] \cdot u$ and $H^2(\{\infty\}/\Gm;\bQ)[\epsilon] \simeq \bQ[\epsilon] \cdot u$.

\begin{defn} \label{defn:positive_bundle}
We say that $\ell \in H^2(\X;\bQ)[\epsilon]$ is \emph{numerically positive} (respectively nef) if $\phi^\ast(\ell)_\infty - \phi^\ast(\ell)|_0 > 0$ (respectively $\geq 0$) for any diagram of the form \eqref{eqn:positive}.
\end{defn}

\begin{ex}
In singular cohomology or operational Chow cohomology, any class in $H^2$ of a curve has a degree, so we can define a class $\ell \in H^2(\X;A)$ to be numerically positive if its restriction to any curve in any fiber of $\pi$ has positive degree. A numerically positive (respectively nef) class in this more common sense is also numerically positive (respectively nef) in the sense above, because the difference $\phi^\ast(\ell)_\infty - \phi^\ast(\ell)|_0$ gives the degree of the invertible sheaf $\cL$ on $\bP^1 / \Gm$ for which $\phi^\ast(\ell) = c_1(\cL)$. We have chosen our notion of positivity because it is the minimal condition needed for the discussion what follows, and it avoids needing a general notion of ``degree'' for classes in $H^2$ of a curve.
\end{ex}

Recall that given a morphism of stacks $\pi : \X \to \Y$ and a numerical invariant $\mu$ on $\Y$, we denote the restriction of $\mu$ to $\X$ by $\mu \circ \pi_\ast$, see \Cref{defn:restricted_invariant}.

\begin{prop} \label{prop:strat_proper_over_reductive}
Let $\Y$ be a stack satisfying \ref{hyp3} over a locally noetherian base $B$, let $\X$ be another algebraic stack, and let $\X \to \Y$ be a proper morphism that is relatively representable by Deligne-Mumford stacks. Let $b \in H^4(\Y;\bQ)$ be positive definite and $\ell \in H^2(\X;\bQ)[\epsilon]$ be numerically positive relative to $\Y$, and let $\mu : \cU \subset \iComp(\Y) \to \bR[\epsilon]$ be a numerical invariant on $\Y$ that is strictly $\Theta$-monotone and has degree $<k$ in $\epsilon$. Then the numerical invariant
\begin{equation} \label{E:perturbed_invariant}
\mu' := \mu \circ \pi_\ast + \epsilon^k \frac{\hat{\ell}}{\sqrt{\widehat{\pi^\ast(b)}}}
\end{equation}
is strictly $\Theta$-monotone on $\X$.

Furthermore, if $\mu'$ defines a weak $\Theta$-stratification on $\X$, then composition with $\pi$ maps HN filtrations in $\X$ to HN filtrations in $\Y$, and $\mu$ defines a weak $\Theta$-stratification of the closed substack $\pi(\Y) \subset \Y$ such that the preimage of each stratum of $\Y$ is a union of strata in $\X$.
\end{prop}
\begin{proof}
Let $R$ be a discrete valuation ring and consider a map $f : \Theta_R \setminus \{(0,0)\} \to \X$. Then by hypothesis we have a proper birational morphism $\cW \to \Theta_R$ that is representable by Deligne-Mumford stacks, and a morphism $\cW \to \Y$ extending the composition of $f$ with $\X \to \Y$. The morphism $f$ determines a section of the map from the fiber product $\cW \times_\Y \X \to \cW \to \Theta_R$ over $\Theta_R \setminus \{(0,0)\}$, and we let $\cW' \hookrightarrow \cW \times_\Y \X$ be the closure of this section. Then $\cW' \to \Theta_R$ is a proper birational morphism that is representable by Deligne-Mumford stacks, and by construction it comes equipped with a morphism $\tilde{f} : \cW' \to \X$ extending the original map $f$ over $\Theta_R \setminus \{(0,0)\}$.

To verify strict monotonicity, consider a diagram as in the condition (3) of \Cref{defn:theta_monotone}:
\[
\xymatrix{\bP^1_k / \Gm \ar[rr] \ar[d] & & (B\Gm)_k \ar[d] \\
\cW' \ar[r] & \cW \ar[r] & \Theta_R}
\]
If the corresponding morphism $\bP^1_k / \Gm \to (B\Gm)_k \times_{\Theta_R} \cW$ is finite, then because $\mu$ is strictly $\Theta$-monotone we have $\mu \circ \pi_\ast(\{0\}/\Gm \to \X) < \mu \circ \pi_\ast(\{\infty\}/\Gm \to \X)$ by hypothesis. The highest power of $\epsilon$ appearing in $\mu$ is lower than $\epsilon^k$, so this inequality is preserved after adding the $\epsilon^k$-term in \eqref{E:perturbed_invariant}, regardless of its value at $0$ and $\infty$. 

Otherwise, if $\bP^1_k/\Gm \to (B\Gm)_k \times_{\Theta_R} \cW$ is not finite, then it factors through the map $\bP^1_k / \Gm \to (B\Gm)_k$, so it fits into a commutative diagram
\[
\xymatrix{ \bP^1_k / \Gm \ar[r] \ar[d] & \cW' \ar[r] \ar[d] & \X \ar[d] \\
(B\Gm)_k \ar[r] & \cW \ar[r] & \Y}.
\]
We have that $\mu \circ \pi_\ast(\{0\}/\Gm \to \X) = \mu \circ \pi_\ast(\{\infty\}/\Gm \to \X)$. Also, the value of $b$ is the same at both graded points $\{0\}/\Gm \to \X$ and $\{\infty\}/\Gm \to \X$, because $b$ is pulled back from $\Y$. The strict inequality in \Cref{defn:positive_bundle} thus implies strict inequality for the $\epsilon^k$-term in \eqref{E:perturbed_invariant}, and we have $\mu'(\{0\}/\Gm \to \X) < \mu'(\{\infty\}/\Gm \to \X)$.

For the last claim, we may replace $\Y$ with the image of $\pi$, and therefore assume $\pi$ is surjective. If $f$ is a filtration in $\X$, the $\mu'(f)>0$ if and only if either $\mu(f \circ \pi)>0$ or $\mu(f \circ \pi)=0$ and $\hat{\ell}(f) > 0$, because the $\epsilon^k$ term in \eqref{E:perturbed_invariant} is smaller than any positive value achieved by $\mu$. If $x \in \X(k)$, then any filtration of $\pi(x) \in \Y$ lifts uniquely to a filtration of $x$, because $\pi$ is proper. These observations lead to the following claims, for $x \in \X$:
\begin{enumerate}
\item If $x$ is semistable, then so is $\pi(x)$, because any filtration $f$ of $\pi(x)$ with $\mu(f)>0$ can be lifted to a filtration $\tilde{f}$ of $x$ with $\mu'(\tilde{f})>0$;
\item If $x$ is unstable with HN filtration $f$, then $\pi(x)$ is semistable if and only if $\mu(f\circ \pi)=0$, i.e., $\mu'(f) \in \epsilon^k \bR[\epsilon]$;
\item If $x$ is unstable with HN filtration $f$ and $\pi(x)$ is unstable, then there can be no filtration $f'$ of $\pi(x)$ with $\mu(f') > \mu(f \circ \pi)$, so $f \circ \pi$ is an HN filtration of $\pi(x)$.
\end{enumerate}
Putting these together, we can check condition \ref{princ:B2} for $\mu$ on $\Y$: For a finite type space $T$ and a map $T \to \Y$, for the purposes of maximizing $\mu$ over filtrations of finite type points of $T$, it suffices to consider filtrations whose associated graded point is the image of the associated graded point of a HN filtration in $\X$ of some point in the image of $T \times_\Y \X \to \X$, and such points of $\Y$ form a bounded family. Therefore because $\mu$ is strictly $\Theta$-monotone on $\Y$, \Cref{T:monotone_stratifications} implies that $\mu$ defines a weak $\Theta$-stratification of $\Y$. Note that we are not assuming that $\mu$ satisfies \ref{princ:R} or either of the conditions (a) or (b) of \Cref{thm:main_improved} -- these hypotheses were only used in \Cref{thm:main_improved} to deduce the existence of HN filtrations of unstable points, but in this case we have demonstrated this directly.

The statement that the preimage of each stratum in $\Y$ is a union of strata in $\X$, and that $\pi$ preserves HN filtrations, follows from the explicit description of HN filtrations above. 
\end{proof}

The most common application of \Cref{prop:strat_proper_over_reductive} is in the case where $\Y$ is $\Theta$-reductive, so that any numerical invariant is strictly $\Theta$-monotone (\Cref{E:reductive_implies_monotone}), and in particular one can take $\mu=0$.

\begin{ex} [Projective-over-affine geometric invariant theory]
When $G$ is a reductive group, $X$ is a scheme that is projective over its affinization $Y:= \Spec(\Gamma(X,\cO_X))$, and $\cO_X(1)$ is a $G$-equivariant ample invertible sheaf, then we consider $\ell = c_1(\cO_X(1)) \in H^2(X/G;\bQ)$. A choice of Weyl group invariant norm on the cocharacter lattice of the maximal torus of $G$ defines a positive definite class in $H^4(BG;\bQ)$ whose pullback we denote by $b \in H^4(X/G;\bQ)$. Then the weak $\Theta$-stratification induced by the numerical invariant associated to $\ell$ and $b$ is the Hesselink-Kempf-Kirwan-Ness stratification of geometric invariant theory \cite{Ki84}.
\end{ex}

\begin{ex} [Projective-over-$\Theta$-reductive] \label{E:projective_over_reductive}
As a slight generalization of the previous example, consider a $\Theta$-reductive stack $\Y$ and a morphism $\X \to \Y$ that is relatively representable by Deligne-Mumford stacks and proper. We say that a line bundle $L$ on $\X$ is ample relative to $\Y$ if for any morphism $\Spec(A) \to \Y$, some power of the restriction of $L$ to $\X_A$ descends to an ample bundle on the coarse moduli space of $\X_A$ (see \cite{kresch}). Then \Cref{prop:strat_proper_over_reductive} implies that the numerical invariant associated to $\ell = c_1(L) \in H^2(\X;\bQ)$ and any positive definite $b \in H^4(\X;\bQ)$ defines a $\Theta$-stratification of $\X$.
\end{ex}

\begin{rem}
By definition, the action of a group scheme $G$ on a DM stack $X$ is a morphism of algebraic stacks $\X \to BG$ along with an isomorphism $\pt \times_{BG} \X \cong X$. Thus, \Cref{E:projective_over_reductive} extends geometric invariant theory to actions of reductive groups on projective-over-affine DM stacks.
\end{rem}

\begin{ex}[Generalization of Kempf's theorem]
In the language of this paper, Kempf's theorem in \cite{Ke78} states that for a stack of the form $\X = \Spec(A)/G$ with $G$ reductive and a closed substack $\Z \subset \X$, for any point $x \in \X(k)$ whose closure meets $\Z$, there is a canonical filtration $f$ of $x$ such that $\ev_0(f)=f(0) \in \Z$. We can generalize this to any $\Theta$-reductive stack that admits a positive definite class $b \in H^4(\X;\bQ)$ as follows:
\begin{itemize}
\item For any closed substack $\Z \subset \X$ and any filtration $f : \Theta_k \to \X$ with $f(1) \notin \Z$ but $f(0) \in \Z$, there is a canonical filtration of $f(1)$ with $f(0) \in \Z$.
\end{itemize}
This follows from \Cref{E:projective_over_reductive} applied to the projection $\X' = \op{Bl}_\Z(\X) \to \X$ and $\ell = c_1(\cO_{\X'}(1))$. One can check that for any filtration $f$ of $p \notin \Z$, $\mu(f) = 0$ if $f(0) \notin \Z$ and $\mu(f)>0$ otherwise, so we can take the canonical filtration of any $p \in \X \setminus \Z$ that admits a filtration whose associated graded lies in $\X$ to be the HN filtration of the unique lift of $p$ to $\X'$. To deduce Kempf's theorem from this more general statement, one must use the Hilbert-Mumford criterion, which says that for $x \in \Spec(A)/G$, if the closure of $x$ meets a closed substack $\Z \subset \Spec(A)/G$ then there is some filtration of $x$ with $f(0) \in \Z$.

\end{ex}

\begin{ex} [Bia{\l}ynicki-Birula stratification]
Let $X$ be a proper algebraic space over a scheme $S$ with an action of a reductive $S$-group $G$, and let $A \subset G$ be a central split torus. Then for any cocharacter $\lambda : \Gm \to A$, the Bia{\l}ynicki-Birula stratification of $X$ with respect to $\lambda$ is $G$-equivariant, and it is a $\Theta$-stratification of $X/G$. We can construct this stratification using \Cref{prop:strat_proper_over_reductive} as follows: We choose an invariant integral inner product on the cocharacter lattice of a maximal torus of $G$, and regard it as a positive definite class $b \in H^4(\pt/G;\bQ)$. We let $\chi$ be the rational character of $G$ dual to $\lambda$ with respect to this form. Then the class $\ell = c_1(\cO_X \otimes \chi) \in H^2(X/G; \bQ)$ is relatively NEF with respect to the projection $X/G \to \pt/G$. Assume that $X/G$ admits a class $\ell' \in H^2(X/G;\bQ)$ that is positive relative to $\pt/G$.

Then the Bia{\l}ynicki-Birula stratification of $X/G$ with respect to $\lambda$ is the $\Theta$-stratification defined by the numerical invariant associated to $\ell + \epsilon \ell'$ and $b$, where $\epsilon$ is a positive infinitesimal formal parameter. The HN filtration of every point $x \in X(k)$ is given by the one parameter subgroup $\lambda : \Gm \to G$, and it does not depend on $\ell'$, but $\ell'$ is used to order the strata. The semistable locus is empty.
\end{ex}

\begin{rem}
In the previous example, the existence of a positive class is necessary in order to have a $\Theta$-stratification. Consider a non-trivial action of $\Gm$ on $\bP^1_k$. Let $X$ be the $\Gm$-scheme obtained by taking two copies of $\bP^1$ and identify $0$ in the first copy with $\infty$ in the second, and vice versa. Then let $\ell$ be the first Chern class of the equivariant invertible sheaf that is $\cO_X$ twisted by a non-trivial character of $\Gm$. Then $\ell$ is NEF relative to $\pt/\Gm$, which is enough to ensure that every point has an HN filtration, but the Bia{\l}ynicki-Birula strata of $X$ can not be ordered by closure containment and thus do not form a $\Theta$-stratification.
\end{rem}

\begin{rem}
In \Cref{prop:strat_proper_over_reductive}, the stratification induced by $\ell \in H^2(\X;\bQ)[\epsilon]$ does not agree with the stratification induced by the class $\ell(r) \in H^2(\X;\bQ)$ substituting a small rational $0 < r \ll 1$ for $\epsilon$. But, we expect that if $\X$ is noetherian, then the $\ell(r)$ stratification is set-theoretically constant for sufficiently small $r$ and refines the $\ell$-stratification set-theoretically, because this is the case in classical situations (see \cite{Te00}*{Lem.~1.2}).
\end{rem}


\subsection{The Recognition Theorem revisited}

Recall from \Cref{thm:recognition} that for a locally quasi-concave numerical invariant $\mu$ on a stack $\X$ such that $\cU^{\mu >0}$ does not contain a pair of antipodal points, $\agr(f)$ is graded-semistable (\Cref{defn:graded_semistable}) for any HN filtration $f$. \Cref{thm:recognition} also includes a converse statement that applies to strictly quasi-concave numerical invariants, e.g., when the stack $\X$ is $\Theta$-reductive (by \Cref{lem:theta_reductive_convex}). In this section we strengthen this converse so that it applies in more common situations, such as reductive GIT.

We will use the following variant of \Cref{defn:theta_monotone}.

\begin{defn}[$\Theta^2$-monotone numerical invariant] \label{defn:theta2_monotone}
Let $\X$ be a stack over $B$, and let $\mu : \cU \subset \iComp(\X) \to \Gamma$ be a numerical invariant. We say that $\mu$ is \emph{(strictly) $\Theta^2$-monotone} if for any field $k$ over $B$ and non-degenerate map $\sigma : \Theta^2_k \setminus \{(0,0)\} \to \X$ over $B$ such that both filtrations $\sigma|_{\Theta_k \times \{1\}}$ and $\sigma|_{\{1\} \times \Theta_k}$ lie in $\cU$ and $\mu(\sigma|_{\{1\} \times \Theta_k}) \geq 0$, there is a proper morphism $p:\cW \to \Theta_k^2$ that is relatively representable by DM stacks and a morphism $\tilde{f} : \cW \to \X$ extending $f$ such that if we denote
\[
\cW' := \left(\Spec(k[\![t]\!]) \times \Theta_k \right) \times_{\Theta^2_k} \cW,
\]
then the projection $p' : \cW' \to \Theta_{k[\![t]\!]}$ and the composition $\tilde{f}' : \cW' \to \cW \to \X$ satisfy conditions (1), (2), and (3) of \Cref{defn:theta_monotone}.
\end{defn}

\begin{ex}
Any numerical invariant on a $\Theta$-reductive stack $\X$ is strictly $\Theta^2$-monotone, because any morphism $\sigma : \Theta_k^2 \setminus \{(0,0)\} \to \X$ extends uniquely to a morphism $\tilde{\sigma} : \Theta_k^2 \to \X$ (see \Cref{lem:theta_reductive_convex} and its proof), so one can take $\cW = \Theta_k^2$ in the definition above.
\end{ex}

\begin{rem}
It is possible that $\mu$ being (strictly) $\Theta$-monotone always implies that $\mu$ is (strictly) $\Theta^2$-monotone. In examples they often follow from identical arguments.
\end{rem}

\begin{thm}[Recognition theorem] \label{thm:recognition_revisit}
Let $\X$ be a stack satisfying \ref{hyp3}, let $\mu : \cU \subset \iComp(\X) \to \Gamma$ be a \gls{numerical_invariant}, and let $f:\Theta_k \to \X$ be a filtration in $\cU$ with $\mu(f)>0$.
\begin{enumerate}
\item If $\mu$ is locally quasi-concave, $\cU^{\mu>0}$ does not contain a pair of antipodal points, and $f$ is an HN filtration, then $\agr(f) \in \Grad(\X)$ is graded-semistable (\Cref{defn:graded_semistable}).
\item If $\agr(f)$ is graded-semistable, then $f$ is an HN filtration of $f(1)$ if either of the following hold:
\begin{enumerate}
\item $\mu$ is $\Theta^2$-monotone; or
\item $\mu$ defines a weak $\Theta$-stratification of $\X$,
\end{enumerate}
\end{enumerate}
\end{thm}

The key observation is the following:

\begin{lem} \label{lem:theta2_monotone_graded}
Let $\X$ be a stack satisfying \ref{hyp3}, let $\mu : \cU \subset \iComp(\X) \to \Gamma$ be a numerical invariant that is $\Theta^2$-monotone, and let $f : \Theta_k \to \X$ be a filtration in $\cU$. If $g : \Theta_k \to \X$ is another filtration in $\cU$ with $g(1) \cong f(1)$ and $\mu(g)>0$, then after replacing the field $k$ with a finite extension, there is a filtration $g' : \Theta_k \to \Grad(\X)$ with $g'(1) = \agr(f) \in |\Grad(\X)|$ and such that $\mu \circ u_\ast (g') \geq \mu(g)$.
\end{lem}
\begin{proof}

There is a unique morphism $f \cup g : \Theta_k^2 \setminus \{(0,0)\} \to \X$ whose restriction to $\{1\} \times \Theta_k$ is $g$ and whose restriction to $\Theta_k \times \{1\}$ is $f$. The condition of $\Theta^2$-monotonicity provides a proper birational morphism $\cW \to \Theta^2_k$ such that $g \cup f$ extends to a morphism $\tilde{g} : \cW \to \X$. If one could lift the map $\{0\}/\Gm \times \Theta_k \to \Theta_k^2$ to a morphism  $\{0\}/\Gm \times \Theta_k \to \cW$, then one could regard the composition $\{0\}/\Gm \times \Theta_k \to \cW \to \X$ as a map $g' : \Theta_k \to \Grad(\X)$ with $g'(1) = \agr(f) \in |\Grad(\X)|$. Combining \Cref{defn:theta2_monotone} and \Cref{lem:stacky_surface_structure}(4) would then show that $\mu(u \circ g') \geq \mu(g)$.

So to complete the proof, it suffices to show that
\[
\cZ := \cW \times_{\Theta_k^2} (\{0\}/\Gm \times \Theta_k) \to B\Gm \times \Theta_k
\]
admits a section after replacing $k$ with a finite field extension and composing with a ramified cover $(\bullet)^m : \Theta_k \to \Theta_k$. It suffices by Artin approximation to show this after base change to $\bar{k}$, so we assume that $k$ is algebraically closed.

We first construct a section under the assumption that $\cZ \cong Z/\Gm^2$ for some proper $\Gm^2$-equivariant morphism of algebraic spaces $p : Z \to \bA^1_k$, where the first factor of $\Gm$ acts trivially on $\bA^1_k$. To start, we construct a $k$-point in $Z_1:=p^{-1}(1)$ that is fixed by the first factor of $\Gm$. For any $k$-point $z \in Z_1(k)$, there is a unique $\Gm$-equivariant morphism $\bA^1_k \setminus 0 \to Z_1$ taking $1 \mapsto z$. This must extend uniquely to a $\Gm$-equivariant morphism $\bA^1_k \to Z_1$ because $Z_1$ is proper, and the image of $0 \in \bA^1_k$ under this map is a fixed $k$-point -- call it $z'$. Next, using the second factor of $\Gm$, the point $z'$ extends uniquely to a $\Gm^2$-equivariant section of the morphism $Z \to \bA^1_k$ over $\bA^1_k \setminus 0$, and the properness of $Z \to \bA^1_k$ guarantees that this section extends uniquely and $\Gm^2$-equivariantly over $\bA^1_k$.

Returning to the general case, we have $\cZ \cong \cY/\Gm^2$ where $\cY$ is a DM stack that is proper over $\bA^1_k$. If $Y$ denotes the coarse moduli space of $\cY$, then we have already shown that $Y / \Gm^2 \to B\Gm \times \Theta_k$ admits a section, and by pulling back along this section we can reduce the claim to the case where $\cY \to \bA^1_k$ is a coarse moduli space morphism. This is proven in \cite{GHLFH}*{Lem.~B.4}, but we summarize the argument here for convenience: After composing with a ramified cover $(\bullet)^m : \Theta_k \to \Theta_k$, one can apply the version of Sumihiro's theorem for DM stacks in \cite{alper2015luna}*{Thm.~4.1} to construct an $\Gm^2$-equivariant \'etale cover $\Spec(A) \to Z$. One can show that some connected component of $\Spec(A)$ is finite over $\bA^1_k$, and then replace $\Spec(A)$ with this component. The morphism $\Spec(A)/\Gm^2 \to B\Gm \times \Theta_k$ then admits a section by the previous paragraph, and hence so does $\cZ \to B\Gm \times \Theta_k$.

\end{proof}

\begin{proof}
\Cref{thm:recognition} already establishes (1) and (2b), so we will show (2a). Assume $\mu$ is $\Theta^2$-monotone, and suppose that there were a filtration $g : \Theta_k \to \X$ with $f(1) \cong g(1)$ and $\mu(g)>\mu(f)$. Then \Cref{lem:theta2_monotone_graded} implies that there is a filtration $g'$ of $\agr(f) \in |\Grad(\X)|$ such that $\mu \circ u_\ast(g') \geq \mu(g) > \mu(f)$, the last of which is the value of $\mu \circ u_\ast$ on the canonical filtration of $\agr(f)$. It follows that $\agr(f)$ is not graded-semistable.

\end{proof}


\subsection{Good moduli spaces and \texorpdfstring{$\Theta$}{Theta}-stratifications}

In this section we reformulate the main theorem of geometric invariant theory from the perspective of $\Theta$-stability. Our main results are \Cref{thm:monotone_moduli_spaces} and \Cref{thm:beyond_git}. Recall from \cite{alper2013good} the following
\begin{defn} \label{defn:gms}
A \emph{good moduli space} for an algebraic stack $\X$ is a quasi-compact and quasi-separated map $q : \X \to Y$, where $Y$ is an algebraic space, such that $q_\ast : \QCoh(\X) \to \QCoh(Y)$ is exact and the canonical map is an equivalence $\cO_Y \simeq q_\ast \cO_{\X}$
\end{defn}
The basic example of a good moduli space morphism is the GIT quotient map
\[
\Spec(R) / G \to \Spec(R) \gitmod G := \Spec(R^G),
\]
where $G$ is a linearly reductive group over a field $k$ and $R$ is a $k$-algebra with a $G$-action. The main results of \cite{alper2013good}*{Sec.~1.2} show that many of the useful properties of GIT quotients are consequences of the simple \Cref{defn:gms}. For instance, if $\X$ is finite type over an excellent noetherian scheme $S$, then $Y$ is finite type over $S$ and universal for maps to algebraic spaces (the noetherian hypotheses were removed from these properties in \cite{ahr2}*{Thm.~13.1\&17.2}). In fact, good moduli spaces are \'etale locally modeled by GIT quotients:

\begin{thm} \cite{alper2015luna}*{Thm.~2.9} \cite{ahr2}*{Thm.~13.1(1)}\label{thm:luna} 
Let $\X$ be an algebraic stack with affine stabilizers, separated inertia, and of finite presentation over a quasi-compact quasi-separated algebraic space $S$. If $\X$ admits a good moduli space $\X \to X$, then there exists an affine scheme $\Spec(A)$ with $\GL_n$-action and a cartesian diagram
\[
\xymatrix{
\Spec(A)/\GL_n \ar[d] \ar[r] & \X \ar[d] \\
\Spec(A^{\GL_{n}}) \ar[r] & X
}
\]
such that $\Spec(A^{\GL_n}) = \Spec(A) \gitmod \GL_n \to X$ is a Nisnevich (and hence \'etale) cover.\footnote{\cite{alper2015luna} works over an algebraically closed field, and in this case one can replace $\GL_n$ with a product of groups arising as the (linearly reductive) stabilizers of closed points of $\X$.}

\end{thm}

This has many implications. For instance, it was previously known that every geometric fiber of the map $\X \to Y$ contains a unique closed point, but we can be more precise:

\begin{cor}[S-equivalence] \label{cor:S_equivalence}
Under the hypotheses of \Cref{thm:luna}, the topological space of the good moduli space $X$ is the quotient of the topological space of $\X$ by the smallest equivalence relation that identifies $f(1)$ and $f(0)$ for any filtration $f: \Theta_k \to \X$. We call this $S$-equivalence.
\end{cor}
\begin{proof}
Using \Cref{thm:luna}, one can immediately reduce to the case of $\Spec(A)/G \to \Spec(A^G)$. The statement is then a consequence of the Hilbert-Mumford criterion, which says that any point in the fiber of this map admits a filtration whose associated graded object is the unique closed point of the fiber.
\end{proof}

Using the main theorem of \cite{ahr2}, the paper \cite{AHLH} identifies necessary and sufficient conditions for a stack to have a good moduli space, which we now recall (sticking to the characteristic $0$ case, for simplicity):

For any discrete valuation ring $R$ with maximal ideal $(\pi)$, we consider the algebraic stack
\[ST_R := \Spec(R[s,t] / (st-\pi)) / \Gm,\]
where the $\Gm$ action is determined by giving $t$ weight $-1$ and $s$ weight $1$ for $\Gm$-actions. We let $(0,0) \in \ST_R$ denote the closed point $\{s=t=0\}$, whose automorphism group is $\Gm$.

\begin{defn}[{\textsf{S}}-completeness]
A locally noetherian algebraic stack $\X$ is {\textsf{S}}-complete over a base stack $B$ if for any discrete valuation ring $R$, any diagram of solid arrows
\[
\xymatrix{\ST_R \setminus \{(0,0)\} \ar[r] \ar[d] &  \ST_R \ar@{-->}[dl] \ar[r] & \Spec(R) \ar[d] \\ \X \ar[rr] & & B }
\] 
admits a unique dotted arrow making the diagram commute.
\end{defn}

Note that $\X$ is {\textsf{S}}-complete as a $B$-stack if and only if $\X_R$ is {\textsf{S}}-complete in the absolute sense for any map $\Spec(R) \to B$.

\begin{warning} \label{warn:s_complete}
As in the case of $\Theta$-reductivity (see \Cref{warn:theta_reductive}), the notion of {\textsf{S}}-completeness of a $B$-stack differs than the relative notion of {\textsf{S}}-completeness for the morphism $\X \to B$ introduced in \cite{AHLH}*{Defn.~3.38}. The two notions agree when $B$ has quasi-finite inertia, because then any morphism $\ST_R \to B$ factors uniquely through a morphism $\Spec(R) \to B$.
\end{warning}

The condition of {\textsf{S}}-completeness is formally similar to that of $\Theta$-reductivity, and often the methods used to verify one applies to the other as well. The main existence result for good moduli spaces in characteristic $0$ states:

\begin{thm}[Special case of \cite{AHLH}*{Thm.~5.4}] \label{T:existence_gms_general}
Let $\X$ be an algebraic stack with affine stabilizers that is of finite presentation over a quasi-separated and locally noetherian algebraic space $S$. If $\X$ admits a good moduli space that is separated over $S$, then $\X$ is $\Theta$-reductive and {\textsf{S}}-complete. The converse holds if $S$ is of characteristic $0$.
\end{thm}

In practice, such as the example of the GIT quotient of a projective variety, one starts with a stack that does not admit a good moduli space, but the semistable locus with respect to some numerical invariant does admit a good moduli space. In order to formulate an intrinsic explanation for this phenomenon, we introduce another variant of \Cref{defn:theta_monotone}:

\begin{defn}[{\textsf{S}}-monotone numerical invariant]\label{D:S_monotone}
Let $\X$ be a stack over $B$, and let $\mu : \cU \subset \iComp(\X) \to \Gamma$ be a numerical invariant. We say that $\mu$ is (strictly) \emph{{\textsf{S}}-monotone} if for any discrete valuation ring $R$ over $B$ and map $f : \ST_R \setminus \{(0,0)\} \to \X$ over $B$, there is a proper birational morphism $\cW \to \ST_R$ that is relatively representable by DM stacks and a morphism $\tilde{f} : \cW \to \X$ extending $f$ such that the conditions (1), (2), and (3) of \Cref{defn:theta_monotone} hold verbatim, but with $\ST_R$ in place of $\Theta_R$.

\end{defn}

\begin{thm} \label{thm:monotone_moduli_spaces}
Let $\X$ be an algebraic stack satisfying \ref{hyp3}, and let $\mu : \cU \subset \iComp(\X) \to \Gamma$ be a locally strictly quasi-concave \gls{numerical_invariant}. Assume that for any pair $f,f'$ of non-degenerate antipodal filtrations, $f$ is in $\cU$ if and only if $f'$ is, and in this case $\mu(f) > 0$ if and only if $\mu(f') < 0$.\footnote{This condition implies that $\cU^{\mu>0}$ does not contain a pair of antipodal points, and hence $\mu$ is standard (\Cref{defn:standard_invariant}). The condition holds if $\Gamma$ is a totally ordered abelian group and $\mu(f) = -\mu(f')$, which is the case for any numerical invariant induced by a class in $H^2(\X;\bQ)$ and a norm on filtrations.} Let $\X^{\rm ss} \subset \X$ denote the substack of points that are semistable with respect to $\mu$, i.e. $M^\mu(p)\leq 0$. Consider the following conditions on $\mu$:
\begin{enumerate}[label=\roman*)]
\item $\X^{\rm ss}$ is the open piece of a weak $\Theta$-stratification for which each HN-filtration has $\mu(f)>0$.
\item $\mu$ is $\Theta^2$-monotone (\Cref{defn:theta2_monotone}).
\end{enumerate}
If either (i) or (ii) holds, then $\X^{\rm ss}$ is $\Theta$-reductive (respectively, {\textsf{S}}-complete) if $\mu$ is strictly $\Theta$-monotone (respectively, strictly {\textsf{S}}-monotone).
\end{thm}

\begin{rem}
One does not need a numerical invariant to define a notion of semistability for points in $\X$. For instance, if $\Gamma$ is a totally ordered real vector space one can define semistability with respect to a function $\ell : |\Comp(\X)_\bullet| \to V$ whose restriction to any $2$-cone is a linear function.  Any class in $H^2(\X;\Gamma)$ induces such a function (see \Cref{lem:cohomology_functions}). One defines a point to be semistable if it admits no filtration with $\ell(f)>0$. Note, however, the HN problem is ill-posed without a norm on graded points.

There is an analogous notion of $\Theta$-monotonicity and {\textsf{S}}-monotonicity in which one replaces \eqref{E:monotonicity} with the corresponding inequalities for the value of $\ell$. Our proof applies, essentially verbatim, to show if $\X$ is an algebraic stack satisfying \ref{hyp3} and $\ell : \iComp(\X) \to V$ is such a function, then under the analogous version of hypothesis (ii) above, which does not require a weak $\Theta$-stratification, $\X^{\rm ss}$ is $\Theta$-reductive or {\textsf{S}}-complete if the corresponding analog of monotonicity holds. See \cite{AHLH}*{Prop.~6.14} for a similar result in this context, but under the stronger hypothesis that $\X$ is $\Theta$-reductive and {\textsf{S}}-complete.
\end{rem}

\Cref{thm:monotone_moduli_spaces} is the last part of our generalization of the main theorem of geometric invariant theory. For convenience, we summarize it here in a form in which it is frequently used. We will say that a $\Theta$ stratification of a stack $\X$ indexed by a totally ordered set $\Gamma$ is \emph{well-ordered} if for every $x \in |\X|$, the subset $\{c \in \Gamma | \overline{\{x\}} \cap \S_c \neq \emptyset \}$ is well-ordered.

\begin{thm}[Intrinsic GIT] \label{thm:beyond_git}
Let $\X$ be an algebraic stack satisfying \ref{hyp3} over a locally noetherian base stack $B$. Let $\mu : \cU \subset \iComp(\X) \to \Gamma$ be a numerical invariant satisfying the conditions of \Cref{thm:main_improved}, and assume in addition that for any pair $f,f'$ of non-degenerate antipodal filtrations in $\cU$, $\mu(f) > 0$ if and only if $\mu(f') < 0$.
\begin{enumerate}
\item If $\mu$ is strictly $\Theta$-monotone and satisfies condition \ref{princ:B2}, then $\mu$ defines a weak $\Theta$-stratification of $\X$, which is a $\Theta$-stratification if $B$ has characteristic $0$.\\
\item In this case, the centers of the unstable strata are (a complete set of $\bN^\times$-orbit representatives for) the open substacks $\cZ_{c}^{\rm ss} \subset \Grad(\X)$ of graded-semistable points (\Cref{defn:graded_semistable}).\\
\item If furthermore $B$ has characteristic $0$, $\mu$ is strictly {\textsf{S}}-monotone, and the connected components of $\X^{\rm ss}$ are quasi-compact over $B$, then $\X^{\rm ss}$ admits a good moduli space that is separated and locally of finite presentation over $B$, and the connected components of this moduli space are proper over $B$ if the $\Theta$-stratification of $\X$ is well-ordered and $\X \to B$ satisfies the existence part\footnote{By this we mean the criterion for universally closed morphisms in \Cref{lem:valuative_variant}.} of the valuative criterion for properness.\\
\end{enumerate}

\end{thm}

\begin{proof}
\Cref{T:monotone_stratifications} implies that $\mu$ defines a weak $\Theta$-stratification of $\X$, and it is a $\Theta$-stratification if $B$ has characteristic $0$ (\Cref{cor:char_0_theta_strat}). The description of the centers of the unstable strata follows from \Cref{thm:recognition}. The existence of a good moduli space is smooth-local over $B$, so we may assume $B$ is affine, in which case the claim follows from \Cref{thm:monotone_moduli_spaces} and \Cref{T:existence_gms_general}. The properness of the good moduli space follows from the semistable reduction theorem \cite{AHLH}*{Cor.~6.12}, and the fact that if $\X^{\rm ss} \to B$ is universally closed, then so is its good moduli space \cite{alper2013good}*{Thm.~4.16}.
\end{proof}

\subsubsection{Completing the proof of \Cref{thm:monotone_moduli_spaces}}

\begin{lem}\label{L:semistable_weight_1}

Under the hypotheses of \Cref{thm:monotone_moduli_spaces}, if $f : \Theta_k \to \X^{\rm ss}$ is a non-degenerate filtration in the semistable locus, then $\mu(f)=0$.
\end{lem}
\begin{proof}
Let $f' = \sigma (\agr(f))$ denote the canonical split filtration of $f(0) = \ev_0(f)$, and let $f''$ be the split filtration coming from the opposite $\bG_m$ action on $f(0)$. Then $\mu(f)=\mu(f') \leq 0$ because $f(1)$ is semistable. On the other hand $\mu(f'') \leq 0$ because of the hypothesis that $f(0)$ is semistable. Neither $\mu(f')$ nor $\mu(f'')$ can be $<0$ without implying the other is $>0$, so $\mu(f) = \mu(f')=\mu(f'')=0$.
\end{proof}

\begin{lem} \label{L:semistable_weight_2}
Under the hypotheses of \Cref{thm:monotone_moduli_spaces}, if either of the conditions (i) or (ii) hold, and if $f : \Theta_k \to \X$ is a non-degenerate filtration in $\cU$ such that $f(1)$ is semistable and $\mu(f)=0$, then $f(0)$ is semistable as well.
\end{lem}
\begin{proof}
We may assume that $k$ is algebraically closed. We shall assume there exists a destabilizing filtration $g : \Theta_k \to \X$ of $f(0)$ and derive a contradiction. First we claim that either condition (i) or (ii) imply that one can choose $g$ to be equivariant with respect to the canonical $\Gm$-action on $f(0)$, i.e., it descends to morphism $g : \Theta_k \times (B\Gm)_k \to \X$:

\medskip
\noindent \textit{Proof of $\Gm$-equivariance under condition (i):}
\medskip

In this case, we can choose $g$ to be the $HN$ filtration of $f(0)$. It is $\Gm$-invariant by part (2) of \Cref{lem:HN_filtrations}.

\medskip
\noindent \textit{Proof of $\Gm$-equivariance under condition (ii):}
\medskip

Applying \Cref{lem:theta2_monotone_graded} to the split filtration $\sigma(\agr(f)) : \Theta_k \to \X$ of $f(0)$, induced by the canonical $\Gm$-action on $f(0)$ allows one to replace $g$ with the filtration $g'$ of $\agr(f) \cong \agr(\sigma(\agr(f)))$ supplied by that lemma.

\medskip

From this point forward, we regard $g$ as a filtration of the point $\agr(f) \in \Grad(\X)(k)$. Let $\mu \circ u_\ast$ denote the numerical invariant on $\Grad(\X)$ induced by the forgetful map $u : \Grad(\X) \to \X$. Note that because $\mu \circ u_\ast(g)>0$ but $\mu \circ u_\ast(\canon)=\mu(f)=0$, the hypothesis that $\mu$ has strictly opposite signs on antipodal filtrations guarantees that $g$ is not antipodal to $\canon$. The remainder of the proof is identical to the proof in \Cref{thm:recognition} that $\agr$ of an HN filtration is graded-semistable:

\Cref{thm:perturbation_light} provides a canonical rational $1$-simplex $\sigma : \simp{\sigma}^1 \to \iDeg(\agr(f))$ with $v_0(\sigma) = \canon$ and $v_1(\sigma)=g$. Because $\mu \circ u_\ast$ is locally quasi-concave and $\mu \circ u_\ast(\canon)=0$ and $\mu \circ u_\ast(g)>0$, the function $\mu \circ u_\ast$ is positive on the interior of $\simp{\sigma}^1$. But \Cref{thm:perturbation_light} also identifies sufficiently small subsimplex of $\simp{\sigma}^1$ that contains $v_0$ with a rational $1$-simplex in $\iDeg(\X,f(1))$ mapping $v_0$ to our original filtration $f$. \Cref{lem:embedded_numerical_invariant} implies that under this identification, the function $\mu \circ u_\ast$ on $\iDeg(\agr(f))$ corresponds to the function $\mu$ on $\iDeg(f(1))$, hence there are points in $\iDeg(f(1))$ where $\mu>0$. This contradicts the assumption that $f(1)$ is a semistable point.

\end{proof}

\begin{proof}[Proof of \Cref{thm:monotone_moduli_spaces}]
Say that one of (i) or (ii) hold and $\mu$ is strictly $\Theta$-monotone. Let $R$ be a discrete valuation ring with fraction field $K$ and residue field $k$, let $\xi : \Spec(R) \to \X^{\rm ss}$ be a morphism, and let $f_K : \Theta_K \to \X$ with $f_K(1) \simeq \xi_K$ be a filtration of $\xi_K$ over $B$. This data defines a morphism $f : \Theta_R \setminus 0 \to \X^{\rm ss}$ over $B$. By \Cref{defn:theta_monotone} there is a normal algebraic stack $\cW$, a proper birational morphism $\cW \to \Theta_R$ with finite relative inertia, and a morphism $\tilde{f} : \cW \to \X$ satisfying certain conditions. After replacing $k$ with a finite extension and composing with a ramified covering $(\bullet)^m : \Theta_k \to \Theta_k$, one can lift the map uniquely $\Theta_k \to \Theta_R$ to $\cW$. Composition with $\tilde{f}$ defines a filtration $f' : \Theta_k \to \X$ of $\xi_k$. In the proof of \Cref{T:monotone_stratifications} (which did not require $B$ to be locally noetherian) we showed that strict $\Theta$-monotonicity of $\mu$ implies that $\mu(f_K) \leq \mu(f') \leq M^\mu(\xi_k)$, and if equality holds then $\cW \to \Theta_R$ is an isomorphism.

By hypothesis $M^\mu(\xi_k)=0$, and $\mu(f_K)=0$ by \Cref{L:semistable_weight_1}, which therefore implies $\cW \cong \Theta_R$. Hence our morphism $\Theta_R \setminus 0 \to \X$ extends to a morphism $\tilde{f} : \Theta_R \to \X$ over $B$. We know that $\mu(f_k)=0$, so by \Cref{L:semistable_weight_2} $\tilde{f}(0) \in |\X^{\rm ss}|$ as well, which shows that $\tilde{f}$ maps $\Theta_R$ to $\X^{\rm ss}$. This shows that $\X^{\rm ss}$ is $\Theta$-reductive.

Next, we prove that if either (i) or (ii) holds and $\mu$ is strictly {\textsf{S}}-monotone, then $\X^{\rm ss}$ is {\textsf{S}}-complete. With $R$, $K$, and $k$ as above, let $f : \ST_R \setminus \{(0,0)\} \to \X$ be a morphism over $B$, and let $\cW$ be the normal algebraic stack, $\cW \to \ST_R$ the proper morphism with finite relative inertia, and $\tilde{f} : \cW \to \X$ the morphism given by \Cref{D:S_monotone}. Let $C = C_1 \cup \cdots \cup C_n$ be the chain of rational curves with $\Gm$-action provided by \Cref{lem:stacky_surface_structure}, and let $(C_1 \cup \cdots \cup C_n)^{\Gm} = \{x_0,\ldots,x_n\}$ be the ordering in part (1) of that lemma.

Recall that $\ST_R = \Spec(R[s,t]/(st-\pi))$ with $s$ having weight $1$ and $t$ having weight $-1$. The divisors $\{s=0\}$ and $\{t=0\}$ are isomorphic to $\Theta_k$, and properness of $\cW \to \ST_R$ implies that after replacing $k$ with an extension field and composing with $(\bullet)^m : \Theta_k \to \Theta_k$ one can lift both maps $\Theta_k \to \ST_R$ uniquely to $\cW$. Let $f_{out} : \Theta_k \to \X$ be the composition of $\{t=0\}/\Gm \to \cW \to \X$, and we let $f_{in}$ denote the composition $\{s=0\}/\Gm \to \cW \to \X$. It follows from parts (1) and (4) of \Cref{lem:stacky_surface_structure}, combined with part (3) of \Cref{D:S_monotone} that
\[
\mu(f_{out}) = \mu(\{x_0\}/\Gm \to \X) \leq \mu(\{x_n\}/\Gm \to \X) = \mu(f_{in}),
\]
and when equality holds $x_0=x_n$. The same argument used in the proof of \Cref{T:monotone_stratifications} now shows that $\cW \to \ST_R$ is an isomorphism. (Briefly: The base change $W := \Spec(R[s,t]/(st-\pi)) \times_{\ST_R} \cW$ is a normal stack with finite inertia such that the map $W \to \Spec(R[s,t]/(st-\pi))$ is a coarse moduli space morphism and an isomorphism away from the closed point $\{s=t=0\}$. The argument in the proof of \Cref{T:monotone_stratifications} applies verbatim to show that it must be an isomorphism.)

Because $\ST_R \setminus 0$ maps to $\X^{\rm ss}$, $\mu(f_{out})\leq 0$ and $\mu(f_{in}) \leq 0$. Note that the canonical filtrations of the associated graded points $\sigma(\agr (f_{in}))$ and $\sigma(\agr(f_{out}))$ are antipodal, so by hypothesis if either of $\mu(f_{in}) = \mu(\sigma(\agr(f_{in})))$ or $\mu(f_{out}) = \mu(\sigma (\agr(f_{out})))$ are negative, then the other must be positive, so both must be $0$. Therefore \Cref{L:semistable_weight_2} implies that $f_{out}(0)=f_{in}(0) \in \X$ is semistable. It follows that the extension $\ST_R \to \X$ lies in the substack $\X^{\rm ss}$.

\end{proof}

\subsection{Variation of good moduli space}

As a final illustration of the interaction between the notion of good moduli spaces and $\Theta$-stratifications, we now formulate an intrinsic version of the classical theory of variation of GIT quotient \cite{DH98}. This is inspired by \Cref{thm:monotone_moduli_spaces}, but we will see that it admits a direct proof based on \Cref{thm:luna}:

\begin{thm}[Variation of good moduli space] \label{thm:main_GIT}
Let $\Y$ be an algebraic stack with affine stabilizers, separated inertia, and of finite presentation over a quasi-separated algebraic space $S$, and let $q : \Y \to Y$ be a good moduli space. Let $\pi : \X \to \Y$ be a projective representable map.
\begin{enumerate}
\item For any $\ell \in H^2(\X;\bQ)$ that is positive relative to $\pi$ and any rational quadratic norm $\lVert\bullet\rVert$ on graded points of $\Y$, the numerical invariant $\mu$ associated to $\ell$ and $\lVert \bullet \rVert$ in \Cref{defn:induced_invariant} defines a weak $\Theta$-stratification of $\X$. This weak $\Theta$-stratification is compatible with base change along any map $Y' \to Y$ in the sense that the stratification defined by the restriction of $\ell$ and $\lVert\bullet\rVert$ to $\X' := \X \times_Y Y'$ agrees with the induced stratification of \Cref{lem:induced_base_change}.\\

\item If furthermore $\ell = c_1(\cL)$ for some $\cL \in \op{Pic}(\X)$ that is ample relative to $\pi$, then
\begin{enumerate}
  \item A point $p \in |\X|$ is $\mu$-semistable if and only if $\exists$ a map $\Spec(A) \to Y$,  a point $p' \in |\X_A|$, and a section $\sigma \in \Gamma(\X_A, \cL^n|_{\X_A})$ for $n>0$ such that $p'$ maps to $p$ and $\sigma(p') \neq 0$.
  \item The canonical map $\X^{\rm{ss}} \to \underline{\op{Proj}}_Y \left( \bigoplus_{n \geq 0} q_\ast(\pi_\ast(\cL^n)) \right)$ is a good moduli space, and
\end{enumerate}

\item If $S$ is noetherian, $\ell \in H^2(\X;\bQ)$ is positive relative to $\pi$, and $\ell' \in H^2(\X;\bQ)$ is arbitrary, then the semistable locus with respect to $\ell + r \ell'$ is independent of $r$ and contained in the semistable locus with respect to $\ell$ for all $0<r\ll 1$.

\end{enumerate}
\end{thm}

\begin{proof}
The claim in (1) that $\mu$ defines a $\Theta$ stratification follows from \Cref{thm:theta_reductive_stratifications}, because $\X$ is $\Theta$-reductive by \Cref{T:existence_gms_general}. The fact that the $\Theta$-stratification is compatible with base change along a map $Y' \to Y$ follows from the claim that for any $p' \in \X'(k)$ with image $p \in \X(k)$, the induced map of formal fans is an isomorphism
\begin{equation}\label{eqn:base_change_fan}
\Deg(\X',p')_\bullet \to \Deg(\X,p)_\bullet.
\end{equation}
This would imply that $p'$ is semistable if and only if $p$ is semistable, and the HN filtration of $p'$ is the HN filtration of $p$.

For \emph{any} map from a stack to an algebraic space $\pi : \X \to Y$, and any $\bZ^n$-weighted filtration $\Theta_k^n \to \X$ of $p \in \X(k)$, the composition $\Theta_k^n \to \X \to Y$ factors uniquely through the projection $\Theta^n_k \to \Spec(k)$, by \Cref{lem:DM_filtrations}. If we let $q = \pi(p) \in Y(k)$ and $\X_q := \pi^{-1}(q) \to \X$, then $p : \Spec(k) \to \X$ factors canonically through $\X_q$, and \Cref{cor:representable_base_change} implies that the both commutative squares in the following diagram are cartesian:
\[
\xymatrix{
\Flag^n(p) \ar[r] \ar[d] & \Filt^n(\X_q) \ar[r] \ar[d]^{\ev_1} & \Filt^n(\X) \ar[d]^{\ev_1} \\ \Spec(k) \ar[r]^{p} & \X_q \ar[r] & \X
}
\]
In particular $\Deg(\X,p)_\bullet \simeq \Deg(\X_q,p)_\bullet$ canonically. In the setting of our map $\X' \to \X$, if we let $q \in Y(k)$ denote the image of $p$ and $q' \in Y'(k)$ denote the image of $p'$, then the canonical map $\X_{q'} \to \X_q$ is an equivalence, hence \eqref{eqn:base_change_fan} is an equivalence.

By (1) we know that $\Theta$-stability of a point can be evaluated \'etale locally over $Y$, so it suffices to verify the claims in (2) after passing to an \'etale cover of $Y$. \Cref{thm:luna} therefore implies that we can assume $\Y \cong \Spec(A)/\GL_n$ and is cohomologically affine, and $\X \cong X/\GL_n$ for some $\GL_n$-equivariant projective morphism $X \to \Spec(A)$. In this setting $\Theta$-stability agrees with the usual Hilbert-Mumford numerical criterion for semistability, and thus it agrees with the notion of semistability in GIT as the complement of the stable base locus of $\cL$ (see \cite{Seshadri} for a discussion). So the claims of (2) follow from the main theorem on the construction of GIT quotients, formulated in this level of generality in \cite{alper2013good}*{Thm.~11.5}.

The claim (3) also can be checked \'etale locally over $Y$, and thus reduces to the classical statement for the GIT quotient of $X/\GL_n$, where $X$ is projective over an affine $\GL_n$-scheme.
\end{proof}

Variation of GIT quotient corresponds to the case where $\pi : \X \to \Y$ is the identity $\Y \to \Y$. Using this same method of passing to an \'etale cover of the good moduli space of $\X$, one may transport many of the other structures from this theory into our intrinsic setting. For instance, one may study the decomposition of $NS(\X)_\bR$ into cells corresponding to equivalences class, where $\ell \sim \ell'$ if $\X^{\ell-\rm{ss}} = \X^{\ell'-\rm{ss}}$ \cite{DH98}*{Defn.~3.4.5}.


\section{Moduli of complexes of coherent sheaves} \label{sect:moduli_derived}

The stack of coherent sheaves on a projective scheme admits a canonical constructible stratification \cite{Sh77}, where each stratum consists of sheaves whose Harder-Narasimhan filtration has associated graded pieces with given Hilbert polynomials. More recently it was observed in \cite{nitsure2011schematic}, without using the language of this paper, that this stratification is a $\Theta$-stratification. However, the stack of coherent sheaves can be approximated by quotients of larger and larger Quot schemes by a general linear group, and it is shown in \cite{hoskins2012quotients} that in a suitable sense the Harder-Narasimhan stratification of the stack of coherent sheaves is a colimit of the usual GIT stratifications on these Quot schemes. So the methods developed in this paper are not strictly necessary in this case.

In this section we use our methods to construct a $\Theta$-stratification on a stack that is not known to be a local quotient stack, the stack of torsion-free objects in the heart of a Bridgeland stability condition (\Cref{thm:theta_stratification_torsion_free}). We show that slope stability and Harder-Narasimhan theory in abelian categories, concepts that are familiar to experts, can be reformulated as $\Theta$-stability on the stack of objects in the heart of a $t$-structure (\Cref{thm:HN_Bridgeland}). In addition to constructing new examples of $\Theta$-stratifications in interesting moduli problems, we hope it is instructive to connect our theory to this well-studied example.

Throughout this section, $k$ will denote a fixed ground field and $X$ a projective scheme over $k$. We denote the bounded derived category of coherent sheaves on $X$ by $\D^b(X)$, and the unbounded derived category of quasi-coherent sheaves $\D_{qc}(X)$.

\subsection{Recollections on constant families of \texorpdfstring{$t$}{t}-structures}
\label{sect:moduli_of_objects}

If $\cA \subset \D^b(X)$ is the heart of a $t$-structure, Polishchuk and Abramovich define in \cites{abramovich2006sheaves,polishchuk2007constant} a moduli functor for each $v \in \Knum(X)$ that assigns to any finite type $k$-scheme $S$ the groupoid of ``flat families in $\cA$,''
\begin{equation} \label{eqn:def_families}
\left\{ F \in \D^b(X \times S) \left|  \begin{array}{c} \forall \text{ closed points } s\in S, \\ (X_s \to X)_\ast (F|_{X_s}) \in \cA \end{array} \right. \right\}.
\end{equation}
We will discuss a moduli functor that is defined on all $k$-algebras and is equivalent to \eqref{eqn:def_families} for finite type $k$-algebras. In order to do this, we recall, in a slightly more general context, the constructions and results of \cite{polishchuk2007constant}.

For $M \in \D_{qc}(R)$ and $E \in \D_{qc}(X)$, we let $M \boxtimes E \in \D_{qc}(\Spec(R) \times_k X)$ denote the exterior tensor product.

\begin{defn} \label{D:induced_t_structure}
Given a $t$-structure on $\D^b(X)$ and a $k$-algebra $R$, the \emph{induced} $t$-structure on $\D_{qc}(X_R)$ is the unique $t$-structure for which
\[
\D_{qc}(X_R)^{\leq 0} := \left( \begin{array}{c}\text{smallest full subcategory of }\D_{qc}(X_R) \\ \text{containing }R \boxtimes E, \forall E \in \D^b(X)^{\leq 0} \text{ and that is closed} \\ \text{under small colimits and extensions} \end{array} \right),
\]
and $\D_{qc}(X_R)^{\geq 0}$ is the category of $E \in \D_{qc}(X_R)$ such that $\Hom(F,E) = 0, \forall F \in \D_{qc}(X_R)^{\leq 0}$. The truncation functors $\tau^{\leq n}$ and $\tau^{\geq n}$ commute with filtered colimits.
\end{defn}
\begin{proof}[Proof that this defines a $t$-structure]
The category $\D_{qc}(X_R)$ is a presentable stable $\infty$-category and the category $\D^b(X)^{\leq 0}$ is essentially small, so by \cite{lurie2012higher}*{Prop.~1.4.4.11} there is an accessible $t$-structure on $\D_{qc}(X_R)$ whose subcategory of connective objects is $\D_{qc}(X_R)^{\leq 0}$. The truncation functors preserve filtered colimits because the $t$-structure is accessible \cite{lurie2012higher}*{Prop.~1.4.4.13}.
\end{proof}

Note that $R$ generates $\D_{qc}(R\Mod)^{\leq 0}$ under colimits and $\D^b(X)^{\leq 0}$ generates $\D_{qc}(X)^{\leq 0}$ under extensions and colimits, so $M \otimes E \in \D_{qc}(X_R)^{\leq 0}$ for any $M \in \D_{qc}(R\Mod)^{\leq 0}$ and $E \in \D_{qc}(X)^{\leq 0}$. \Cref{D:induced_t_structure} is simply an $\infty$-categorical version of the construction in \cite{polishchuk2007constant}, and when $R$ is finitely generated the induced $t$-structure on $X_R$ constructed above agrees with the previous constructions of induced $t$-structures by \cite{abramovich2006sheaves}*{Thm.~2.7.2} and \cite{polishchuk2007constant}*{Thm.~3.3.6}.

\begin{lem} \label{lem:t_structure_properties}
For any ring map $R \to S$, the induced map $\phi : X_S \to X_R$ has the following properties with respect to the $t$-structure we have constructed:
\begin{enumerate}
\item $\phi^\ast : \D_{qc}(X_R) \to \D_{qc}(X_S)$ is right $t$-exact,
\item $\phi_\ast : \D_{qc}(X_S) \to \D_{qc}(X_R)$ is $t$-exact,
\item any $E \in \D_{qc}(X_S)$ lies in $\D_{qc}(X_S)^\heartsuit$ (respectively $\D_{qc}(X_S)^{\leq 0}$ or $\D_{qc}(X_S)^{\geq 0}$) if and only if $\phi_\ast(E)$ does,
\item if $R \to S$ is flat then $\phi^\ast$ is $t$-exact, and
\item if $\{R \to S_\alpha\}_{\alpha \in I}$ is a flat cover of $\Spec(R)$ then $E \in \D_{qc}(X_R)$ lies in the heart if and only if $\phi_\alpha^\ast(E) \in \D_{qc}(X_{S_\alpha})^\heartsuit$ for all $\alpha \in I$.
\end{enumerate}
\end{lem}
\begin{proof}
The first claim is immediate from the definitions, because $\phi^\ast(R \boxtimes E) = S \boxtimes E$ and $\phi^\ast$ commutes with colimits because it has a right adjoint $\phi_\ast$. Furthermore, the adjunction between $\phi^\ast$ and $\phi_\ast$ implies that $\phi_\ast$ must map the right orthogonal complement $\D_{qc}(X_S)^{>0}$ to $\D_{qc}(X_R)^{>0}$. Furthermore, under the equivalence of stable $\infty$-categories $\D_{qc}(X_S) \simeq \op{Mod}_{\cO_{X_R} \boxtimes S}(\D_{qc}(X_R)$, colimits in the former are colimits in $\D_{qc}(X_R)$ along with their induced $\cO_{X_R} \boxtimes R$-module structure. It follows that $\phi_\ast$ commutes with all small colimits. Note that for $E \in \D_{qc}(X)^{\leq 0}$ we have $\phi_\ast(S \boxtimes E) \in \D$ because $S \in \D_{qc}(R\Mod)^{\leq 0}$. This combined with the fact that $\phi_\ast$ commutes with colimits implies that $\D_{qc}(X_S)^{\leq 0} \subset (\phi_\ast)^{-1}(\D_{qc}(X_R)^{\leq 0})$, which implies (2). (3) follows formally from the fact that $\phi_\ast$ is $t$-exact and conservative. (4) follows from (2) and the fact that $\phi_\ast(\phi^\ast(E)) \simeq S \otimes^L_R E$, so if $R \to S$ is flat then $S$ is a filtered colimit of free $R$-modules. (5) follows from (4) and the fact that $\prod_\alpha \phi_\alpha^\ast : \D_{qc}(X_R) \to \prod_\alpha \D_{qc}(X_{S_\alpha})$ is conservative.
\end{proof}

Note as a consequence that if $p : X_R \to X$ is the projection, then we have
\[
\D_{qc}(X_R)^{[a,b]} = \{ E \in \D_{qc}(X_R) | p_\ast(E) \in \D_{qc}(X)^{[a,b]} \},
\]
where $\D_{qc}(X_R)^{[a,b]} = \D_{qc}(X_R)^{\leq b} \cap \D_{qc}(X_R)^{\geq a}$.

\begin{cor} \label{cor:families}
For any algebraic $k$-stack $\Y$, there is a canonical $t$-structure induced on $\D_{qc}(X_{\Y}) = \D_{qc}(X \times_{\Spec(k)} \Y)$ in which $\D_{qc}(X_{\Y})^{\leq 0}$ (respectively $\D_{qc}(X_{\Y})^{\geq 0}$) is the full subcategory of complexes $E$ such that for any smooth map $\Spec(R) \to \Y$ we have $E|_{X_R} \in \D_{qc}(X_R)^{\leq 0}$ (respectively $\D_{qc}(X_R)^{\geq 0}$). It suffices to check if $E \in \D_{qc}(X_\Y)^{\leq 0}$ or $E \in \D_{qc}(X_\Y)^{\geq 0}$ after restricting to a smooth cover of $\Y$ by affine schemes.
\end{cor}
\begin{proof}
Given a smooth hypercover $Y_\bullet \to \Y$ by a simplicial scheme $Y_\bullet$ where each $Y_n$ is a disjoint union of affines, hyperdescent for quasi-coherent sheaves gives an equivalence with the totalization of the associated cosimplicial diagram of stable $\infty$-categories
\[
\D_{qc}(X_{\Y}) = \varprojlim \D_{qc}(X_{Y_\bullet}),
\]
where the limit is taken over the semisimplex category $\mathbf{\Delta}_{inj} \subset \mathbf{\Delta}$ consisting of face maps only. This is equivalent to taking the limit over $\mathbf{\Delta}$, because $\mathbf{\Delta}_{inf} \hookrightarrow \mathbf{\Delta}$ is cofinal. \Cref{lem:t_structure_properties}(4) implies that all of the functors in this diagram are $t$-exact for the induced $t$-structure on each level.

The existence of a $t$-structure on $\D_{qc}(X_\Y)$ as described above now follows from the following claim: given any diagram of stable $\infty$-categories $\{\cC_i\}_{i \in I}$ indexed by some category $I$ along with $t$-structures $(\cC_i^{>0},\cC_i^{\leq 0})$ on each $\cC_i$ such that for each arrow $i \to j$ in $I$ the corresponding functor $\cC_i \to \cC_j$ is $t$-exact, the pair of full subcategories $(\varprojlim (\cC_i^{> 0}), \varprojlim (\cC_i^{\leq 0}))$ is a $t$-structure on $\varprojlim \cC_i$. This is proved in \cite{halpernleistner2020derived}*{Lem.~2.2.7} for semiorthogonal decompositions, but the proof applies verbatim for $t$-structures. (A semiorthogonal decomposition is just a special case of a $t$-structure where $\cC^{\leq 0}$ and $\cC^{>0}$ are closed under homological shift in either direction.) The idea is that one can replace each $\cC_i$ with an equivalent $\infty$-category $\cC'_i$ consisting of diagrams $A \to B \to C$ that are exact triangles in $\cC_i$ with $A \in \cC_i^{\leq 0}$ and $C \in \cC_i^{>0}$. The $t$-exactness of the functors in the original diagram implies that one gets a diagram $\{\cC'_i\}_{i \in I}$ with an equivalence $\varprojlim_i \cC'_i \cong \varprojlim_i \cC_i$. $\varprojlim_i \cC'_i$ admits clear projection functors to $\varprojlim_i (\cC_i^{\leq 0})$ and $\varprojlim_i (\cC_i^{>0})$, which are the truncation functors for the $t$-structure on $\varprojlim_i \cC_i$.

\end{proof}

The following depends on the key technical result of \cites{abramovich2006sheaves,polishchuk2007constant}, and is a slight generalization of what is stated in those papers. Following \cite{polishchuk2007constant}, we say that a $t$-structure on $\D^b(X)$ is nondegenerate if $\bigcap_n \D^b(X)^{\leq n} = \bigcap_n \D^b(X)^{\geq n} = 0$.

\begin{prop} \label{prop:t_structure_coherent}
Assume the $t$-structure on $\D^b(X)$ is noetherian and nondegenerate, and let $R$ be an algebra that is essentially of finite type over $k$. Then the truncation functors on $\D_{qc}(X_R)$ preserve $\D^b(X_R)$, and the induced $t$-structure on $\D^b(X_R)$ is noetherian.
\end{prop}

\begin{lem} \label{lem:localization_surjectivity}
Let $R \to S$ be a localization of noetherian $k$-algebras. Then $\phi^\ast : \D^b(X_S) \to \D^b(X_R)$ is essentially surjective.
\end{lem}
\begin{proof}
We may work with the usual $t$-structure on $\D^b(-)$. The key observation is that because $R \to S$ is a localization, for any $E \in \QCoh(X_S)$ we have $\phi^\ast(\phi_\ast(E)) = \phi_\ast(E) \otimes_R S \simeq E$, where $\phi_\ast(E)$ is just $E$ regarded as a quasi-coherent sheaf on $X_R$. It follows that for any bounded complex of quasi-coherent sheaves $E^\bullet$ on $X_S$ with coherent homology sheaves, the pushforward $\phi_\ast(E^\bullet)$ is a complex of quasi-coherent sheaves on $X_R$ whose restriction to $X_S$ is $E^\bullet$.

From this point, a descending induction argument allows one to replace $\phi_\ast(E^\bullet)$ with a bounded coherent subcomplex $F^\bullet \subset \phi_\ast(E^\bullet)$ such that the induced map $F^\bullet \otimes_R S \to \phi_\ast(E^\bullet) \otimes_R S \simeq E^\bullet$ is a quasi-isomorphism. More precisely, the proof of \cite{polishchuk2007constant}*{Lem.~2.3.1} applies verbatim once one establishes the following:
\begin{itemize}
\item[Claim:] For any surjective map $F \to G$ in $\QCoh(X_R)$ such that $G\otimes_R S \in \QCoh(X_S)$ is coherent, there is a coherent subsheaf $F' \subset F$ such that $F'\otimes_R S \to G \otimes_R S$ is surjective.
\end{itemize}
This claim follows by writing $F$ as a filtered union $F = \bigcup_\alpha F_\alpha$ of coherent subsheaves $F_\alpha \subset F$, then observing that because the category of coherent sheaves on $X_R$ is noetherian one of the maps $F_\alpha \otimes_R S \to G \otimes_R S$ must be surjective.
\end{proof}

\begin{proof}[Proof of \Cref{prop:t_structure_coherent}]
The case of a smooth $k$-algebra is \cite{abramovich2006sheaves}*{Thm.~2.6.1\& 2.7.2}, and the case of a general $k$-algebra of finite type in \cite{polishchuk2007constant}*{Thm.~3.3.6}. Here we simply extend this observation to the case where $S$ is a localization $R \to S$ of a $k$-algebra $R$ of finite type. We already know from \Cref{lem:t_structure_properties} that $\phi^\ast : \D_{qc}(X_R) \to \D_{qc}(X_S)$ is $t$-exact. It therefore suffices to show that $\phi^\ast :\D^b(X_R) \to \D^b(X_S)$ is essentially surjective, which is \Cref{lem:localization_surjectivity}. By the same logic, a sequence of surjective maps $E_1 \to E_2 \to \cdots$ in $\D^b(X_S)^\heartsuit$ can be lifted to a descending chain of surjections in $\D^b(X_R)^\heartsuit$, so it must stabilize and hence the $t$-structure on $\D^b(X_S)$ is noetherian.
\end{proof}

\begin{cor}
If the $t$-structure on $\D^b(X)$ is noetherian and nondegenerate, then for any field extension $K/k$, the subcategory $\D^b(X_K) \subset \D_{qc}(X_K)$ is preserved by the truncation functors for the induced $t$-structure and thus inherits a $t$-structure.
\end{cor}
\begin{proof}
Any $E \in \D^b(X_K)$ is relatively perfect (see \Cref{S:moduli_functor} the definition), so by \cite{lieblich2005moduli}*{Prop.~2.2.1} we can find a finitely generated $k$-subalgebra $R \subset K$ and an $E'\in \D^b(X_R)$ such that $E \simeq E' \otimes^L_R K$. The extension $R \subset K$ must be flat because $K$ contains the field of fractions of $R$, so the restriction functor $\D_{qc}(X_R) \to \D_{qc}(X_K)$ is $t$-exact. It follows that all truncations of $E$ are the restriction of truncations of $E'$, hence they lie in $\D^b(X_K)$.
\end{proof}

After this paper was first released, the paper \cite{AHLH} established existence results for moduli spaces of objects in abelian categories using a slightly different formulation of the moduli functor, $\cM_{\cC}^{\rm ab}$, which is defined in terms of a locally noetherian $k$-linear abelian category $\cC$ (see \cite{AHLH}*{Sect.~7.1} for background). For any $k$-algebra $R$, one can define a locally finitely presented abelian category $R\Mod(\cC)$ of $R$-module objects in $\cC$, i.e., objects $E \in \cC$ along with a $k$-algebra homomorphism $R \to \End_k(E)$. In order to relate the results of \cite{AHLH} with the construction here, we observe:

\begin{prop} \label{prop:heart_category}
Equip $\D^b(X)$ with a t-structure that is noetherian and bounded with respect to the usual $t$-structure, and let $\cA_{qc} := \D_{qc}(X)^{\heartsuit}$. Then,
\begin{enumerate}
\item $\cA_{qc}$ is a locally noetherian $k$-linear abelian category whose subcategory of compact objects is $\D^b(X)^\heartsuit$, and 
\item there is a natural equivalence $\D_{qc}(X_R)^{\heartsuit} \cong R\Mod(\cA_{qc})$ for any $k$-algebra $R$ such that for any morphism of $k$-algebras $\phi : R \to S$, the base change functor $S \otimes_R (-) : R \Mod(\cA_{qc}) \to S\Mod(\cA_{qc})$ is identified with $H^0(\phi^\ast(-)) : \D_{qc}(X_R)^\heartsuit \to \D_{qc}(X_S)^\heartsuit$.
\end{enumerate}
\end{prop}
\begin{proof}
The fact that the $t$-truncation functors on $\D_{qc}(X)$ commute with filtered colimits implies that filtered colimits exist and are exact in $\cA_{qc}$. If we temporarily use the notation $\D_{qc}(X)^{\leq_\tau n}$ and $\D_{qc}(X)^{\geq_\tau n}$ for the usual $t$-structure on $\D_{qc}(X)$, then the hypothesis that the $t$-structure on $\D^b(X)$ is bounded with respect to the usual $t$-structure along with \Cref{D:induced_t_structure} implies that
\[
\D_{qc}(X)^{\leq -n} \subset \D_{qc}(X)^{\leq_{\tau} 0} \subset \D_{qc}(X)^{\leq n}
\]
for some $n\geq 0$. This implies inclusions of the semiorthogonal complements $\D_{qc}(X)^{> n} \subset \D_{qc}(X)^{>_{\tau} 0} \subset \D_{qc}(X)^{> -n}$. For any fixed $n$, the category $\D_{qc}(X)^{\geq_\tau n}$ is compactly generated by $\D^b(X)^{\geq_{\tau} n}$ (see for example \cite{HLPreygel}*{Thm.~A.2.1}), which implies that any $F \in \D^b(X)^{\heartsuit}$ is compact in $\cA_{qc}$. Furthermore, one can write any $E \in \cA_{qc} \subset \D_{qc}(X)^{\geq_{\tau} -n}$ as a filtered colimit $E \cong \colim E_\alpha$ with $E_\alpha \in \D^b(X)^{\geq_{\tau} -n}$. Then the fact that the truncation functors commute with filtered colimits implies that $E \cong H^0(E) \cong \colim H^0(E_\alpha)$, so $\cA_{qc}$ is compactly generated by $\D^b(X)^\heartsuit$. Because $\D^b(X)^\heartsuit$ is closed under retracts, $\D^b(X)^\heartsuit \subset \cA_{qc}$ coincides with the category of compact objects.

To show that $\cA_{qc}$ is locally noetherian, it remains to show that any $E \in \D^b(X)^\heartsuit$ is noetherian, meaning every filtered system of subobjects $\{F_\alpha \subset E\}$ stabilizes. It suffices to show that if $F = \bigcup_\alpha F_\alpha \subset E$, then $F \in \D^b(X)^\heartsuit$, because in this case $F$ would be finitely presented by the previous paragraph. By the hypothesis that $\D^b(X)^\heartsuit$ is noetherian, there exists a maximal $F' \subset F$ such that $F' \in \D^b(X)^\heartsuit$. The inclusion $\D^b(X) \subset \D_{qc}(X)$ is $t$-exact, so the image of any morphism in $\D^b(X)^\heartsuit$ agrees with the image of the same regarded as a morphism in $\cA_{qc}$. This implies that for any $G \in \D^b(X)^{\heartsuit}$, any morphism $G \to F$ factors uniquely through $F' \subset F$. Because $F$ can be written as a filtered colimit of objects of $\D^b(X)^\heartsuit$, it follows that $F = F'$, hence $F \in \D^b(X)^\heartsuit$.

To construct the canonical equivalence $\D_{qc}(X_R)^{\heartsuit} \cong R\Mod(\cA_{qc})$ for a $k$-algebra $R$, let $p : X_R \to X$ be the projection, and observe that the pushforward functor lifts to an ``enhanced" pushforward functor
\[
\xymatrix{
\D_{qc}(X_R)^{\heartsuit} \ar@{-->}[rr]^{p_\ast^{en}} \ar[dr]^{p_\ast} & &  R\Mod(\cA_{qc}) \ar[dl]^{forget} \\
& \cA_{qc} & 
}
\]
where for $E \in \D_{qc}(X_R)^\heartsuit$, $p_\ast^{en}(E)$ regards the underlying complex $p_\ast(E) \in \cA_{qc} = \D_{qc}(X)^\heartsuit$ as an $R$-module object via the $k$-algebra map
\[
R \to H^0(X_R,\cO_{X_R}) \to \End(E) \to \End(p_\ast(E)).
\]
The functor $p_\ast^{en}$ is exact, because both forgetful functors to $\cA_{qc}$ are exact and faithful. We leave it to the reader to verify that the formation of $p_\ast^{en}$ commutes with base change along a map of $k$-algebras $R \to S$.

Observe that after identifying $p_\ast^{en}(p^\ast(E)) \cong R\otimes_k E$, the homomorphism induced by $p_\ast^{en}$ on Hom groups agrees with the composition of equivalences coming from adjunction
\begin{align*}
\Hom_{\D_{qc}(X_R)}(p^\ast(E),G) &\cong \Hom_{\cA_{qc}}(E,p_\ast(G)) \\
&\cong \Hom_{R\Mod(\cA_{qc})}(R \otimes_k E, p_\ast^{en}(G)).
\end{align*}
We also observe for any $F \in \D_{qc}(X_R)^\heartsuit$, the canonical morphism $p^\ast p_\ast(F) \to F$ is surjective, because after applying $p_\ast$ it admits a section. It follows that for any $E \in \cA_{qc}$ one can write $E \cong \coker(p^\ast(F_1) \to p^\ast(F_0))$ with $F_0,F_1 \in \cA_{qc}$. Then for any other $G \in \D_{qc}(X_R)^\heartsuit$, the exactness of $p_\ast^{en}$ shows that
\[
\begin{array}{l}
\Hom_{\D_{qc}(X_R)^\heartsuit}(E,G) \cong \\
\quad \cong \ker(\Hom_{\D_{qc}(X_R)^\heartsuit}(p^\ast(F_0),G) \to \Hom_{\D_{qc}(X_R)^\heartsuit}(p^\ast(F_1),G)) \\
\quad \cong \ker(\Hom_{R\Mod(\cA_{qc})}(R\otimes_k F_0,p_\ast^{en} G) \to \Hom_{R\Mod(\cA_{qc})}(R \otimes_k F_1,p_\ast^{en} G)) \\
\quad \cong \Hom_{R\Mod(\cA_{qc})}(p_\ast^{en} E,p_\ast^{en} G).
\end{array}
\]
Hence, $p_\ast^{en}$ is fully faithful.

An argument analogous to that above shows that any $G \in R\Mod(\cA_{qc})$ is the cokernel of a map $R \otimes_k E_1 \to R \otimes_k E_0$ with $E_0,E_1 \in \cA_{qc}$. It follows that $p_\ast^{en}$ is essentially surjective, because its essential image contains $R \otimes_k E$ for any $E \in \cA_{qc}$ and is closed under the formation of cokernels because $p_\ast^{en}$ is fully faithful and exact.
\end{proof}


\subsection{The moduli functor \texorpdfstring{$\cM$}{M}}
\label{S:moduli_functor}

In \cite{lieblich2005moduli}*{Defn.~2.1.8}, Lieblich defines the stack $\cD^b_{pug}(X)$ of ``universally gluable relatively perfect complexes on $X$'' and shows that it is an algebraic stack locally of finite type over $k$. Recall that for any $k$-algebra $R$, a complex $E \in \D_{qc}(X_R)$ is perfect relative to $R$ if it is pseudo-coherent and the functor $E \otimes^L_R (-) : \D_{qc}(R\Mod) \to \D_{qc}(X_R)$ has bounded cohomological amplitude in the usual $t$-structure (see \cite{stacks-project}*{Tag 0DHZ, Tag0DI9}).

\begin{defn} [Moduli functor] \label{defn:moduli_functor}
Given a $t$-structure on $\D^b(X)$ and a $k$-algebra $R$, we say that a complex $E \in \D_{qc}(X_R)$ is \emph{$R$-flat} if $E \otimes^L_R M \in \D_{qc}(X_R)^\heartsuit, \forall M \in R\Mod$. We define the moduli of \emph{flat families of objects in $\D^b(X)^\heartsuit$} to be the functor $\cM$ that assigns to an affine $k$-scheme $\Spec(R)$ the groupoid
\[
\cM(R) = \left \{ E \in \D^b(X_R) \text{ perfect relative to } R \text{ and } R\text{-flat}\right\}
\]
\end{defn}

\begin{lem}[Flatness criteria]\label{lem:flatness_criteria}
For $E \in \D_{qc}(X_R)$, the following are equivalent:
\begin{enumerate}
\item $E$ is $R$-flat, i.e., $E \otimes^L_R M \in \D_{qc}(X_R)^\heartsuit, \forall M \in R\Mod$,
\item $\phi^\ast(E) \in \D_{qc}(X_S)^\heartsuit$ for any map $\phi : X_S \to X_R$ induced by a map of $k$-algebras $R \to S$, and
\item $E \otimes^L_R (R/I) \in \D_{qc}(X)^\heartsuit$ for all finitely generated ideals $I \subset R$.
\end{enumerate}
If the $t$-structure on $\D^b(X)$ is non-degenerate, these are equivalent to
\begin{enumerate}[resume]
\item $E \in \D_{qc}(X_R)^\heartsuit$ and the functor $E \otimes_R (-) = H^0(E \otimes^L_R(-)) : R\Mod \to \D_{qc}(X_R)^\heartsuit$ is exact.
\end{enumerate}
Furthermore, if $R$ is Noetherian and $E \in \D_{qc}(X_R)^{\leq 0}$ and is pseudo-coherent, then these are equivalent to
\begin{enumerate}[resume]
\item $E|_{R/\mathfrak{m}} \in \D_{qc}(X_{R/\mathfrak{m}})^\heartsuit$ for all maximal ideals $\mathfrak{m} \subset R$.
\end{enumerate}
\end{lem}
\begin{proof}
By \Cref{lem:t_structure_properties}, $\phi^\ast(E) \in \D_{qc}(X_S)$ lies in the heart if and only if $\phi_\ast(\phi^\ast(E)) = E \otimes^L_R S \in \D_{qc}(X_R)^\heartsuit$, so $(1) \Rightarrow (2) \Rightarrow (3)$ tautologically. To show that $(3) \Rightarrow (1)$, it suffices to show that $E \otimes^L_R M \in \D_{qc}(X_R)^\heartsuit$ for finitely presented $M$, because $\D_{qc}(X_R)^\heartsuit$ is closed under filtered colimits and $R\Mod$ is compactly generated by finitely presented modules. One can show that any finitely generated module $M = R^{\oplus n}/L$ admits a finite filtration whose associated graded pieces are of the form $R/I$ for some ideal $I \subset R$ using induction on $n$ (see the latter part of the proof of \cite{stacks-project}*{Tag 00HD}). $\D_{qc}(X_R)^\heartsuit$ is closed under extension, so it suffices to prove that $E \otimes^L_R (R/I) \in \D_{qc}(X_R)^\heartsuit$. Finally, $R/I = \colim_\alpha R/I_\alpha$, where $\alpha$ indexes the filtered system of finitely generated submodules $I_\alpha \subset I$, so it suffices to prove the claim when $I$ is finitely generated.

One can show directly that $E \otimes^L_R(-) : \D_{qc}(R) \to \D_{qc}(X_R)$ is right $t$-exact for any $E \in \D_{qc}(X_R)^\heartsuit$, by presenting any connective complex of $R$-modules as a complex of free modules $M^\bullet$ in cohomological degree $\leq 0$, then observing that $E \otimes^L_R M^\bullet$ lies in the category generated by $E$ under extensions, left shifts, and filtered colimits. The implication $(1) \Rightarrow (4)$ follows from this observation and the long exact sequence for the homology of an exact triangle.

To show $(4) \Rightarrow (1)$, one considers for any $M \in R\Mod$ a presentation $0 \to K \to R^S \to M \to 0$. The exactness of $H^0(E \otimes^L_R(-))$ and the long exact homology sequence implies that $H^{-1}(E \otimes^L_R M) = 0$, and $H^{-i}(E \otimes^L_R K) \cong H^{-i-1}(E \otimes^L_R M)$ for all $i>0$. Because this holds for all $R$-modules simultaneously, it follows that $H^{-i}(E \otimes^L_R M) = 0$ for all $i>0$. Assuming the $t$-structure is non-degenerate, i.e., $\bigcap_{n \geq 0} \D_{qc}(X)^{\leq -n} = 0$, this implies that $E \otimes^L_R M \in \D_{qc}(X_R)^\heartsuit$.

Finally, we show that when $R$ is Noetherian and $E \in \D^b(X_R)$ then (5) implies (3): \Cref{lem:t_structure_properties}(5) implies that if $T \to \Spec(R)$ is a faithfully flat map of schemes and $E|_{X_T}$ is $T$-flat, then $E$ is $R$-flat. If we let $\hat{R}_\mathfrak{m}$ denote the completion of $R$ at a maximal ideal $\mathfrak{m}$, then $T = \bigsqcup_{\mathfrak{m} \subset R} \Spec(\hat{R}_\mathfrak{m})$ is faithfully flat over $\Spec(R)$ because $R$ is noetherian. It therefore suffices to assume $R$ is a complete noetherian local ring with maximal ideal $\mathfrak{m} \subset R$.

Note that (3) is equivalent to the condition that $E|_{\Spec(R/I)} \in \D_{qc}(X_{R/I})^\heartsuit$ for all ideals $I \subset R$ by \Cref{lem:t_structure_properties}(3). $R/I$ is again a complete local noetherian ring with maximal ideal $\mathfrak{m} / I$.  By replacing $R$ with $R/I$, it suffices to prove the weaker claim that if $E \in \D_{qc}(X_R)^{\leq 0} \cap \D^-_{coh}(X_R)$ and $E|_{R/\mathfrak{m}} \in \D_{qc}(X_{R/\mathfrak{m}})^\heartsuit$, then $E \in \D_{qc}(X_{R})^\heartsuit$.

Considering the exact triangle
\[
E \otimes^L_R \mathfrak{m}^n / \mathfrak{m}^{n+1} \to E \otimes^L_R (R/\mathfrak{m}^{n+1}) \to E \otimes^L_R (R/\mathfrak{m}^{n}) \to
\]
shows that $E \otimes^L_R (R/\mathfrak{m}^n) \in \D_{qc}(X_R)^\heartsuit$ for all $n\geq 1$. The Grothendieck existence theorem for the stable $\infty$-category $\D^-_{coh}(X_R)$ implies that $E = \op{holim}_{n} E \otimes^L_R (R/\mathfrak{m}^n) \in \D_{qc}(X_R)$. Cofiltered homotopy limits are left $t$-exact with respect to any accessible $t$-structure on a presentable stable $\infty$-category, so the fact that $E \otimes^L_R (R/\mathfrak{m}^n) \in \D_{qc}(X_R)^\heartsuit$ for all $n$ implies $E \in \D_{qc}(X_R)^{\geq 0}$ and hence $E \in \D_{qc}(X_R)^\heartsuit$.

\end{proof}

\begin{cor}[Open heart property] \label{cor:open_heart}
Let $R$ be a finite type $k$-algebra and let $E \in \D^b(X_R)$. The set of prime ideals
\[
U :=\left\{ \fp \in \Spec(R) \left| E|_{R_{\fp}} \in \D^b(X_{R_{\fp}})^\heartsuit \right.\right\}
\]
is open, and it contains those primes for which $E|_{\kappa(\fp)} \in \D_{qc}(X_{\kappa(\fp)})^\heartsuit$, where $\kappa(\fp)$ denotes the field of fractions of $R/\fp$.
\end{cor}
\begin{proof}
Because restriction along the map $X_{R_{\fp}} \to X_R$ is $t$-exact, the subset $U$ is the complement of the image under the projection $X_R \to \Spec(R)$ of the closed subsets $\op{Supp}(\tau^{<0}(E))$ and $\op{Supp}(\tau^{>0}(E))$. Therefore $\Spec(R) \setminus U$ is closed because the projection $X_R \to \Spec(R)$ is proper. The fact that if $E|_{\kappa(\fp)} \in \D_{qc}(X_{\kappa(\fp)})^\heartsuit$ then $E|_{R_\fp}\in \D^b(X_{R_\fp})^\heartsuit$ is \Cref{lem:flatness_criteria}(5).
\end{proof}

\subsubsection{Algebraicity of the moduli functor}

\begin{defn}[Generic Flatness \cite{abramovich2006sheaves}] A $t$-structure on $\D^b(X)$ has the \emph{Generic Flatness} property if given a domain $R$ of finite type over $k$ with fraction field $K$ and an object $E \in \D^b(X_R)$ such that $E_K \in \D^b(X_K)^\heartsuit$, there is an $f \in R$ such that $E|_{\Spec(R_f)} \in \D^b(X_{R_f})$ is flat, where $R_f=R[f^{-1}]$.
\end{defn}

This is not exactly how Generic Flatness was formulated in \cite{abramovich2006sheaves}*{Prob.3.5.1}, but over a field of characteristic $0$ this is equivalent to the definition in \cite{abramovich2006sheaves} by the following:
\begin{lem}
If $\op{char}(k)=0$, then Generic Flatness is equivalent to the following condition: for every smooth $k$-algebra $R$ and every $E \in \D^b(X_R)^\heartsuit$, there is a dense open subset $U \subset \Spec(R)$ such that $E|_U$ is flat.
\end{lem}
\begin{proof}
Note that the restriction $\D^b(X_R) \to \D^b(X_K)$ is $t$-exact, so if $E \in \D^b(X_R)^\heartsuit$ then so is $E_K$, and Generic Flatness implies the condition in the lemma. For the other direction assume the condition of the lemma, and consider an integral $k$-algebra $R$ and $E \in \D^b(X_R)$. If $E_K \in \D^b(X_K)^\heartsuit$, then by \Cref{cor:open_heart} we can find a dense open $U \subset \Spec(R)$ such that $E|_U \in \D^b(X_U)^\heartsuit$. By generic smoothness for reduced $k$-algebras we can pass to a smaller open subset $U'' \subset U \subset \Spec(R)$ that is smooth over $k$. The condition of the lemma implies that there is an open $\Spec(R_f) \subset U''$ such that $E|_{\Spec(R_f)}$ is flat, hence we have Generic Flatness.
\end{proof}

\begin{ex}
The usual $t$-structure on $\D^b(X)$ satisfies Generic Flatness (see \cite{stacks-project}*{\href{https://stacks.math.columbia.edu/tag/051R}{Tag 051R}} or \cite{stacks-project}*{\href{https://stacks.math.columbia.edu/tag/052A}{Tag 052A}}).
\end{ex}

Our main result is the following:

\begin{prop} \label{prop:moduli_heart}
Given a noetherian $t$-structure on $\D^b(X)$ that satisfies the Generic Flatness condition, the functor $\cM$ of \Cref{defn:moduli_functor} is an open substack of $\cD_{pug}^b(X)$, hence it is an algebraic stack locally of finite type over $k$ with affine diagonal.\footnote{The fact that $\cD_{pug}^b(X)$ has an affine diagonal is not established in \cite{lieblich2005moduli}, but this follows from a formal argument that applies to any moduli of objects in a $k$-linear category with finite dimensional Hom spaces. See \cite{stacks-project}*{Tag 0DPW}.} Furthermore, this moduli functor agrees with other similar descriptions of moduli functors in the following ways:
\begin{enumerate}
\item On finite type $k$-algebras $R$, $\cM(R)$ is naturally equivalent to the moduli functor of \cite{abramovich2006sheaves}*{Defn.~3.3.1}.
\item If the $t$-structure on $\D^b(X)$ is bounded with respect to the usual $t$ structure, then for any $k$-algebra $R$, $\cM(R)$ is naturally equivalent to
\[
\cM'(R) := \left\{ E \in \D_{qc}(X_R) \text{ that are pseudo-coherent and } R\text{-flat}\right\}.
\]
\item If the $t$-structure on $\D^b(X)$ is noetherian and bounded with respect to the usual $t$-structure, then $\cM(R)$ is naturally equivalent to the moduli functor associated to $\cA_{qc} := \D_{qc}(X)^\heartsuit$ in \cite{AHLH}*{Defn.~7.8},
\[
\cM^{ab}_{\cA_{qc}}(R) := \left\{ E \in R\Mod(\cA_{qc}) \text{ that are finitely presented and } R\text{-flat} \right\},
\]
where $R$-flatness here means that the (non-derived) tensor product functor $E \otimes_R (-) : R\Mod \to R\Mod(\cA_{qc})$ is exact.
\end{enumerate}
\end{prop}

\begin{proof}
Because any object in the heart of a $t$-structure is gluable, we know from \Cref{lem:flatness_criteria}(2) that any complex $E \in \cM(R)$ is universally gluable, hence $\cM$ is a full subfunctor of $\cD_{pug}^b(X)$. We must show that for any $k$-algebra $R$ and any $E \in \cD_{pug}^b(X)(R)$, there is an open subset $U \subset \Spec(R)$ such that for any homomorphism of $k$-algebras $\phi : R \to S$, $\phi^\ast(E) \in \cM(S) \subset \cD_{pug}^b(X)(S)$ if and only if the image of $\Spec(S) \to \Spec(R)$ lies in $U$. Because $E$ is relatively perfect, there is a subalgebra $R' \subset R$ of finite type over $k$ and a relatively perfect complex $E' \in \D^b(X_{R'})$ such that $E = E' \otimes_{R'} R$ \cite{lieblich2005moduli}*{Prop.~2.2.1}. If $U' \subset \Spec(R')$ represents the corresponding functor for $E'$, then the preimage $U$ of $U'$ under $\Spec(R) \to \Spec(R')$ satisfies the desired condition, so we may assume that $R'$ is finite type.

We now claim that if $R$ is finite type and the Generic Flatness property holds, then the set of prime ideals
\[
U :=\left\{ \fp \in \Spec(R) \left| E|_{R/\fp} \in \D_{qc}(X_{R/\fp})^\heartsuit\right.\right\}
\]
is open and satisfies the desired condition. We will show by noetherian induction that $Z := \Spec(R) \setminus U$ is closed. A simple inductive argument reduces one to the case where $R$ is integral, so we will assume this. The property $E|_{R/\fp} \notin \D_{qc}(X_{R/\fp})^\heartsuit$ is closed under specialization by \Cref{cor:open_heart}, so if $K$ is the field of fractions of $R$ and $E_K \notin \D^b(X_K)^\heartsuit$ then $Z = \emptyset$. On the other hand, if $E|_K \in \D^b(X_K)^\heartsuit$ then by Generic Flatness we know that there is an $f$ such that $E|_{R_f}$ is flat, and hence $Z \subset \Spec(R / (f)) \subset \Spec(R)$, which is closed by the inductive hypothesis. By \Cref{lem:flatness_criteria}(5), the restriction $E|_{U} \in \D_{qc}(X_U)$ is $U$-flat and relatively perfect, so $\phi^\ast(E)$ is $S$-flat for any morphism $\phi : \Spec(S) \to \Spec(R)$ landing in $U$. Conversely for any morphism $\phi : \Spec(S) \to \Spec(R)$ such that there is some point $p \in S$ lying over $Z$, $\phi^\ast(E)|_p \notin \D_{qc}(X_{k(p)})^\heartsuit$, so $\phi^\ast(F)$ is not flat by \Cref{lem:flatness_criteria}(2).

\medskip
\noindent{\textit{Comparison with other moduli functors:}}
\medskip

The equivalence of $\cM(R)$ with the moduli functor of \cite{abramovich2006sheaves}*{Defn.~3.3.1} for finite type $k$-algebras $R$ follows immediately from \Cref{lem:flatness_criteria}(5).

If $\D^b(X)$ is bounded with respect to the usual $t$ structure, then we saw in the proof of \Cref{prop:heart_category} that the induced $t$-structure on $\D_{qc}(X)$ is bounded with respect to the usual $t$-structure, and \Cref{lem:t_structure_properties}(2) implies that the induced $t$-structure on $\D_{qc}(X_R)$ is also bounded with respect to the usual $t$-structure for any $k$-algebra $R$. Thus any $R$-flat object $E \in \D_{qc}(X_R)$ also has finite tor amplitude in the usual $t$-structure, hence a pseudo-coherent and $R$-flat complex is relatively perfect, giving the equivalence $\cM'(R) \cong \cM(R)$.

Finally, to prove the natural equivalence $\cM'(R) \cong \cM^{ab}(R)$ under the hypotheses in (3), we use the natural equivalence $R\Mod(\cA_{qc}) \cong \D_{qc}(X_R)^\heartsuit$ from \Cref{prop:heart_category}(2). \Cref{lem:flatness_criteria}(4) implies that this equivalence identifies the given notions of $R$-flatness  in each category. So one must show that an $R$-flat object $E \in \D_{qc}(X_R)^\heartsuit$ is pseudo-coherent if and only if it is finitely presented as an object of the abelian category $\D_{qc}(X_R)^\heartsuit$.

First assume that $R$ is a finitely generated $k$-algebra. Because $R$ is noetherian and the induced $t$-structure on $\D_{qc}(X_R)$ is bounded with respect to the usual $t$-structure, $E \in \D_{qc}(X_R)^\heartsuit$ is pseudo-coherent if and only if $E \in \D^b(X_R)$. On the other hand, in the proof that $\cA_{qc}$ is locally noetherian in \Cref{prop:heart_category}, we showed that $\D^b(X)^\heartsuit \subset \D_{qc}(X)^\heartsuit$ is the subcategory of finitely presented objects. Because $X_R$ is noetherian and the $t$-structure on $\D_{qc}(X_R)$ preserves $\D^b(X_R)$ by \Cref{prop:t_structure_coherent}, the same proof applies verbatim to show that $\D^b(X_R)^\heartsuit \subset \D_{qc}(X_R)^\heartsuit$ is the subcategory of finitely presented objects. The claim follows.

For general $R$, we use noetherian approximation: If $E \in \D_{qc}(X_R)^\heartsuit$ is $R$-flat and pseudo-coherent, then because the stack $\cM'$ is algebraic and locally of finite presentation, $E$ is obtained by base change from a flat and pseudo-coherent $E' \in \D_{qc}(X_{R'})^\heartsuit$ for some finitely generated subalgebra $R' \subset R$. From the previous paragraph, we know that $E'$ is finitely presented as an object of the abelian category $\D_{qc}(X_{R'})^\heartsuit$, and thus $E = E' \otimes_{R'} R \in \D_{qc}(X_R)^\heartsuit$ is finitely presented.

Conversely, assume $E \in R\Mod(\cA_{qc})$ is finitely presented. Using the fact that $R\Mod(\cA_{qc})$ is compactly generated by objects of the form $R \otimes_k E'$ for $E' \in \D^b(X)^\heartsuit$, one can show that $E \cong \coker(R \otimes_k F_1 \to R \otimes_k F_0)$ for some $F_0,F_1 \in \D^b(X)^\heartsuit$. The morphism $R \otimes_k F_1 \to R \otimes_k F_0$ corresponds to a map $F_1 \to R \otimes_k F_0$ in $\D_{qc}(X)^\heartsuit$. Because $F_1$ is finitely presented, this morphism must factor through $R' \otimes_k F_0 \to R \otimes_k F_0$ for some finitely generated subalgebra $R' \subset R$. In particular, $E \cong H^0(R \otimes_{R'}^L E')$ for some finitely generated $R' \subset R$ and $E' \in \D^b(X_{R'})^\heartsuit$. The morphism $\Spec(R) \to \Spec(R')$ factors through the maximal open subscheme $U \subset \Spec(R')$ over which $E'$ is flat, which implies that $R \otimes_{R'}^L E' \in \D_{qc}(X_R)^\heartsuit$, so $E \cong R \otimes_{R'}^L E'$ is pseudo-coherent.

\end{proof}

\begin{ex} \label{ex:toda_1}
In \cite{toda2008moduli}*{Lem.~4.7, Prop.~3.18}, Toda shows that the Generic Flatness property holds for a set of algebraic Bridgeland stability conditions (meaning $Z(\cA) \subset \bQ+i\bQ$) that is dense (up to the action of $\tilde{\GL}_2^+(2,\bR)$) in the connected component of the space of Bridgeland stability conditions on $K3$ surfaces over $\bC$ constructed and studied in \cite{bridgeland2008stability}. More recently, Piyaratne and Toda establish the Generic Flatness property for certain Bridgeland stability conditions on $3$-folds in \cite{piyaratne2015moduli}.
\end{ex}

\begin{rem}
For a Bridgeland stability condition $(Z,\cA)$, the methods of \cites{toda2008moduli,piyaratne2015moduli} show that the stack of semistable objects of a given numerical $K$-theory class is an open substack of $\cD^b_{pug}(X)$ whenever $(Z,\cA)$ is sufficiently close to a stability condition with noetherian heart that satisfies the Generic Flatness property.
\end{rem}

\subsubsection{$\Theta$-reductivity and {\textsf{S}}-completeness}

In an earlier version of this paper, we showed that for any noetherian $t$-structure on $\D^b(X)$, the morphism $\ev_1 : \Filt(\cM) \to \cM$ satisfies the valuative criterion for properness for any discrete valuation ring that is essentially finite type over $k$. Subsequently, this material was included and expanded to treat {\textsf{S}}-completeness in \cite{AHLH}*{Sect.~7.2}, so for the sake of brevity here we simply recall those results in our context:

\begin{prop} \label{prop:theta_reductive_S_complete}
Fix a $t$-structure on $\D^b(X)$ that is noetherian, bounded with respect to the usual $t$-structure, and satisfies Generic Flatness. Then the stack $\cM$ of \Cref{defn:moduli_functor} is $\Theta$-reductive and {\textsf{S}}-complete with respect to DVR's essentially of finite type over $k$.
\end{prop}

\begin{proof}
Because we are assuming Generic Flatness, the moduli functor $\cM$ agrees with the moduli functor $\cM^{ab}$ of \cite{AHLH}*{Defn.~7.8} by \Cref{prop:moduli_heart}(3), so this follows from \cite{AHLH}*{Lem.~7.15\&7.16}.
\end{proof}


\subsection{Filtered points of \texorpdfstring{$\cM$}{M}}

Let $R$ be a $k$-algebra. We regard $\bZ^n$ as the objects of a category in which there is a unique morphism $(w_1,\ldots,w_n) \to (w'_1,\ldots,w'_n)$ if $w_i \geq w'_i$ for all $i$. The Reese construction (see \Cref{prop:quasicoh_theta}) gives an equivalence of $\infty$-categories
\[
\D_{qc}(\Theta^n_R \times X) \simeq \op{Fun}(\bZ^n,\D_{qc}(X_R)).
\]
More explicitly, given a quasi-coherent complex on $X\times \Theta^n_R$, we may regard it as a complex of $\bZ^n$-graded quasi-coherent modules over $\cO_{X_R}[t_1,\ldots,t_n]$, where $t_i$ has multi-degree $(0,\ldots,0,-1,0,\ldots,0)$ with $-1$ in the $i^{th}$ coordinate. Pushing forward along the map $X \times \Theta_R^n \to X \times (\pt/\Gm)_R^n$, which corresponds to forgetting the $\cO_{X_R}[t_1,\ldots,t_n]$-module structure, we have an isomorphism $E \simeq \bigoplus_{w\in \bZ^n} E_w$, where the quasi-coherent complex $E_w$ is the summand in weight $w$. Multiplication by the monomial $t^v = t_1^{v_1} \cdots t_v^{v_n}$ gives a map $E_{w+v} \to E_w$ for all $w \in \bZ^n$.

Now consider a functor $\bZ^n \to \D_{qc}(X_R)^\heartsuit$, corresponding to objects $E_w \in \D_{qc}(X_R)^\heartsuit$ and maps $E_w \to E_{w'}$ for $w \geq w'$. We call such a functor a $\bZ^n$-weighted filtration if all of the maps $E_w \to E_{w'}$ are injective. We denote the $i^{th}$ standard basis vector in $\bZ^n$ by $e_i$. By definition $\gr(E_\bullet)$ is the $\bZ^n$-graded object for which
\[
\gr_w(E_\bullet) := \op{coker}(\bigoplus_i E_{w+e_i} \to E_w).
\]
For any $w$ and any ordered tuple $i_1 < \ldots < i_m$, let us denote the map by
\[
\delta^p_{i_1<\ldots<i_m} : E_{w+e_{i_1} + \cdots + e_{i_m}} \to E_{w+e_{i_1} + \cdots+e_{i_{p-1}} + e_{i_{p+1}} + \cdots + e_{i_m}}.
\]
Using this notation we define the augmented Koszul complex in $\D_{qc}(X_R)^\heartsuit$
\begin{multline}\label{eqn:koszul_complex}
0 \to E_{w+e_1+\cdots+e_n} \to \cdots \to \bigoplus_{i_1<i_2} E_{w+e_{i_1}+e_{i_2}} \to \bigoplus_{i_1} E_{w+e_{i_1}} \to E_w \to \\ \op{gr}_w(E_\bullet) \to 0,
\end{multline}
where the differential is given by $\sum_{i_1 < \cdots < i_m} \sum_{p=1}^n (-1)^{p+1} \delta_{i_1<\cdots<i_m}^p$. We define the Koszul complex $\cK^\bullet_E$ to be the complex above without the $\gr_w(E_\bullet)$ term on the right. So there is a canonical map from the totalization $\op{Tot}^{\oplus}(\cK_E^\bullet) \to \gr_w(E)$ in $\D_{qc}(X_R)$, and the totalization of the augmented Koszul complex is the cone of this map.

\begin{lem} \label{lem:filtrations_derived}
Let $R$ be a $k$-algebra. Under the Reese construction above, $\Filt^n(\cM)(R) = \Map(\Theta^n_R,\cM)$ is naturally isomorphic to the groupoid of $\bZ^n$-weighted filtrations $\{E_w\}_{w \in \bZ^n}$ in $\D_{qc}(X_R)^\heartsuit$ such that
\begin{enumerate}
\item for all $w \in \bZ^n$, $\gr_w(E_\bullet)$ is $R$-flat and relatively perfect over $R$,
\item $\exists N \gg0$ such that $E_{w_1,\ldots,w_n} = 0$ if $w_i > N$ for any $i$,
\item $\exists N \ll 0$ such that $E_{w_1,\ldots,w_n} \to E_{w'_1,\ldots,w'_n}$ is an isomorphism if $w_i < N$ for any $i$, and
\item for all $w \in \bZ^n$ the augmented Koszul complex \eqref{eqn:koszul_complex} is exact.
\end{enumerate}
The fiber of a map $f : \Theta_R^n \to \cM$ at $(1,\ldots,1)$ corresponds to $\op{colim}_w {E_w} = E_{w'} \in \D^b(X_R)$ for $w' \ll 0$, and the fiber at $(0,\ldots,0)$ corresponds to $\bigoplus_w \gr_w(E_\bullet)$ in $\D^b(X_R)$.
\end{lem}

\begin{rem} \label{rem:characterize_filtration_complexes}
Conditions (1) and (4) together should be interpreted as the condition that if $E \in \D_{qc}(X \times \Theta_R^n)$ is the complex corresponding to the $\bZ^n$ weighted filtration $\{E_w\}_{w \in \bZ^n}$, then each weight summand of
\[
E \otimes_{\Theta_R^n} \cO_{\Spec(R) \times \{0,\ldots,0\}} \in \D_{qc}(X \times (\pt/\Gm^n)_R)
\]
is flat and relatively perfect over $R$. This graded object is then $\bigoplus_{w \in \bZ^n} \gr_w(E)$. Therefore (1) and (4) are equivalent to the condition that the (derived) restriction of $E$ along the closed immersion $(\{0\}/\Gm)_R^n \hookrightarrow \Theta_R^n$ is flat, and its weight summands with respect to $(\Gm)_R^n$ are relatively perfect over $R$.
\end{rem}

\begin{proof}

\medskip
\noindent{\textit{First we show that given a filtration $\{E_w\}_{w\in \bZ^n}$ satisfying these conditions, the corresponding object $E \in \D_{qc}(X \times \Theta_R^n)^\heartsuit$ is flat and relatively perfect over $\Theta_R^n$:}}
\medskip

Consider the restriction of $E$ along the closed immersion $i_n : \Theta_R^{n-1} \times (\pt/\Gm) = \{t_n=0\} \hookrightarrow \Theta_R^n$, and let $F_r$ for $r \in \bZ$ denote the summand of $(i_n)^\ast (E)$ concentrated in weight $r$ with respect to the last factor of $\pt/\Gm$. Regarding $i_n^\ast(E)$ as a $\bZ^n$-graded object in $\D_{qc}(X_R)$, we weight $w$ piece is
\[
\op{Cone}(E_{w_1,\ldots,w_{n-1},w_n+1} \xrightarrow{t_n} E_{w_1,\ldots,w_n}),
\]
and conditions (2) and (3) imply that this vanishes for $w_n \ll 0$ or $w_n \gg 0$. It follows that $F_r \neq 0$ for only finitely many $r \in \bZ$. If $r_{max}$ denotes the largest $r$ for which $F_r \neq 0$, then we have an exact triangle
\[
F_{r_{max}} \otimes k[t_n] \twist{r_{max}} \to E \to E' \to
\]
where the left object denotes the pullback along the projection $\Theta_R^n \to \Theta_R^{n-1}$ along the last copy of $\Theta$ followed by a twist by a character of $\Gm$ so that $(i_n)^\ast(F_{r_{max}} \otimes k[t_n] \twist{r_{max}})$ is concentrated in weight $r_{max}$ with respect to the last $\Gm$. This map is uniquely determined by requiring that it is an isomorphism on highest weight pieces (with respect to the last factor of $\pt/\Gm$) after restriction along $i_n$. In particular if $F'_r$ denote the weight summands of $i_n^\ast(E')$, then $F'_r = F_r$ for $r \neq r_{max}$ and $F'_{r_{max}} = 0$. It follows that $E$ can be constructed as an iterated extension of objects of the form $F \otimes k[t_n] \twist{r}$ for certain $F \in \D_{qc}(X \times \Theta_R^{n-1})$ and $r \in \bZ$.

We argue by induction on $n$ that $E$ is flat and relatively perfect over $\Theta_R^n$. Because flat and relatively perfect complexes are closed under extension, pullback, and tensoring by an invertible sheaf on the base (which is what the shift by a character of $\Gm$ is doing), it suffices to show that the $F_r$ constructed above are flat and relatively perfect. The base case is when $n=1$, in this case (4) implies that the restriction $(i_1)^\ast(E)$ lies in $\D_{qc}(X \times (\Gm)_R)^\heartsuit$ and then the weight $r$ pieces are just $\gr_r(E)$, which are flat and relatively perfect by (1).

For $n>1$, consider $\bigoplus_r F_r \twist{r} \simeq i_n^\ast(E)$. As discussed in \Cref{rem:characterize_filtration_complexes}, conditions (1) and (4) are equivalent to the condition that $E|_{(\pt/\Gm)_R^n} \in \D_{qc}(X \times (\pt/\Gm)_R^n)$ is flat and its weight summands are relatively perfect. It follows that (1) and (4) still hold for each $F_r$ because the closed immersion $(\pt/\Gm)_R^n \hookrightarrow \Theta_R^n$ factors through $\Theta_R^{n-1} \times (\pt/\Gm)$. The weight $(w_1,\ldots,w_{n-1})$ component of $F_r$, regarded as a graded object in $\D_{qc}(X_R)$ is the cone of the map $E_{w_1,\ldots,w_{n-1},r+1} \to E_{w_1,\ldots,w_{n-1},r}$, so conditions (2) and (3) imply that this cone is $0$ if $w_i \ll 0$ or $w_i \gg 0$ for any $i$. It follows that $F_r$ satisfies (1)-(4), and hence by the inductive hypothesis $F_r$ is flat and relatively perfect over $\Theta_R^{n-1}$.

\medskip
\noindent{\textit{Second, we show that for $E \in \D_{qc}(X \times \Theta_R^n)$ that is flat and relatively perfect over $\Theta_R^n$ the corresponding functor $\{E_w\}_{w\in \bZ^n}$ is a filtration satisfying (1)-(4):}}
\medskip

As remarked above, conditions (1) and (4) together are equivalent to requiring that the restriction of $E$ along the inclusion $(\pt/\Gm)_R^n \hookrightarrow \Theta_R^n$ is flat and has relatively perfect weight summands with respect to $\Gm^n$. This is clearly satisfied if $E$ is relatively perfect and flat.\footnote{Note that in fact $E|_{(\pt/\Gm)_R^n}$ will have finitely many non-zero graded piece for only finitely many $w$. This is close to but slightly weaker than (2) and (3).} Next we claim that (2) and (3) hold for \emph{any} relatively perfect and universally gluable $E \in \D^b(X \times \Theta_R^n)$. First by Artin approximation for relatively perfect and universally gluable complexes \cite{lieblich2005moduli}*{Prop.~2.2.1}, we can find a finitely generated sub-algebra $R' \subset R$ and a universally gluable and relatively perfect $E' \in \D^b(X \times \Theta_{R'}^n)$ along with an isomorphism $E \simeq E'\otimes_{R'}R$. We may therefore assume that $R$ is noetherian.

We now claim that for a noetherian $k$-algebra $R$, conditions (2) and (3) apply for any $E \in \D^b(X \times \Theta_R^n)$. Indeed because pushforward to $X \times (\pt/\Gm)_R^n$ commutes with taking homology in the usual $t$-structure, we can identify the graded summands $H^i(E)_{w} \simeq H^i(E_w)$ in the usual $t$-structure for $w \in \bZ^n$. So it suffices to prove the claim for $E \in \Coh(X \times \Theta_R^n)$. In this case there are graded coherent sheaves $F_0,F_1 \in \Coh(X \times (\pt/\Gm)_R^n)$ such that
\[
E = \op{coker}(F_1 \otimes k[t_1,\ldots,t_n] \to F_0 \otimes k[t_1,\ldots,t_n]),
\]
where $F_i \otimes k[t_1,\ldots,t_n]$ denotes the pullback of $F_i$ along the projection $X \times \Theta_R^n \to X \times (\pt/\Gm)_R^n$. Conditions (2) and (3) hold for the weight $w$ summands of $F_i \otimes k[t_1,\ldots,t_n]$ for $i=0,1$, and it follows that these conditions hold for the $E_w$ as well.
\end{proof}

\begin{ex}
When $n=2$ in the previous lemma, the conditions (1) and (4) amount to the requirement that all maps $E_{v+1,w} \to E_{v,w}$ and $E_{v,w+1} \to E_{v,w}$ are injective, and that every square
$$\xymatrix{E_{v+1,w+1} \ar[r] \ar[d] & E_{v+1,w} \ar[d] \\ E_{v,w+1} \ar[r] & E_{v,w} }$$
is a fiber product.
\end{ex}

For clarity of notation we will use underscored letters to denote vectors $\underline{w}_1,\underline{w}_2,\ldots \in \bZ^n$. As above we say $\underline{w}_1 \leq \underline{w}_2$ if the inequality holds for every entry. Given a flat and relatively perfect $E \in D_{qc}(\Theta_R^n \times X)$, let us denote $E_{-\infty} = E_w$ for $w \ll 0$. Then forgetting the $\cO_{X_R}[t_1,\ldots,t_n]$ module structure, we have an injective map of $\bZ^n$-graded objects in $\D_{qc}(X_R)$
\[
E \hookrightarrow E_{-\infty} \otimes^L_R R[t_1^\pm,\ldots,t_n^\pm].
\]
So a $\bZ^n$-weighted filtration is determined by an object $E_{-\infty} \in \cM(\Spec(R))$ and an $\cO_{X_R}[t_1,\ldots,t_n]$-submodule of $\cO_{X_R}[t_1^\pm,\ldots,t_n^\pm] \otimes E_{-\infty}$ of the form
\begin{equation} \label{eqn:filtered_object}
E = \sum \cO_{X_K} [t_1,\ldots,t_n] \cdot t^{-\underline{w}_i} \otimes E_i \subset E_{-\infty} \otimes_K K[t_1^\pm, \ldots,t_n^\pm]
\end{equation}
for some finite collection of subobjects $E_i \subset E_{-\infty}$ with the property that $E_i \subset E_j$ whenever $\underline{w}_i \geq \underline{w}_j$, and which satisfy the conditions of \Cref{lem:filtrations_derived}. Over a field, a filtration is non-degenerate if and only if the weights $w \in \bZ^n$ for which $\gr_w(E_\bullet) \neq 0$ span $\bQ^n$.

Any map $\phi : \Theta_K^p \to \Theta_K^n$ is encoded by an $n \times p$ matrix with nonnegative integer entries $A$, and we have
\[
\phi^\ast (E) = \sum \cO_X[t_1,\ldots,t_p] \cdot t^{-A^\tau \underline{w}_i} \otimes E_i \subset E \otimes_K  K[t_1^\pm, \ldots, t_p^\pm]
\]
where $A^\tau$ denotes the transpose matrix, and we are regarding $\underline{w}_i$ as a column vector. In particular for any point in $a \in \bR^n_{\geq 0}$ with integer coefficients, regarded as a $n \times 1$ matrix, the pullback $\phi^\ast(E)$ under the corresponding map $\phi : \Theta_K \to \Theta_K^n$ is the $\bZ$-weighted filtration
\[
E \supset \cdots \supset E'_v \supset E'_{v+1} \supset \cdots, \text{ where } E'_v = \sum_{i | a^\tau \underline{w}_i \geq v} E_i.
\]

The object \eqref{eqn:filtered_object} automatically satisfies conditions (2) and (3) of \Cref{lem:filtrations_derived}, but it does not necessarily define a flat $\cO_{X_R}[t_1,\ldots,t_n]$ module if $n>1$. However for any descending non-weighted filtration $E_1 \supsetneq E_2 \supsetneq E_3 \supsetneq \cdots \supsetneq E_n \neq 0$ in $\D^b(X_R)^\heartsuit$ whose associated graded complexes are flat and relatively perfect over $\Spec(R)$, and for any weight sequence $\underline{w}_1 \leq \underline{w}_2 \leq \cdots \leq \underline{w}_n$ we define the canonical $\bZ^n$-weighted filtration $\sigma(E_\bullet,\underline{w}_\bullet) \in \D^b(X \times \Theta_R^n)^\heartsuit$ as the object \eqref{eqn:filtered_object}. We leave it as an exercise to the reader to check conditions (1) and (4), i.e., that the totalization of the Koszul complex is flat and relatively perfect over $\Spec(R)$ for each weight $\underline{w}$. More precisely, $\op{Tot}^{\oplus}(\cK^\bullet_{\sigma(E_\bullet,\underline{w}_\bullet)})$ in weight $\underline{w}_i$ is quasi-isomorphic to $\gr_{\underline{w}_i}(\sigma(E_\bullet,\underline{w}_\bullet)) = E_i / E_{i+1}$ for $i = 1,\ldots,n$, and $\op{Tot}^{\oplus}(\cK^\bullet_{\sigma(E_\bullet,\underline{w}_\bullet)})$ is exact in every other weight. In particular the filtration is non-degenerate if and only if $\underline{w}_1,\ldots,\underline{w}_n$ span $\bQ^n$. We use the same notation to denote the corresponding non-degenerate map $\sigma(E_\bullet,\underline{w}_\bullet) : \Theta_R^n \to \cM$.

\begin{prop} \label{prop:degeneration_space}
Let $K/k$ be a field extension, and let $E \in \cM(\Spec(K))$. The rational simplices corresponding to the weighted filtrations $\sigma(E_\bullet,\underline{w}_\bullet)$, for all finite length filtrations $E = E_1 \supsetneq \cdots \supsetneq E_n \supsetneq 0$ and weight sequences $\underline{w}_1 \leq \cdots \leq \underline{w}_n \in \bZ^n$, cover $\iDeg(\cM,[E])$.
\end{prop}
\begin{proof}
We will show that for any rational cone $\bR_{\geq 0}^n \to |\Deg(\cM,E)_\bullet|$, corresponding to a non-degenerate filtration $f : \Theta_K^n \to \cM$, there is a decomposition of $\bR^{n}_{\geq 0}$ into rational polyhedral sub-cones such that the restriction of $f$ along any map $\phi : \Theta_K^p \to \Theta_K^n$ for which $\phi(\bR^{p}_{\geq 0}) \subset \bR^n_{\geq 0}$ is contained in one of the subcones, the restriction of $f \circ \phi$ factors through one of our canonical maps $\sigma(E_\bullet,\underline{w}_\bullet) : \Theta^N_K \to \cM$. In the language of \Cref{sect:structures}, if $F_\bullet \subset \Deg(\cM,E)_\bullet$ denotes the smallest sub-fan containing cones of the form $\sigma(E_\bullet,\underline{w}_\bullet)$, then we will show that the inclusion $F_\bullet \subset \Deg(\cM,E)_\bullet$ is a bounded inclusion. Because $F_1 = \Deg(\cM,E)_1$, it then follows from \Cref{lem:fans_bounded_surjectity} and \Cref{cor:fans_closed_embedding} that $\bP(F_\bullet) \to \iDeg(\cM,E)$ is a homeomorphism, which implies the claim of the proposition.

Consider a presentation of $E$ of the form \eqref{eqn:filtered_object}. Note that the objects $E_i$ need not be totally ordered by inclusion. However, if the weights $\underline{w}_i$ are totally ordered with respect to the partial order on $\bZ^n$, then we can re-express the object $E = \sum \cO_{X_K} [t_1,\ldots,t_n] \cdot t^{-\underline{w}_i} \otimes E'_i$, where $E'_i = \sum_{\underline{w}_j \geq \underline{w}_i} E_j$. If $E'_i=E'_j$ for some $i\neq j$, then we only keep one term in the sum with $E_i$ corresponding to the largest value of $\underline{w}_j$ for which $E_j = E_i$. The result, after reindexing, is an equality $E = \sum \cO_{X_K} [t_1,\ldots,t_n] \cdot t^{-\underline{w}_i} \otimes E'_i$ where $E'_1 \supsetneq E'_2 \supsetneq \cdots \supsetneq E'_{N} \supsetneq 0$ and $\underline{w}_i \geq \underline{w}_{i+1}$ (and not equal) for all $i$. Let $A$ be the $N \times n$ matrix whose $i^{th}$ row vector is $\underline{w}_i$ for $i=1,\ldots,N$. Then $A$ classifies a map $\Theta_K^n \to \Theta_K^N$ such that the given $\bZ^n$-weighted filtration $f : \Theta_K^n \to \cM$ is the composition of this map with the standard map $\sigma(E'_\bullet,\underline{w}_\bullet) : \Theta_K^N \to \cM$. Thus the filtration corresponding to $E$ lies in $F_\bullet$.

Now for a general filtration $f : \Theta_K^n \to \cM$, consider the subdivision of the positive cone $(\bR_{\geq 0})^n$ into rational polyhedral cones obtained by intersecting with the arrangement of hyperplanes $H_{i,j}=\{ r \in \bR_{\geq 0}^n | r \cdot (\underline{w}_i - \underline{w}_j) = 0\}$ for all $\underline{w}_i$ and $\underline{w}_j$ in \eqref{eqn:filtered_object}. Let $A : \bR^p \to \bR^n$ be a linear map with integer coefficients mapping $\bR^p_{\geq 0}$ to a single cone in this decomposition of $\bR^n_{\geq 0}$. Then by definition $r \cdot A^\tau \underline{w}_i - r \cdot A^\tau \underline{w}_j$ has the same sign (or vanishes) for all $r \in \bR_{\geq 0}^p$. In particular the vectors $A^\tau \underline{w}_i$ are totally ordered with respect to the partial order on $\bZ^p$ discussed above. It follows from the previous paragraph that the restriction to $f \circ \phi$ factors through a map of the form $\sigma(E'_\bullet,\underline{w}'_\bullet)$.
\end{proof}

\begin{rem} \label{rem:real_weights}
The canonical filtration $\sigma(E_\bullet,w_\bullet)$ defines a cone $\bR_{\geq 0}^n \to |\Deg(\cM,E)_\bullet|$. It is natural to regard vectors $a \in \bR_{\geq 0}^n$ as assigning the real weights $a \cdot w_1 \leq a \cdot w_2 \leq \cdots \leq a \cdot w_n$ to the filtration $E_1 \supset E_2 \supset \cdots \supset E_n$. We therefore interpret \Cref{prop:degeneration_space} as identifying $|\Deg(\cM,E)_\bullet|$ with the space of finite real weighted filtrations of $E$, and identifying $\iDeg(\cM,E)$ with the space of finite real weighted filtrations of $E$ up to positive rescaling of weights. 
\end{rem}

\subsection{Slope semistability as \texorpdfstring{$\Theta$}{Theta}-stability}

We let $\Knum(X)$ denote the numerical $K$-group of coherent sheaves on $X$, which we define as the image of the map $K_0(\D^b(X)) \to \Hom(K^0(\Perf(X)),\bZ)$ induced by the Euler pairing. $\Knum(X)$ is the quotient of $K_0(\D^b(X))$ by the subgroup of classes $[F]$ for which $\chi(X,E \otimes F) = 0, \forall E \in \Perf(X)$. Likewise we define $\Knum^{perf}(X)$ to be the quotient of $K^0(\Perf(X))$ by the subgroup of classes $[E]$ for which $\chi(E \otimes F)=0, \forall F \in \D^b(X)$. By construction $\Knum(X)$ and $\Knum^{perf}(X)$ are torsion-free, and $\chi$ descends to a bilinear pairing $\chi : \Knum^{perf}(X) \otimes \Knum(X) \to \bZ$. When $k = \bC$, then Riemann-Roch shows that $\Knum^{perf}(X)$ and $\Knum(X)$ are finitely generated, and hence $\chi$ is a perfect pairing after tensoring with $\bQ$.

\begin{lem}
Let $T$ be a connected $k$ scheme of finite type, and let $E \in \D^b(X_T)$ be relatively perfect. For any finite type point $t \in T$, consider the class
\[
v := \frac{1}{\deg(\kappa(t)/k)} [E_{\kappa(t)}] \in \Knum(X)_\bQ,
\]
where we regard $E_{\kappa(t)}$ as a complex on $X$ via pushforward along the map $X_{\kappa(t)} \to X$. Then $v$ is independent of the choice of $t$, and we may write $\cD_{pug}^b(X)$ (respectively $\cM$) as a disjoint union of open and closed substacks $\cD_{pug}^b(X)_v$ (respectively $\cM_v$) parameterizing families of a fixed class $v \in \Knum(X)$.
\end{lem}
\begin{proof}
If $p :X_T \to T$ denotes the projection, then the semicontinuity theorem implies that $p_\ast((\cO_T \boxtimes F) \otimes E) \in \Perf(T)$ for any $F \in \Perf(X)$. Therefore the Euler characteristic
\[
\chi(X_{\kappa(t)},E_{\kappa(t)} \otimes (\kappa(t)\otimes_k F)) = \deg(\kappa(t)/k) \chi(X,E_{\kappa(t)} \otimes_k F)
\]
does not depend on the finite type point $t \in T$, and the claim follows. This splits $\cD^b_{pug}(X)$ into open and closed substacks because $\cD^b_{pug}(X)$ is algebraic and locally of finite type over $k$, and it splits $\cM$ into open and closed substacks (even when $\cM$ is not algebraic), because $\cM$ is a substack of $\cD^b_{pug}(X)$.
\end{proof}

\begin{defn} \label{def:slope_stability}
A \emph{slope-stability condition} on $X$ is a triple $(Z,\cA,\D^b(X)^\heartsuit)$, where: 1) $\D^b(X)^\heartsuit$ is the heart of a $t$-structure on $\D^b(X)$; 2) $\cA \subset \D^b(X)^\heartsuit$ is a numerical subcategory in the sense that there is a finitely generated subgroup $\Lambda \subset \Knum(X)$ such that $\cA = \{E \in \D^b(X)^\heartsuit | [E] \in \Lambda\}$; and 3) $Z$ is a group homomorphism $Z : \Lambda \to \bC$, known as the central charge, such that $Z(\cA) \subset {\mathbb H} \cup \bR_{\leq 0}$.
\end{defn}

Note that for any short exact sequence in $\D^b(X)^\heartsuit$, if any two of the objects lie in the category $\cA$ then so does the third, so $\cA$ is abelian and the inclusion $\cA \subset \D^b(X)^\heartsuit$ is exact. The previous lemma allows us to regard $\cM_v$ for $v \in \Lambda$ as the moduli of objects in $\cA$, and in particular filtrations in $\cM_v$ correspond under \Cref{lem:filtrations_derived} to weighted filtrations in the abelian category $\D^b(X)^\heartsuit$ whose associated graded object lies in $\cA$.

\begin{ex} \label{ex:Bridgeland_curves}
The simplest example is where $X$ is a smooth curve, and we let $\cA$ be the category of coherent sheaves on $X$. The central charge is $Z(E) = - \deg(E) + i \op{rank}(E)$.
\end{ex}

The value of \Cref{def:slope_stability} is that it simultaneously covers the following commonly studied generalizations of \Cref{ex:Bridgeland_curves}.

\begin{ex}
If $(Z,\cA)$ is numerical Bridgeland stability condition on a smooth projective variety $X$, then $\cA \subset \D^b(X)$ is the heart of a $t$-structure by definition, and $(Z,\cA,\cA)$ is a slope stability condition. In addition to the Harder-Narasimhan property, which we discuss below, a Bridgeland stability condition must have $Z(E) \neq 0$ for nonzero $E \in \cA$. We will not need this nondegeneracy hypothesis until \Cref{sect:theta_strat_Bridgeland}.
\end{ex}

\begin{ex} \label{ex:moduli_coherent_sheaves}
We can fix a projective variety $X$ with ample invertible sheaf $L$, and choose an integer $0 < d \leq \dim X$. We consider the usual $t$-structure on $\D^b(X)$ and let $\Lambda \subset \Knum(X)$ be the subgroup of classes whose Hilbert polynomial has degree $\leq d$. Then $\cA = \op{Coh}(X)_{\leq d}$ is the category of coherent sheaves whose support has dimension $\leq d$. This is a full abelian subcategory of $\op{Coh}(X)$. The assignment
$$E \mapsto p_E(n) := \chi(E \otimes L^n) = \sum \frac{p_{E,k}}{k!} n^k$$
defines a group homomorphism $K_0(X) \to \bZ[n]$. We define our central charge to be $Z(E) := i p_{E,d} - p_{E,d-1}$. This has the property that $Z(E) \subset \mathbb{H} \cup \bR_{\leq 0}$ for any $E \in \cA$: if $\dim (\op{supp}(E)) < d-1$ then $Z(E) = 0$, if $\dim(\op{supp} E) = d-1$ then $p_{E,d} = 0$ and $p_{E,d-1} > 0$, and if $\dim(\op{supp} E) = d$ then $p_{E,d} > 0$.
\end{ex}

For a nonzero $E \in \cA$, the phase $\phi(E) \in (0,1]$ is the unique number such that $Z(E) = |Z(E)| e^{2 \pi i \phi(E)}$, and by convention $\phi = 1$ if $Z(E) = 0$. The following is a common notion of semi-stability in abelian categories:

\begin{defn} \label{defn:slope_semistability}
Let $(Z,\cA,\D^b(X)^\heartsuit)$ be a slope stability condition on $X$. An object $E \in \cA$ is \emph{slope semistable} if there are no sub-objects $F \subset E$ in the abelian category $\cA$ such that $\phi(F) > \phi(E)$. We say that $(Z,\cA,\D^b(X)^\heartsuit)$ has the \emph{Harder-Narasimhan} property if every object $E$ admits a filtration by subobjects $E = E_1 \supset \cdots \supset E_p \supset E_{p+1} = 0$ in $\cA$ with each subquotient $\gr_j E_\bullet$ semistable and $\phi(\gr_j E_\bullet)$ increasing in $j$.
\end{defn}

Our goal in this subsection is to construct a numerical invariant on $\cM_v$ for $v \in \Lambda$ such that $\Theta$-stability and the HN problem correspond to the notions of \Cref{defn:slope_semistability}.

It follows from the fact that $\Lambda$ is finitely generated that the Euler characteristic induces a perfect pairing $\Lambda_\bQ \otimes_\bQ (\Knum^{perf}(X) / \Lambda^\perp)_\bQ \to \bQ$, where $\Lambda^\perp := \{v \in \Knum(X)^{perf} | \chi(v\cdot \lambda)=0, \forall \lambda \in \Lambda\}$. Therefore we can find a unique element $\omega_Z \in (\Knum^{perf}(X) / \Lambda^\perp)_\bC$ such that
$$Z(x) = \chi(\omega_Z \otimes x), \forall x \in \Lambda_\bC$$
We will regard $\omega_Z$ as a class $\omega_Z \in K^0(\Perf(X))\otimes \bC$ by choosing a lift under the surjective map $K_0(\Perf(X)) \otimes \bC \twoheadrightarrow (\Knum^{perf}(X) / \Lambda^\perp)_\bC$. Let $\cE$ denote the universal object in the derived category of $\cM \times X$, and consider the diagram
$$\xymatrix{ \cM_v & \cM_v \times X \ar[r]^{p_2} \ar[l]_{p_1} & X }$$
We define the cohomology classes
\begin{equation} \label{eqn:classes_Bridgeland}
\renewcommand\arraystretch{1.5}
\begin{array}{l}
l := |Z(v)|^2 \ch_1 \left( (p_1)_\ast \left( \cE \otimes p_2^\ast \Im \left(\frac{-\omega_Z}{Z(v)} \right)  \right) \right)\\
b := 2 \ch_2 \left( (p_1)_\ast \left( \cE \otimes p_2^\ast \Im (\omega_Z) \right) \right)
\end{array}
\end{equation}
Where $\Im : \Knum(X)_\bC \to \Knum(X)_\bR$ denotes the imaginary part. Note that the cohomology classes \eqref{eqn:classes_Bridgeland} are real, but if $Z(\cA) \subset \bQ+i\bQ$, then we may assume $\omega_Z \in K^0(\Perf(X)) \otimes (\bQ+i\bQ)$, and the cohomology classes $l$ and $b$ are rational.

\begin{notn}
We define the \emph{rank-degree-weight} sequence associated to a $\bZ$-weighted filtration $\cdots \supset E_w \supset E_{w+1} \supset \cdots$ as the order sequence of triples
\begin{equation} \label{eqn:rank_deg_weight_seq}
\alpha = \left\{(r_j,d_j,w_j) = (\Im Z(\gr_j E_\bullet), - \Re Z(\gr_j E_\bullet), w_j) \right\}_{j = 1,\ldots,p}
\end{equation}
where $w_j \in \bZ$ are the finite set of integers $w$ for which $\gr_w(E_\bullet) \neq 0$. Because $Z(E)$ only depends on the numerical class of $E$, $\alpha$ is locally constant on $\Filt(\cM_v)$. We sometimes use the notation $E_j := E_{w_j}$, and regard the data of a $\bZ$-weighted filtration as a finite filtration $E = E_0 \supsetneq E_1 \supsetneq \cdots \supsetneq E_p \supsetneq E_{p+1} = 0$ along with the strictly increasing weights $w_0 < w_1 < \cdots < w_p$ in $\bZ$. It is natural to allow the weights $w_i$ to lie in $\bQ$ or even $\bR$, and we refer to such data as a \emph{descending weighted filtration}.
\end{notn}

\begin{lem} \label{lem:compute_classes_Bridgeland}
Let $k'/k$ be a field extension and let $f: \Theta_{k'} \to \cM_v$ correspond to a descending weighted filtration $E_\bullet$. Let $\{(r_j,d_j,w_j)\}$ be the rank-degree-weight sequence associated to $E_\bullet$ as in \eqref{eqn:rank_deg_weight_seq} and let $Z(v) = -D + i R$. Then we have
\begin{equation} \label{eqn:compute_classes_Bridgeland}
\renewcommand\arraystretch{2}
\begin{array}{ccc}
\frac{1}{q} f^\ast l = \sum \limits_{j=1}^p w_j (R d_j - D r_j), & \text{and} & \frac{1}{q^2} f^\ast b = \sum \limits_{j=1}^p w_j^2 r_j
\end{array}
\end{equation}
\end{lem}
\begin{proof}
If $\pi : \Theta_{k'} \times X \to \Theta$ and $\cE_f$ is the object classified by $f$, then in $K$-theory, we have $\cE_f \simeq \sum_j u^{-w_j} [\gr_j E_\bullet]$ under the decomposition $K_0(X \times \Theta_{k'}) \simeq K_0(X_{k'}) \otimes \bZ[u^\pm]$. We can thus compute
\begin{align*}
f^\ast l &= |Z(v)|^2 \op{ch}_1 \left(\pi_\ast (\cE_f \otimes p_2^\ast \Im(\frac{-\omega_Z}{Z(v)})) \right) \\
&= |Z(v)|^2 \op{ch}_1 \left(  \sum u^{-w_j} \Im \left( \frac{-\chi(X,\omega_Z \otimes \gr_j E_\bullet )}{Z(v)} \right) \right) \\
&= \sum w_j q \Im \left( \chi(X,\omega_Z \otimes \gr_j E_\bullet ) \cdot \overline{Z(v)} \right)
\end{align*}
By the defining property of $\omega_Z$, we have $\chi(\omega_Z \otimes \gr_j E_\bullet) = (-d_j + i r_j)$, and the claim follows. The computation for $f^\ast b$ is almost identical, so we omit it.
\end{proof}

\begin{defn} \label{defn:torsion_theory}
Given a slope-stability condition $(Z,\cA,\D^b(X)^\heartsuit)$, we refer to $E \in \cA$ as \emph{torsion} if $Z(E) \in \bR_{\leq 0}$, and we call an object \emph{torsion-free} if it has no torsion subobjects (in particular a torsion-free object has $\Im Z(E)>0$). We define $\cT \subset \cA$ (resp. $\cF \subset \cA$) to be the full subcategory consisting of torsion objects (resp. torsion-free objects). Note that torsion objects are automatically semistable.
\end{defn}

\begin{ex}
Continuing \Cref{ex:moduli_coherent_sheaves}, let $\cA$ be the category of sheaves whose support has dimension $\leq d$. The torsion subcategory $\cT$ consists of sheaves whose support has dimension $\leq d-1$, and the $\cF$ consists of sheaves that are pure of dimension $d$.
\end{ex}

We shall consider the numerical invariant $\mu : \cU \subset \iComp(\cM_v) \to \bR$ associated to the cohomology classes \eqref{eqn:classes_Bridgeland}. Observe from \eqref{eqn:compute_classes_Bridgeland} that for any filtration $f: \Theta_{k'} \to \cM_v$, $f^\ast(b) = 0$ if and only if $r_i=0$ for any $i$ for which $w_i \neq 0$. It follows that $b$ is not positive definite on $\cM_v$, and $\mu$ is not strictly quasi-concave, if $\cA$ has torsion objects. The rational points of the subset $\cU_{E} \subset \iDeg(\cM_v,E)$ on which $\mu$ is defined corresponds to the set of all $\bZ$-weighted filtrations that are \emph{not} of the form
\[
0 = 0 = \cdots = E_1 \subset E_0 \subset E_{-1} \subset E_{-2} \subset \cdots \subset E
\]
with $\gr_{i}(E_\bullet)$ torsion for all $i<0$. The main result of this subsection states that the HN problem for this numerical invariant on $\cM_v$ is equivalent to the Harder-Narasimhan property for $(Z,\cA,\D^b(X)^\heartsuit)$.

\begin{thm} \label{thm:HN_Bridgeland}
Let $(Z,\cA,\D^b(X)^\heartsuit)$ be a slope-semistability condition on a projective $k$-scheme $X$. Let $\mu$ be the numerical invariant on $\cM_v$ associated to the cohomology classes \eqref{eqn:classes_Bridgeland} as in \Cref{defn:induced_invariant}, and let $M^\mu$ the corresponding stability function of \Cref{defn:numerical_invariant}. Then $E \in \cF$ is slope semistable if and only if $M^\mu([E]) \leq 0$, and the following are equivalent:
\begin{enumerate}
\item $(Z,\cA,\D^b(X)^\heartsuit)$ has the Harder-Narasimhan property; and
\item Every object in $\cA$ has a maximal torsion subobject, and for every unstable $E \in \cF$, the function $\mu : \cU_E \subset \iDeg(\cM_v,E) \to \bR$ obtains a maximum.
\end{enumerate}
If $\cA$ is Noetherian, this is equivalent to the condition:
\begin{enumerate}
\setcounter{enumi}{2}
\item For any $E \in \cA$, the set $\{\phi(F) | F \subset E\}$ has a maximal element.
\end{enumerate}
Furthermore, if (2) holds then the maximum of $\mu$ occurs at the point in $\iDeg(\cM_v,E)$ corresponding to the Harder-Narasimhan filtration $E = E_1 \supset E_2 \supset \cdots \supset E_n$ with weights proportional to the slopes (defined below) of the associated graded pieces $E_i/E_{i+1}$, and this is the unique maximizer corresponding to a filtration whose associated graded pieces are torsion-free. If $Z$ is such that $Z(E)=0 \Rightarrow E=0$ for all $E \in \cA$, then this is the unique maximizer of $\mu$.
\end{thm}

We shall prove \Cref{thm:HN_Bridgeland} at the end of this subsection, after collecting several intermediate results. The proof bears a strong formal resemblance to the analysis of semistability for filtered isocrystals in $p$-adic Hodge theory as presented in \cite{dat2010period}.

As in the context of Bridgeland stability conditions, we let $\cP(\phi) \subset \cA$ be the full subcategory of semistable objects of phase $\phi$. In addition, for any subset $I \subset (0,1]$, we let $\cP(I)$ be the full subcategory generated under extensions by semistable objects whose phase lies in $I$.

\begin{lem} \label{lem:truncations}
For any $\epsilon \in (0,1]$, $\op{Hom}(E,F)= 0$ for all $E \in \cP([\epsilon,1])$ and $F \in \cP((0,\epsilon))$.
\end{lem}
\begin{proof}
It suffices to show that $\op{Hom}(E_1,E_2) = 0$ for semistable objects with $\phi(E_1) > \phi(E_2)$. Letting $f : E_1 \to E_2$, if $\op{im}(f)$ is nonzero, then the semistability of $E_1$ and $E_2$ implies that $\phi(E_1) \leq \phi(\op{im}(f)) \leq \phi(E_2)$. Hence if this inequality is violated we must have $f = 0$. 
\end{proof}

If $(Z,\cA,\D^b(X)^\heartsuit)$ has the Harder-Narasimhan property, we let $E^{\phi \geq \epsilon} \subset E$ denote the largest sub-object in the Harder-Narasimhan filtration whose associated graded pieces have phase $\geq \epsilon$.

\begin{cor} \label{cor:truncations}
If $(Z,\cA,\D^b(X)^\heartsuit)$ has the Harder-Narasimhan property, then the Harder-Narasimhan filtration is functorial in the sense that for any $\epsilon \in (0,1]$, any homomorphism $E \to F$ maps $E^{\phi \geq \epsilon}$ to $F^{\phi\geq \epsilon}$.
\end{cor}
\begin{proof}
The map $E^{\phi \geq \epsilon} \to F^{\phi < \epsilon} := F / F^{\phi \geq \epsilon}$ must vanish by the previous lemma.
\end{proof}

For any ordered sequence of complex numbers $(z_1,\ldots,z_p)$ with $z_j \in \mathbb{H} \cup \bR_{\leq 0}$, we associate the convex polyhedron
\[
\op{Pol}(\{z_j\}) = \left\{ \left. \sum \lambda_j z_j + c \right| 0 \leq \lambda_1 \leq \cdots \leq \lambda_p \leq 1 \text{ and } c \in \bR_{\geq 0} \right\} \subset {\mathbb H} \cup \bR_{\geq 0}
\]
This is the convex hull of the points $z_p,z_p+z_{p-1},\ldots, z_p+\cdots+z_1$ plus an arbitrary shift by $c \in \bR_{\geq 0}$. For any filtration $E_\bullet$ in $\cA$, we let $\op{Pol}(E_\bullet) = \op{Pol}(\{Z(\gr_j E_\bullet)\})$, and for any object $E \in \cA$ we let $\op{Pol}^{\rm HN}(E) := \op{Pol}(E_\bullet)$, where $E_\bullet$ is the HN filtration of $E$.

\begin{lem} \label{lem:containment}
If $(Z,\cA,\D^b(X)^\heartsuit)$ has the Harder-Narasimhan property, then $\op{Pol}(E_\bullet) \subset \op{Pol}^{\rm HN}(E)$ for any weighted descending filtration $E_\bullet$ of $E$.
\end{lem}

\begin{proof}
It suffices to show that for any $E$ and any subobject $F \subset E$, one has $Z(F) \in \op{Pol}^{\rm HN}(E)$, and we prove this by induction on the length of the HN filtration of $E$. The base case is when $E$ is semistable, in which case the statement of the lemma is precisely the definition of semistability, combined with the fact that $\Im Z(F) \leq \Im Z(E)$ for any subobject.

Let $E_p \subset E$ be the first subobject in the HN filtration, so by hypothesis $E_p$ is semistable, say with phase $\phi_p$. \Cref{cor:truncations} guarantees that any subobject of $E$ must have phase $\leq \phi_p$. Let $E^\prime = E / E_p$ and let $F^\prime$ be the image of $F \to E^\prime$. Then can pullback the defining sequence for $E^\prime$ to $F^\prime$:
\[
\xymatrix@R=10pt{ 0 \ar[r] & E_p \ar@{=}[d] \ar[r] & E_p + F \ar@{^{(}->}[d] \ar[r] & F^\prime \ar@{^{(}->}[d] \ar[r] & 0 \\ 0 \ar[r] & E_p \ar[r] & E \ar[r] & E^\prime \ar[r] & 0}
\]
The length of the HN filtration of $E^\prime$ is one less than $E$, and $\op{Pol}^{\rm HN}(E)$ decomposes as a union of $\op{Pol}^{\rm HN}(E_p)$ and a shift of $\op{Pol}^{\rm HN}(E^\prime)$ by $Z(E_p)$. Thus by the inductive hypothesis $Z(E_p + F) = Z(F^\prime) + Z(E_p)$ lies in $Z(E_p) + \op{Pol}^{\rm HN}(E^\prime)$, which we illustrate graphically:

\medskip
\begin{center}
\begin{tikzpicture}
\draw (-2,1) -- (-2.5,2) -- (-2,3) -- (6,3);
\draw (6,0) -- (0,0) -- (-2,1) -- (6,1);

\draw[fill] (-2,1) circle [radius = .1];
\node[below left] at (-2,1) {$Z(E_p)$};

\node at (.5,.5) {$\op{Pol}^{\rm HN}(E_p)$};
\node at (-.5,2) {$\op{Pol}^{\rm HN}(E^\prime) + Z(E_p)$};

\draw[fill=gray] (6,.5) -- (4,.5) -- (2,1.5) -- (6,1.5);
\draw[fill] (2,1.5) circle [radius = .1];
\node[above right] at (2,1.5) {$Z(E_p + F)$};
\end{tikzpicture}
\end{center}
\medskip

We also have the short exact sequence $0 \to E_p \cap F \to E_p \oplus F \to E_p + F \to 0$, which implies that
\begin{align*}
Z(F) &= Z(E_p + F) + Z(E_p \cap F) - Z(E_p) \\
&= Z(E_p + F) - Z(E_p / E_p \cap F)
\end{align*}
We have $\Im Z(E_p / E_p \cap F) \leq \Im Z(E_p)$, and the phase of $Z(E_p / E_p \cap F)$ is $\geq \phi_p$ because $E_p$ is semistable. It follows that $Z(F)$ must lie in the shaded region above, which is a translate of $\op{Pol}^{\rm HN}(E_p)$, and hence $Z(F) \in \op{Pol}^{\rm HN}(E)$. Even though we have drawn the diagram assuming that $\Im Z(E_p) > 0$, the argument works without modification in the case where $E_p$ is torsion as well.
\end{proof}

Let $\alpha = \{(r_j,d_j,w_j)\}$ be a rank-degree-weight sequence. For $j=1,\ldots, p$, we define the $j^{th}$ phase $\phi_j \in (0,1]$ by the property that $i r_j - d_j \in \bR_{>0} e^{i \phi_j}$. Note that by \Cref{lem:compute_classes_Bridgeland}, the numerical invariant can be expressed formally in terms of the rank-degree-weight sequence $\mu = \mu(\alpha)$. Furthermore it extends continuously to rank-degree-weight sequences where the $w_j$ are arbitrary real numbers rather than integers.

\begin{lem}[Insertion and deletion] \label{lem:convex_rdw}
Let $\alpha$ be a rank-degree-weight sequence such that $r_k+r_{k+1}>0$ for some $k$, and let $\alpha^\prime$ be the rank-degree-weight sequence obtained from $\alpha$ by discarding the $(k+1)^{st}$ element and relabeling $r^\prime_k = r_k + r_{k+1}$ and $d^\prime_k = d_k + d_{k+1}$, and
$$w_j^\prime := \left\{ \begin{array}{ll} w_j, & \text{if } j<k \\ \frac{w_j r_j + w_{j+1} r_{j+1}}{r_j+r_{j+1}}, & \text{if } j=k \\ w_{j+1}, & \text{if } j>k \end{array} \right. .$$
Then we have:
\begin{itemize}
\item If $\mu(\alpha) \geq 0$ and $\phi_{k} \geq \phi_{k+1}$, then $\mu(\alpha^\prime) \geq \mu(\alpha)$ with equality if and only if $\phi_k = \phi_{k+1}$ and $\mu(\alpha) = 0$.
\item If $\mu(\alpha^\prime) \geq 0$ and $\phi_k < \phi_{k+1}$, then $\mu(\alpha) > \mu(\alpha^\prime)$.
\end{itemize}
\end{lem}
\begin{proof}

We denote $\mu = L / \sqrt{B}$. Substituting $\alpha^\prime$ for $\alpha$, the numerator and denominator change by
\begin{align*}
\Delta L &= w_k^\prime (d_k + d_{k+1}) - w_k d_k - w_{k+1} d_{k+1} \\
&= \frac{w_{k+1}-w_k}{r_k + r_{k+1}} (d_k r_{k+1} - r_k d_{k+1}) \\
\Delta B &= (w^\prime_k)^2 (r_{k+1}+r_k) - w_k^2 r_k - w_{k+1}^2 r_{k+1} \\
&= - \frac{r_k r_{k+1}}{r_k + r_{k+1}} (w_k - w_{k+1})^2
\end{align*}
Note that $\Delta B \leq 0$. Also the sign of $d_k r_{k+1} - d_{k+1} r_k$, and hence the sign of $\Delta L$, is the same as the sign of $\phi_k - \phi_{k+1}$, and $\Delta L \neq 0$ if $\phi_k \neq \phi_{k+1}$. The claim follows from these observations.
\end{proof}

In terms of descending weighted filtrations, the modification $\alpha \mapsto \alpha^\prime$ in \Cref{lem:convex_rdw} corresponds to deleting the $(k+1)^{st}$ subobject from the filtration $E = E_1 \supset \cdots \supset E_p \supset E_{p+1} = 0$ and adjusting the weights appropriately. We also observe that if $\phi_k \geq \phi_{k+1}$ then this does not change $\op{Pol}(E_\bullet)$. Hence we have the following:

\begin{cor} \label{cor:deletion}
For any descending weighted filtration $(E_\bullet,w_\bullet)$, there is a sequence of deletions resulting in a descending weighted filtration $(E_\bullet^\prime,w'_\bullet)$ that is convex, in the sense that $\phi_1^\prime < \cdots < \phi_{p^\prime}^\prime$, such that $\op{Pol}(E_\bullet) = \op{Pol}(E^\prime_\bullet)$ and $\mu(E_\bullet^\prime) \geq \mu(E_\bullet)$, with strict inequality if $E_\bullet$ was not convex to begin with.
\end{cor}

Given a sequence $(z_1,\ldots,z_p)$ in $\mathbb{H} \cup \bR_{\leq 0}$, we let $R = \sum r_j$ and define a continuous piecewise linear function $h_{\{z_j\}} : [0,R] \to \bR$,
$$h_{\{z_j\}}(r) := \sup \left\{ x \left| i r - x \in \op{Pol}(\{z_j\}) \right. \right\}.$$
Note that this only depends on $\op{Pol}(\{z_j\})$.
\begin{lem} \label{lem:explicit_maximum}
Let $i r_j - d_j \in {\mathbb{H}}$ be sequence of points such that $\phi_1 < \cdots < \phi_p < 1$. Let $\nu_j := d_j / r_j$, $D := \sum d_j$, $R := \sum r_j$, and $\nu := D/R$, and let $h(x) := h_{\{i r_j - d_j\}}(x)$. Then $\mu$ is maximized by assigning the weights $w_j \propto \nu_j - \nu$. The maximum is
\begin{equation} \label{eqn:explicit_maximum}
\frac{\mu}{R} = \sqrt{ (\sum \nu_j^2 r_j) - \nu^2 R} = \sqrt{\int_0^R (h^\prime(x))^2 dx - \nu^2 R}.
\end{equation}
\end{lem}

\begin{proof}
We can think of the numbers $\ev_1,\cdots,r_p$ as defining an inner product $\vec{a} \cdot \vec{b} = \sum a_j b_j r_j$. Then given a choice of weights $\vec{w} = (w_1,\cdots,w_p)$, the numerical invariant can be expressed as
$$\mu = \frac{R}{|\vec{w}|} \vec{w} \cdot (\vec{\nu} - \nu \vec{1})$$
where $\vec{\nu} = (\nu_1, \cdots, \nu_p)$ and $\vec{1} = (1,\ldots,1)$. From linear algebra we know that this quantity is maximized when $\vec{w} \propto \vec{\nu} - \nu \vec{1}$, and the maximum value is $ R |\vec{\nu} - \nu \vec{1}|$. In the case when $\nu_1 < \cdots < \nu_p$ the assignment $\vec{w} \propto \vec{\nu} - \nu \vec{1}$ satisfies the constraints $w_1 < w_2 < \cdots < w_p$. The integral formula simply expresses $\nu_j^2 r_j$ as the integral of $h^\prime(x)$ on an interval of length $r_j$ along which it is constant of value $\nu_j$.
\end{proof}

\begin{lem} \label{lem:monotonicity}
The expression \eqref{eqn:explicit_maximum} is strictly monotone increasing with respect to inclusion of polyhedra.
\end{lem}
\begin{proof}
Let $h_1,h_2 : [0,R] \to \bR$ be two continuous piecewise linear functions with $h_i^\prime(x)$ decreasing and with $h_1(x) \leq h_2(x)$ with equality at the endpoints of the interval. We must show that $\int_0^R (h_1^\prime(x))^2 < \int_0^R (h_2^\prime(x))^2$. First by suitable approximation with respect to a Sobolev norm it suffices to prove this when $h_i$ are smooth functions with $h^{\prime \prime} < 0$. Then we can use integration by parts
\[
\int_0^R (h_2^\prime)^2 - (h_1^\prime)^2 dx = (h_2^\prime + h_1^\prime)(h_2 - h_1)|_0^R - \int_0^R (h_2 - h_1)(h_2^{\prime \prime} + h_1^{\prime \prime})dx
\]
The first term vanishes because $h_1 = h_2$ at the endpoints, and the second term is strictly positive unless $h_1 = h_2$.
\end{proof}

\begin{proof} [Proof of Theorem \Cref{thm:HN_Bridgeland}]

First we show that $E \in \cF$ is unstable if and only if $M^\mu([E]) > 0$. If $F \subset E$ with $\phi(F) > \phi(E)$, then consider the two step filtration $\gr_2 E_\bullet = F$ and $\gr_1 E_\bullet = E / F$ and $w_2 > w_1$ arbitrary. Because $E$ is torsion-free we know $Z(F) \neq 0$, so by \Cref{lem:compute_classes_Bridgeland}
\begin{align*}
\frac{1}{q} f^\ast l &= w_1 \Im \left((Z(E) - Z(F)) \overline{Z(E)} \right) + w_2 \Im \left(Z(F) \overline{Z(E)} \right) \\
&= (w_2 - w_1) \Im \left(Z(F) \overline{Z(E)} \right) > 0,
\end{align*}
and hence $M^\mu(E)>0$. Conversely, if $f : \Theta \to \cM_v$ corresponds to a weighted filtration $E_\bullet$ of $E$ such that $\mu(f) > 0$, then \Cref{cor:deletion} provides a convex filtration $E'_\bullet$ with $\mu(E^\prime_\bullet) \geq \mu(E_\bullet) > 0$. In particular $E'_\bullet$ is a non-trivial filtration, and so the first subobject $E^\prime_{p^\prime} \subset E$ destabilizes $E$.

\medskip
{\noindent \textit {Proof that $(1) \Rightarrow (2)$:}}
\medskip

If $(Z,\cA,\D^b(X)^\heartsuit)$ has the Harder-Narasimhan property, then $E^{\phi \geq 1} \subset E$ is the maximal torsion subsheaf of any $E \in \cA$. Now consider $E \in \cF$. By \Cref{cor:deletion}, it thus suffices to maximize $\mu$ over the set of convex filtrations, and because $E$ is torsion-free any convex filtration will have $\phi<1$ for all graded pieces. For a fixed convex filtration, \Cref{lem:explicit_maximum} computes the weights that maximize $\mu$ for a fixed filtration, and this expression is strictly monotone increasing with respect to inclusion of polyhedra by \Cref{lem:monotonicity}. By \Cref{lem:containment} we have $\op{Pol}(E_\bullet) \subset \op{Pol}^{\rm HN}(E)$ for any filtration, so if $\mu(E_\bullet) >0$ for some weighted descending filtration of $E$, then $\mu(E_\bullet)$ is strictly less than the value of $\mu$ obtained by assigning weights to the HN filtration as in \Cref{lem:explicit_maximum}.

\medskip
{\noindent \textit {Proof that $(2) \Rightarrow (1)$:}}
\medskip

First if any $E$ admits a maximal torsion subobject $T \subset E$, then $T$ is semistable, and it suffices to show that the torsion-free quotient $E / T$ admits an HN filtration with phases $<1$. Assume $E\in \cF$ is unstable and choose a point in $\iDeg(\cM_v,E)$ that maximizes $\mu$, corresponding to a descending filtration $E_1 \supset \cdots \supset E_n$ with real weights $w_1\leq \cdots \leq w_n$ (see \Cref{rem:real_weights}). \Cref{cor:deletion} implies that $\phi_1 < \cdots < \phi_p < 1$, where $\phi_j := \phi(E_j /E_{j+1})$.

Consider a subobject $F \subset E_j / E_{j+1}$. We refine the filtration $E_\bullet$ to a new filtration
\[
E'_\bullet = (E_n \subset \cdots \subset E_{j+1} \subset \tilde{F} \subset E_j \subset \cdots \subset E_1),
\]
where $\tilde{F}$ is the preimage of $F$ under the map $E_j \to E_j/E_{j+1}$. One automatically has $\op{Pol}(E_\bullet) \subset \op{Pol}(E'_\bullet)$, and the inclusion is strict if $Z(F) \neq 0$ and $\phi(F) > \phi(E_j/E_{j+1})$. Because $\mu(E_\bullet)$ is maximal, \Cref{cor:deletion} and \Cref{lem:monotonicity} imply that $\op{Pol}(E_\bullet) = \op{Pol}(E'_\bullet)$, so it follows that either $Z(F) = 0$ or $\phi(F) \leq \phi_j$.

This shows that each object $E_j/E_{j+1}$ is either torsion-free, in which case it is semistable, or the maximal torsion subobject $F \subset E_j / E_{j+1}$ has $Z(F)=0$. By a simple inductive procedure, re-defining $E_{j+1}$ to be the preimage of the maximal torsion subobject in $E_j / E_{j+1}$, we can construct a new filtration of the same length $E'_1 \supset \cdots \supset E'_n$ with $Z(E'_j) = Z(E_j)$ for all $j$ and $E'_{j} / E'_{j+1}$ torsion-free. It follows that $E'_\bullet$ maximizes $\mu$, so the analysis of the previous paragraph implies that $E'_{j}/E'_{j+1}$ is semistable of increasing phase, hence $E'_\bullet$ is a Harder-Narasimhan filtration for $E$. Note that this analysis shows that \emph{any} maximizer of $\mu$ for which $E_j/E_{j+1}$ is torsion-free is actually the Harder-Narasimhan filtration with weights given by \Cref{lem:explicit_maximum}.

\medskip
{\noindent \textit {Proof that $(3) \Leftrightarrow (1)$ when $\cA$ is noetherian:}}
\medskip

The proof is well-known in closely related contexts, see for instance \cite{bridgeland2007stability}*{Prop.~2.4}, but for completeness we summarize it here: For any $E$, let $\phi$ be the maximal phase of a subobject. Any subobject $F \subset E$ of phase $\phi$ must be semistable, and for any two subobjects $F,F' \subset E$ of phase $\phi$ the image $F+F' = \op{im}(F\oplus F' \to E)$ has phase $\phi$ as well, or else the kernel would be a subobject of $F \oplus F'$ with phase $>\phi$. Because $\cA$ is noetherian it follows that there is a maximal subobject of phase $\phi$. Using this fact, one can build an ascending sequence $0 \subsetneq E_0 \subsetneq E_1 \subset \cdots \subset E$ such that $\phi(\gr_i(E_\bullet))$ is increasing in $i$, and $E_i$ is the preimage of the maximal subobject of minimal phase in $E/E_{i-1}$. Again because $\cA$ is noetherian this sequence must stabilize, and it is a HN filtration by construction.

\end{proof}

\subsection{\texorpdfstring{$\Theta$}{Theta}-stratification of the moduli of torsion-free objects}
\label{sect:theta_strat_Bridgeland}

We now specialize to the setting of a (numerical) Bridgeland stability condition $(Z,\cA)$ on $\D^b(X)$ for a projective $k$-scheme $X$. This is defined to be a slope stability condition (\Cref{def:slope_stability}) that satisfies the Harder-Narasimhan property (\Cref{defn:slope_semistability}) and in which $\cA = \D^b(X)^\heartsuit$ and $Z(E) \neq 0$ for any $E \in \cA$.\footnote{In some sources, stability conditions are required to satisfy an additional "support property," and the notion we use here is referred to as a prestability condition.}

If $\cA$ is noetherian, then for any finite extension $k'/k$, one can define a Bridgeland stability condition $(Z_{k'},\cA_{k'})$ on $\D^b(X_{k'})$ as follows: $\cA_{k'} \subset \D^b(X_{k'})$ is the heart of the induced $t$-structure (\Cref{D:induced_t_structure}), and $Z_{k'}(E) = Z(p_\ast(E))$, where $p_\ast : \D^b(X_{k'}) \to \D^b(X)$ is the pushforward. The Harder-Narasimhan property follows from \Cref{thm:HN_Bridgeland}(3), because $\cA_{k'}$ is noetherian by \Cref{prop:t_structure_coherent}.

We use this induced stability condition on $\D^b(X_{k'})$ to define semistable objects in $\cA_{k'}$ (\Cref{defn:slope_semistability}), as well as the category of torsion free objects $\cF_{k'} \subset \cA_{k'}$ (\Cref{defn:torsion_theory}). Using the fact that for any object $E \in \cA_{k'}$ the canonical map $k' \otimes_k p_\ast(E) \to E$ is surjective, one can show that $E \in \cF_{k'}$ if and only if $p_\ast(E) \in \cF$ as well.

\begin{defn}[Semistable and torsion free moduli functors] \label{defn:moduli_torsion_free}
Let $(Z,\cA)$ be a numerical Bridgeland stability condition on $\D^b(X)$ for which $\cA$ is noetherian and satisfies the Generic Flatness condition. For any finite type $k$-scheme $T$, the groupoid of \emph{flat families of torsion-free objects} of class $v \in \Lambda \subset \Knum(X)$ in $\cA$ is
\[
\cM^{\cF}_v(T) := \left\{ E \in \cM_v(T) | E_p \in \cF_{\kappa(p)} \subset \cA_{\kappa(p)}, \forall \text{ closed }p \in T \right\}.
\]
The groupoid of \emph{flat families of semistable objects} of class $v$ is
\[
\cM^{\rm ss}_v(T) := \left\{ E \in \cM_v(T) | E_p  \in \cA_{\kappa(p)} \text{ is semistable}, \forall \text{ closed }p \in T \right\}.
\]
Both conditions are local for the smooth topology, so these define substacks of the restriction of $\cM_v$ to the \'etale site of finite type $k$-schemes.
\end{defn}

\begin{defn}[Boundedness of Quotients] \label{defn:boundedness_of_quotients}
We say that a Bridgeland stability condition $(Z,\cA)$ on $\D^b(X)$ satisfies the \emph{Boundedness of Quotients} condition if for any $E \in \cA$ and any $\phi \in (0,1)$, the set of points of $\cM$ that parameterize a torsion-free object $E' \in \cF_{k'}$ over a finite extension $k'/k$ of phase $\leq \phi$ that admits a surjection $E \otimes_k k' \to E'$ is bounded.
\end{defn}

\begin{thm} \label{thm:theta_stratification_torsion_free}
Let $X$ be a projective scheme over a field $k$, and let $(Z,\cA)$ be a Bridgeland stability condition on $\D^b(X)$ such that $Z : \Knum(X) \to \bC$ factors through a finite rank quotient of $\Knum(X)$.  Assume that:
\begin{enumerate}
\item $Z(\cA) \subset \bQ + i \bQ$;
\item $\cA$ satisfies the Generic Flatness condition; and
\item $(Z,\cA)$ satisfies the Boundedness of Quotients condition.
\end{enumerate}
Then for any $v \in \Knum(X)$ with $\Im(Z(v)) > 0$, the moduli functors in \Cref{defn:moduli_torsion_free} extend uniquely to open substacks $\cM^{\rm ss}_v \subset \cM_v^\cF \subset \cM_v$ on the \'etale site of all $k$-schemes.

Furthermore, the numerical invariant $\mu$ associated to the classes of \eqref{eqn:classes_Bridgeland} defines a $\Theta$-stratification of the algebraic stack $\cM^\cF_v$, where: i) $\cM^{\rm ss}_v \subset \cM^\cF_v$ is the $\mu$-semistable locus; and ii) the $\mu$-HN filtration of any unstable point in the sense of \Cref{defn:HN_filtration} corresponds under \Cref{lem:filtrations_derived} to a canonical weighting of the Harder-Narasimhan filtration in the sense of \Cref{defn:slope_semistability}.
\end{thm}

\begin{proof}
The hypothesis (1) and the assumption that $Z$ factors through a finite rank quotient of $\Knum(X)$ implies that $\cA$ is noetherian \cite{abramovich2006sheaves}*{Prop.~5.0.1}. We first verify that $\cM_v^\cF$ extends uniquely to an open subfunctor of $\cM_v$. The analogous claim for $\cM^{\rm ss}_v$ follows from the claim about the $\Theta$-stratification of $\cM^\cF_v$, which we prove below.

We first show that the Boundedness of Quotients condition holds for families, i.e., for a finite type $k$-scheme $T$, a flat family $E \in \cM_v(T)$, and a $\phi \in (0,1)$, the set of finite type points of $\cM$ parameterizing objects $E' \in \cF_{\kappa(t)}$ of phase $\leq \phi$ that arise as quotients $E_t \twoheadrightarrow E'$ for some finite type point $t \in T$ is bounded. It suffices to prove this when $T=\Spec(R)$ is affine. In this case there is a surjection $R \otimes_k F \twoheadrightarrow E$ for some $F \in \cA$, because $E \in \D_{qc}(X_R)^{\heartsuit}$ is compact and $\D_{qc}(X_R)^\heartsuit$ is compactly generated by objects of the form $R \otimes_k F$ for $F \in \cA$. This reduces the claim to the boundedness of quotients of $F\otimes_k k'$ for the single object $F$, which is provided by \Cref{defn:boundedness_of_quotients}.

Now consider the open and closed substack $\cY \subset \Filt(\cM_v)$ classifying weighted descending filtrations of the form $\cdots = 0 \subset F_1 \subset F_0 = \cdots$ where the class of $F_0$ is $v$ and the class of $F_1$ is torsion and nonzero. Observe that the value of $Z(F_1) \in \bR_{\leq 0}$ is locally constant on $\cY$ and its image is discrete in $\bR_{\leq 0}$ by hypothesis (1). We regard $\cY$ as a stack over $\cM_v$ via the morphism $\ev_1 : \cY \to \cM_v$.

For any finite type $k$-scheme $T$ and morphism $\xi : T \to \cM_v$, the family version of the Boundedness of Quotients condition implies that there is a minimal value of $Z(F_1)$ occurring among all filtrations in $T \times_{\cM_v} \cY$ for which $F_1/F_0$ is torsion-free. Let $Y' \subset T \times_{\cM_v} \cY$ be the connected component on which $Z(F_1)$ is minimal. Then any finite type point in $Y'$ parameterizes a filtration such that $F_1/F_0$ is torsion free (or else $Z(F_1)$ would not be minimal), so the family Boundedness of Quotients condition implies that $Y'$ is quasi-compact. Also, $Y' \to T$ satisfies the valuative criterion for properness with respect to essentially finite type DVR's by \Cref{prop:theta_reductive_S_complete}. It follows that the image of $\ev_1 : Y' \to T$ is closed.

We now replace $T$ with $T \setminus \ev_1(Y')$ and iterate the construction of the previous paragraph. Because the values of $Z(F_1)$ are discrete in $\bR_{\leq 0}$, after finitely many steps this results in an open subscheme $U \subset T$ whose closed points are precisely those $t \in T$ such that $\xi(t)$ classifies a torsion free object in $\cA_{\kappa(t)}$. Because $\cM_v$ is algebraic and locally of finite type over $k$ (\Cref{prop:moduli_heart}), this implies that $\cM^\cF_v$ extends uniquely to an open substack of $\cM_v$.

\medskip
\noindent{\textit{Proof that $\mu$ defines a $\Theta$-stratification of $\cM_v^\cF$}:}
\medskip

Note that by \Cref{lem:compute_classes_Bridgeland} the class $b \in H^4(\cM^\cF_v)$ in \eqref{eqn:classes_Bridgeland} is positive definite, so it defines a norm on graded points (\Cref{ex:norm_from_H4}). \Cref{defn:induced_invariant} associates an $\bR$-valued numerical invariant $\mu$ on $\cM^\cF_v$ to this norm and the class $\ell$ in \eqref{eqn:classes_Bridgeland}, where $\cU = \iComp(\cM^\cF_v)$. We note that the existence and uniqueness of $\mu$-HN filtrations for finite type points of $\cM^\cF_v$, as well as the conclusions (i) and (ii), follow from \Cref{thm:HN_Bridgeland}.

We will show that $\mu$ defines a weak $\Theta$-stratification by verifying the hypotheses of \Cref{thm:main_improved}. We have already shown above that $\cM^\cF_v$ is algebraic and locally finite type over $k$ with affine diagonal. $\mu$ is standard by \Cref{lem:induced_invariant} and satisfies \ref{princ:R} by \Cref{lem:induced_invariants_R}.

Next we verify the HN-boundedness property \ref{princ:B2}: Let $\xi : T \to \cM^\cF_v$ be a family over a finite type $k$-scheme $T$. \Cref{thm:HN_Bridgeland} implies that the HN filtration of any finite type point of $\cM^\cF_v$ has a convex rank-degree sequence. It therefore suffices to show that the set of points of the form $\gr(f) \in \cM^\cF_v$, where $f$ is a finite type point of $\Flag(\xi)$ classifying a filtration with convex rank-degree sequence, is bounded. The family version of the Boundedness of Quotients property implies that there is a minimal slope appearing among all quotients of the fiber $\xi(t)$ for any closed point $t \in T$. This implies that the image under $Z$ of any subobject appearing in a convex filtration of $\xi(t)$ for any closed point $t \in T$ is contained in some bounded region of $\bC$. Because $Z(\cA) \subset \bC$ is discrete, there is a maximal length for all convex filtrations of all fibers $\xi(t)$ over closed points $t \in T$. From here, one can show using induction that the family version of the Boundedness of Quotients property implies that the set of finite type points of $\Flag(\xi)$ classifying convex filtrations is bounded, and hence so is its image under $\gr : \Flag(\xi) \to \cM^\cF_v$.

All that remains is to verify the HN-specialization property \ref{princ:S}: Consider a family $E \in \cM^\cF_v(R)$ over an essentially finite type DVR $R$ with fraction field $K$ and residue field $\kappa$ and let $(E_K)_\bullet$ be an HN filtration in $\cM^\cF_v$ of $E_K$. Note that $(E_K)_\bullet$ is also an HN filtration in $\cM_v$ by \Cref{thm:HN_Bridgeland}. The stack $\cM_v$ is $\Theta$-reductive by \Cref{prop:theta_reductive_S_complete}, so the filtration extends uniquely to a filtration in $\cM_v$ over $R$. The filtration over the special fiber lies in $\cU_{E_\kappa} \subset \iDeg(\cM_v,E_\kappa)$ and has the same value of $\mu$ as the filtration of $E_K$. This implies that $M^\mu(E_K) \leq M^\mu(E_\kappa)$ for the stability function on $\cM_v$. If equality holds, then the filtration over the special fiber is a $\mu$-HN filtration of $E_\kappa$ in $\cM_v$, and it follows from \Cref{thm:HN_Bridgeland} that the filtration of $E$ lies in $\cM^\cF_v$.

\medskip
\noindent \textit{Upgrading from a weak $\Theta$-stratification to a $\Theta$-stratification:}
\medskip

In order to show that we have a $\Theta$-stratification, it suffices, by \Cref{lem:tangent_spaces}, to show that for any finite extension $k'/k$ and any point $E_\bullet \in \Filt(\cM_v^\cF)(k')$ corresponding to a \emph{split} Harder-Narasimhan filtration (with its canonical weights), the map of vector spaces $\op{Lie}(\op{Aut}_{\Filt(\cM_v)}(E_\bullet)) \to \op{Lie}(\op{Aut}_{\cM_v}(E))$, which is simply the restriction map
\[
\Hom_{X\times \Theta_{k'}}(E_\bullet,E_\bullet)\to \Hom_{X_k'}(E,E),
\]
is surjective (see \cite{behrend1995semi} for a version of this argument). Let us write $E = F_1 \oplus \cdots \oplus F_n$, where $F_i$ is semistable and the phase $\phi(F_i)$ is increasing in $i$. Then the HN filtration $E_\bullet$ of $E$ has $E_i = \bigoplus_{j \geq i} F_j$ with $\bZ$-weights $w_1<\cdots<w_n$ that are positively proportional to the slopes $\nu_i = d_i/r_i$ of $F_i$. The $\bZ$-weighted filtration $E_\bullet$ corresponds to the complex
\[
\bigoplus_i F_i \otimes_{k'} k'[t] \cdot t^{-w_i} \subset E \otimes_{k'} k'[t^\pm] \in \D^b(X \times \Theta_{k'}).
\]
Then $\Hom_{X\times \Theta_{k'}}(E_\bullet,E_\bullet) \simeq \bigoplus_{i\leq j} \Hom_{X_{k'}}(F_i,F_j)$, and the surjectivity of the map above is equivalent to the claim that $\Hom_{X_{k'}}(F_i,F_j) = 0$ for $i>j$. This follows from \Cref{lem:truncations}.

\end{proof}

\begin{ex} We observed in \Cref{ex:toda_1} that the connected component of the space of Bridgeland stability conditions on a $K3$ surface $X$ over $\bC$ constructed in \cite{bridgeland2008stability} contains a dense set of stability conditions for which $Z(\cA) \subset \bQ+i\bQ$ and $\cA$ satisfies the Generic Flatness property. The results of \cite{toda2008moduli}*{Sect.~4} imply that these Bridgeland stability conditions also satisfy the Boundedness of Quotients property, by the following:

\begin{lem}[\cite{toda2008moduli}*{Lem.~3.15}] \label{lem:alternate_boundedness}
Let $(Z,\cA)$ be a Bridgeland stability condition on $\D^b(X)$. Assume that the functor $\cM_v^{\rm ss}$ is bounded for every class $v \in \Knum(X)$. Then $(Z,\cA)$ satisfies the Boundedness of Quotients property.
\end{lem}

Therefore, a dense set of Bridgeland stability conditions in this connected component satisfy the hypotheses of \Cref{thm:theta_stratification_torsion_free}. This gives examples of $\Theta$-stratifications on the non-finite-type algebraic $k$-stacks $\cM^\cF_v$, which are not known to admit Zariski covers by quotient stacks. These stacks therefore lie beyond the scope of geometric invariant theory.
\end{ex}

\begin{rem}
One can extend the $\Theta$-stratification of \Cref{thm:theta_stratification_torsion_free} to a $\Theta$-stratification of $\cM_v$, where the HN filtration of every unstable point corresponds to some choice of weights for the Harder-Narasimhan filtration in the sense of slope semistability. For every Harder-Narasimhan polytope whose last subobject is torsion, simply choose a weight for this maximal torsion subobject that is greater than the weights appearing in the HN filtration for the torsion-free quotient. Using \Cref{thm:main_stratification} one can show that this defines a $\Theta$-stratification of $\cM_v$, but it is not clear if there is a numerical invariant on $\cM_v$ that allows one to assign \emph{canonical} rational weights to the Harder-Narasimhan filtration of any $E \in \cA$. More speculatively, it would be interesting if a Bridgeland stability condition allowed one to define a structure analogous to a $\Theta$-stratification on the higher derived stack of all relatively perfect complexes on $X$, as opposed to just those in $\D^b(X)^\heartsuit$.
\end{rem}

\begin{rem}
If $Z(\cA)$ does not lie in $\bQ + i \bQ$, the the cohomology classes \eqref{eqn:classes_Bridgeland} are not rational, and the unique maximizer of the numerical invariant $\mu$ in \Cref{thm:HN_Bridgeland} need not lie at a rational point. One could perturb the canonical assignment of real weights of \Cref{lem:explicit_maximum} for each HN polytope so that they are rational, but a more satisfying solution would be to develop an algebro-geometric interpretation for filtrations weighted by real numbers.

A similar issue arises when trying to recover the full Gieseker-Harder-Narasimhan filtration of an unstable coherent sheaf on a projective scheme $X$. Once could consider a central charge valued in polynomials $Z : K_0(\D^b(X)) \to \bC[n]$ as in \cite{bayer2009polynomial}, leading to cohomology classes on $\cM = \underline{\op{Coh}}(X)$ with values in $\bQ[n]$. The natural generalization of the second class in \eqref{eqn:classes_Bridgeland} is
\[
b = 2 \op{ch}_2((p_1)_\ast(\cE \boxtimes L^n)) \in H^4(\cM;\bQ[n]),
\]
where $p_1 : \underline{\op{Coh}}(X) \times X \to \underline{\op{Coh}}(X)$ is the projection and $\cE$ is the universal coherent sheaf on the product $\underline{\op{Coh}}(X) \times X$. This class is positive definite in the sense that for any filtration $f : \Theta_k \to \cM$, $f^\ast(b) \in \bQ[n]$ will be a polynomial that takes positive values for $n\gg0$, because the higher cohomology of $R\Gamma(X, E \otimes L^n)$ vanishes for $n\gg 0$. The numerical invariant $\mu(f) = f^\ast(l) / \sqrt{f^\ast(b)}$ can then be regarded as taking values in the fraction field $F$ of the ring of germs of analytic functions in one real variable $n$ in a neighborhood of $n=\infty$. A similar approach is carried out in \cite{halpernleistner2021moduli}, but the resulting filtrations are coarser than the full Gieseker-Harder-Narasimhan filtration. To get the full filtration, one should allow the weights to take values in $F$ as well, which would require developing a geometric theory of filtrations weighted by $F$ and a generalization of the theory of formal fans developed here.

\end{rem}


\appendix

\section{Proof of \texorpdfstring{\Cref{thm:describe_strata_global_quotient}}{Theorem}}
\label{appendix:filtrations_global_quotient}
The piece of the proof of \Cref{thm:describe_strata_quotient_GLN} which reduces the general case of $X/G$ to the case of $\pt/G$ applies verbatim. It therefore suffices to show that if $G$ is a smooth algebraic $k$ group and $T \subset G$ a split maximal torus, then the canonical $P_\psi$ bundles $E_\psi$ on $\Theta_k$ for any homomorphism $\psi : \Gm^n \to T$ induce equivalences
\begin{gather*}
\bigsqcup_{\psi \in \Hom(\Gm^n,T) / W} \pt / P_\psi \xrightarrow{\simeq} \Filt^n(\pt/G) \quad \text{and} \\
\bigsqcup_{\psi \in \Hom(\Gm^n,T)/W} \pt / L_\psi \xrightarrow{\simeq} \Grad^n(\pt/G).
\end{gather*}
We will prove this in the case where $n=1$. The case for general $n$ follows from a similar argument. Alternatively, using the fact that $T \subset P_\psi$ will be a split maximal torus for all $\psi \in \Hom(\Gm^n,T)$, one can prove the general statement inductively using the fact that
\[
\filt[n]{\pt/G} \simeq \filt{\filt{\ldots,\filt{\pt / G}}}\ldots),
\]
and likewise for $\grad[n]{\pt/G}$.

\begin{prop} \label{prop:parabolic_reduction}
Let $k$ be a perfect field, let $S$ be a connected finite type $k$-scheme, and let $E$ be a $G$-bundle over $\Theta_S := \Theta \times S$. Let $s \in S(k)$, thought of as the point $(0,s) \in \bA^1_S$, and assume that $\op{Aut}(E_s) \simeq G$. Let $\lambda : \Gm \to G$ be a 1PS conjugate to the one parameter subgroup $\Gm \to \op{Aut} (E_s)$. Then
\begin{enumerate}
\item There is a unique reduction of structure group $E^\prime \subset E$ to a $P_\lambda$-torsor such that $\Gm \to \op{Aut}(E^\prime_s)\simeq P_\lambda$ is conjugate in $P_\lambda$ to $\lambda$, and
\item the restriction of $E^\prime$ to $\{1\}\times S$ is canonically isomorphic to the sheaf on the \'{e}tale site of $S$ mapping $T / S \mapsto \op{Iso}((E_\lambda)_{\Theta_T}, E|_{\Theta_T})$.
\end{enumerate}
\end{prop}

\begin{proof}
$(E_\lambda)_{\Theta_S} = E_\lambda \times S / \Gm$ is a $G$-bundle over $\Theta_S$, and $\inner{\op{Iso}}((E_\lambda)_{\Theta_S},E)$ is a sheaf over $\Theta_S$ representable by a (relative) scheme over $\Theta_S$. In fact, if we define a twisted action of $\Gm$ on $E$ given by $t \star e := t \cdot e \cdot \lambda(t)^{-1}$, then
\begin{equation} \label{eqn:parabolic_reduction_isomorphism}
\inner{\op{Iso}}((E_\lambda)_{\Theta_S},E) \simeq E/ \Gm \text{ w.r.t. the }\star \text{-action}
\end{equation}
as sheaves over $\Theta_S$.\footnote{To see this, note that a map $T \to \Theta_S$ corresponds to a $\Gm$-bundle $P \to T$ along with a $\Gm$ equivariant map $f : P \to \bA^1 \times S$. Then the restrictions $(E_\lambda \times S)_T$ and $E|_T$ correspond (via descent for $G$-bundles) to the $\Gm$-equivariant bundles $f^{-1}(E_\lambda \times S)$ and $f^{-1}E$ over $P$. Forgetting the $\Gm$-equivariant structure, the $G$-bundle $E_\lambda \times S$ is trivial, so an isomorphism $f^{-1}(E_\lambda \times S) \to f^{-1}E$ as $G$-bundles corresponds to a section of $f^{-1} E$, or equivalently to a lifting $$\xymatrix{ & E \ar[d] \\ P \ar[r]^f \ar@{-->}[ur]^{\tilde{f}} & \bA^1 \times S}$$ to a map $\tilde{f} : P \to E \to \bA^1 \times S$. The isomorphism of $G$-bundles defined by the lifting $\tilde{f}$ descends to an isomorphism of $\Gm$-equivariant $G$-bundles $f^{-1} (E_\lambda \times S) \to f^{-1}E$ if and only if the lift $\tilde{f}$ is equivariant with respect to the twisted $\Gm$ action on $E$.}

The twisted $\Gm$ action on $E$ is compatible with base change. Let $T \to S$ be an $S$-scheme. From the isomorphism of sheaves \eqref{eqn:parabolic_reduction_isomorphism}, there is a natural bijection between the set of isomorphisms $(E_\lambda)_{\bA^1_T} \xrightarrow{\simeq} E|_{\bA^1_T}$ as $\Gm$-equivariant $G$-bundles and the set of $\Gm$-equivariant sections of $E|_{\bA^1_T} \to \bA^1_T$ with respect to the twisted $\Gm$ action.

The morphism $E|_{\bA^1_T} \to \bA^1_T$ is separated, so a twisted equivariant section is uniquely determined by its restriction to $\Gm \times T$, and by equivariance this is uniquely determined by its restriction to $\{1\}\times T$. Thus we can identify $\Gm$-equivariant sections with the set of maps $T \to E$ such that $\lim \limits_{t \to 0} t \star e$ exists and $T \to E \to \bA^1_S$ factors as the given morphism $T \to \{1\} \times S \to \bA^1_S$.

If we define the subsheaf of $E$ over $\bA^1_S$
$$E^\prime (T) := \left\{ e \in E(T) | \Gm \times T \xrightarrow{t\star e(x)} E \text{ extends to } \bA^1 \times T \right\} \subset E(T),$$
then we have shown that $E^\prime|_{\{1\}\times S} (T) \simeq \op{Iso} ( (E_\lambda)_{\Theta_T}, E|_{\Theta_T} )$. Next we show in several steps that the subsheaf $E^\prime \subset E$ over $\bA^1_S$ is a torsor for the subgroup $P_\lambda \subset G$, so $E^\prime$ is a reduction of structure group to $P_\lambda$.
\begin{itemize}
\setlength{\itemindent}{12pt}
\setlength{\itemsep}{5pt}
\item[Step 1:] \textit{$E^\prime$ is representable:} The functor $E^\prime$ is exactly the functor of \Cref{cor:contraction_functor}, so \Cref{prop:Hesselink_concentration} implies that $E^\prime$ is representable by a disjoint union of $\Gm$ equivariant locally closed subschemes of $E$.

\item[Step 2:] \textit{$P_\lambda \subset G$ acts simply transitively on $E^\prime \subset E$:} Because $E$ is a $G$-bundle over $\bA^1_S$, right multiplication $(e,g) \mapsto (e,e\cdot g)$ defines an isomorphism $E \times G \to E \times_{\bA^1_S} E$. The latter has a $\Gm$ action, which we can transfer to $E \times G$ using this isomorphism.

For $g \in G(T)$, $e \in E(T)$, and $t \in \Gm(T)$ we have $t \star (e \cdot g) = (t \star e) \cdot (\lambda(t) g \lambda(t)^{-1})$. This implies that the $\Gm$ action on $E \times G$ corresponding to the diagonal action on $E \times_{\bA^1_S} E$ is given by
$$t\cdot (e,g) = (t \star e, \lambda(t) g \lambda(t)^{-1})$$
The subfunctor of $E \times G$ corresponding to $E^\prime \times_{\bA^1_S} E^\prime \subset E \times_{\bA^1_S} E$ consists of those points for which $\lim_{t\to 0} t \cdot (e,g)$ exists. This is exactly the subfunctor represented by $E^\prime \times P_\lambda \subset E \times G$. We have thus shown that $E^\prime$ is equivariant for the action of $P_\lambda$, and $E^\prime \times P_\lambda \to E^\prime \times_{\bA^1_S} E^\prime$ is an isomorphism of sheaves.

\item[Step 3:] \textit{$p : E^\prime \to \bA^1_S$ is smooth:} \Cref{prop:stratum_regularity} implies that $E^\prime$ and $E^{\Gm} \subset E^\prime$ are both smooth over $S$. The restriction of the tangent bundle $T_{E / S}|_{E^{\Gm}}$ is an equivariant locally free sheaf on a scheme with trivial $\Gm$ action, hence it splits into a direct sum of vector bundles of fixed weight with respect to $\Gm$. The tangent sheaf $T_{E^\prime / S}|_{E^{\Gm}} \subset T_{E/S}|_{E^{\Gm}}$ is precisely the subsheaf with weight $\geq 0$. By hypothesis $T_{E/S} \to p^\ast T_{\bA^1_S / S}$ is surjective, and $p^\ast T_{\bA^1_S / S}|_{E^{\Gm}}$ is concentrated in nonnegative weights, therefore the map
$$T_{E^\prime / S}|_{E^{\Gm}} = ( T_{E / S}|_{E^{\Gm}} )_{\geq 0} \to p^\ast T_{\bA^1_S / S}|_{E^{\Gm}} = (p^\ast T_{\bA^1_S / S}|_{E^{\Gm}})){\geq 0}$$
is surjective as well.

Thus we have shown that $T_{E^\prime/S} \to p^\ast T_{\bA^1_S / S}$ is surjective when restricted to $E^{\Gm} \subset E^\prime$, and by Nakayama's Lemma it is also surjective in a Zariski neighborhood of $E^{\Gm}$. On the other hand, the only equivariant open subscheme of $E^\prime$ containing $E^{\Gm}$ is $E^\prime$ itself. It follows that $T_{E^\prime / S} \to p^\ast T_{\bA^1_S / S}$ is surjective, and therefore that the morphism $p$ is smooth.

\item[Step 4:] \textit{$p : E^\prime \to \bA^1_S$ admits sections \'{e}tale locally:} We consider the $\Gm$ equivariant $G$-bundle $E|_{\{0\}\times S}$. After \'{e}tale base change we can assume that $E|_{\{0\}\times S}$ admits a non-equivariant section, hence the $\Gm$-equivariant structure is given by a homomorphism $(\Gm)_{S^\prime} \to G_{S^\prime}$. \Cref{lem:rigidity_1PS} implies that after further \'{e}tale base change this homomorphism is conjugate to a constant homomorphism. Thus $E|_{ \{0\}\times S^\prime}$ is isomorphic to the trivial equivariant $\Gm$-bundle $(E_\lambda)_{\bA^1_S} = \bA^1_S \times G \to \bA^1_S$ with $\Gm$ acting by left multiplication by $\lambda(t)$.

It follows that $E|_{ \{0\}\times S^\prime}$ admits an invariant section with respect to the twisted $\Gm$ action. In other words $(E_{\bA^1_{S^\prime}})^{\Gm} \to \{0\} \times S^\prime$ admits a section, and $E^{\Gm} \subset E^\prime$, so we have shown that $E^\prime \to \bA^1_{S^\prime}$ admits a section over $\{0\}\times S^\prime$. On the other hand, because $p : E^\prime \to \bA^1_{S^\prime}$ is smooth and $\Gm$-equivariant, the locus over which $p$ admits an \'{e}tale local section is open and $\Gm$-equivariant. It follows that $p$ admits an \'{e}tale local section over every point of $\bA^1_{S^\prime}$.
\end{itemize}

\end{proof}

\begin{rem}
In fact we have shown something slightly stronger than the existence of \'{e}tale local sections of $E^\prime \to \bA^1_S$ in Step 4. We have shown that there is an \'{e}tale map $S^\prime \to S$ such that $E^\prime|_{S^\prime} \to \bA^1_{S^\prime}$ admits a $\Gm$-equivariant section.
\end{rem}

The following fact, that families of one parameter subgroups of $G$ are \'{e}tale locally constant up to conjugation, was the key input to the proof of \Cref{prop:parabolic_reduction}.
\begin{lem} \label{lem:rigidity_1PS}
Let $S$ be a connected scheme of finite type over a perfect field $k$, let $G$ be a smooth affine $k$-group, and let $\phi : (\Gm)_S \to G_S$ be a homomorphism of group schemes over $S$. Let $\lambda : \Gm \to G$ be a 1PS conjugate to $\phi_s$ for some $s \in S(k)$. Then the subsheaf
$$F(T) = \left\{ g \in G(T) | \phi_{T} = g \cdot (\op{id}_T, \lambda) \cdot g^{-1} : (\Gm)_T \to G_T \right\} \subset G_S(T)$$
is an $L_\lambda$-torsor. In particular $\phi$ is \'{e}tale-locally conjugate to a constant homomorphism.
\end{lem}
\begin{proof}
Verifying that $F \times L_\lambda \to F \times_S F$ given by $(g,l) \mapsto (g,gl)$ is an isomorphism of sheaves is straightforward. The more important question is whether $F(T) \neq \emptyset$ \'{e}tale locally.

As in the proof of \Cref{prop:parabolic_reduction} we introduce a twisted $\Gm$ action on $G\times S$ by $t \star (g,s) = \phi_s(t) \cdot g \cdot \lambda(t)^{-1}$. Then $G \times S \to S$ is $\Gm$ invariant, and the functor $F(T)$ is represented by the map of schemes $(G\times S)^{\Gm} \to S$. By \Cref{prop:stratum_regularity}, $(G \times S)^{\Gm} \to S$ is smooth, and in particular it admits a section after \'{e}tale base change in a neighborhood of a point $s \in S(k)$ for which $(G \times S)^{\Gm}_s = (G \times \{s\})^{\Gm} \neq \emptyset$. By construction this set is nonempty precisely when $\phi_s$ is conjugate to $\lambda$, so by hypothesis it is nonempty for some $s \in S(k)$.

By the same reasoning, for every finite separable extension $k^\prime / k$, every $k^\prime$-point has an \'{e}tale neighborhood on which $\phi$ is conjugate to a constant homomorphism determined by some one parameter subgroup defined over $k^\prime$. Because $k$ is perfect and $S$ is locally finite type, we have a cover of $S$ by Zariski opens over each of which $\phi$ is \'{e}tale-locally conjugate to a constant homomorphism determined by a one parameter subgroup defined over some finite separable extension of $k$. Because $S$ is connected, all of these 1PS's are conjugate to $\lambda$, possibly after further finite separable field extensions. Thus $(G \times S)^{\Gm} \to S$ admits a global section after \'{e}tale base change.
\end{proof}

Because $T$ is split, the set of one parameter subgroups up to conjugacy by $W$ is unaffected by passing to the algebraic closure, so it suffices to assume that $k = \bar{k}$. \Cref{prop:parabolic_reduction} implies that the map $\bigsqcup_{\Hom(\Gm,T)/W} \pt/P_\lambda \to \Filt(\pt/G)$ is an equivalence over the sub-site of $k$-schemes of finite type. The functor $\Filt(\pt/G)$ is limit preserving by the formal observation
$$\Map(\varprojlim_i T_i , \filt{\pt/G}) \simeq \varprojlim_i \Map(T_i \times \Theta,\pt/G)$$
where the last equality holds because $\pt/G$ is an algebraic stack locally of finite presentation. The stack $\bigsqcup_{I} \pt / P_\lambda$ is locally of finite presentation and thus limit preserving as well. Every affine scheme over $k$ can be written as a limit of finite type $k$ schemes, so the isomorphism for finite type $k$-schemes implies the isomorphism for all $k$-schemes.

To show that $\bigsqcup_{\lambda \in \Hom(\Gm,T)/W} \pt / L_\lambda \to \Grad(\pt / G)$ is an equivalence, note that for $S$ of finite type over $k = \bar{k}$, the proof of \Cref{prop:parabolic_reduction} carries over unchanged for $G$-bundles over $(\pt / \Gm) \times S$, showing that \'{e}tale locally in $S$ they are isomorphic to $S \times G$ with $\Gm$ acting by left multiplication by $\lambda(t)$ on $G$. In fact, we had to essentially prove this when we considered the $\Gm$-equivariant bundle $E|_{\{0\} \times S}$ in Step 4 of that proof. The amplification of the statement from finite type $\bar{k}$-schemes to all $\bar{k}$-schemes, and from $\bar{k}$-schemes to $k$-schemes, is identical to the argument for $\Filt(\pt/G)$ above. This completes the proof of \Cref{thm:describe_strata_global_quotient}.

\begin{rem}
The argument above does not a priori use the fact that $\Filt(\pt/G)$ is algebraic. It thus provides an alternative proof for the algebraicity of $\Filt(X/G)$ when $G$ is a split algebraic $k$-group and $X$ is a $G$-quasi-projective scheme which is independent of the general results of \cite{HLPreygel}.
\end{rem}

\bibliographystyle{plain}
\bibliography{refs}{}

\end{document}